\documentclass[12pt,leqno]{article}
\usepackage{amssymb}
\usepackage[mathscr]{eucal}
\usepackage{amsmath,amssymb,latexsym,theorem,bbm}
\newcommand{\DS}{\displaystyle}
\newcommand{\TS}{\textstyle}

\newcommand{\CC}{\mathbb{C}}

\newcommand{\EE}{\mathsf{E}}

\newcommand{\NN}{\mathbb{N}}
\newcommand{\PP}{\mathsf{P}}

\newcommand{\RR}{\mathbb{R}}

\newcommand{\ZZ}{\mathbb{Z}}

\newcommand{\cA}{{\mathcal A}}

\newcommand{\cF}{{\mathcal F}}

\newcommand{\cL}{{\mathcal L}}
\newcommand{\cN}{{\mathcal N}}

\newcommand{\dd}{\mathrm{d}}

\newcommand{\sums}{\sideset{}{^{(s)}}\sum_{k=1}^n}
\newcommand{\sumsf}{\sideset{}{^{(s+1)}}\sum_{k=1}^n}
\newcommand{\sumss}{\sideset{}{^{(s,s+1)}}\sum_{k=1}^n}
\newcommand{\sumssd}{\sideset{}{^{(s_1,s_2)}}\sum_{k=1}^n}

\newcommand{\var}{\operatorname{Var}}
\newcommand{\cov}{\operatorname{Cov}}

\newcommand{\hmuen}{\widehat\mu_{\vare,n}}
\newcommand{\halpha}{\widehat{\alpha}}
\newcommand{\htheta}{\widehat{\theta}}
\newcommand{\talpha}{\widetilde{\alpha}}
\newcommand{\ttheta}{\widetilde{\theta}}

\renewcommand{\mid}{\,|\,}

\renewcommand{\leq}{\leqslant}
\renewcommand{\geq}{\geqslant}

\newcommand{\stoch}{\stackrel{\PP}{\longrightarrow}}
\newcommand{\distr}{\stackrel{\cL}{\longrightarrow}}

\newcommand{\mean}{\stackrel{L_1}{\longrightarrow}}

\newcommand{\lpmean}{\stackrel{L_p}{\longrightarrow}}
\newcommand{\bone}{\mathbbm{1}}
\newcommand{\vare}{\varepsilon}
\newcommand{\proofend}{\hfill\mbox{$\Box$}}

\numberwithin{equation}{subsection}

\theoremstyle{change} \theorembodyfont{\em}
\newtheorem{Lem}{Lemma.}[subsection]
\newtheorem{Thm}{Theorem.}[subsection]
\newtheorem{Pro}{Proposition.}[subsection]

\newtheorem{Def}{Definition.}[subsection]
\newtheorem{Lem2}{Lemma.}[section]

\theorembodyfont{\rm}
\newtheorem{Rem}{Remark.}[subsection]

\begin{document}

\begin{center}
 {\bfseries\Large Outliers in INAR(1) models} \\[5mm]

 {\sc\large M\'aty\'as $\text{Barczy}^{*,\diamond}$, \ M\'arton $\text{Isp\'any}^*$, Gyula $\text{Pap}^{**}$, \\
    Manuel $\text{Scotto}^{***}$} {\large and} {\sc\large Maria Eduarda $\text{Silva}^{****}$}
\end{center}

\vskip0.2cm

\noindent * University of Debrecen, Faculty of Informatics, Pf.~12, H--4010 Debrecen, Hungary;\\
** University of Szeged, Bolyai Institute, H-6720 Szeged, Aradi v\'ertan\'uk tere 1, Hungary;
*** Universidade de Aveiro, Departamento de Matem\'atica, Campus Universit\'ario de Santiago, 3810-193 Aveiro,
    Portugal;\\
**** Universidade do Porto, Faculdade de Economia, Rua Dr.~Roberto Frias s/n, 4200 464 Porto, Portugal.

\noindent e--mails: barczy.matyas@inf.unideb.hu (M. Barczy), ispany.marton@inf.unideb.hu (M. Isp\'any),
          papgy@math.u-szeged.hu (G. Pap), mscotto@ua.pt (M. Scotto), mesilva@fep.up.pt (M. E. Silva).

$\diamond$ Corresponding author.



\renewcommand{\thefootnote}{}
\vspace*{0.2cm}
\footnote{The authors have been supported by the Hungarian
 Portuguese Intergovernmental S \& T Cooperation Programme for 2008-2009 under Grant No.\ PT-07/2007.
 M\'aty\'as Barczy and Gyula Pap were partially supported by the Hungarian Scientific
 Research Fund under Grants No.~OTKA-T-079128.}

\vspace*{-10mm}

\begin{abstract}
In this paper the integer-valued autoregressive model of order one,
 contaminated with additive or innovational outliers is studied in some detail.
Moreover, parameter estimation is also addressed.
Supposing that the time points of the outliers are known but their sizes are
 unknown, we prove that the Conditional Least Squares (CLS) estimators of the
 offspring and innovation means are strongly consistent.
In contrast, however, the CLS estimators of the outliers' sizes are
 not strongly consistent, although they converge to a random limit with
 probability 1.
This random limit depends on the values of the process at the outliers' time points and on the values
 at the preceding time points and in case of additive outliers also on the values at the
 following time points.
We also prove that the joint CLS estimator of the offspring and innovation
 means is asymptotically normal.
Conditionally on the above described values of the process, the joint CLS estimator of the sizes of
 the outliers is also asymptotically normal.
\end{abstract}

\vfill

\eject

\tableofcontents

\eject

\section{Introduction}

Recently, there has been considerable interest in integer-valued
 time series models and a sizeable volume of work is now available
 in specialized monographs (e.g., MacDonald and Zucchini \cite{MacZuc},
 Cameron and Trivedi \cite{CamTri}, and Steutel and van Harn \cite{SteHar2}) and
 review papers (e.g., McKenzie \cite{McK2}, Jung and Tremayne \cite{JunTre}, and
 Wei{\ss} \cite{Wei}).
Motivation to include discrete data models comes from the need to account for the
 discrete nature of certain data sets, often counts of events, objects or individuals.
Examples of applications can be found in the analysis of time series of count
 data that are generated from stock transactions (Quoreshi \cite{Quo}),
 where each transaction refers to a trade between a buyer and a seller in a
 volume of stocks for a given price, in optimal alarm systems
 to predict whether a count process will upcross a certain level and give an
 alarm whenever the upcrossing level is predicted (Monteiro, Pereira and Scotto
 \cite{MonPerSco}), international tourism demand (Br\"{a}nn\"{a}s and
 Nordstr\"{o}m \cite{BraNor}), experimental biology (Zhou and Basawa \cite{ZhoBas}),
 social science (McCabe and Martin \cite{McCMar}), and queueing systems
 (Ahn, Gyemin and Jongwoo \cite{AhnGyeJon}).

Several integer-valued time series models were proposed in the
 literature, we mention the INteger-valued AutoRegressive model of order $p$ \
 (INAR($p$)) and the INteger-valued Moving Average model of order $q$ (INMA($q$)).
The former was first introduced by McKenzie \cite{McK} and Al-Osh and Alzaid \cite{AloAlz}
 for the case $p=1$.
The INAR(1) and INAR($p$) models have been investigated by several authors, see, e.g.,
 Silva and Oliveira \cite{SilOli1}, \cite{SilOli2}, Silva and Silva \cite{SilSil},
 Isp\'any, Pap and van Zuijlen \cite{IspPapZui1}, \cite{IspPapZui2},
 \cite{IspPapZui3}, Drost, van den Akker and Werker \cite{DroAkkWer4} (local asymptotic normality
 for INAR($p$) models), Gy\"orfi, Isp\'any, Pap and Varga \cite{GyorIspPapVar},
 Isp\'any \cite{Isp}, Drost, van den Akker and Werker \cite{DroAkkWer2}, \cite{DroAkkWer3}
 (nearly unstable INAR(1) models and semiparametric INAR($p$) models),
 Bu and McCabe \cite{BuMcCab} (model selection) and Andersson and Karlis \cite{AndKar} (missing values).
Empirical relevant extensions have been suggested by Br\"{a}nn\"{a}s \cite{Bra}
 (explanatory variables), Blundell \cite{BluGriWin} (panel data), Br\"{a}nn\"{a}s and Hellstr\"{o}m
 \cite{BraHel} (extended dependence structure), and more recently by Silva, Silva, Pereira
  and Silva \cite{SilSilPerSil} (replicated data) and by Wei{\ss} \cite{Wei3} (combined INAR($p$) models).
Extensions and generalizations were proposed by Du and Li \cite{DuLi} and Latour \cite{Lat}.
The INMA($q$) model was proposed by Al-Osh and Alzaid \cite{AloAlz2}, and subsequently studied by
 Br\"{a}nn\"{a}s and Hall \cite{BraHal} and Wei{\ss} \cite{Wei2}.
Related models were introduced by Aly and Bouzar \cite{AlyBou}, \cite{AlyBou2}, Zhu and Joe \cite{ZhuJoe}
 and more recently by Wei{\ss} \cite{Wei2}.
Extensions for random coefficients integer-valued autoregressive
 models have been proposed by Zheng, Basawa and Datta \cite{ZheBasDat}, \cite{ZheBasDat2}
 who investigated basic probabilistic and statistical properties of these models.
Zheng and co-authors illustrated their performance in the analysis of epileptic seizure counts
 (e.g., Latour \cite{Lat}) and in the analysis of the monthly number cases of poliomyelitis in
 the US for the period 1970-1983.
Doukhan, Latour and Oraichi \cite{DouLatOra} introduced the class of non-negative integer-valued
 bilinear time series, some of their results concerning the existence of stationary solutions
 were extended by Drost, van den Akker and Werker \cite{DroAkkWer1}.
Recently, the so called $p$-order Rounded INteger-valued AutoRegressive (RINAR($p$)) time series model
 was introduced and studied by Kachour and Yao \cite{KacYao} and Kachour \cite{Kac}.

 Moreover, topics of major current interest in time series modeling are to detect outliers in sample
 data and to investigate the impact of outliers on the estimation of conventional ARIMA models.
Motivation comes from the need to assess for data quality and to the robustness of subsequent
 statistical analysis in the presence of discordant observations.
Fox \cite{Fox} introduced the notion of additive and innovational outliers and proposed the use
 of maximum likelihood ratio test to detect them.
Chang and Chen \cite{ChaChe} extended Fox's results to ARIMA models and proposed a likelihood ratio test
 and an iterative procedure for detecting outliers and estimating the model parameters.
Some generalizations were obtained by Tsay \cite{Tsa} for the detection of level shifts and
 temporary changes.
Random level shifts were studied by Chen and Tiao \cite{CheTia}.
Extensions of Tsay's results can be found in Balke \cite{Bal}.
Abraham and Chuang \cite{AbrChu} applied the EM algorithm to the estimation of outliers.
Other useful references for outlier detection and estimation in time series models are Guttman
 and Tiao \cite{GutTia}, Bustos and Yohai \cite{BusYoh}, McCulloch and Tsay \cite{McCTsa},
 Pe\~{n}a \cite{Pen}, S\'{a}nchez and Pe\~{n}a \cite{SanPen}, Perron and Rodriguez \cite{PerRod}
 and Burridge and Taylor \cite{BurTay}.

It is worth mentioning that all references given in the previous
 paragraph deal with the case of continuous-valued processes.
A general motivation for studying outliers for integer-valued time series can be the fact that
 it may often difficult to remove outliers in the integer-valued case, and hence
 an important and interesting problem, which has not yet been addressed, is to investigate the impact
 of outliers on the parameter estimation of series of counts which are represented through
 integer-valued autoregressive models.
This paper aims at giving a contribution towards this direction.
A more specialized motivation is the possibility of potential applications, for example
 in the field of statistical process control (a good description of this topic can be found
 in Montgomery \cite[Chapter 4, Section 3.7]{Mon}).
In this paper we consider the problem of Conditional Least Squares (CLS) estimation of
 some parameters of the INAR(1) model contaminated with additive or innovational outliers
 starting from a general initial distribution (having finite second or third moments).
We suppose that the time points of the outliers are known, but their sizes are unknown.
Under the assumption that the second moment of the innovation distribution is finite, we prove that the
 CLS estimators of the means of the offspring and innovation distributions are
 strongly consistent, but the joint CLS estimator of the sizes of the outliers is not strongly consistent; nevertheless, it converges to a random limit with probability 1.
This random limit depends on the values of the process at the outliers' time points and on the values
 at the preceding time points and in case of additive outliers also on the values at the
 following time points.
Under the assumption that the third moment of the innovation distribution is finite, we prove that the
 joint CLS estimator of the means of the offspring and innovation distributions
 is asymptotically normal with the same asymptotic variance as in the case when there are no outliers.
Conditionally on the above described values of the process, the joint CLS estimator of the sizes of
 the outliers is also asymptotically normal.
We calculate its asymptotic covariance matrix as well.
In this paper we present results in the case of one or two additive or innovational outliers for INAR(1) models,
 the general case of finitely many additive or innovational outliers may be handled in a similar way,
 but we renounce to consider it.

The rest of the paper is organized as follows.
Section \ref{INAR1_section} provides a background description of basic theoretical results
 related with the asymptotic behavior of CLS estimator for the INAR(1) model.
In Sections \ref{INAR1_additive} and \ref{INAR1_innovational} we consider
 INAR(1) models contaminated with one or two additive or innovational outliers, respectively.
The cases of one outlier and two outliers are handled separately.
Section \ref{Appendix} is an appendix containing the proofs of some auxiliary results.

\section{The INAR(1) model}\label{INAR1_section}

\subsection{The model and some preliminaries}

Let \ $\ZZ_+$ \ and \ $\NN$ \ denote the set of non-negative integers and
 positive integers, respectively.
Every random variable will be defined on a fixed probability space \ $(\Omega,\cA,\PP)$.

One way to obtain models for integer-valued data is replacing multiplication in the conventional
 ARMA models in order to ensure the integer discreteness of the process and to
 adopt the terms of self-decomposability and stability for integer-valued time series.

\begin{Def}
Let \ $(\vare_k)_{k\in\NN}$ \ be an independent and identically distributed (i.i.d.) sequence of
non-negative integer-valued random variables.
An INAR(1) time series model is a stochastic process \ $(X_n)_{n\in\ZZ_+}$ \ satisfying
 the recursive equation
 \begin{align}\label{INAR1}
   X_k = \sum\limits_{j=1}^{X_{k-1}}\xi_{k,j}+\vare_k, \qquad k\in\NN,
 \end{align}
 where for all \ $k\in\NN$, \ $(\xi_{k,j})_{j\in\NN}$ \ is a sequence of i.i.d.\ Bernoulli random
 variables with mean \ $\alpha\in[0,1]$ \ such that these sequences are mutually independent and
 independent of the sequence \ $(\vare_\ell)_{\ell\in\NN}$, \ and \ $X_0$ \ is a non-negative
 integer-valued random variable independent of the sequences \ $(\xi_{k,j})_{j\in\NN}$, \ $k\in\NN$,
 \ and \ $(\vare_\ell)_{\ell\in\NN}$.
\end{Def}

\begin{Rem}
The INAR(1) model in \eqref{INAR1} can be written in another way using
 the binomial thinning operator \ $\alpha\,\circ$ \ (due to Steutel and van Harn~\cite{SteHar})
 which we recall now.
Let \ $X$ \ be a non-negative integer-valued random variable.
Let \ $(\xi_j)_{j\in\NN}$ \ be a sequence of i.i.d.\ Bernoulli random variables with
 mean \ $\alpha\in[0,1]$.
\ We assume that the sequence \ $(\xi_j)_{j\in\NN}$ \ is independent of \ $X$.
\ The non-negative integer-valued random variable \ $\alpha\,\circ X$ \ is defined by
 \[
   \alpha\circ X
     :=\begin{cases}
        \sum\limits_{j=1}^X\xi_j, & \quad \text{if \ $X>0$},\\[2mm]
         0, & \quad \text{if \ $X=0$}.
       \end{cases}
 \]
The sequence \ $(\xi_j)_{j\in\NN}$ \ is called a counting sequence.
The INAR(1) model in \eqref{INAR1} takes the form
 \[
    X_k=\alpha\circ X_{k-1}+\vare_k,\qquad k\in\NN.
 \]
\proofend
\end{Rem}

In the sequel we always assume that \ $\EE X_0^2<\infty$ \ and that
 \ $\EE\vare_1^2<\infty$, \ $\PP(\vare_1\ne 0)>0$.
\ Let us denote the mean and variance of \ $\vare_1$ \ by \ $\mu_\vare$ \ and \ $\sigma_\vare^2$,
 \ respectively.
\ Clearly, \ $0<\mu_\vare<\infty$.

It is easy to show that
 \begin{align}\label{Seged_stationary1}
   \lim_{k\to\infty}\EE X_k=\frac{\TS\mu_\vare}{\TS1-\alpha},\qquad
     \lim_{k\to\infty}\var X_k=\frac{\TS\sigma_\vare^2+\alpha\mu_\vare}{\TS1-\alpha^2},\qquad
      \text{if \ $\alpha\in(0,1)$,}
 \end{align}
 and that \ $\lim_{k\to\infty}\EE X_k=\lim_{k\to\infty}\var X_k=\infty$ \ if \ $\alpha=1$
 \ (e.g., Isp\'any, Pap and van Zuijlen \cite[page 751]{IspPapZui1}).
\ The case \ $\alpha\in(0,1)$ \ is called \emph{stable} or \emph{asymptotically
 stationary}, whereas the case \ $\alpha=1$ \ is called \emph{unstable}.
For the stable case, there exists a unique stationary distribution of the INAR(1)
 model in \eqref{INAR1}, see Lemma \ref{LEMMA_UNIQUE_STATDISTR} in the Appendix.

In the sequel we assume that \ $\alpha\in(0,1)$, \ and we denote by \ $\cF_k^X$ \ the \ $\sigma$--algebra
 generated by the random variables \ $X_0,X_1,\ldots,X_k$.

\subsection{Estimation of the mean of the offspring distribution}

In this section we concentrate on the CLS estimation of the parameter
 \ $\alpha$.
\ Clearly, \ $\EE(X_k\mid \cF^X_{k-1})=\alpha X_{k-1}+\mu_\vare$, \ $k\in\NN$, \ and thus
 \begin{align*}
  \sum_{k=1}^n\big(X_k-\EE(X_k\mid \cF^X_{k-1})\big)^2
     =\sum_{k=1}^n\big(X_k-\alpha X_{k-1}-\mu_\vare\big)^2, \qquad n\in\NN.
 \end{align*}
For all \ $n\in\NN$, \ we define the function \ $Q_n:\RR^{n+1}\times\RR\to\RR$, \ as
 \begin{align*}
    Q_n(x_0,x_1,\ldots,x_n;\alpha')
      & :=\sum_{k=1}^n\big(x_k-\alpha' x_{k-1}-\mu_\vare\big)^2,
       \qquad x_0,x_1,\ldots,x_n,\alpha'\in\RR.
 \end{align*}
By definition, for all \ $n\in\NN$, \ a CLS estimator for
 the parameter \ $\alpha\in(0,1)$ \ is a measurable function
 \ $\talpha_n:\RR^{n+1}\to\RR$ \ such that
 \begin{align*}
   Q_n(x_0,x_1,\ldots,x_n;\,&\talpha_n(x_0,x_1,\ldots,x_n))\\
       &= \inf_{\alpha'\in\RR}Q_n(x_0,x_1,\ldots,x_n;\alpha')
       \qquad \forall\;\;  (x_0,x_1,\ldots,x_n)\in\RR^{n+1}.
 \end{align*}
It is well-known that
 \begin{align}\label{SEGED91}
    \widetilde\alpha_n(X_0,X_1,\ldots,X_n)
       =\frac{\sum_{k=1}^n(X_k-\mu_\vare)X_{k-1}}{\sum_{k=1}^nX_{k-1}^2}
 \end{align}
 holds asymptotically as \ $n\to\infty$ \ with probability one.
Hereafter by the expression {\sl `a property holds asymptotically as \ $n\to\infty$ \
 with probability one'} we mean that there exists an event \ $S\in\cA$ \ such that \ $\PP(S)=1$
 \ and for all \ $\omega\in S$ \ there exists an \ $n(\omega)\in\NN$ \ such that
 the property in question holds for all \ $n\geq n(\omega)$.
\ The reason why \eqref{SEGED91} holds only asymptotically as \ $n\to\infty$ \ with
 probability one and not for all \ $n\in\NN$ \ and \ $\omega\in\Omega$ \ is that for all
 \ $n\in\NN$, \ the probability that the denominator \ $\sum_{k=1}^nX_{k-1}^2$ \ equals zero
 is positive (provided that \ $\PP(X_0=0)>0$ \ and \ $\PP(\vare_1=0)>0$), \
 but \ $\PP(\lim_{n\to\infty}\sum_{k=1}^nX_{k-1}^2=\infty)=1$ \ (which follows by the later
 formula \eqref{Ergodic2}).
In what follows we simply denote \ $\widetilde\alpha_n(X_0,X_1,\ldots,X_n)$ \ by
 \ $\widetilde\alpha_n$.
\ Using the same arguments as in Hall and Heyde \cite[Section 6.3]{HalHey},
 one can easily check that \ $\widetilde\alpha_n$ \ is a strongly consistent
 estimator of \ $\alpha$ \ as \ $n\to\infty$ \ for all \ $\alpha\in(0,1)$,
 \ i.e.,
 \begin{align}\label{Strong_consistency_INAR}
   \PP\left(\lim_{n\to\infty}\frac{\sum_{k=1}^n(X_k-\mu_\vare)X_{k-1}}{\sum_{k=1}^nX_{k-1}^2}
             =\alpha\right) =1, \qquad \forall\;\alpha\in(0,1).
 \end{align}
Namely, if \ $\widetilde X$ \ denotes a random variable with the unique
 stationary distribution of the INAR(1) model in \eqref{INAR1}, then
 \begin{align}\label{STAC_MOMENT1}
    \EE \widetilde X
       & = \frac{\TS\mu_\vare}{\TS1-\alpha},\\ \label{STAC_MOMENT2}
    \EE \widetilde X^2
       & = \frac{\TS\sigma_\vare^2+\alpha\mu_\vare}{\TS1-\alpha^2}
         + \frac{\mu_\vare^2}{(1-\alpha)^2}.
 \end{align}
For the proofs of \eqref{STAC_MOMENT1} and \eqref{STAC_MOMENT2}, see the Appendix.
By the existence of a unique stationary distribution,
 we obtain that \ $\{i\in\ZZ_+: i\geq i_{\min}\}$ \ with
 \[
    i_{\min}:=\min\Big\{ i\in\ZZ_+ : \PP(\vare_1=i)>0\Big\}
 \]
 is a  positive recurrent class of the Markov chain \ $(X_k)_{k\in\ZZ_+}$
 (see, e.g., Bhattacharya and Waymire \cite[Section II, Theorem 9.4 (c)]{BhaWay}
 or Chung \cite[Section I.6, Theorem 4 and Section I.7, Theorem 2]{Chu}).
 By ergodic theorems (see, e.g., Bhattacharya and Waymire
 \cite[Section II, Theorem 9.4 (d)]{BhaWay} or Chung \cite[Section I.15, Theorem 2]{Chu}),
 we get
 \begin{align}
    & \PP\left(\lim_{n\to\infty}\frac{1}{n}\sum_{k=1}^nX_k
                  = \EE \widetilde X
        \right)=1,\label{Ergodic1} \\
     & \PP\left(\lim_{n\to\infty}\frac{1}{n}\sum_{k=1}^nX_k^2
                  = \EE \widetilde X^2
        \right)=1, \label{Ergodic2} \\
      & \PP\left(\lim_{n\to\infty}\frac{1}{n}\sum_{k=1}^nX_{k-1}X_k
                  = \EE(\widetilde X(\alpha\circ \widetilde X + \vare))
                  = \alpha\EE \widetilde X^2 + \mu_\vare\EE\widetilde X
        \right)=1 \label{Ergodic3},
 \end{align}
 where \ $\vare$ \ is a random variable independent of \ $\widetilde X$ \ with
 the same distribution as \ $\vare_1$.
\ (For \eqref{Ergodic3}, one uses that the distribution of  \ $(\widetilde X,\alpha\circ\widetilde X+\vare)$
 \ is the unique stationary distribution of the Markov chain \ $(X_k,X_{k+1})_{k\in\ZZ_+}$.)
\ By \eqref{Ergodic1}--\eqref{Ergodic3},
 \[
   \PP\left(\lim_{n\to\infty} \widetilde\alpha_n =
               \frac{\alpha \EE \widetilde X^2 + \mu_\vare \EE\widetilde X
                          - \mu_\vare \EE\widetilde X}{\EE\widetilde X^2}
           =\alpha \right)=1.
 \]
Furthermore, if \ $\EE X_0^3<\infty$ \ and \ $\EE\vare_1^3<\infty$, \ then using
 the same arguments as in Hall and Heyde \cite[Section 6.3]{HalHey}, it follows easily that
 \begin{align}\label{SEGED11}
   \sqrt{n}(\widetilde\alpha_n-\alpha)\distr \cN(0,\sigma_{\alpha,\,\vare}^2)
      \qquad \text{as \ $n\to\infty$,}
 \end{align}
 where \ $\distr$ \ denotes convergence in distribution and
 \begin{align}\label{SEGED_SZIGMA_ALPHA}
   \sigma_{\alpha,\,\vare}^2
    :=\frac{\alpha(1-\alpha)\EE \widetilde X^3+\sigma_\vare^2\EE \widetilde X^2}
        {(\EE\widetilde X^2)^2},
  \end{align}
 with
 \begin{align}
  \begin{split}\label{STAC_MOMENT3}
    \EE \widetilde X^3
       & = \frac{\EE\vare^3-3\sigma_\vare^2(1+\mu_\vare)-\mu_\vare^3+2\mu_\vare}{1-\alpha^3}
         +3\frac{\sigma_\vare^2+\alpha\mu_\vare}{1-\alpha^2}
         -2\frac{\mu_\vare}{1-\alpha}\\
      &\phantom{=\;\;}
         +3\frac{\mu_\vare(\sigma_\vare^2+\alpha\mu_\vare)}{(1-\alpha)(1-\alpha^2)}
         +\frac{\mu_\vare^3}{(1-\alpha)^3}.
   \end{split}
 \end{align}
For the proof of \eqref{STAC_MOMENT3}, see the Appendix.

We remark that one uses in fact Corollary 3.1 in Hall and Heyde \cite{HalHey}
 to derive \eqref{SEGED11}.
It is important to point out that the moment conditions \ $\EE X_0^3<\infty$ \ and
 \ $\EE\vare_1^3<\infty$ \ are needed to check the conditions of this corollary
 (the so called conditional Lindeberg condition and an analogous condition
 on the conditional variance).

\subsection{Estimation of the mean of the offspring and innovation distributions}

Now we consider the joint CLS estimation of \ $\alpha$ \ and \ $\mu_\vare$.
\ For all \ $n\in\NN$, \ we define the function \ $Q_n:\RR^{n+1}\times\RR^2\to\RR$, \ as
 \begin{align*}
    Q_n(x_0,x_1,\ldots,x_n;\alpha',\mu_\vare')
      & :=\sum_{k=1}^n\big(x_k-\alpha' x_{k-1}-\mu_\vare'\big)^2,
       \qquad x_0,x_1,\ldots,x_n,\alpha',\mu_\vare'\in\RR.
 \end{align*}
By definition, for all \ $n\in\NN$, \ a CLS estimator for
 the parameter \ $(\alpha,\mu_\vare)\in(0,1)\times(0,\infty)$ \ is a measurable function
 \ $(\halpha_n,\hmuen):\RR^{n+1}\to\RR^2$ \ such that
 \begin{align*}
   Q_n(x_0,x_1,\ldots,x_n;\,&\halpha_n(x_0,x_1,\ldots,x_n),\hmuen(x_0,x_1,\ldots,x_n))\\
       &= \inf_{(\alpha',\mu_\vare')\in\RR^2}Q_n(x_0,x_1,\ldots,x_n;\alpha',\mu_\vare')
       \qquad \forall\;\;  (x_0,x_1,\ldots,x_n)\in\RR^{n+1}.
 \end{align*}
It is well-known that
 \begin{align*}
   &\sum_{k=1}^n (X_k-\halpha_n X_{k-1} - \hmuen)X_{k-1} = 0,\\
   &\sum_{k=1}^n (X_k-\halpha_n X_{k-1} - \hmuen) = 0,
 \end{align*}
 hold asymptotically as \ $n\to\infty$ \ with probability one, or equivalently
  \begin{align*}
     \begin{bmatrix}
       \sum_{k=1}^n X_{k-1}^2 & \sum_{k=1}^n X_{k-1} \\
       \sum_{k=1}^n X_{k-1} & n \\
     \end{bmatrix}
     \begin{bmatrix}
       \halpha_n \\
        \hmuen \\
     \end{bmatrix}
    = \begin{bmatrix}
        \sum_{k=1}^n X_{k-1}X_k \\
        \sum_{k=1}^n X_k \\
      \end{bmatrix}
  \end{align*}
  holds asymptotically as \ $n\to\infty$ \ with probability one.
Using that, by \eqref{Ergodic1} and \eqref{Ergodic2},
 \begin{align*}
  \PP\left(\lim_{n\to\infty}\frac{1}{n^2}
      \left( n\sum_{k=1}^n X_{k-1}^2-\left(\sum_{k=1}^n X_{k-1}\right)^2\right)
       =\EE\widetilde X^2-(\EE\widetilde X)^2
       =\var \widetilde X >0
     \right)
    =1,
 \end{align*}
we get
 \begin{align*}
   \halpha_n(X_0,X_1,\ldots,X_n)
     &=\frac{n\sum_{k=1}^nX_{k-1}X_k-\left(\sum_{k=1}^nX_{k-1}\right)\left(\sum_{k=1}^nX_k\right)}
         {n\sum_{k=1}^nX_{k-1}^2-\left(\sum_{k=1}^nX_{k-1}\right)^2},\\[2mm]
   \widehat\mu_{\vare,n}(X_0,X_1,\ldots,X_n)
     &=\frac{\left(\sum_{k=1}^nX_{k-1}^2\right)\left(\sum_{k=1}^nX_k\right)
             - \left(\sum_{k=1}^nX_{k-1}\right)\left(\sum_{k=1}^nX_{k-1}X_k\right)}
           {n\sum_{k=1}^nX_{k-1}^2-\left(\sum_{k=1}^nX_{k-1}\right)^2} \\
     &=\frac{1}{n}\left(\sum_{k=1}^nX_k-\halpha_n\sum_{k=1}^nX_{k-1}\right),
 \end{align*}
  hold asymptotically as \ $n\to\infty$ \ with probability one, see, e.g.,
  Hall and Heyde \cite[formulae (6.36) and (6.37)]{HalHey}.
In the sequel we simply denote \ $\widehat\alpha_n(X_0,X_1,\ldots,X_n)$ \ and
 \ $\widehat\mu_{\vare,n}(X_0,X_1,\ldots,X_n)$ \ by \ $\widehat\alpha_n$ \ and
 \ $\widehat\mu_{\vare,n}$, \ respectively.
It is well-known that \ $(\halpha_n,\,\widehat\mu_{\vare,\,n})$ \ is a strongly consistent
 estimator of \ $(\alpha,\mu_\vare)$ \ as \ $n\to\infty$ \ for all
 \ $(\alpha,\mu_\vare)\in(0,1)\times(0,\infty)$, \ see, e.g., Hall and Heyde \cite[Section 6.3]{HalHey}.
Moreover, if \ $\EE X_0^3<\infty$ \ and \ $\EE\vare_1^3<\infty$, \ by Hall and Heyde
 \cite[formula (6.44)]{HalHey},
 \begin{align}\label{ALPHA_MU_CONVERGENCE_INAR}
      \begin{bmatrix}
        \sqrt{n}(\halpha_n-\alpha) \\
        \sqrt{n}(\hmuen-\mu_\vare) \\
      \end{bmatrix}
     \distr
      \cN\left( \begin{bmatrix}
                  0 \\
                  0 \\
                \end{bmatrix}
               , \;B_{\alpha,\vare}
      \right)
   \qquad \text{as \ $n\to\infty$,}
 \end{align}
 where
 \begin{align}\label{SEGED_BALPHA}
   \begin{split}
   &B_{\alpha,\vare}:=
      \begin{bmatrix}
             \EE \widetilde X^2 & \EE \widetilde X \\
             \EE \widetilde X & 1 \\
           \end{bmatrix}^{-1}
           A_{\alpha,\vare}
           \begin{bmatrix}
             \EE \widetilde X^2 & \EE \widetilde X \\
             \EE \widetilde X & 1 \\
           \end{bmatrix}^{-1} \\
   &\phantom{B_{\alpha,\vare}\,}
       = \frac{1}{(\var \widetilde X)^2}
           \begin{bmatrix}
             1 & -\EE \widetilde X \\
             -\EE \widetilde X & \EE \widetilde X^2 \\
           \end{bmatrix}
           A_{\alpha,\vare}
           \begin{bmatrix}
             1 & -\EE \widetilde X \\
             -\EE \widetilde X & \EE \widetilde X^2 \\
           \end{bmatrix},
   \end{split}        \\[2mm] \label{SEGED_AALPHA}
  &A_{\alpha,\vare}:=
     \alpha(1-\alpha)
          \begin{bmatrix}
            \EE \widetilde X^3 & \EE \widetilde X^2 \\
            \EE \widetilde X^2 & \EE \widetilde X \\
          \end{bmatrix}
          +\sigma_\vare^2
           \begin{bmatrix}
            \EE \widetilde X^2 & \EE \widetilde X \\
            \EE \widetilde X & 1 \\
          \end{bmatrix},
 \end{align}
 and \ $\widetilde X$ \ denotes a random variable with the unique stationary distribution of
 the INAR(1) model in \eqref{INAR1}.
For our later purposes, we sketch a proof of \eqref{ALPHA_MU_CONVERGENCE_INAR}.
Using that
 \[
     \begin{bmatrix}
        \halpha_n \\
        \hmuen \\
      \end{bmatrix}
      =
       \begin{bmatrix}
          \sum_{k=1}^n X_{k-1}^2 & \sum_{k=1}^n X_{k-1} \\
          \sum_{k=1}^n X_{k-1} & n \\
       \end{bmatrix}^{-1}
       \begin{bmatrix}
          \sum_{k=1}^n X_{k-1}X_k \\
          \sum_{k=1}^n X_k \\
       \end{bmatrix}
 \]
 holds asymptotically as \ $n\to\infty$ \ with probability one, we obtain
 \begin{align*}
   &\begin{bmatrix}
     \halpha_n - \alpha\\
     \hmuen - \mu_\vare\\
   \end{bmatrix}=\\
   &\quad
    =
      \begin{bmatrix}
          \sum_{k=1}^n X_{k-1}^2 & \sum_{k=1}^n X_{k-1} \\
          \sum_{k=1}^n X_{k-1} & n \\
       \end{bmatrix}^{-1}
     \left(
       \begin{bmatrix}
          \sum_{k=1}^n X_{k-1}X_k \\
          \sum_{k=1}^n X_k \\
       \end{bmatrix}
       - \begin{bmatrix}
          \sum_{k=1}^n X_{k-1}^2 & \sum_{k=1}^n X_{k-1} \\
          \sum_{k=1}^n X_{k-1} & n \\
       \end{bmatrix}
       \begin{bmatrix}
         \alpha\\
         \mu_\vare\\
    \end{bmatrix}
       \right)  \\
   &\quad
    = \begin{bmatrix}
          \sum_{k=1}^n X_{k-1}^2 & \sum_{k=1}^n X_{k-1} \\
          \sum_{k=1}^n X_{k-1} & n \\
       \end{bmatrix}^{-1}
      \begin{bmatrix}
          \sum_{k=1}^n (X_k-\alpha X_{k-1}-\mu_\vare)X_{k-1} \\
          \sum_{k=1}^n (X_k-\alpha X_{k-1}-\mu_\vare) \\
       \end{bmatrix}
 \end{align*}
 holds asymptotically as \ $n\to\infty$ \ with probability one.
By \eqref{Ergodic1} and \eqref{Ergodic2}, we have
 \begin{align*}
   \frac{1}{n}
     \begin{bmatrix}
       \sum_{k=1}^n X_{k-1}^2 & \sum_{k=1}^n X_{k-1} \\
       \sum_{k=1}^n X_{k-1} & n \\
     \end{bmatrix}
     \to
     \begin{bmatrix}
       \EE\widetilde X^2 & \EE\widetilde  X \\
       \EE\widetilde  X & 1 \\
     \end{bmatrix}
  \qquad \text{as \ $n\to\infty$ \ with probability one},
 \end{align*}
 and, by Hall and Heyde \cite[Section 6.3, formula (6.43)]{HalHey},
 \begin{align*}
   \begin{bmatrix}
     \frac{1}{\sqrt{n}}\sum_{k=1}^n (X_k-\alpha X_{k-1}-\mu_\vare)X_{k-1} \\
     \frac{1}{\sqrt{n}}\sum_{k=1}^n (X_k-\alpha X_{k-1}-\mu_\vare) \\
   \end{bmatrix}
    \distr \cN\left(\begin{bmatrix}
                      0 \\
                      0 \\
                    \end{bmatrix},
                    A_{\alpha,\vare}
    \right)
  \qquad \text{as \ $n\to\infty$.}
 \end{align*}
Using Slutsky's lemma, we get \eqref{ALPHA_MU_CONVERGENCE_INAR}.

Let us introduce some notations which will be used throughout the paper.
For all \ $k,\ell\in\ZZ_+$, \  let
 \[
   \delta_{k,\ell}:=
       \begin{cases}
         1 & \text{if \ $k=\ell$},\\
         0 & \text{if \ $k\ne \ell$.}
       \end{cases}
 \]
For a sequence of random variables \ $(\zeta_k)_{k\in\NN}$
 \ and for \ $s_1,\ldots,s_I\in\NN$, \ $I\in\NN$, \ we define
 \begin{align*}
   \sideset{}{^{(s_1,\ldots,s_I)}}\sum_{k=1}^n \zeta_k
         := \DS\sum_{\substack{k=1 \\ k\ne s_1,\ldots,k\ne s_I}}^n\zeta_k.
 \end{align*}

\section{The INAR(1) model with additive outliers}\label{INAR1_additive}

\subsection{The model}

In this section we only introduce INAR(1) models contaminated with additive outliers.

\begin{Def}
A stochastic process \ $(Y_k)_{k\in\ZZ_+}$ \ is called an INAR(1) model with finitely many
 additive outliers if
 \[
   Y_k=
        X_k+\sum_{i=1}^I\delta_{k,s_i}\theta_i, \qquad k\in\ZZ_+,
 \]
 where \ $(X_k)_{k\in\ZZ_+}$ \ is an INAR(1) process given by \eqref{INAR1}
 with \ $\alpha\in(0,1)$, \ $\EE X_0^2<\infty$, \ $\EE\vare_1^2<\infty$, \ $\PP(\vare_1\ne 0)>0$,
 \ and \ $I\in\NN$, \ $s_i,\theta_i\in\NN$, \ $i=1,\ldots,I$ \ such that \ $s_i\ne s_j$ \ if
 \ $i\ne j$, \ $i,j=1,\ldots,I$.
\end{Def}

Notice that \ $\theta_i$, \ $i=1,\ldots,I$, \ represents the \ $i$th additive outlier's size
 and \ $\delta_{k,s_i}$ \ is an impulse taking the value \ $1$ \ if \ $k=s_i$ \ and \ $0$
 \ otherwise.
Roughly speaking, an additive outlier can be interpreted as a measurement error at time
 \ $s_i$, \ $i=1,\ldots,I$, \ or as an impulse due to some unspecified exogenous source.
Note also that \ $Y_0=X_0$.
\ Let \ $\cF_k^Y$ \ be the \ $\sigma$--algebra generated by the random variables \ $Y_0,Y_1,\ldots,Y_k$.
For all \ $n\in\NN$, \ $y_0,\ldots,y_n\in\RR$ \ and \ $\omega\in\Omega$, \ let us introduce
 the notations
 \begin{align*}
    {\bf Y}_n(\omega):=(Y_0(\omega),Y_1(\omega),\ldots,Y_n(\omega)),\qquad\;
    {\bf Y}_n:=(Y_0,Y_1,\ldots,Y_n), \qquad\;
    {\bf y}_n:=(y_0,y_1,\ldots,y_n).
 \end{align*}

\subsection{One outlier, estimation of the mean of the offspring distribution and the outlier's size}

First we assume that \ $I=1$ \ and that the relevant time point \ $s_1:=s$ \ is known.
We concentrate on the CLS estimation of the parameter \ $(\alpha,\theta)$ \ with
 \ $\theta:=\theta_1$.
\ An easy calculation shows that
 \begin{align*}
   \EE(Y_k\mid\cF^Y_{k-1})
     & =\alpha X_{k-1}+\mu_\vare+\delta_{k,s}\theta
       = \alpha(Y_{k-1}-\delta_{k-1,s}\theta)+\mu_\vare+ \delta_{k,s}\theta\\
     &= \alpha Y_{k-1} + \mu_\vare + (-\alpha\delta_{k-1,s} + \delta_{k,s})\theta
      =\begin{cases}
        \alpha Y_{k-1}+\mu_\vare & \quad \text{if \ $k=1,\ldots,s-1$,}\\
        \alpha Y_{k-1}+\mu_\vare+\theta & \quad \text{if \ $k=s$,}\\
        \alpha Y_{k-1}+\mu_\vare-\alpha\theta & \quad \text{if \ $k=s+1$,}\\
        \alpha Y_{k-1}+\mu_\vare & \quad \text{if \ $k\geq s+2$.}
      \end{cases}
 \end{align*}
Hence
 \begin{align} \label{SEGED32}
   \begin{split}
    \sum_{k=1}^n\big(Y_k-\EE(Y_k\mid \cF^Y_{k-1})\big)^2
         & =\sumss \big(Y_k-\alpha Y_{k-1}-\mu_\vare\big)^2
           + \big(Y_s-\alpha Y_{s-1}-\mu_\vare-\theta\big)^2 \\
         & \phantom{=\;\;} + \big(Y_{s+1}-\alpha Y_{s}-\mu_\vare+\alpha\theta\big)^2.
   \end{split}
 \end{align}
For all \ $n\geq s+1$, \ $n\in\NN$, \ we define the function \ $Q_n:\RR^{n+1}\times\RR^2\to\RR$, \ as
 \begin{align*}
    Q_n({\bf y}_n;\alpha',\theta')
      & :=\sumss \big(y_k-\alpha' y_{k-1}-\mu_\vare\big)^2
           + \big(y_s-\alpha' y_{s-1}-\mu_\vare-\theta'\big)^2 \\
      & \phantom{=\;\;} + \big(y_{s+1}-\alpha' y_{s}-\mu_\vare+\alpha'\theta'\big)^2,
       \qquad {\bf y}_n\in\RR^{n+1},\;\alpha',\theta'\in\RR.
 \end{align*}
By definition, for all \ $n\geq s+1$, \ a CLS estimator for
 the parameter \ $(\alpha,\theta)\in(0,1)\times\NN$ \ is a measurable function
 \ $(\talpha_n,\ttheta_n):S_n\to\RR^2$ \ such that
 \begin{align*}
   Q_n({\bf y}_n;\,&\talpha_n({\bf y}_n),\ttheta_n({\bf y}_n))
        = \inf_{(\alpha',\theta')\in\RR^2}Q_n({\bf y}_n;\alpha',\theta')
       \qquad \forall\;\;  {\bf y}_n\in S_n,
 \end{align*}
 where \ $S_n$ \ is a suitable subset of \ $\RR^{n+1}$ (defined in the proof of Lemma \ref{LEMMA2}).
We note that we do not define the CLS estimator
 \ $(\talpha_n,\ttheta_n)$ \ for all samples \ ${\bf y}_n\in \RR^{n+1}$.
\ We have for all \ $({\bf y}_n;\alpha',\theta')\in\RR^{n+1}\times\RR^2$,
 \begin{align*}
   &\frac{\partial Q_n}{\partial \alpha'}({\bf y}_n;\alpha',\theta')
      = \sumss 2\big(y_k-\alpha' y_{k-1}-\mu_\vare\big)(-y_{k-1})
         + 2\big(y_s-\alpha' y_{s-1}-\mu_\vare-\theta'\big)(-y_{s-1})\\
     &\phantom{\frac{\partial Q_n}{\partial \alpha'}({\bf y}_n;\alpha',\theta') =\;\;}
         + 2\big(y_{s+1}-\alpha' y_{s}-\mu_\vare+\alpha'\theta'\big)(-y_s+\theta')\\
     &=\sum_{k=1}^n 2\big(y_k-\alpha' y_{k-1}-\mu_\vare\big)(-y_{k-1})
         - 2\theta'(-y_{s-1}) + 2\alpha'\theta'(-y_s+\theta')
         + 2\big(y_{s+1}-\alpha' y_{s}-\mu_\vare\big)\theta',
 \end{align*}
 and
 \begin{align*}
   \frac{\partial Q_n}{\partial\theta'}({\bf y}_n;\alpha',\theta')
      =-2\big(y_s-\alpha' y_{s-1}-\mu_\vare-\theta'\big)
       +2\big(y_{s+1}-\alpha' y_{s}-\mu_\vare+\alpha'\theta'\big)\alpha'.
  \end{align*}

The next lemma is about the existence and uniqueness of the CLS estimator of \ $(\alpha,\theta)$.

\begin{Lem}\label{LEMMA2}
There exist subsets \ $S_n\subset\RR^{n+1}$, $n\geq s+1$ \ with the following properties:
 \begin{enumerate}
  \item[\upshape{(i)}]
   there exists a unique CLS estimator
   \ $(\talpha_n,\ttheta_n):S_n\to\RR^2$,
  \item[\upshape{(ii)}]
   for all \ ${\bf y}_n\in S_n$,
   \ $(\talpha_n({\bf y}_n),\ttheta_n({\bf y}_n))$
   \ is the unique solution of the system of equations
 \begin{align}\label{Additive_CLSE_EQ}
    \frac{\partial Q_n}{\partial \alpha'}({\bf y}_n;\alpha',\theta')=0,\qquad
    \frac{\partial Q_n}{\partial \theta'}({\bf y}_n;\alpha',\theta')=0,
 \end{align}
  \item [\upshape{(iii)}]
  ${\bf Y}_n\in S_n$ \ holds asymptotically as \ $n\to\infty$ \ with probability one.
 \end{enumerate}
\end{Lem}

\noindent{\bf Proof.}
For any fixed \ ${\bf y}_n \in \RR^{n+1}$ \ and \ $\alpha' \in \RR$, \ the quadratic function
 \ $\RR\ni \theta' \mapsto Q_n({\bf y}_n;\alpha',\theta')$
 \ can be written in the form
 \begin{align*}
  Q_n({\bf y}_n;\alpha',\theta')
   = A_n(\alpha') \! \left( \theta' - A_n(\alpha')^{-1} t_n({\bf y}_n;\alpha') \right)^2
     \!\! + \widetilde{Q}_n({\bf y}_n;\alpha'),
 \end{align*}
 where
 \begin{align*}
  A_n(\alpha')
   &:= 1 + (\alpha')^2,\\[1mm]
  t_n({\bf y}_n;\alpha')
   &:= ( 1 + (\alpha')^2 ) y_s - \alpha' ( y_{s-1} + y_{s+1} ) - (1-\alpha')\mu_\vare,\\[2mm]
  \widetilde{Q}_n({\bf y}_n;\alpha')
   &:= \sum_{k=1}^n \big( y_k-\alpha' y_{k-1} \big)^2
      - A_n(\alpha')^{-1} t_n({\bf y}_n;\alpha')^2 .
 \end{align*}
We have \ $\widetilde{Q}_n({\bf y}_n;\alpha') = R_n({\bf y}_n;\alpha') / A_n(\alpha')$, \ where
 \ $\RR\ni\alpha' \mapsto R_n({\bf y}_n;\alpha')$ \ is a polynomial of order 4 with leading coefficient
 \begin{align*}
   c_n({\bf y}_n) &:= \sum_{k=1}^n y_{k-1}^2 - y_s^2.
 \end{align*}
Let
 \[
   \widetilde{S}_n := \left\{{\bf y}_n\in\RR^{n+1} : c_n({\bf y}_n) > 0 \right\}.
 \]
For \ ${\bf y}_n \in \widetilde{S}_n$, \ we have
 \ $\lim_{|\alpha'|\to\infty} \widetilde{Q}_n({\bf y}_n;\alpha') = \infty$ \ and the continuous function
 \ $\RR \ni \alpha' \mapsto \widetilde{Q}_n({\bf y}_n;\alpha')$ \ attains its infimum.
Consequently, for all \ $n\geq s+1$ \ there exists a CLS estimator
 \ $(\talpha_n,\ttheta_{n}):\widetilde{S}_n\to\RR^2$, \ where
 \begin{align}\nonumber
   \widetilde{Q}_n({\bf y}_n;\talpha_n({\bf y}_n))
   &=\inf_{\alpha'\in\RR} \widetilde{Q}_n({\bf y}_n;\alpha')
   \qquad \forall\;{\bf y}_n\in \widetilde{S}_n ,\\[2mm]\label{SEGED_UJ2}
   \ttheta_n({\bf y}_n)
      & = A_n(\talpha_n({\bf y}_n))^{-1} t_n({\bf y}_n;\talpha_n({\bf y}_n)) ,
  \qquad {\bf y}_n\in \widetilde{S}_n,
 \end{align}
 and for all \ ${\bf y}_n\in \widetilde{S}_n$,
 \ $(\talpha_n({\bf y}_n),\ttheta_n({\bf y}_n))$
 \ is a  solution of the system of equations \eqref{Additive_CLSE_EQ}.

By \eqref{Ergodic1} and \eqref{Ergodic2}, we get
 \ $\PP\left(\lim_{n\to\infty} n^{-1} c_n({\bf Y}_n) = \EE\widetilde X^2 \right)=1$, \ where \ $\widetilde X$
 \ denotes a random variable with the unique stationary distribution of the INAR(1) model in \eqref{INAR1}.
Hence \ ${\bf Y}_n\in \widetilde{S}_n$ \ holds asymptotically as \ $n\to\infty$ \ with probability one.

Now we turn to find sets \ $S_n \subset \widetilde{S}_n$, $n\geq s+1$ \ such that the system of equations
 \eqref{Additive_CLSE_EQ} has a unique solution with respect to
 \ $(\alpha',\theta')$ \ for all \ ${\bf y}_n\in S_n$.
\ Let us introduce the \ $(2\times 2)$ \ Hessian matrix
 \[
    H_n({\bf y}_n;\alpha',\theta'):=
        \begin{bmatrix}
       \frac{\partial^2 Q_n}{\partial (\alpha')^2}
      & \frac{\partial^2 Q_n}{\partial \theta' \partial \alpha'} \\
       \frac{\partial^2 Q_n}{\partial \alpha'\partial\theta'}
       & \frac{\partial^2 Q_n}{\partial (\theta')^2} \\
     \end{bmatrix}
       ({\bf y}_n;\alpha',\theta'),
 \]
 and let us denote by \ $\Delta_{i,n}({\bf y}_n;\alpha',\theta')$ \
 its \ $i$-th order leading principal minor, \ $i=1,2$.
\ Further, for all \ $n\geq s+1$, \ let
 \[
   S_n:=\Big\{{\bf y}_n\in\widetilde{S}_n : \Delta_{i,n}({\bf y}_n;\alpha',\theta')>0,
                                   \;\, i=1,2,\, \forall\;(\alpha',\theta')
                                   \in\RR^2 \Big\}.
 \]
By Berkovitz \cite[Theorem 3.3, Chapter III]{Ber},
 the function \ $\RR^2\ni(\alpha',\theta')\mapsto Q_n({\bf y}_n;\alpha',\theta')$ \
 is strictly convex for all \ ${\bf y}_n\in S_n$.
\ Since it was already proved that the system of equations \eqref{Additive_CLSE_EQ} has a solution for all
 \ ${\bf y}_n\in \widetilde{S}_n$, \ we obtain that this solution is unique for all \ ${\bf y}_n\in S_n$.

Next we check that \ ${\bf Y}_n\in S_n$ \ holds asymptotically as \ $n\to\infty$ \ with probability one.
For all \ $(\alpha',\theta')\in\RR^2$,
 \begin{align*}
   \frac{\partial^2 Q_n}{\partial (\alpha')^2} &({\bf Y}_n;\alpha',\theta')
       = 2\sum_{k=1}^n Y_{k-1}^2 + 2\theta'(-Y_s+\theta')-2Y_s\theta'
       =2\left(\sum_{k=1}^n Y_{k-1}^2 -2Y_s\theta'+(\theta')^2\right)\\
       &=2\left( \sumsf Y_{k-1}^2 +(Y_s-\theta')^2\right)
        =2\left( \sumsf  X_{k-1}^2 + (X_s+\theta-\theta')^2 \right),
 \end{align*}
 and
 \begin{align*}
  &\frac{\partial^2 Q_n}{\partial\alpha'\partial\theta'}({\bf Y}_n;\alpha',\theta')
    = \frac{\partial^2 Q_n}{\partial\theta'\partial\alpha'}({\bf Y}_n;\alpha',\theta')
    = 2(Y_{s-1}+Y_{s+1}-2\alpha' Y_s-\mu_\vare+2\alpha'\theta')\\
  &\phantom{\frac{\partial^2 Q_n}{\partial\alpha'\partial\theta'}({\bf Y}_n;\alpha',\theta')}
    =2(X_{s-1}+X_{s+1}-2\alpha' X_s-\mu_\vare-2\alpha'(\theta-\theta')),\\
  &\frac{\partial^2 Q_n}{\partial(\theta')^2}({\bf Y}_n;\alpha',\theta')
    = 2((\alpha')^2+1).
 \end{align*}
Then
 \begin{align*}
    &H_n({\bf Y}_n;\alpha',\theta')\\
    &=2\begin{bmatrix}
        \DS\sumsf X_{k-1}^2 + (X_s+\theta-\theta')^2 &
                  X_{s-1}+X_{s+1}-2\alpha' X_s-\mu_\vare-2\alpha'(\theta-\theta') \\
        X_{s-1}+X_{s+1}-2\alpha' X_s-\mu_\vare-2\alpha'(\theta-\theta') & (\alpha')^2+1 \\
      \end{bmatrix}
 \end{align*}
 has leading principal minors
  \ $\Delta_{1,n}({\bf Y}_n;\alpha',\theta')=2\left( \DS\sumsf  X_{k-1}^2 + (X_s+\theta-\theta')^2 \right)$
 \ and
 \begin{align*}
    \Delta_{2,n}({\bf Y}_n;\alpha',\theta')=\det H_n({\bf Y}_n;\alpha',\theta')
      & = 4((\alpha')^2+1)\left(\sumsf  X_{k-1}^2 + (X_s+\theta-\theta')^2\right) \\
      &\phantom{= \;\;}
                - 4\Big(X_{s-1}+X_{s+1}-2\alpha' X_s-\mu_\vare-2\alpha'(\theta-\theta')\Big)^2.
 \end{align*}
By \eqref{Ergodic2},
 \begin{align*}
    & \PP\left(\lim_{n\to\infty}\frac{1}{n}\Delta_{1,n}({\bf Y}_n;\alpha',\theta')
               = 2\EE\widetilde X^2,\;\;
               \forall\;\; (\alpha',\theta')\in\RR^2 \right)=1,\\
    &\PP\left(\lim_{n\to\infty}\frac{1}{n}\Delta_{2,n}({\bf Y}_n;\alpha',\theta')
               = 4((\alpha')^2+1)\EE\widetilde X^2,\;\;
               \forall\;\; (\alpha',\theta')\in\RR^2 \right)=1,
 \end{align*}
 where \ $\widetilde X$ \ denotes a random variable with the unique stationary
 distribution of the INAR(1) model in \eqref{INAR1}.
Hence
 \begin{align*}
   &\PP\big(\lim_{n\to\infty}\Delta_{1,n}({\bf Y}_n;\alpha',\theta') =\infty,\;\;
            \forall\;\; (\alpha',\theta')\in\RR^2 \big)=1,\\
   &\PP\big(\lim_{n\to\infty}\Delta_{2,n}({\bf Y}_n;\alpha',\theta') =\infty,\;\;
             \forall\;\; (\alpha',\theta')\in\RR^2\big)=1,
 \end{align*}
which yields that \ ${\bf Y}_n\in S_n$ \ asymptotically as \ $n\to\infty$
 \ with probability one, since we have already proved that \ ${\bf Y}_n\in \widetilde{S}_n$
 \ asymptotically as \ $n\to\infty$ \ with probability one.
\proofend

By Lemma \ref{LEMMA2}, \ $(\talpha_n({\bf Y}_n),\ttheta_n({\bf Y}_n))$
 \ exists uniquely asymptotically as \ $n\to\infty$ \ with probability one.
In the sequel we will simply denote it by \ $(\talpha_n,\ttheta_n)$.

The next result shows that \ $\talpha_n$ \ is a strongly consistent estimator of \ $\alpha$, \
 whereas \ $\ttheta_n$ \ fails to be also a strongly consistent estimator of \ $\theta$.

\begin{Thm}\label{THEOREM1}
For the CLS estimators \ $(\talpha_n,\ttheta_n)_{n\in \NN}$ \ of the parameter
 \ $(\alpha,\theta)\in(0,1)\times\NN$,
\ the sequence \ $(\talpha_n)_{n\in \NN}$ \ is strongly consistent for all
 \ $(\alpha,\theta)\in(0,1)\times\NN$, \ i.e.,
 \begin{align}\label{Strong_consistency1}
   \PP(\lim_{n\to\infty}\talpha_n=\alpha)=1, \qquad \forall\;(\alpha,\theta)\in(0,1)\times\NN.,
 \end{align}
 whereas the sequence \ $(\ttheta_n)_{n\in \NN}$ \ is not strongly consistent
 for any \ $(\alpha,\theta)\in(0,1)\times\NN$, \ namely,
 \begin{align}\label{Strong_consistency2}
   \PP\left(\lim_{n\to\infty}\ttheta_n
         =Y_s-\frac{\alpha}{1+\alpha^2}(Y_{s-1}+Y_{s+1})
              -\frac{1-\alpha}{1+\alpha^2}\mu_\vare
         \right)=1,\qquad \forall\;(\alpha,\theta)\in(0,1)\times\NN.
 \end{align}
\end{Thm}

\noindent{\bf Proof.}
An easy calculation shows that
 \begin{align}\label{SEGED1}
   \talpha_n
     =\frac{\sum_{k=1}^n(Y_k-\mu_\vare)Y_{k-1}-\ttheta_n(Y_{s-1}+Y_{s+1}-\mu_\vare)}
           {\sum_{k=1}^nY_{k-1}^2-2\ttheta_n Y_s+(\ttheta_n)^2},\\[2mm]\label{SEGED2}
   \ttheta_n
      = Y_s-\frac{\talpha_n}{1+(\talpha_n)^2}(Y_{s-1}+Y_{s+1})
                  -\frac{1-\talpha_n}{1+(\talpha_n)^2}\mu_\vare,
 \end{align}
 hold asymptotically as \ $n\to\infty$ \ with probability one.
Since \ $Y_k=X_k+\delta_{k,s}\theta$, \ $k\in\ZZ_+$, \ we get
 \begin{align*}
   &\sum_{k=1}^n(Y_k-\mu_\vare)Y_{k-1} -\ttheta_n(Y_{s-1}+Y_{s+1}-\mu_\vare) \\
     & = \DS\sumss(X_k-\mu_\vare)X_{k-1} + (X_s+\theta-\mu_\vare)X_{s-1}+(X_{s+1}-\mu_\vare)(X_s+\theta)\\
     & \phantom{=\;\;}
        -\ttheta_n(X_{s-1}+X_{s+1}-\mu_\vare),
 \end{align*}
 and
 \begin{align*}
   \sum_{k=1}^nY_{k-1}^2-2\ttheta_n Y_s+(\ttheta_n)^2
      = \DS\sideset{}{^{(s+1)}}\sum_{k=1}^nX_{k-1}^2 + (X_s+\theta)^2 -2\ttheta_n(X_s+\theta)+(\ttheta_n)^2,
 \end{align*}
 hold asymptotically as \ $n\to\infty$ \ with probability one.
Hence
 \begin{align}\label{SEGED30}
     \talpha_n
      & =\frac{\sum_{k=1}^n(X_k-\mu_\vare)X_{k-1} + (\theta-\ttheta_n)(X_{s-1}+X_{s+1}-\mu_\vare)}
           {\sum_{k=1}^nX_{k-1}^2+(\theta-\ttheta_n)(\theta-\ttheta_n+2X_s)},
 \end{align}
 holds asymptotically as \ $n\to\infty$ \ with probability one.
We check that in proving \eqref{Strong_consistency1} it is enough to verify that
 \begin{align}\label{SEGED34}
   &\PP\left(\lim_{n\to\infty}
       \frac{(\theta-\ttheta_n)(X_{s-1}+X_{s+1}-\mu_\vare)}{n}
                = 0\right)=1,\\[1mm]\label{SEGED35}
   &\PP\left(\lim_{n\to\infty}
        \frac{(\theta-\ttheta_n)(\theta-\ttheta_n+2X_s)}{n} = 0\right)=1.
 \end{align}
Indeed, using \eqref{Ergodic1}, \eqref{Ergodic2} and \eqref{Ergodic3}, we get
 \eqref{SEGED34} and \eqref{SEGED35} yield that
 \begin{align*}
    \PP\left(\lim_{n\to\infty}\talpha_n
              =\frac{\alpha\EE\widetilde X^2+\mu_\vare\EE\widetilde X-\mu_\vare\EE\widetilde X}
                    {\EE\widetilde X^2}=\alpha\right)=1.
 \end{align*}
Now we turn to prove \eqref{SEGED34} and \eqref{SEGED35}.
By \eqref{SEGED2} and using again the decomposition \ $Y_k=X_k+\delta_{k,s}\theta$,
 \ $k\in\ZZ_+$, \ we obtain
 \begin{align*}
   \ttheta_n
      = X_s+\theta-\frac{\talpha_n}{1+(\talpha_n)^2}(X_{s-1}+X_{s+1})
                  -\frac{1-\talpha_n}{1+(\talpha_n)^2}\mu_\vare,
 \end{align*}
 and hence
 \begin{align}\label{THETA_BOUND}
  \vert \ttheta_n-\theta \vert \leq X_s+\frac{1}{2}(X_{s-1}+X_{s+1})+\frac{3}{2}\mu_\vare,
 \end{align}
 i.e., the sequences \ $(\ttheta_n)_{n\in\NN}$ \ and \ $(\ttheta_n-\theta)_{n\in\NN}$
  \ are bounded with probability one.
This implies \eqref{SEGED34} and \eqref{SEGED35}.
By \eqref{SEGED2} and \eqref{Strong_consistency1}, we get \eqref{Strong_consistency2}.
\proofend

\begin{Rem}\label{REMARK2}
We check that \ $\EE\big(\lim_{n\to\infty}\ttheta_n\big)= \theta$,
 \ $\forall$ $(\alpha,\theta)\in(0,1)\times\NN$, \ and
 \begin{align}\label{SEGED36}
  \var\big(\lim_{n\to\infty}\ttheta_n\big)
    =\frac{\mu_\vare(\alpha+\alpha^3-\alpha^s-\alpha^{s+3})+\sigma_\vare^2(1+\alpha^2)
           +(1-\alpha)(\alpha^s+\alpha^{s+3})\EE X_0}{(1+\alpha^2)^2}.
 \end{align}
Note that, by \eqref{Strong_consistency2}, with probability one it holds that
 \begin{align*}
    \lim_{n\to\infty}\ttheta_n
        & = X_s+\theta -\frac{\alpha}{1+\alpha^2}(X_{s-1}+X_{s+1})-\frac{1-\alpha}{1+\alpha^2}\mu_\vare \\
        & = \theta +\frac{X_s-\alpha X_{s-1}-\mu_\vare-\alpha(X_{s+1}-\alpha X_s-\mu_\vare)}{1+\alpha^2}
          = \theta +\frac{M_s - \alpha M_{s+1}}{1+\alpha^2},
 \end{align*}
 where \ $M_k:=X_k-\alpha X_{k-1}-\mu_\vare$, \ $k\in\NN$.
\ Notice that \ $\EE M_k=0$, \ $k\in\NN$, \ and
 \ $\cov(M_k,M_\ell)=\delta_{k,\ell}\var M_k$, \ $k,\ell\in\NN$, \ where
 \ $\var M_k=\alpha\mu_\vare(1-\alpha^{k-1})+\alpha^k(1-\alpha)\EE X_0+\sigma_\vare^2$, \ $k\in\NN$.
\ Indeed, by the recursion \ $\EE X_\ell=\alpha\EE X_{\ell-1}+\mu_\vare$, \ $\ell\in\NN$,
 \ we get \ $\EE M_k=0$, $k\in\NN$, \ and
 \ we get
 \begin{align}\label{SEGED_POT}
   \EE X_\ell=\alpha^\ell\EE X_0+(1+\alpha+\cdots+\alpha^{\ell-1})\mu_\vare
     =\alpha^\ell\EE X_0+\frac{1-\alpha^\ell}{1-\alpha}\mu_\vare, \qquad \ell\in\NN,
 \end{align}
 and hence
 \begin{align}\label{SEGED64}
  \begin{split}
   \var M_k
     & =\var(X_k-\alpha X_{k-1}-\mu_\vare)
       = \var \left( \sum_{j=1}^{X_{k-1}}(\xi_{k,j}-\alpha)+(\vare_k-\mu_\vare)\right) \\
     & = \alpha(1-\alpha)\EE X_{k-1}+\sigma_\vare^2
       = \alpha(1-\alpha)\mu_\vare\frac{1-\alpha^{k-1}}{1-\alpha}
         + \alpha^k(1-\alpha)\EE X_0 + \sigma_\vare^2\\
     & = \alpha\mu_\vare (1-\alpha^{k-1}) + \alpha^k(1-\alpha)\EE X_0 + \sigma_\vare^2,
     \qquad k\in\NN.
   \end{split}
 \end{align}
Hence \ $\EE(\lim_{n\to\infty}\ttheta_n)=\theta$ \ and
 \begin{align*}
   \var\big(\lim_{n\to\infty}\ttheta_n\big)
    =\frac{1}{(1+\alpha^2)^2}
      &\Big(\alpha\mu_\vare(1-\alpha^{s-1})+\alpha^s(1-\alpha)\EE X_0+\sigma_\vare^2 \\
      & \;\,+\alpha^2\big(\alpha\mu_\vare(1-\alpha^s) +\alpha^{s+1}(1-\alpha)\EE X_0 +\sigma_\vare^2\big)\Big),
 \end{align*}
 which implies \eqref{SEGED36}.

We also check that \ $\ttheta_n$ \ is an asymptotically unbiased estimator of \ $\theta$ \
 as \ $n\to\infty$ \ for all \ $(\alpha,\theta)\in(0,1)\times\NN$.
\ By \eqref{Strong_consistency2}, the sequence \ $\ttheta_n-\theta$, \ $n\in\NN$,
 \ converges with probability one, and, by \eqref{THETA_BOUND}, the dominated
 convergence theorem yields that
 \[
   \lim_{n\to\infty}\EE(\ttheta_n-\theta)
   =\EE\big[\lim_{n\to\infty}(\ttheta_n-\theta)\big]
   =0.
  \]
Finally, we note that \ $\lim_{n\to\infty}\ttheta_n$ \ can be negative with positive probability,
 despite the fact that \ $\theta\in\NN$.
\proofend
\end{Rem}

\begin{Def}
Let \ $(\zeta_n)_{n\in\NN}$, \ $\zeta$ \ and \ $\eta$ \ be random variables on
 \ $(\Omega,\cA,\PP)$ \ such that \ $\eta$ \ is non-negative and integer-valued.
By the expression {\sf "conditionally on the values \ $\eta$, \ the weak convergence
 \ $\zeta_n\distr\zeta$ \ as \ $n\to\infty$ \ holds"} we mean that for all
 non-negative integers \ $m\in\NN$ \ such that \ $\PP(\eta=m)>0$, \ we have
  \[
      \lim_{n\to\infty} F_{\zeta_n \mid \{\eta=m\}}(y)
         =  F_{\zeta \mid \{\eta=m\}}(y)
  \]
 for all \ $y\in\RR$ \ being continuity points of \ $F_{\zeta \mid \{\eta=m\}}$, \ where
 \ $F_{\zeta_n \mid \{\eta=m\}}$ \ and \ $F_{\zeta \mid \{\eta=m\}}$ \ denote the conditional
 distribution function of \ $\zeta_n$ \ and \ $\zeta$ \ with respect to the event
 \ $ \{\eta=m\}$, \ respectively.
\end{Def}

The asymptotic distribution of the CLS estimation is given in the next theorem.

\begin{Thm}\label{Proposition3}
Under the additional assumptions \ $\EE X_0^3<\infty$ \ and \ $\EE\vare_1^3<\infty$, \ we have
 \begin{align}\label{CONVERGENCE6}
   \sqrt{n}(\talpha_n-\alpha)\distr \cN(0,\sigma_{\alpha,\,\vare}^2)
      \qquad \text{as \ $n\to\infty$,}
  \end{align}
 where \ $\sigma_{\alpha,\,\vare}^2$ \ is defined in \eqref{SEGED_SZIGMA_ALPHA}.
Furthermore, conditionally on the values \ $Y_{s-1}$ \ and \ $Y_{s+1}$,
 \begin{align}\label{CONVERGENCE7}
   \sqrt{n}\left(\ttheta_n - \lim_{k\to\infty}\ttheta_k\right)
         \distr \cN\left(0,c_{\alpha,\,\vare}^2\right)
         \qquad \text{as \ $n\to\infty$,}
 \end{align}
 where
 \[
   c_{\alpha,\,\vare}^2:=\frac{\sigma_{\alpha,\,\vare}^2}{(1+\alpha^2)^4}
           \big((\alpha^2-1)(Y_{s-1}+Y_{s+1})+(1+2\alpha-\alpha^2)\mu_\vare\big)^2.
 \]
\end{Thm}

\noindent{\bf Proof.}
By \eqref{SEGED30}, we have \ $\sqrt{n}(\talpha_n-\alpha)=\frac{A_n}{B_n}$ \
 holds asymptotically as \ $n\to\infty$ \ with probability one, where
 \begin{align*}
    &A_n:= \frac{1}{\sqrt n}\sum_{k=1}^n(X_k-\alpha X_{k-1}-\mu_\vare)X_{k-1}
          + \frac{1}{\sqrt n}(\theta-\ttheta_n)(X_{s-1}+X_{s+1}-\mu_\vare
             -\alpha(\theta-\ttheta_n)-2\alpha X_s),\\
    &B_n:= \frac{1}{n}\sum_{k=1}^nX_{k-1}^2+\frac{1}{n}(\theta-\ttheta_n)(\theta-\ttheta_n+2X_s).
 \end{align*}
By \eqref{SEGED11}, we have
 \[
     \frac{\frac{1}{\sqrt n}\sum_{k=1}^n(X_k-\alpha X_{k-1}-\mu_\vare)X_{k-1}}
          {\frac{1}{n}\sum_{k=1}^nX_{k-1}^2}
        \distr \cN(0,\sigma_{\alpha,\,\vare}^2)
       \qquad \text{as \ $n\to\infty$.}
 \]
By \eqref{THETA_BOUND},
 \begin{align*}
   &\PP\left(\lim_{n\to\infty}\frac{1}{\sqrt n}(\theta-\ttheta_n)
              (X_{s-1}+X_{s+1}-\mu_\vare-\alpha(\theta-\ttheta_n)-2\alpha X_s)=0\right)=1,\\[1mm]
   &\PP\left(\lim_{n\to\infty} \frac{1}{n}(\theta-\ttheta_n)(\theta-\ttheta_n+2X_s)=0\right)=1.
 \end{align*}
Hence Slutsky's lemma yields \eqref{CONVERGENCE6}.
Using \eqref{Strong_consistency2} and that
\begin{align}\label{SEGED31}
   \ttheta_n
      & =  Y_s + \frac{-(\talpha_n-\alpha)(Y_{s-1}+Y_{s+1})+(\talpha_n-\alpha)\mu_\vare
                    -\alpha(Y_{s-1}+Y_{s+1})-(1-\alpha)\mu_\vare}
                    {1+(\talpha_n)^2},
 \end{align}
  holds asymptotically as \ $n\to\infty$ \ with probability one, we get
 \begin{align*}
    \sqrt{n}&\left(\ttheta_n - \lim_{k\to\infty}\ttheta_k\right)
             =\sqrt{n}\left(\ttheta_n - \left(Y_s-\frac{\alpha}{1+\alpha^2}(Y_{s-1}+Y_{s+1})
              -\frac{1-\alpha}{1+\alpha^2}\mu_\vare\right)\right)\\
            &=\sqrt{n}
               \left(\frac{-(\talpha_n-\alpha)(Y_{s-1}+Y_{s+1})
                   +(\talpha_n-\alpha)\mu_\vare
                   -\alpha(Y_{s-1}+Y_{s+1})
                   -(1-\alpha)\mu_\vare}
                  {1+(\talpha_n)^2}\right.\\
            &\phantom{=\sqrt{n} \;\;\;} \left.
                +\frac{\alpha}{1+\alpha^2}(Y_{s-1}+Y_{s+1})
                +\frac{1-\alpha}{1+\alpha^2}\mu_\vare
             \right)\\
            &=\sqrt{n}(\talpha_n-\alpha)
              \left[\frac{-Y_{s-1}-Y_{s+1}+\mu_\vare}{1+(\talpha_n)^2}+\alpha(Y_{s-1}+Y_{s+1})
                    \frac{\alpha+\talpha_n} {(1+(\talpha_n)^2)(1+\alpha^2)} \right.\\
            &\phantom{=\sqrt{n}(\talpha_n-\alpha)\Big[\;\;}\left.
                   + \mu_\vare\frac{(1-\alpha)(\alpha+\talpha_n)} {(1+(\talpha_n)^2)(1+\alpha^2)}
                    \right].
 \end{align*}
Using \eqref{CONVERGENCE6} and \eqref{Strong_consistency1}, Slutsky's lemma yields \eqref{CONVERGENCE7}
 with
 \begin{align*}
   c_{\alpha,\,\vare}^2
     &= \sigma_{\alpha,\,\vare}^2
         \left(\frac{-Y_{s-1}-Y_{s+1}+\mu_\vare}{1+\alpha^2}
                + \frac{2\alpha^2}{(1+\alpha^2)^2}(Y_{s-1}+Y_{s+1})
                +\frac{2\alpha(1-\alpha)}{(1+\alpha^2)^2}\mu_\vare\right)^2 \\
     &=\frac{\sigma_{\alpha,\,\vare}^2}{(1+\alpha^2)^2}
         \left(\frac{\alpha^2-1}{1+\alpha^2}(Y_{s-1}+Y_{s+1})
               +\frac{1+2\alpha-\alpha^2}{1+\alpha^2}\mu_\vare\right)^2.
 \end{align*}
\proofend

\subsection{One outlier, estimation of the mean of the offspring and innovation distributions
            and the outlier's size}

We assume \ $I=1$ \ and that the relevant time point \ $s\in\NN$ \ is known.
We concentrate on the CLS estimation of \ $\alpha$, \ $\mu_\vare$ \ and
 \ $\theta:=\theta_1$.
\ For all \ $n\geq s+1$, \ $n\in\NN$, \ we define the function \ $Q_n:\RR^{n+1}\times\RR^3\to\RR$, \ as
 \begin{align*}
    Q_n({\bf y}_n;\alpha',\mu_\vare',\theta')
      & :=\sumss \big(y_k-\alpha' y_{k-1}-\mu_\vare'\big)^2
           + \big(y_s-\alpha' y_{s-1}-\mu_\vare'-\theta'\big)^2 \\
      & \phantom{=\;\;} + \big(y_{s+1}-\alpha' y_{s}-\mu_\vare'+\alpha'\theta'\big)^2,
       \qquad {\bf y}_n\in\RR^{n+1},\;\alpha',\mu_\vare',\theta'\in\RR.
 \end{align*}
By definition, for all \ $n\geq s+1$, \ a CLS estimator for
 the parameter \ $(\alpha,\mu_\vare,\theta)\in(0,1)\times(0,\infty)\times\NN$ \ is
 a measurable function \ $(\halpha_n,\hmuen,\htheta_n):S_n\to\RR^3$ \ such that
 \begin{align*}
   Q_n({\bf y}_n;\,&\halpha_n({\bf y}_n),\hmuen({\bf y}_n),
            \htheta_n({\bf y}_n))
       = \inf_{(\alpha',\mu_\vare',\theta')\in\RR^3}Q_n({\bf y}_n;\alpha',\mu_\vare',\theta')
       \qquad \forall\;\;  {\bf y}_n\in S_n,
 \end{align*}
 where \ $S_n$ \ is suitable subset of \ $\RR^{n+1}$ \ (defined in the proof of
 Lemma \ref{LEMMA5}).
We note that we do not define the CLS estimator
 \ $(\halpha_n,\hmuen,\htheta_n)$ \ for all samples \ ${\bf y}_n\in \RR^{n+1}$.
We get for all \ $({\bf y}_n;\alpha',\mu_\vare',\theta')\in\RR^{n+1}\times\RR^3$,
 \begin{align*}
   &\frac{\partial Q_n}{\partial \alpha'}({\bf y}_n;\alpha',\mu_\vare',\theta')\\
   &= \sumss 2\big(y_k-\alpha' y_{k-1}-\mu_\vare'\big)(-y_{k-1})
         + 2\big(y_s-\alpha' y_{s-1}-\mu_\vare'-\theta'\big)(-y_{s-1})\\
   &\phantom{=\;\,}
         + 2\big(y_{s+1}-\alpha' y_{s}-\mu_\vare'+\alpha'\theta'\big)(-y_s+\theta')\\
   &= 2\alpha'\left(\sumsf y_{k-1}^2 + (y_s-\theta')^2\right)
      +2\mu_\vare'\left(\sum_{k=1}^n y_{k-1}-\theta'\right)\\
   &\phantom{=\;\,}
      - 2\sumss y_{k-1}y_k - 2(y_s-\theta')y_{s-1} - 2y_{s+1}(y_s-\theta') \\
   &=\sum_{k=1}^n 2\big(y_k-\alpha' y_{k-1}-\mu_\vare'\big)(-y_{k-1})
         - 2\theta'(-y_{s-1}) + 2\alpha'\theta'(-y_s+\theta')
         + 2\big(y_{s+1}-\alpha' y_{s}-\mu_\vare'\big)\theta',
 \end{align*}
 \begin{align*}
   \frac{\partial Q_n}{\partial \mu_\vare'}&({\bf y}_n;\alpha',\mu_\vare',\theta')\\
   &= \sumss \big(y_k-\alpha' y_{k-1}-\mu_\vare'\big)(-2)
         - 2\big(y_s-\alpha' y_{s-1}-\mu_\vare'-\theta'\big)
         - 2\big(y_{s+1}-\alpha' y_s-\mu_\vare'+\alpha'\theta'\big)\\
   &=2\alpha'\left(\sum_{k=1}^n y_{k-1} -\theta'\right) + 2n\mu_\vare' - 2\sum_{k=1}^n y_k +2\theta',
 \end{align*}
 and
 \begin{align*}
   \frac{\partial Q_n}{\partial \theta'}({\bf y}_n;\alpha',\mu_\vare',\theta')
      =-2\big(y_s-\alpha' y_{s-1}-\mu_\vare'-\theta'\big)
       +2\big(y_{s+1}-\alpha' y_{s}-\mu_\vare'+\alpha'\theta'\big)\alpha'.
  \end{align*}

The next lemma is about the existence and uniqueness of the CLS estimator of \ $(\alpha,\mu_\vare,\theta)$.

\begin{Lem}\label{LEMMA5}
There exist subsets \ $S_n\subset\RR^{n+1}$, $n\geq \max(3,s+1)$ \ with the following properties:
 \begin{enumerate}
  \item[\upshape{(i)}]
   there exists a unique CLS estimator
   \ $(\halpha_n,\hmuen,\htheta_n):S_n\to\RR^3$,
  \item[\upshape{(ii)}]
   for all \ ${\bf y}_n\in S_n$,
   \ $(\halpha_n({\bf y}_n),\hmuen({\bf y}_n),\htheta_n({\bf y}_n))$
   \ is the unique solution of the system of equations
 \begin{align}\label{Additive_CLSE_EQ4}
   \begin{split}
    \frac{\partial Q_n}{\partial \alpha'}({\bf y}_n;\alpha',\mu_\vare',\theta')=0,\;\; 
    \frac{\partial Q_n}{\partial \mu_\vare'}({\bf y}_n;\alpha',\mu_\vare',\theta')=0,\;\; 
    \frac{\partial Q_n}{\partial \theta'}({\bf y}_n;\alpha',\mu_\vare',\theta')=0,
   \end{split}
 \end{align}
  \item [\upshape{(iii)}]
  ${\bf Y}_n\in S_n$ \ holds asymptotically as \ $n\to\infty$ \ with probability one.
 \end{enumerate}
\end{Lem}

\noindent{\bf Proof.}
For any fixed \ ${\bf y}_n \in \RR^{n+1}$, \ $n\geq \max(3,s+1)$ \ and \ $\alpha' \in \RR$,
 \ the quadratic function
 \ $\RR^2 \ni (\mu_\vare',\theta') \mapsto Q_n({\bf y}_n;\alpha',\mu_\vare',\theta')$
 \ can be written in the form
 \begin{align*}
  &Q_n({\bf y}_n;\alpha',\mu_\vare',\theta') \\[2mm]
  &=  \!\!\left(\! \begin{bmatrix}
             \mu_\vare' \smallskip \\
             \theta'
            \end{bmatrix}
           - A_n(\alpha')^{-1} t_n({\bf y}_n;\alpha') \!\! \right)^{\hspace*{-1mm}\top}
      \!\!\! A_n(\alpha') \! \left( \! \begin{bmatrix}
                          \mu_\vare' \smallskip \\
                          \theta'
                         \end{bmatrix}
                         - A_n(\alpha')^{-1} t_n({\bf y}_n;\alpha') \! \right)
    \!\! + \widehat{Q}_n({\bf y}_n;\alpha'),
 \end{align*}
 where
 \begin{align*}
  t_n({\bf y}_n;\alpha')
  &:= \begin{bmatrix}
       \sum_{k=1}^n ( y_k - \alpha' y_{k-1}) \smallskip \\
       ( 1 + (\alpha')^2 ) y_s - \alpha' ( y_{s-1} + y_{s+1} )
      \end{bmatrix} , \\[2mm]
  \widehat{Q}_n({\bf y}_n;\alpha')
  &:= \sum_{k=1}^n \big( y_k-\alpha' y_{k-1} \big)^2
      - t_n({\bf y}_n;\alpha')^\top A_n(\alpha')^{-1} t_n({\bf y}_n;\alpha') ,
 \end{align*}
 and the matrix
 \begin{align*}
  A_n(\alpha')
  := \begin{bmatrix}
      n & 1-\alpha' \smallskip \\
      1-\alpha' & 1+(\alpha')^2
     \end{bmatrix}
 \end{align*}
 is strictly positive definite for all \ $n\geq3$ \ and \ $\alpha' \in \RR$.
\ Indeed, the leading principal minors of \ $A_n(\alpha')$ \ take the following forms:
 \ $n,$
 \begin{align*}
  D_n(\alpha'):= n(1+(\alpha')^2) - (1-\alpha')^2 = (n-1)(\alpha')^2 + 2\alpha' + n-1,
 \end{align*}
 and for all \ $n\geq 3$, \ the discriminant \ $4-4(n-1)^2$ \ of the equation
 \ $(n-1)x^2 + 2x +n-1 = 0$ \ is negative.

The inverse matrix \ $A_n(\alpha')^{-1}$ \ takes the form
 \begin{align*}
  \frac{1}{D_n(\alpha')}\!\!
   \begin{bmatrix}
    1+(\alpha')^2 \, & \, -(1-\alpha') \\
    -(1-\alpha') \, & \, n
  \end{bmatrix}.
 \end{align*}
The polynomial \ $\RR\ni\alpha' \mapsto D_n(\alpha')$ \ is of order 2 with leading coefficient \ $n-1$.
\ We have \ $\widehat{Q}_n({\bf y}_n;\alpha') = R_n({\bf y}_n;\alpha') / D_n(\alpha')$, \ where
 \ $\RR\ni\alpha' \mapsto R_n({\bf y}_n;\alpha')$ \ is a polynomial of order 4 with leading coefficient
 \begin{align*}
   c_n({\bf y}_n) &:= (n-1) \sum_{k=1}^n y_{k-1}^2 - \left( \sum_{k=1}^n y_{k-1} \right)^2
                      + 2\left( \sum_{k=1}^n y_{k-1} \right)y_s - n y_s^2.
 \end{align*}
Let
 \[
   \widehat{S}_n := \left\{ {\bf y}_n\in\RR^{n+1} : c_n({\bf y}_n) > 0 \right\}.
 \]
For \ ${\bf y}_n \in \widehat{S}_n$, \ we have
 \ $\lim_{|\alpha'|\to\infty} \widehat{Q}_n({\bf y}_n;\alpha') = \infty$ \ and the continuous function
 \ $\RR \ni \alpha' \mapsto \widehat{Q}_n({\bf y}_n;\alpha')$ \ attains its infimum.
Consequently, for all \ $n\geq\max(3,s+1)$ \ there exists a CLS estimator
 \ $(\halpha_n, \hmuen, \htheta_n):\widehat{S}_n\to\RR^3$,
\ where
 \begin{align}\nonumber
   \widehat{Q}_n({\bf y}_n;\halpha_n({\bf y}_n))
   &=\inf_{\alpha'\in\RR} \widehat{Q}_n({\bf y}_n;\alpha')
   \qquad \forall\;{\bf y}_n\in \widehat{S}_n ,\\[2mm] \label{SEGED_UJ4}
  \begin{bmatrix}
   \hmuen({\bf y}_n) \smallskip \\
   \htheta_n({\bf y}_n)
  \end{bmatrix}
  & = A_n(\halpha_n({\bf y}_n))^{-1} t_n({\bf y}_n;\halpha_n({\bf y}_n)) ,
  \qquad {\bf y}_n\in \widehat{S}_n,
 \end{align}
 and for all \ ${\bf y}_n\in \widehat{S}_n$,
   \ $(\halpha_n({\bf y}_n),\hmuen({\bf y}_n),\htheta_n({\bf y}_n))$
   \ is a  solution of the system of equations \eqref{Additive_CLSE_EQ4}.

By \eqref{Ergodic1} and \eqref{Ergodic2}, we get
 \ $\PP\left(\lim_{n\to\infty} n^{-2} c_n({\bf Y}_n) = \var\widetilde X \right)=1$, \ where \ $\widetilde X$
 \ denotes a random variable with the unique stationary distribution of the INAR(1) model in \eqref{INAR1}.
Hence \ ${\bf Y}_n\in \widehat{S}_n$ \ holds asymptotically as \ $n\to\infty$ \ with probability one.

Now we turn to find sets \ $S_n \subset \widehat{S}_n$, $n\geq \max(3,s+1)$ \ such that the system of equations
 \eqref{Additive_CLSE_EQ4} has a unique solution with respect to
 \ $(\alpha',\mu_\vare',\theta')$ \ for all \ ${\bf y}_n\in S_n$.
\ Let us introduce the \ $(3\times 3)$ \ Hessian matrix
 \[
    H_n({\bf y}_n;\alpha',\mu_\vare',\theta'):=
         \begin{bmatrix}
          \frac{\partial^2 Q_n}{\partial (\alpha')^2}
           & \frac{\partial^2 Q_n}{\partial\mu_\vare' \partial \alpha'}
           & \frac{\partial^2 Q_n}{\partial\theta' \partial \alpha'} \\
           \frac{\partial^2 Q_n}{\partial \alpha'\partial\mu_\vare'}
           & \frac{\partial^2 Q_n}{\partial(\mu_\vare')^2 }
           & \frac{\partial^2 Q_n}{\partial\theta' \partial \mu_\vare'} \\
           \frac{\partial^2 Q_n}{\partial \alpha'\partial\theta'}
           & \frac{\partial^2 Q_n}{\partial\mu_\vare' \partial\theta'}
           & \frac{\partial^2 Q_n}{ \partial(\theta')^2} \\
        \end{bmatrix}({\bf y}_n;\alpha',\mu_\vare',\theta'),
 \]
 and let us denote by \ $\Delta_{i,n}({\bf y}_n;\alpha',\mu_\vare',\theta')$ \
 its \ $i$-th order leading principal minor, \ $i=1,2,3$.
\ Further, for all \ $n\geq \max(3,s+1)$, \ let
 \[
   S_n:=\Big\{{\bf y}_n\in \widehat{S}_n : \Delta_{i,n}({\bf y}_n;\alpha',\mu_\vare',\theta')>0,
                                   \;\, i=1,2,3,\, \forall\;(\alpha',\mu_\vare',\theta')
                                   \in\RR^3 \Big\}.
 \]
By Berkovitz \cite[Theorem 3.3, Chapter III]{Ber}, the function
 \ $\RR^3 \ni (\alpha',\mu_\vare',\theta')
     \mapsto Q_n({\bf y}_n;\alpha',\mu_\vare',\theta')$
 \ is strictly convex for all \ ${\bf y}_n\in S_n$.
\ Since it was already proved that the system of equations \eqref{Additive_CLSE_EQ4} has a solution for all
 \ ${\bf y}_n\in \widehat{S}_n$,
 \ we obtain that this solution is unique for all \ ${\bf y}_n\in S_n$.

Next we check that \ ${\bf Y}_n\in S_n$ \ holds asymptotically as \ $n\to\infty$ \ with probability one.
Using also the proof of Lemma \ref{LEMMA2}, for all \ $(\alpha',\mu_\vare',\theta')\in\RR^3$, \ we get
 \begin{align*}
   & \frac{\partial^2 Q_n}{\partial (\alpha')^2}({\bf Y}_n;\alpha',\mu_\vare',\theta')
       = 2\sum_{k=1}^n Y_{k-1}^2 + 2\theta'(-Y_s+\theta')-2Y_s\theta' \\
   &\phantom{\frac{\partial^2 Q_n}{\partial (\alpha')^2}({\bf Y}_n;\alpha',\mu_\vare',\theta') \;}
     = 2\left(\sumsf X_{k-1}^2+(X_s+\theta-\theta')^2\right),\\
   &\frac{\partial^2 Q_n}{\partial\alpha'\partial\theta'}({\bf Y}_n;\alpha',\mu_\vare',\theta')
    = \frac{\partial^2 Q_n}{\partial\theta'\partial\alpha'}({\bf Y}_n;\alpha',\mu_\vare',\theta')
    = 2(Y_{s-1}+Y_{s+1}-2\alpha' Y_s-\mu_\vare'+2\alpha'\theta')\\
  &\phantom{\frac{\partial^2 Q_n}{\partial\alpha'\partial\theta'}({\bf Y}_n;\alpha',\mu_\vare',\theta')}
    =2(X_{s-1}+X_{s+1}-2\alpha' X_s-\mu_\vare'-2\alpha'(\theta-\theta')),\\[1mm]
  \end{align*}
 and
 \begin{align*}
    & \frac{\partial^2 Q_n}{\partial\alpha'\partial\mu_\vare'}({\bf Y}_n;\alpha',\mu_\vare',\theta')
    = \frac{\partial^2 Q_n}{\partial\mu_\vare'\partial\alpha'}({\bf Y}_n;\alpha',\mu_\vare',\theta')
    = 2 \sum_{k=1}^n Y_{k-1} - 2\theta'
    = 2\left(\sum_{k=1}^n X_{k-1} + \theta - \theta'\right),\\[1mm]
  & \frac{\partial^2 Q_n}{\partial\theta'\partial\mu_\vare'}({\bf Y}_n;\alpha',\mu_\vare',\theta')
    = \frac{\partial^2 Q_n}{\partial\mu_\vare'\partial\theta'}({\bf Y}_n;\alpha',\mu_\vare',\theta')
    = 2(1-\alpha'),\\
  & \frac{\partial^2 Q_n}{\partial(\theta')^2}({\bf Y}_n;\alpha',\mu_\vare',\theta')
    = 2((\alpha')^2+1),
    \qquad \qquad
  \frac{\partial^2 Q_n}{\partial(\mu_\vare')^2}({\bf Y}_n;\alpha',\mu_\vare',\theta')
    = 2n.
 \end{align*}
Then
 \begin{align*}
   &H_n({\bf Y}_n;\alpha',\mu_\vare',\theta')\\
   &=2
      \begin{bmatrix}
        \DS\sumsf X_{k-1}^2+(X_s+\theta-\theta')^2
         & \sum_{k=1}^n X_{k-1} + \theta - \theta'  & a  \\
         \sum_{k=1}^n X_{k-1} + \theta - \theta' & n & 1-\alpha'  \\
         a & 1-\alpha' & (\alpha')^2+1
      \end{bmatrix},
 \end{align*}
 where \ $a:=X_{s-1}+X_{s+1}-2\alpha' X_s-\mu_\vare'-2\alpha'(\theta-\theta')$.
\ Then \ $H_n({\bf Y}_n;\alpha',\mu_\vare',\theta')$ \ has leading principal minors
 \ $\Delta_{1,n}({\bf Y}_n;\alpha',\mu_\vare',\theta')
       :=2\left(\DS\sumsf X_{k-1}^2+(X_s+\theta-\theta')^2\right)$,
 \begin{align*}
    \Delta_{2,n}({\bf Y}_n;\alpha',\mu_\vare',\theta'):=
      4\left( n\left(\sumsf X_{k-1}^2+(X_s+\theta-\theta')^2\right)
                - \left(\sum_{k=1}^n X_{k-1} +\theta-\theta' \right)^2 \right),
 \end{align*}
 and
 \begin{align*}
   \Delta_{3,n}({\bf Y}_n;\alpha',\mu_\vare',\theta')
         & = \det H_n({\bf Y}_n;\alpha',\mu_\vare',\theta')\\
         & = 8\Bigg[ n \Bigg(((\alpha')^2+1)\left(\sumsf X_{k-1}^2+(X_s+\theta-\theta')^2\right)
                         - a^2 \Bigg)\\
   &\phantom{=\Bigg[\;}
    -((\alpha')^2+1)\left(\sum_{k=1}^n X_{k-1} + \theta - \theta' \right)^2
     +2(1-\alpha')a \left(\sum_{k=1}^n X_{k-1} + \theta - \theta'\right) \\
   & \phantom{=\Bigg[\;}
     -(1-\alpha')^2\left( \sumsf X_{k-1}^2+(X_s+\theta-\theta')^2 \right)
         \Bigg].
 \end{align*}
By \eqref{Ergodic1} and \eqref{Ergodic2}, we have
 \begin{align*}
   &\PP\left(\lim_{n\to\infty}
              \frac{1}{2n}\Delta_{1,n}({\bf Y}_n;\alpha',\mu_\vare',\theta')
              =\EE\widetilde X^2,\;\; \forall\;\; (\alpha',\mu_\vare',\theta')\in\RR^3
        \right)=1,\\
   &\PP\left(\lim_{n\to\infty}
              \frac{1}{4n^2}\Delta_{2,n}({\bf Y}_n;\alpha',\mu_\vare',\theta')
              =\var\widetilde X,\;\; \forall\;\; (\alpha',\mu_\vare',\theta')\in\RR^3
        \right)=1,
 \end{align*}
 and
 \begin{align*}
   \PP\left(\lim_{n\to\infty}
              \frac{1}{8n^2}\Delta_{3,n}({\bf Y}_n;\alpha',\mu_\vare',\theta')
              = ((\alpha')^2+1)\var \widetilde X,
              \;\; \forall\;\; (\alpha',\mu_\vare',\theta')\in\RR^3 \right)=1,
 \end{align*}
 for all \ $(\alpha',\mu_\vare',\theta')\in\RR^3$, \ where \ $\widetilde X$ \ denotes
 a random variable with the unique stationary distribution of the INAR(1) model in \eqref{INAR1}.
Hence
 \begin{align*}
   &\PP\big(\lim_{n\to\infty}\Delta_{i,n}({\bf Y}_n;\alpha',\mu_\vare',\theta') =\infty,
     \;\; \forall\;\; (\alpha',\mu_\vare',\theta')\in\RR^3 \big)=1,
    \qquad i=1,2,3,
 \end{align*}
  which yields that \ ${\bf Y}_n\in S_n$ \ asymptotically as \ $n\to\infty$
 \ with probability one, since we have already proved that \ ${\bf Y}_n\in \widehat{S}_n$
 \ asymptotically as \ $n\to\infty$ \ with probability one.
\proofend

By Lemma \ref{LEMMA5}, \ $(\halpha_n({\bf Y}_n), \hmuen({\bf Y}_n), \htheta_n({\bf Y}_n))$
 \ exists uniquely asymptotically as \ $n\to\infty$ \ with probability one.
In the sequel we will simply denote it
 by \ $(\halpha_n,\hmuen,\htheta_n)$.

The next result shows that \ $\halpha_n$ \ and \ $\hmuen$ \ are strongly consistent estimators
 of \ $\alpha$ \ and \ $\mu_\vare$, \ respectively, whereas \ $\htheta_n$ \ fails to be a strongly
 consistent estimator of \ $\theta$.

\begin{Thm}\label{THEOREM5}
For the CLS estimators \ $(\halpha_n,\hmuen,\htheta_n)_{n\in\NN}$ \ of
 the parameter \ $(\alpha,\mu_\vare,\theta)\in(0,1)\times(0,\infty)\times\NN$,
\ the sequences \ $(\halpha_n)_{n\in\NN}$ \ and \ $(\hmuen)_{n\in\NN}$ \ are strongly consistent for all
 \ $(\alpha,\mu_\vare,\theta)\in(0,1)\times(0,\infty)\times\NN$, \ i.e.,
 \begin{align}\label{Strong_consistency7}
   &\PP(\lim_{n\to\infty}\halpha_n=\alpha)=1,
         \qquad \forall\;(\alpha,\mu_\vare,\theta)\in(0,1)\times(0,\infty)\times\NN,\\ \label{Strong_consistency12}
   &\PP(\lim_{n\to\infty}\hmuen=\mu_\vare)=1,
      \qquad \forall\; (\alpha,\mu_\vare,\theta)\in(0,1)\times(0,\infty)\times\NN,
 \end{align}
 whereas the sequence \ $(\htheta_n)_{n\in\NN}$ \ is not strongly consistent
 for any \ $(\alpha,\mu_\vare,\theta)\in(0,1)\times(0,\infty)\times\NN$,
 \ namely,
 \begin{align}\label{Strong_consistency13}
   \PP\left(\lim_{n\to\infty}\htheta_n= Y_s-\frac{\alpha}{1+\alpha^2}(Y_{s-1}+Y_{s+1})
              -\frac{1-\alpha}{1+\alpha^2}\mu_\vare\right)=1,
 \end{align}
 for all \ $(\alpha,\mu_\vare,\theta)\in(0,1)\times(0,\infty)\times\NN$.
\end{Thm}

\noindent{\bf Proof.}
An easy calculation shows that
 \begin{align*}
  &\left(\sum_{k=1}^n Y_{k-1}^2 + (Y_s-\htheta_n)^2 - Y_s^2\right)\halpha_n
     + \left(\sum_{k=1}^n Y_{k-1} - \htheta_n\right)\hmuen
          = \sum_{k=1}^n Y_{k-1}Y_k - \htheta_n(Y_{s-1}+Y_{s+1})\\
  &\left(\sum_{k=1}^n Y_{k-1} -\htheta_n\right)\halpha_n + n\hmuen
          = \sum_{k=1}^n Y_k - \htheta_n,
 \end{align*}
 or equivalently
 \begin{align}\label{SEGED40}
  \begin{split}
   &\begin{bmatrix}
     \sum_{k=1}^n Y_{k-1}^2 + (Y_s-\htheta_n)^2-Y_s^2 &  \sum_{k=1}^n Y_{k-1} -\htheta_n \\
     \sum_{k=1}^n Y_{k-1} -\htheta_n & n \\
    \end{bmatrix}
    \begin{bmatrix}
     \halpha_n \\
     \hmuen \\
   \end{bmatrix}\\
  &\qquad\qquad
     =\begin{bmatrix}
       \sum_{k=1}^n Y_{k-1}Y_k - \htheta_n(Y_{s-1}+Y_{s+1}) \\
       \sum_{k=1}^n Y_k - \htheta_n \\
      \end{bmatrix}
   \end{split}
 \end{align}
 holds asymptotically as \ $n\to\infty$ \ with probability one.
Let us introduce the notation
 \begin{align*}
   E_n& := n\left(\sum_{k=1}^n Y_{k-1}^2 + (Y_s-\htheta_n)^2-Y_s^2\right)
          - \left(\sum_{k=1}^n Y_{k-1} -\htheta_n\right)^2\\
       & = n\left(\sumsf X_{k-1}^2 + (X_s+\theta-\htheta_n)^2\right)
          - \left(\sum_{k=1}^n X_{k-1} +\theta-\htheta_n\right)^2,
          \qquad n\geq s+1, \;\;n\in\NN.
 \end{align*}
By \eqref{Ergodic1} and \eqref{Ergodic2},
 \begin{align}\label{SEGED39}
    \PP\left(\lim_{n\to\infty}\frac{E_n}{n^2}=\EE\widetilde X^2 - (\EE\widetilde X)^2
                 = \var\widetilde X >0 \right)=1,
 \end{align}
 which yields that \ $\PP(\lim_{n\to\infty} E_n=\infty)=1$.
\ Hence asymptotically as \ $n\to\infty$ \ with probability one we get
 \begin{align} \label{SEGED92}
   \begin{split}
   \begin{bmatrix}
     \halpha_n \\
     \hmuen \\
    \end{bmatrix}
    &=\frac{1}{E_n}
      \begin{bmatrix}
        n & -\sum_{k=1}^n Y_{k-1} +\htheta_n \\
        -\sum_{k=1}^n Y_{k-1} +\htheta_n & \sum_{k=1}^n Y_{k-1}^2 + (Y_s-\htheta_n)^2-Y_s^2  \\
      \end{bmatrix}\\
    &\phantom{=\frac{1}{E_n}\;}\times
      \begin{bmatrix}
        \sum_{k=1}^n Y_{k-1}Y_k - \htheta_n(Y_{s-1}+Y_{s+1}) \\
        \sum_{k=1}^n Y_k - \htheta_n \\
      \end{bmatrix}\\
    &=\frac{1}{E_n}
      \begin{bmatrix}
        n & -\sum_{k=1}^n X_{k-1} +(\htheta_n-\theta) \\
        -\sum_{k=1}^n X_{k-1} +(\htheta_n-\theta) & \sum_{k=1}^n X_{k-1}^2 + (X_s+\theta-\htheta_n)^2-X_s^2  \\
      \end{bmatrix}\\
    &\phantom{=\frac{1}{E_n}\;}\times
      \begin{bmatrix}
        \sum_{k=1}^n X_{k-1}X_k + (\theta-\htheta_n)(X_{s-1}+X_{s+1}) \\
        \sum_{k=1}^n X_k - (\htheta_n -\theta)\\
      \end{bmatrix}\\
    &=:\frac{1}{E_n}\begin{bmatrix}
                     V^{(1)}_n \\
                     V^{(2)}_n \\
                   \end{bmatrix},
   \end{split}
 \end{align}
 where
 \begin{align*}
  V^{(1)}_n:=&n\sum_{k=1}^n X_{k-1}X_k - \left(\sum_{k=1}^n X_{k-1}\right)\left(\sum_{k=1}^n X_k\right)
        + n(\theta-\htheta_n)(X_{s-1}+X_{s+1})
        + (\htheta_n -\theta)\sum_{k=1}^n X_{k-1}\\
        &+ (\htheta_n -\theta)\sum_{k=1}^n X_k
        - (\htheta_n -\theta)^2,
 \end{align*}
  and
 \begin{align*}
    V^{(2)}_n:=&\left(\sum_{k=1}^n X_{k-1}^2\right)\left(\sum_{k=1}^n X_k\right)
         - \left(\sum_{k=1}^n X_{k-1}\right)\left(\sum_{k=1}^n X_{k-1}X_k\right)
         - (\theta-\htheta_n)(X_{s-1}+X_{s+1})\sum_{k=1}^n X_{k-1} \\
        &+ (\htheta_n -\theta)\sum_{k=1}^n X_{k-1}X_k
         - (\htheta_n-\theta)^2(X_{s-1}+X_{s+1})
         - (\htheta_n-\theta)\sum_{k=1}^n X_{k-1}^2\\
        & + (X_s+\theta-\htheta_n)^2\sum_{k=1}^n X_k
         - (\htheta_n-\theta)(X_s+\theta-\htheta_n)^2
         - X_s^2 \sum_{k=1}^n X_k
         + X_s^2 (\htheta_n-\theta).
 \end{align*}

Similarly, an easy calculation shows that
 \begin{align*}
   & n\hmuen + (1-\halpha_n)\htheta_n = \sum_{k=1}^n Y_k-\halpha_n \sum_{k=1}^n Y_{k-1},\\
   & (1-\halpha_n)\hmuen + (1+(\halpha_n)^2)\htheta_n
         = (1+(\halpha_n)^2)Y_s -\halpha_n(Y_{s-1}+Y_{s+1}),
 \end{align*}
 or equivalently
 \begin{align*}
   \begin{bmatrix}
     n & 1-\halpha_n \\
     1-\halpha_n & 1+(\halpha_n)^2 \\
   \end{bmatrix}
      \begin{bmatrix}
        \hmuen \\
        \htheta_n \\
      \end{bmatrix}
   =
     \begin{bmatrix}
          \sum_{k=1}^n Y_k - \halpha_n \sum_{k=1}^n Y_{k-1}) \\
         (1+(\halpha_n)^2)Y_s -\halpha_n(Y_{s-1}+Y_{s+1}) \\
      \end{bmatrix}
 \end{align*}
 holds asymptotically as \ $n\to\infty$ \ with probability one.
Recalling that \ $D_n(\halpha_n)=n(1+(\halpha_n)^2)-(1-\halpha_n)^2$, \ we have
 asymptotically as \ $n\to\infty$ \ with probability one,
 \begin{align*}
   \begin{bmatrix}
      \hmuen \\
      \htheta_n \\
     \end{bmatrix}
    &=\frac{1}{D_n(\halpha_n)}
      \begin{bmatrix}
        1+(\halpha_n)^2 & -(1-\halpha_n) \\
        -(1-\halpha_n) &  n \\
      \end{bmatrix}
      \begin{bmatrix}
          \sum_{k=1}^n Y_k-\halpha_n \sum_{k=1}^n Y_{k-1} \\
         (1+(\halpha_n)^2)Y_s -\halpha_n(Y_{s-1}+Y_{s+1}) \\
      \end{bmatrix} \\[1mm]
    &=\begin{bmatrix}
        \frac{(1+(\halpha_n)^2)\left(\sum_{k=1}^nY_k-\halpha_n\sum_{k=1}^nY_{k-1}\right)
             - (1-\halpha_n)\left((1+(\halpha_n)^2)Y_s-\halpha_n(Y_{s-1}+Y_{s+1})\right)}
             {D_n(\halpha_n)}
             \\[2mm]
       \frac{ - (1-\halpha_n)\left(\sum_{k=1}^nY_k-\halpha_n\sum_{k=1}^nY_{k-1}\right)
           + n\left((1+(\halpha_n)^2)Y_s-\halpha_n(Y_{s-1}+Y_{s+1})\right)}
           {D_n(\halpha_n)} \\
      \end{bmatrix}.
 \end{align*}

We show that the sequence \ $(\htheta_n - \theta)_{n\in\NN}$ \ is bounded
 with probability one.
\ Using the decomposition \ $Y_k=X_k+\delta_{k,s}\theta$, \ $k\in\ZZ_+$, \ we get
 \begin{align}\label{SEGED33}
    \begin{bmatrix}
       \hmuen - \mu_\vare\\
       \htheta_n -\theta \\
     \end{bmatrix}
    = \frac{1}{D_n(\halpha_n)}
     \begin{bmatrix}
       V^{(3)}_n \\
       V^{(4)}_n \\
     \end{bmatrix},
  \end{align}
 holds asymptotically as \ $n\to\infty$ \ with probability one, where
 \begin{align*}
   & V^{(3)}_n
      := (1+(\halpha_n)^2)\left(\sum_{k=1}^nX_k-\halpha_n\sum_{k=1}^nX_{k-1}
                                  +(1-\halpha_n)\theta\right) \\
   &\phantom{V^{(3)}_n :=}
      -(1-\halpha_n)\Big((1+(\halpha_n)^2)X_s-\halpha_n(X_{s-1}+X_{s+1})+(1+(\halpha_n)^2)\theta\Big) \\
   &\phantom{V^{(3)}_n :=}
      - n(1+(\halpha_n)^2)\mu_\vare +(1-\halpha_n)^2\mu_\vare\\
   &\phantom{V^{(3)}_n }=
      (1+(\halpha_n)^2)\left(\sum_{k=1}^nX_k-\halpha_n\sum_{k=1}^nX_{k-1}
                                  -n\mu_\vare\right)\\
   &\phantom{V^{(3)}_n :=}
             -(1-\halpha_n)\Big((1+(\halpha_n)^2)X_s-\halpha_n(X_{s-1}+X_{s+1})
                  - (1-\halpha_n)\mu_\vare \Big),
 \end{align*}
 and
 \begin{align*}
    V^{(4)}_n
      := - (1-\halpha_n)\left(\sum_{k=1}^nX_k-\halpha_n\sum_{k=1}^nX_{k-1}\right)
             +n\Big( (1+(\halpha_n)^2) X_s -\halpha_n(X_{s-1}+X_{s+1}) \Big).
 \end{align*}
\ By \eqref{SEGED33}, we have asymptotically as \ $n\to\infty$ \ with probability one,
 \begin{align*}
   \vert \htheta_n - \theta \vert
      &\leq \frac{(1+(\halpha_n)^2)n}{(1+(\halpha_n)^2)n-(1-\halpha_n)^2}
           \left[
             \frac{\vert 1-\halpha_n \vert}{1+(\halpha_n)^2}
             \frac{\sum_{k=1}^nX_k}{n}
             + \frac{\vert \halpha_n(1-\halpha_n) \vert}{1+(\halpha_n)^2}
               \frac{\sum_{k=1}^nX_{k-1}}{n} \right.\\
       &\phantom{\leq \frac{(1+(\halpha_n)^2)n}{(1+(\halpha_n)^2)n-(1-\halpha_n)^2}\;\;\;}\left.
             + X_s
             + \frac{\vert \halpha_n \vert}{1+(\halpha_n)^2}(X_{s-1}+X_{s+1})
             \right] \\
       & \leq \frac{1}{1-\frac{(1-\halpha_n)^2}{n(1+(\halpha_n)^2)}}
             \left[\frac{3}{2}\frac{\sum_{k=1}^nX_k}{n}
                  +\frac{3}{2}\frac{\sum_{k=1}^nX_{k-1}}{n}
                  +X_s+X_{s-1}+X_{s+1} \right], \\
       & \leq \frac{1}{1-\frac{(1-\halpha_n)^2}{n(1+(\halpha_n)^2)}}
             \left[3\frac{\sum_{k=0}^nX_k}{n}
                   +X_s+X_{s-1}+X_{s+1} \right],
 \end{align*}
Using \eqref{Ergodic1} and that
 \[
   \frac{(1-\halpha_n)^2}{1+(\halpha_n)^2} < 3, \qquad n\in\NN,
 \]
 we have the sequences \ $(\htheta_n -\theta)_{n\in\NN}$ \ and \ $(\htheta_n)_{n\in\NN}$ \ are bounded
 with probability one.

Similarly to \eqref{SEGED2}, one can check that
 \begin{align}\label{SEGED93}
  \htheta_n
      = Y_s-\frac{\halpha_n}{1+(\halpha_n)^2}(Y_{s-1}+Y_{s+1})
                  -\frac{1-\halpha_n}{1+(\halpha_n)^2}\hmuen
 \end{align}
 holds asymptotically as \ $n\to\infty$ \ with probability one.

Using \eqref{SEGED39} and \eqref{SEGED92}, to prove \eqref{Strong_consistency7} and
 \eqref{Strong_consistency12}, it is enough to check that
 \[
   \PP\left(\lim_{n\to\infty}\frac{V^{(1)}_n}{n^2}=\alpha\var\widetilde X\right)=1
    \qquad \text{and} \qquad \PP\left(\lim_{n\to\infty}\frac{V^{(2)}_n}{n^2}=\mu_\vare\var\widetilde X\right)=1,
  \]
 for all \ $(\alpha,\mu_\vare,\theta)\in(0,1)\times(0,\infty)\times\NN$.
\ Using that the sequence \ $(\htheta_n -\theta)_{n\in\NN}$ \ is bounded with probability one,
 \ by \eqref{Ergodic1}, \eqref{Ergodic2} and \eqref{Ergodic3}, we get with probability one
 \begin{align*}
   &\lim_{n\to\infty}\frac{V^{(1)}_n}{n^2}
      = \alpha\EE\widetilde X^2 + \mu_\vare\EE\widetilde X
        - (\EE\widetilde X)^2
      = \alpha\var\widetilde X + \mu_\vare\EE\widetilde X
        + (\alpha-1) (\EE\widetilde X)^2
      = \alpha\var\widetilde X,\\
   &\lim_{n\to\infty}\frac{V^{(2)}_n}{n^2}
      = \EE\widetilde X^2 \EE\widetilde X
        - \EE\widetilde X(\alpha\EE\widetilde X^2 + \mu_\vare\EE\widetilde X)
      =\mu_\vare\var\widetilde X+((1-\alpha)\EE\widetilde X-\mu_\vare)\EE\widetilde X^2
      =\mu_\vare\var\widetilde X,
 \end{align*}
 where the last equality follows by \eqref{STAC_MOMENT1}.

Finally, \eqref{Strong_consistency13} follows from \eqref{SEGED93}, \eqref{Strong_consistency7}
 and  \eqref{Strong_consistency12}.
\proofend

The asymptotic distribution of the CLS estimation is given in the next theorem.

\begin{Thm}
Under the additional assumptions \ $\EE X_0^3<\infty$ \ and \ $\EE\vare_1^3<\infty$, \ we have
   \begin{align} \label{CONVERGENCE8}
     \begin{bmatrix}
       \sqrt{n}(\halpha_n-\alpha) \\
       \sqrt{n}(\hmuen-\mu_\vare) \\
     \end{bmatrix}
      \distr \cN\left(\begin{bmatrix}
                       0 \\
                       0  \\
                      \end{bmatrix},B_{\alpha,\vare}\right)
      \qquad \text{as \ $n\to\infty$,}
  \end{align}
 where the \ $(2\times 2)$-matrix \ $B_{\alpha,\vare}$ \ is defined in \eqref{SEGED_BALPHA}.
Moreover, conditionally on the values \ $Y_{s-1}$ \ and \ $Y_{s+1}$,
 \begin{align}\label{CONVERGENCE9}
   \sqrt{n}\left(\htheta_n - \lim_{k\to\infty}\htheta_k\right)
         \distr \cN\left(0,d_{\alpha,\,\vare}^\top B_{\alpha,\,\vare} d_{\alpha,\,\vare}\right)
         \qquad \text{as \ $n\to\infty$,}
 \end{align}
 where
 \[
    d_{\alpha,\,\vare}
      :=\frac{1}{(1+\alpha^2)^2}
         \begin{bmatrix}
           (\alpha^2-1)(Y_{s-1}+Y_{s+1}) + (2\alpha+1-\alpha^2)\mu_\vare \\
            -(1+\alpha^2)(1-\alpha) \\
         \end{bmatrix}.
 \]
\end{Thm}

\noindent{\bf Proof.}
 By \eqref{SEGED40}, with the notation
 \[
   B_n:=
     \begin{bmatrix}
      \sum_{k=1}^n Y_{k-1}^2 + (Y_s-\htheta_n)^2-Y_s^2 &  \sum_{k=1}^n Y_{k-1} -\htheta_n \\
      \sum_{k=1}^n Y_{k-1} -\htheta_n & n \\
     \end{bmatrix}, \qquad n\in\NN,
 \]
 we get
 \[
    \begin{bmatrix}
     \halpha_n \\
     \hmuen \\
   \end{bmatrix}
   =B_n^{-1}
      \begin{bmatrix}
       \sum_{k=1}^n Y_{k-1}Y_k - \htheta_n(Y_{s-1}+Y_{s+1}) \\
       \sum_{k=1}^n Y_k - \htheta_n \\
      \end{bmatrix},
 \]
 holds asymptotically as \ $n\to\infty$ \ with probability one.
Hence
 \begin{align*}
   \begin{bmatrix}
     \halpha_n -\alpha\\
     \hmuen -\mu_\vare\\
    \end{bmatrix}
   &=B_n^{-1}
     \left(
       \begin{bmatrix}
         \sum_{k=1}^n Y_{k-1}Y_k - \htheta_n(Y_{s-1}+Y_{s+1}) \\
          \sum_{k=1}^n Y_k - \htheta_n \\
       \end{bmatrix}
         - B_n
         \begin{bmatrix}
           \alpha \\
           \mu_\vare \\
         \end{bmatrix}
     \right)\\
    &=B_n^{-1}
      \left(
       \begin{bmatrix}
         \sum_{k=1}^n Y_{k-1}Y_k - \htheta_n(Y_{s-1}+Y_{s+1}) \\
         \sum_{k=1}^n Y_k - \htheta_n \\
      \end{bmatrix} \right.\\
    &\phantom{\;\;=B_n^{-1}\left(\right.} \left.
      - \begin{bmatrix}
         \alpha\sum_{k=1}^n Y_{k-1}^2 +\mu_\vare\sum_{k=1}^n Y_{k-1}
           + \alpha (Y_s-\htheta_n)^2 -\alpha Y_s^2 -\mu_\vare\htheta_n \\
         \alpha\sum_{k=1}^n Y_{k-1} + n\mu_\vare - \alpha\htheta_n \\
       \end{bmatrix}
      \right).
 \end{align*}
Then
 \[
   \begin{bmatrix}
     \halpha_n -\alpha\\
     \hmuen -\mu_\vare\\
    \end{bmatrix}
    = B_n^{-1}
     \begin{bmatrix}
      V_n^{(5)}\\
      V_n^{(6)}\\
     \end{bmatrix}
 \]
 holds asymptotically as \ $n\to\infty$ \ with probability one, where
 \begin{align*}
     &V_n^{(5)} :=  \sum_{k=1}^n (Y_k-\alpha Y_{k-1}-\mu_\vare)Y_{k-1}
                   - \htheta_n(Y_{s-1}+Y_{s+1}) + 2\alpha Y_s\htheta_n
                   - \alpha (\htheta_n)^2 + \mu_\vare \htheta_n,\\
     &V_n^{(6)} := \sum_{k=1}^n (Y_k-\alpha Y_{k-1}-\mu_\vare) - (1-\alpha)\htheta_n.
 \end{align*}
To prove \eqref{CONVERGENCE8}, it is enough to show that
 \begin{align}\label{SEGED41}
   &\PP\left(\lim_{n\to\infty}\frac{B_n}{n}
              = \begin{bmatrix}
                  \EE\widetilde X^2 & \EE\widetilde X \\
                  \EE\widetilde X & 1 \\
                \end{bmatrix}
              \right)=1,\\ \label{SEGED42}
   &\frac{1}{\sqrt n}
          \begin{bmatrix}
            V_n^{(5)} \\
            V_n^{(6)} \\
          \end{bmatrix}
         \distr \cN\left(
           \begin{bmatrix}
             0 \\
             0 \\
          \end{bmatrix},A_{\alpha,\vare}\right)
         \qquad \text{as \ $n\to\infty$, }
 \end{align}
 where \ $\widetilde X$ \ is a random variable having the unique stationary distribution
 of the INAR(1) model in \eqref{INAR1} and the \ $(2\times 2)$-matrix \ $A_{\alpha,\vare}$
 \ is defined in \eqref{SEGED_AALPHA}.
Using \eqref{Ergodic1}, \eqref{Ergodic2} and that the sequence \ $(\htheta_n)_{n\in\NN}$
 \ is bounded with probability one, we get \eqref{SEGED41}.
Now we turn to prove \eqref{SEGED42}.
An easy calculation shows that
 \begin{align*}
   V_n^{(5)}
    & = \sumss (X_k-\alpha X_{k-1}-\mu_\vare)X_{k-1} + (X_s-\alpha X_{s-1}-\mu_\vare+\theta)X_{s-1} \\
    &\phantom{=\;}
       + (X_{s+1}-\alpha X_s-\mu_\vare-\alpha\theta)(X_s+\theta)
       - \htheta_n(X_{s-1}+X_{s+1}) + 2\alpha (X_s+\theta)\htheta_n
                   - \alpha (\htheta_n)^2 + \mu_\vare \htheta_n  \\
    &= \sum_{k=1}^n (X_k-\alpha X_{k-1}-\mu_\vare)X_{k-1}
      +\theta (X_{s-1}+X_{s+1}) -2\alpha \theta X_s - \alpha\theta^2-\theta\mu_\vare  \\
    &\phantom{=\;}
      - \htheta_n(X_{s-1}+X_{s+1}) + 2\alpha (X_s+\theta)\htheta_n
                   - \alpha (\htheta_n)^2 + \mu_\vare \htheta_n \\
   &= \sum_{k=1}^n (X_k-\alpha X_{k-1}-\mu_\vare)X_{k-1}
     + (\theta-\htheta_n)(X_{s-1}+X_{s+1}-2\alpha X_s-\mu_\vare - \alpha(\theta-\htheta_n)),
 \end{align*}
 and
 \begin{align*}
   V_n^{(6)}
     = \sum_{k=1}^n (X_k-\alpha X_{k-1}-\mu_\vare) + (1-\alpha)(\theta-\htheta_n).
 \end{align*}
By formula (6.43) in Hall and Heyde \cite[Section 6.3]{HalHey},
 \begin{align*}
   \begin{bmatrix}
     \frac{1}{\sqrt n} \sum_{k=1}^n X_{k-1}(X_k-\alpha X_{k-1}-\mu_\vare)   \\
     \frac{1}{\sqrt n} \sum_{k=1}^n (X_k-\alpha X_{k-1}-\mu_\vare) \\
   \end{bmatrix}
         \distr
   \cN\left(\begin{bmatrix}
              0 \\
              0 \\
            \end{bmatrix}
             ,A_{\alpha,\vare}\right)
     \qquad \text{as \ $n\to\infty$.}
 \end{align*}
Using that the sequence \ $(\htheta_n-\theta)_{n\in\NN}$ \ is bounded with probability one,
 by Slutsky's lemma, we get \eqref{SEGED42}.

Now we turn to prove \eqref{CONVERGENCE9}.
Using \eqref{Strong_consistency13} and \eqref{SEGED93}, we have
 \begin{align*}
   &\sqrt{n}(\htheta_n-\lim_{k\to\infty}\htheta_k)
       = \sqrt{n}
        \left(
          \htheta_n
            - \left(
                  Y_s - \frac{\alpha}{1+\alpha^2}(Y_{s-1}+Y_{s+1})
                      - \frac{1-\alpha}{1+\alpha^2}\mu_\vare
              \right)
         \right)\\
    & = \sqrt{n}\left( \left(\frac{\alpha}{1+\alpha^2}-\frac{\halpha_n}{1+(\halpha_n)^2}\right)
                        (Y_{s-1}+Y_{s+1})
                       + \frac{1-\alpha}{1+\alpha^2}\mu_\vare
                       - \frac{1-\halpha_n}{1+(\halpha_n)^2}\hmuen
                        \right) \\
   & = \sqrt{n}\left(\frac{(\halpha_n-\alpha)(\alpha\halpha_n-1)}{(1+\alpha^2)(1+(\halpha_n)^2)}
                       (Y_{s-1}+Y_{s+1})
                    + \frac{1-\alpha}{1+\alpha^2}(\mu_\vare - \hmuen)
                    + \left(\frac{1-\alpha}{1+\alpha^2} - \frac{1-\halpha_n}{1+(\halpha_n)^2}\right)\hmuen
                        \right) \\
   & = \sqrt{n}\left(
                     \frac{(\alpha\halpha_n-1)(Y_{s-1}+Y_{s+1})
                            + (\halpha_n+\alpha+1-\alpha\halpha_n)\hmuen}{(1+\alpha^2)(1+(\halpha_n)^2)}
                       (\halpha_n-\alpha)
                     - \frac{1-\alpha}{1+\alpha^2}(\hmuen-\mu_\vare)
             \right)
 \end{align*}
 holds asymptotically as \ $n\to\infty$ \ with probability one.
Hence
 \begin{align*}
   \sqrt{n}(\htheta_n &- \lim_{k\to\infty}\htheta_k) \\
   &= \begin{bmatrix}
         \DS\frac{(\alpha\halpha_n-1)(Y_{s-1}+Y_{s+1}) + (\halpha_n+\alpha+1-\alpha\halpha_n)\hmuen}
               {(1+\alpha^2)(1+(\halpha_n)^2)}
          &  \DS-\frac{1-\alpha}{1+\alpha^2} \\
       \end{bmatrix}
       \begin{bmatrix}
         \sqrt{n}(\halpha_n-\alpha) \\
         \sqrt{n}(\hmuen-\mu_\vare)  \\
       \end{bmatrix}
 \end{align*}
 holds asymptotically as \ $n\to\infty$ \ with probability one.
Using Slutsky's lemma, by \eqref{Strong_consistency7}, \eqref{Strong_consistency12} and
 \eqref{CONVERGENCE8}, we have  \eqref{CONVERGENCE9}.
\proofend

It can be checked that the asymptotic variances of
 \ $\sqrt{n}\left(\widetilde\theta_n - \lim_{k\to\infty}\widetilde\theta_k\right)$ \ and
  \ $\sqrt{n}\left(\htheta_n - \lim_{k\to\infty}\htheta_k\right)$ \ are not equal.

\subsection{Two not neighbouring outliers, estimation of the mean of the offspring distribution
            and the outliers' sizes}

In this section we assume that \ $I=2$ \ and that the relevant time points
 \ $s_1$, $s_2\in\NN$ \ are known.
We concentrate on the CLS estimation of \ $\alpha$, \ $\theta_1$ \ and \ $\theta_2$.
\ Since \ $Y_k=X_k+\delta_{k,s_1}\theta_1+\delta_{k,s_2}\theta_2$,
 \ $k\in\ZZ_+$, \ we get for all \ $s_1,s_2\in\NN$,
 \begin{align}
   \begin{split}\label{SEGED94}
   \EE(Y_k\mid\cF^Y_{k-1})
    & =\alpha X_{k-1}+\mu_\vare + \delta_{k,s_1}\theta_1+\delta_{k,s_2}\theta_2 \\
    & =\alpha(Y_{k-1}-\delta_{k-1,s_1}\theta_1-\delta_{k-1,s_2}\theta_2)
       + \mu_\vare + \delta_{k,s_1}\theta_1 + \delta_{k,s_2}\theta_2  \\
    & = \alpha Y_{k-1} + \mu_\vare + (-\alpha\delta_{k-1,s_1}+\delta_{k,s_1})\theta_1
       + (-\alpha\delta_{k-1,s_2}+\delta_{k,s_2})\theta_2,
    \qquad k\in\NN.
   \end{split}
 \end{align}

In the sequel we also suppose that \ $s_1<s_2-1$, \ i.e., the time points \ $s_1$ \ and
 \ $s_2$ \ are not neighbouring.
Then, by \eqref{SEGED94},
 \begin{align*}
   \EE(Y_k\mid\cF^Y_{k-1})
     = \begin{cases}
        \alpha Y_{k-1}+\mu_\vare & \text{\quad if \ $1\leq k\leq s_1-1$,}\\
        \alpha Y_{k-1}+\mu_\vare + \theta_1 & \text{\quad if \ $k=s_1$,}\\
        \alpha Y_{k-1}+\mu_\vare - \alpha\theta_1 & \text{\quad if \ $k=s_1+1$,}\\
        \alpha Y_{k-1}+\mu_\vare & \text{\quad if \ $s_1+2 \leq k\leq s_2-1$,}\\
        \alpha Y_{k-1}+\mu_\vare + \theta_2 & \text{\quad if \ $k = s_2$,}\\
        \alpha Y_{k-1}+\mu_\vare - \alpha\theta_2 & \text{\quad if \ $k=s_2+1$,}\\
        \alpha Y_{k-1}+\mu_\vare & \text{\quad if \ $k\geq s_2+2$.}
       \end{cases}
 \end{align*}
Hence for all \ $ n\geq s_2+1$, \ $n\in\NN$,
 \begin{align}\label{SEGED53}
   \begin{split}
    \sum_{k=1}^n\big(Y_k-\EE(Y_k\mid \cF^Y_{k-1})\big)^2
         & =\DS\sum_{\substack{k=1 \\ k\not\in \{s_1,s_1+1,s_2,s_2+1\}}}^n
            \big(Y_k-\alpha Y_{k-1}-\mu_\vare\big)^2\\
          &\phantom{=\;}
           + \big(Y_{s_1}-\alpha Y_{s_1-1}-\mu_\vare-\theta_1\big)^2
           + \big(Y_{s_1+1}-\alpha Y_{s_1}-\mu_\vare+\alpha\theta_1\big)^2  \\
          &\phantom{=\;}
           + \big(Y_{s_2}-\alpha Y_{s_2-1}-\mu_\vare-\theta_2\big)^2
           + \big(Y_{s_2+1}-\alpha Y_{s_2}-\mu_\vare+\alpha\theta_2\big)^2.
   \end{split}
 \end{align}
For all \ $n\geq s_2+1$, \ $n\in\NN$, \ we define the function \ $Q_n^{\dag}:\RR^{n+1}\times\RR^3\to\RR$, \ as
 \begin{align*}
    &Q_n^{\dag}({\bf y}_n;\alpha',\theta_1',\theta_2')\\
      & :=\DS\sum_{\substack{k=1 \\ k\not\in \{s_1,s_1+1,s_2,s_2+1\}}}^n
            \big(y_k-\alpha' y_{k-1}-\mu_\vare\big)^2
           + \big(y_{s_1}-\alpha' y_{s_1-1}-\mu_\vare-\theta_1'\big)^2 \\
           &\phantom{=\;}
           + \big(y_{s_1+1}-\alpha' y_{s_1}-\mu_\vare+\alpha'\theta_1'\big)^2
           + \big(y_{s_2}-\alpha' y_{s_2-1}-\mu_\vare-\theta_2'\big)^2
           + \big(y_{s_2+1}-\alpha' y_{s_2}-\mu_\vare+\alpha'\theta_2'\big)^2,
 \end{align*}
 for all \ ${\bf y}_n\in\RR^{n+1}$, $\alpha',\theta_1',\theta_2'\in\RR$.
\ By definition, for all \ $n\geq s_2+1$, \ a CLS estimator for
 the parameter \ $(\alpha,\theta_1,\theta_2)\in(0,1)\times\NN^2$ \ is
 a measurable function \ $(\talpha_n^{\,\dag},\ttheta_{1,n}^{\,\dag},\ttheta_{2,n}^{\,\dag}):S_n\to\RR^3$
 \ such that
 \begin{align*}
   Q_n^{\dag}({\bf y}_n;\,&\talpha_n^{\,\dag}({\bf y}_n),\ttheta_{1,n}^{\,\dag}({\bf y}_n),
             \ttheta_{2,n}^{\,\dag}({\bf y}_n))
       = \inf_{(\alpha',\theta_1',\theta_2')\in\RR^3}Q_n^{\dag}({\bf y}_n;\alpha',\theta_1',\theta_2')
       \qquad \forall\;\;  {\bf y}_n\in S_n,
 \end{align*}
 where \ $S_n$ \ is suitable subset of \ $\RR^{n+1}$ \ (defined in the proof of
 Lemma \ref{LEMMA8}).
We note that we do not define the CLS estimator
 \ $(\talpha_n^{\,\dag},\ttheta_{1,n}^{\,\dag},\ttheta_{2,n}^{\,\dag})$ \
 for all samples \ ${\bf y}_n\in \RR^{n+1}$.
\ For all \ ${\bf y}_n\in\RR^{n+1}$ \ and \ $(\alpha',\theta_1',\theta_2')\in\RR^3$,
 \begin{align*}
    \frac{\partial Q_n^{\dag}}{\partial \alpha'}&({\bf y}_n;\alpha',\theta_1',\theta_2')\\
        &= \DS\sum_{\substack{k=1 \\ k\not\in \{s_1,s_1+1,s_2,s_2+1\}}}^n
            \big(y_k-\alpha' y_{k-1}-\mu_\vare\big)(-2y_{k-1})
           -2\big(y_{s_1}-\alpha' y_{s_1-1}-\mu_\vare-\theta_1'\big)y_{s_1-1} \\
        &\;\;\; + 2\big(y_{s_1+1}-\alpha' y_{s_1}-\mu_\vare+\alpha'\theta_1'\big)(-y_{s_1}+\theta_1')
                - 2\big(y_{s_2}-\alpha' y_{s_2-1}-\mu_\vare-\theta_2'\big)y_{s_2-1} \\
        &\;\;\; + 2\big(y_{s_2+1}-\alpha' y_{s_2}-\mu_\vare+\alpha'\theta_2'\big)(-y_{s_2}+\theta_2'),
 \end{align*}
 and
 \begin{align*}
   & \frac{\partial Q_n^{\dag}}{\partial \theta_1'}({\bf y}_n;\alpha',\theta_1',\theta_2')
       = -2(y_{s_1}-\alpha' y_{s_1-1}-\mu_\vare-\theta_1')
         +2\alpha'(y_{s_1+1}-\alpha' y_{s_1}-\mu_\vare+\alpha'\theta_1'),\\
   & \frac{\partial Q_n^{\dag}}{\partial \theta_2'}({\bf y}_n;\alpha',\theta_1',\theta_2')
       = -2(y_{s_2}-\alpha' y_{s_2-1}-\mu_\vare-\theta_2')
         +2\alpha'(y_{s_2+1}-\alpha' y_{s_2}-\mu_\vare+\alpha'\theta_2').
 \end{align*}

The next lemma is about the existence and uniqueness of the CLS estimator of \ $(\alpha,\theta_1,\theta_2)$.

\begin{Lem}\label{LEMMA8}
There exist subsets \ $S_n\subset\RR^{n+1}$, $n\geq s_2+1$ \ with the following properties:
 \begin{enumerate}
  \item[\upshape{(i)}]
   there exists a unique CLS estimator
   \ $(\talpha_n^{\,\dag},\ttheta_{1,n}^{\,\dag},\ttheta_{2,n}^{\,\dag}):S_n\to\RR^3$,
  \item[\upshape{(ii)}]
   for all \ ${\bf y}_n\in S_n$,
   \ $(\talpha_n^{\,\dag}({\bf y}_n),\ttheta_{1,n}^{\,\dag}({\bf y}_n),\ttheta_{2,n}^{\,\dag}({\bf y}_n))$
   \ is the unique solution of the system of equations
   \begin{align}\label{Additive_CLSE_EQ5}
   \begin{split}
    \frac{\partial Q_n^{\dag}}{\partial \alpha'}
     ({\bf y}_n;\alpha',\theta_1',\theta_2')=0,\\
    \frac{\partial Q_n^{\dag}}{\partial \theta_1'}
     ({\bf y}_n;\alpha',\theta_1',\theta_2')=0,\\
    \frac{\partial Q_n^{\dag}}{\partial \theta_2'}
      ({\bf y}_n;\alpha',\theta_1',\theta_2')=0,
   \end{split}
 \end{align}
  \item [\upshape{(iii)}]
  ${\bf Y}_n\in S_n$ \ holds asymptotically as \ $n\to\infty$ \ with probability one.
 \end{enumerate}
\end{Lem}

\noindent{\bf Proof.}
For any fixed \ ${\bf y}_n \in \RR^{n+1}$ \ and \ $\alpha' \in \RR$, \ the quadratic function
 \ $\RR^2 \ni (\theta_1',\theta_2') \mapsto Q_n^{\dag}({\bf y}_n;\alpha',\theta_1',\theta_2')$
 \ can be written in the form
 \begin{align*}
  &Q_n^{\dag}({\bf y}_n;\alpha',\theta_1',\theta_2') \\[2mm]
  &=  \!\!\left(\! \begin{bmatrix}
             \theta_1' \smallskip \\
             \theta_2'
            \end{bmatrix}
           - A_n(\alpha')^{-1} t_n({\bf y}_n;\alpha') \!\! \right)^{\hspace*{-1mm}\top}
      \!\!\! A_n(\alpha') \! \left( \! \begin{bmatrix}
                          \theta_1' \smallskip \\
                          \theta_2'
                         \end{bmatrix}
                         - A_n(\alpha')^{-1} t_n({\bf y}_n;\alpha') \! \right)
    \!\! + \widetilde{Q}_n^{\dag}({\bf y}_n;\alpha'),
 \end{align*}
 where
 \begin{align*}
  t_n({\bf y}_n;\alpha')
  &:= \begin{bmatrix}
       ( 1 + (\alpha')^2 ) y_{s_1} - \alpha' ( y_{s_1-1} + y_{s_1+1} )  - (1-\alpha')\mu_\vare \smallskip \\
       ( 1 + (\alpha')^2 ) y_{s_2} - \alpha' ( y_{s_2-1} + y_{s_2+1}) - (1-\alpha')\mu_\vare
      \end{bmatrix} , \\[2mm]
  \widetilde{Q}_n^{\dag}({\bf y}_n;\alpha')
  &:= \sum_{k=1}^n \big( y_k-\alpha' y_{k-1} \big)^2
      - t_n({\bf y}_n;\alpha')^\top A_n(\alpha')^{-1} t_n({\bf y}_n;\alpha') ,\\[2mm]
 A_n(\alpha')
   & := \begin{bmatrix}
       1+(\alpha')^2 & 0 \smallskip \\
                    0 &  1+(\alpha')^2 \\
     \end{bmatrix}.
 \end{align*}
Then \ $\widetilde{Q}_n^{\dag}({\bf y}_n;\alpha') = R_n({\bf y}_n;\alpha') / D_n(\alpha')$, \ where
 \ $D_n(\alpha'):=(1+(\alpha')^2)^2$ \ and
 \ $\RR\ni\alpha' \mapsto R_n({\bf y}_n;\alpha')$ \ is a polynomial of order 6 with leading coefficient
 \begin{align*}
   c_n({\bf y}_n) &:= \sum_{k=1}^n y_{k-1}^2  - (y_{s_1}^2 + y_{s_2}^2).
 \end{align*}
Let
 \[
   \widetilde{S}^{\,\dag}_n := \left\{{\bf y}_n\in\RR^{n+1} : c_n({\bf y}_n) > 0 \right\}.
 \]
For \ ${\bf y}_n \in \widetilde{S}^{\,\dag}_n$, \ we have
 \ $\lim_{|\alpha'|\to\infty} \widetilde{Q}_n^{\,\dag}({\bf y}_n;\alpha') = \infty$ \ and
 the continuous function
 \ $\RR \ni \alpha' \mapsto \widetilde{Q}_n^{\,\dag}({\bf y}_n;\alpha')$ \ attains its infimum.
Consequently, for all \ $n\geq s_2+1$ \ there exists a CLS estimator
 \ $(\talpha_n^{\,\dag},
    \ttheta_{1,n}^{\,\dag},
    \ttheta_{2,n}^{\,\dag}):\widetilde{S}^{\,\dag}_n\to \RR^3$,
\ where
 \begin{align}\nonumber
   \widetilde{Q}_n^{\,\dag}({\bf y}_n;\talpha_n^{\,\dag}({\bf y}_n))
   &=\inf_{\alpha'\in\RR} \widetilde{Q}_n^{\,\dag}({\bf y}_n;\alpha')
   \qquad \forall\;{\bf y}_n\in \widetilde{S}^{\,\dag}_n ,\\[2mm] \label{SEGED_UJ6}
  \begin{bmatrix}
   \ttheta_{1,n}^{\,\dag}({\bf y}_n) \smallskip \\
   \ttheta_{2,n}^{\,\dag}({\bf y}_n)
  \end{bmatrix}
  &= A_n(\talpha_n^{\,\dag}({\bf y}_n))^{-1} t_n({\bf y}_n;\talpha_n^{\,\dag}({\bf y}_n)) ,
  \qquad {\bf y}_n\in \widetilde{S}^{\,\dag}_n,
 \end{align}
 and for all \ ${\bf y}_n\in \widetilde{S}^{\,\dag}_n$,
   \ $(\talpha_n^{\,\dag}({\bf y}_n),
       \ttheta_{1,n}^{\,\dag}({\bf y}_n),
       \ttheta_{2,n}^{\,\dag}({\bf y}_n))$
   \ is a  solution of the system of equations \eqref{Additive_CLSE_EQ5}.

By \eqref{Ergodic1} and \eqref{Ergodic2}, we get
 \ $\PP\left(\lim_{n\to\infty} n^{-1} c_n({\bf Y}_n) = \EE\widetilde X^2 \right)=1$, \ where \ $\widetilde X$
 \ denotes a random variable with the unique stationary distribution of the INAR(1) model in \eqref{INAR1}.
Hence \ ${\bf Y}_n\in \widetilde{S}^{\,\dag}_n$ \ holds asymptotically as \ $n\to\infty$ \ with probability one.

Now we turn to find sets \ $S_n \subset \widetilde{S}^{\,\dag}_n$, $n\geq s_2+1$ \ such that
 the system of equations \eqref{Additive_CLSE_EQ5} has a unique solution with respect to
 \ $(\alpha',\theta_1',\theta_2')$ \ for all \ ${\bf y}_n\in S_n$.
\ Let us introduce the \ $(3\times 3)$ \ Hessian matrix
 \[
    H_n({\bf y}_n;\alpha',\theta_1',\theta_2'):=
     \begin{bmatrix}
       \frac{\partial^2 Q_n^{\dag}}{\partial (\alpha')^2}
       & \frac{\partial^2 Q_n^{\dag}}{\partial \theta_1' \partial \alpha'}
       & \frac{\partial^2 Q_n^{\dag}}{\partial \theta_2' \partial \alpha'}\\
          \frac{\partial^2 Q_n^{\dag}}{\partial \alpha' \partial \theta_1'}
        & \frac{\partial^2 Q_n^{\dag}}{ \partial (\theta_1')^2}
        &  \frac{\partial^2 Q_n^{\dag}}{\partial \theta_2' \partial \theta_1'} \\
           \frac{\partial^2 Q_n^{\dag}}{\partial \alpha' \partial \theta_2'}
        & \frac{\partial^2 Q_n^{\dag}}{ \partial \theta_1' \partial \theta_2'}
        &  \frac{\partial^2 Q_n^{\dag}}{ \partial (\theta_2')^2 } \\
     \end{bmatrix}
       ({\bf y}_n;\alpha',\theta_1',\theta_2'),
 \]
 and let us denote by \ $\Delta_{i,n}({\bf y}_n;\alpha',\theta_1',\theta_2')$ \
 its \ $i$-th order leading principal minor, \ $i=1,2,3$.
\ Further, for all \ $n\in\NN$, \ let
 \[
   S_n:=\Big\{{\bf y}_n\in  \widetilde{S}^{\,\dag}_n : \Delta_{i,n}({\bf y}_n;\alpha',\theta_1',\theta_2')>0,
                                   \;\, i=1,2,3,\, \forall\;(\alpha',\theta_1',\theta_2')
                                   \in\RR^3 \Big\}.
 \]
By Berkovitz \cite[Theorem 3.3, Chapter III]{Ber}, the function
 \ $\RR^3 \ni (\alpha',\theta_1',\theta_2')
     \mapsto Q_n^{\dag}({\bf y}_n;\alpha',\theta_1',\theta_2')$
 \ is strictly convex for all \ ${\bf y}_n\in S_n$.
\ Since it was already proved that the system of equations \eqref{Additive_CLSE_EQ5} has a solution
 for all \ ${\bf y}_n\in \widetilde{S}^{\,\dag}_n$, \ we obtain that this solution is unique for
 all \ ${\bf y}_n\in S_n$.

Next we check that \ ${\bf Y}_n\in S_n$ \ holds asymptotically as \ $n\to\infty$ \ with probability one.
For all \ $(\alpha',\theta_1',\theta_2')\in\RR^3$,
 \begin{align*}
   &\frac{\partial^2 Q_n^{\dag}}{\partial (\alpha')^2}  ({\bf Y}_n;\alpha',\theta_1',\theta_2')\\
   &\qquad = 2 \DS\sum_{\substack{k=1 \\ k\not\in \{s_1,s_1+1,s_2,s_2+1\}}}^n Y_{k-1}^2
            +2Y_{s_1-1}^2+2(Y_{s_1}-\theta_1')^2
            +2Y_{s_2-1}^2+2(Y_{s_2}-\theta_2')^2 \\[1mm]
   &\qquad  =2\DS\sum_{\substack{k=1 \\ k\not\in \{s_1+1,s_2+1\}}}^n X_{k-1}^2
            + 2(X_{s_1}+\theta_1 - \theta_1')^2
            + 2(X_{s_2}+\theta_2 - \theta_2')^2 ,\\[2mm]
  & \frac{\partial^2
    Q_n^{\dag}}{\partial\theta_1'\partial\alpha'}({\bf Y}_n;\alpha',\theta_1',\theta_2')
      = \frac{\partial^2 Q_n^{\dag}}
        {\partial\alpha'\partial\theta_1'}({\bf Y}_n;\alpha',\theta_1',\theta_2')
        = 2(Y_{s_1-1}+Y_{s_1+1}-2\alpha' Y_{s_1}-\mu_\vare+2\alpha'\theta_1')\\[1mm]
  &\phantom{\frac{\partial^2
      Q_n^{\dag}}{\partial\alpha'\partial\theta_1'}({\bf Y}_n;\alpha',\theta_1',\theta_2')}
      =2(X_{s_1-1}+X_{s_1+1}-2\alpha' X_{s_1}-\mu_\vare-2\alpha'(\theta_1-\theta_1')), \\[2mm]
  &\frac{\partial^2 Q_n^{\dag}}{\partial\theta_2'\partial\alpha'}({\bf Y}_n;\alpha',\theta_1',\theta_2')
    = \frac{\partial^2 Q_n^{\dag}}
     {\partial\alpha'\partial\theta_2'}({\bf Y}_n;\alpha',\theta_1',\theta_2')
    = 2(Y_{s_2-1}+Y_{s_2+1}-2\alpha' Y_{s_2}-\mu_\vare+2\alpha'\theta_2')\\[2mm]
  &\phantom{\frac{\partial^2 Q_n^{\dag}}{\partial\alpha'\partial\theta_2'}
      ({\bf Y}_n;\alpha',\theta_1',\theta_2')}
    =2(X_{s_2-1}+X_{s_2+1}-2\alpha' X_{s_2}-\mu_\vare-2\alpha'(\theta_2-\theta_2')),
 \end{align*}
 and
 \begin{align*}
  &\frac{\partial^2 Q_n^{\dag}}{\partial(\theta_1')^2}({\bf Y}_n;\alpha',\theta_1',\theta_2')
    =\frac{\partial^2 Q_n^{\dag}}{\partial(\theta_2')^2}({\bf Y}_n;\alpha',\theta_1',\theta_2')
    = 2((\alpha')^2+1),\\[2mm]
  &\frac{\partial^2 Q_n^{\dag}}{\partial\theta_1'\partial\theta_2'}({\bf Y}_n;\alpha',\theta_1',\theta_2')
    =\frac{\partial^2 Q_n^{\dag}}{\partial\theta_2'\partial\theta_1'}({\bf Y}_n;\alpha',\theta_1',\theta_2')
     = 0.
 \end{align*}
Then \ $H_n({\bf Y}_n;\alpha',\theta_1',\theta_2')$ \ has the following leading principal minors
 \begin{align*}
    &\Delta_{1,n}({\bf Y}_n;\alpha',\theta_1',\theta_2')
       =\DS\sum_{\substack{k=1 \\ k\not\in \{s_1+1,s_2+1\}}}^n \!\!2X_{k-1}^2
            + 2(X_{s_1}+\theta_1 - \theta_1')^2
            + 2(X_{s_2}+\theta_2 - \theta_2')^2, \\
    &\Delta_{2,n}({\bf Y}_n;\alpha',\theta_1',\theta_2')
      = 4\Bigg( ((\alpha')^2+1)\left(\DS\sum_{\substack{k=1 \\ k\not\in \{s_1+1,s_2+1\}}}^n X_{k-1}^2
                                   + (X_{s_1}+\theta_1 - \theta_1')^2
                                   + (X_{s_2}+\theta_2 - \theta_2')^2\right) \\
    & \phantom{\Delta_{2,n}({\bf Y}_n;\alpha',\theta_1',\theta_2')=\Bigg(\;\,}
       -\big(X_{s_1-1}+X_{s_1+1}-2\alpha' X_{s_1}-\mu_\vare-2\alpha'(\theta_1-\theta_1')\big)^2
       \Bigg),
 \end{align*}
 and
 \begin{align*}
    &\Delta_{3,n}({\bf Y}_n;\alpha',\theta_1',\theta_2')
        = \det H_n({\bf Y}_n;\alpha',\theta_1',\theta_2')\\
    &   = 8\Bigg(
          ((\alpha')^2+1)^2\left(\DS\sum_{\substack{k=1 \\ k\not\in \{s_1+1,s_2+1\}}}^n X_{k-1}^2
                                   + (X_{s_1}+\theta_1 - \theta_1')^2
                                   + (X_{s_2}+\theta_2 - \theta_2')^2\right) \\
     &\phantom{\det H_n({\bf Y}_n;\alpha',\theta_1',\theta_2') = \;}
          - ((\alpha')^2+1)\big(X_{s_1-1}+X_{s_1+1}-2\alpha' X_{s_1}-\mu_\vare-2\alpha'(\theta_1-\theta_1')\big)^2 \\
     &\phantom{\det H_n({\bf Y}_n;\alpha',\theta_1',\theta_2') = \;}
          - ((\alpha')^2+1)\big(X_{s_2-1}+X_{s_2+1}-2\alpha' X_{s_2}-\mu_\vare-2\alpha'(\theta_2-\theta_2')\big)^2
         \Bigg).
 \end{align*}
By \eqref{Ergodic2},
 \begin{align*}
    &\PP\left(\lim_{n\to\infty}\frac{1}{n}\Delta_{1,n}({\bf Y}_n;\alpha',\theta_1',\theta_2')
               = 2 \EE\widetilde X^2,
               \;\; \forall\;\; (\alpha',\theta_1',\theta_2')\in\RR^3 \right)=1, \\
    & \PP\left(\lim_{n\to\infty}\frac{1}{n}\Delta_{2,n}({\bf Y}_n;\alpha',\theta_1',\theta_2')
               = 4((\alpha')^2+1)\EE\widetilde X^2,
              \;\; \forall\;\; (\alpha',\theta_1',\theta_2')\in\RR^3   \right)=1,
 \end{align*}
 and
 \begin{align*}
    \PP\left(\lim_{n\to\infty}\frac{1}{n}\Delta_{3,n}({\bf Y}_n;\alpha',\theta_1',\theta_2')
               = 8((\alpha')^2+1)^2\EE\widetilde X^2,
               \;\; \forall\;\; (\alpha',\theta_1',\theta_2')\in\RR^3  \right)=1,
 \end{align*}
 where \ $\widetilde X$ \ denotes a random variable with the unique stationary
 distribution of the INAR(1) model in \eqref{INAR1}.
Hence
 \begin{align*}
    \PP\big(\lim_{n\to\infty}\Delta_{i,n}({\bf Y}_n;\alpha',\theta_1',\theta_2')=\infty,
            \;\;\forall\;\; (\alpha',\theta_1',\theta_2')\in\RR^3 \big)=1,
            \qquad i=1,2,3,
 \end{align*}
which yields that \ ${\bf Y}_n\in S_n$ \ asymptotically as \ $n\to\infty$
 \ with probability one, since we have already proved that \ ${\bf Y}_n\in \widetilde{S}^{\,\dag}_n$
 \ asymptotically as \ $n\to\infty$ \ with probability one.
\proofend

By Lemma \ref{LEMMA8},
 \ $(\talpha_n^{\,\dag}({\bf Y}_n),
    \ttheta_{1,n}^{\,\dag}({\bf Y}_n),
    \ttheta_{2,n}^{\,\dag}({\bf Y}_n))$
 \ exists uniquely asymptotically as \ $n\to\infty$ \ with probability one.
In the sequel we will simply denote it by
 \ $(\talpha_n^{\,\dag},\ttheta_{1,n}^{\,\dag},\ttheta_{2,n}^{\,\dag})$.

An easy calculation shows that
 \begin{align}\label{SEGED45}
   &\talpha_n^{\,\dag}
       = \frac{\sum_{k=1}^n(Y_k-\mu_\vare)Y_{k-1} - \ttheta_{1,n}^{\,\dag}(Y_{s_1-1}+Y_{s_1+1}-\mu_\vare)
               - \ttheta_{2,n}^{\,\dag}(Y_{s_2-1}+Y_{s_2+1}-\mu_\vare)}
              {\sum_{k=1}^nY_{k-1}^2 - 2\ttheta_{1,n}^{\,\dag} Y_{s_1} + (\ttheta_{1,n}^{\,\dag})^2
               - 2 \ttheta_{2,n}^{\,\dag} Y_{s_2} + (\ttheta_{2,n}^{\,\dag})^2},\\[2mm] \label{SEGED46}
  &\ttheta_{1,n}^{\,\dag} = Y_{s_1}-\frac{\talpha_n^{\,\dag}}{1+(\talpha_n^{\,\dag})^2}(Y_{s_1-1}+Y_{s_1+1})
                    - \frac{1-\talpha_n^{\,\dag}}{1+(\talpha_n^{\,\dag})^2}\mu_\vare,\\[2mm] \label{SEGED47}
  &\ttheta_{2,n}^{\,\dag} = Y_{s_2}-\frac{\talpha_n^{\,\dag}}{1+(\talpha_n^{\,\dag})^2}(Y_{s_2-1}+Y_{s_2+1})
                    - \frac{1-\talpha_n^{\,\dag}}{1+(\talpha_n^{\,\dag})^2}\mu_\vare,
 \end{align}
 hold asymptotically as \ $n\to\infty$ \ with probability one.

The next result shows that \ $\talpha_n^{\,\dag}$ \ is a strongly consistent estimator
 of \ $\alpha$, \ whereas \ $\ttheta_{1,n}^{\,\dag}$ \ and  \ $\ttheta_{2,n}^{\,\dag}$ \ fail to be
 strongly consistent estimators of \ $\theta_1$ \ and \ $\theta_2$, \ respectively.

\begin{Thm}\label{THEOREM2}
For the CLS estimators $(\talpha_n^{\,\dag},\ttheta_{1,n}^{\,\dag},\ttheta_{2,n}^{\,\dag})_{n\in\NN}$  of
 the parameter \ $(\alpha,\theta_1,\theta_2)\in(0,1)\times\NN^2$,
 \ the sequence \ $(\talpha_n^{\,\dag})_{n\in\NN}$ \ is strongly consistent
 for all \ $(\alpha,\theta_1,\theta_2)\in(0,1)\times\NN^2$, \ i.e.,
 \begin{align} \label{Strong_consistency10}
   \PP(\lim_{n\to\infty}\talpha_n^{\,\dag}=\alpha)=1,
     \qquad \forall\;(\alpha,\theta_1,\theta_2)\in(0,1)\times\NN^2,
 \end{align}
 whereas the sequences \ $(\ttheta_{1,n}^{\,\dag})_{n\in\NN}$ \ and \ $(\ttheta_{2,n}^{\,\dag})_{n\in\NN}$ \ are
 not strongly consistent for any
 \ $(\alpha,\theta_1,\theta_2)\in(0,1)\times\NN^2$, \ namely,
 \begin{align}\label{Strong_consistency11}
   & \PP\left(\lim_{n\to\infty}\ttheta_{1,n}^{\,\dag}
         = Y_{s_1}-\frac{\alpha}{1+\alpha^2}(Y_{s_1-1}+Y_{s_1+1})
                    - \frac{1-\alpha}{1+\alpha^2}\mu_\vare \right)=1, \\[2mm]\label{Strong_consistency14}
   & \PP\left(\lim_{n\to\infty}\ttheta_{2,n}^{\,\dag}
         = Y_{s_2}-\frac{\alpha}{1+\alpha^2}(Y_{s_2-1}+Y_{s_2+1})
                    - \frac{1-\alpha}{1+\alpha^2}\mu_\vare  \right)=1,
 \end{align}
 for all \ $(\alpha,\theta_1,\theta_2)\in(0,1)\times\NN^2$.
\end{Thm}

\noindent{\bf Proof.}
Similarly to \eqref{THETA_BOUND}, we obtain
 \begin{align}\label{SEGED43}
  &\vert \ttheta_{1,n}^{\,\dag}-\theta_1 \vert \leq X_{s_1}+\frac{1}{2}(X_{s_1-1}+X_{s_1+1})+\frac{3}{2}\mu_\vare,\\
         \label{SEGED44}
  &\vert \ttheta_{2,n}^{\,\dag}-\theta_2 \vert \leq X_{s_2}+\frac{1}{2}(X_{s_2-1}+X_{s_2+1})+\frac{3}{2}\mu_\vare,
 \end{align}
 which yield that the sequences \ $(\ttheta_{1,n}^{\,\dag}-\theta_1)_{n\in\NN}$
 \ and \ $(\ttheta_{2,n}^{\,\dag}-\theta_2)_{n\in\NN}$ \ are bounded with probability one.
\ Using \eqref{SEGED45}, \eqref{SEGED43} and \eqref{SEGED44}, by the same arguments as in the proof
 of Theorem \ref{THEOREM1}, one can derive \eqref{Strong_consistency10}.
Then \eqref{Strong_consistency10}, \eqref{SEGED46} and \eqref{SEGED47} yield \eqref{Strong_consistency11}
 and \eqref{Strong_consistency14}.
\proofend

The asymptotic distribution of the CLS estimation is given in the next theorem.

\begin{Thm}
Under the additional assumptions \ $\EE X_0^3<\infty$ \ and \ $\EE\vare_1^3<\infty$, \ we have
 \begin{align}\label{CONVERGENCE10}
   \sqrt{n}(\talpha_n^{\,\dag}-\alpha)\distr \cN(0,\sigma_{\alpha,\,\vare}^2)
      \qquad \text{as \ $n\to\infty$,}
  \end{align}
 where \ $\sigma_{\alpha,\,\vare}^2$ \ is defined in \eqref{SEGED_SZIGMA_ALPHA}.
Moreover, conditionally on the values \ $Y_{s_1-1}$, \ $Y_{s_2-1}$ \ and \ $Y_{s_1+1}$,
 \ $Y_{s_2+1}$,
 \begin{align}\label{CONVERGENCE11}
    \begin{bmatrix}
      \sqrt{n}\big(\ttheta_{1,n}^{\,\dag} - \lim_{k\to\infty}\ttheta_{1,k}^{\,\dag}\big) \\
      \sqrt{n}\big(\ttheta_{2,n}^{\,\dag} - \lim_{k\to\infty}\ttheta_{2,k}^{\,\dag}\big) \\
    \end{bmatrix}
      \distr \cN\left(\begin{bmatrix}
                        0 \\
                        0 \\
                      \end{bmatrix},
                      e_{\alpha,\vare} \sigma_{\alpha,\,\vare}^2 e_{\alpha,\vare}^\top
      \right)
      \qquad \text{as \ $n\to\infty$,}
 \end{align}
 where
 \[
   e_{\alpha,\vare}
      := \frac{1}{(1+\alpha^2)^2}
         \begin{bmatrix}
           (\alpha^2-1)(Y_{s_1-1}+Y_{s_1+1}) + (1+2\alpha-\alpha^2)\mu_\vare \\
           (\alpha^2-1)(Y_{s_2-1}+Y_{s_2+1}) + (1+2\alpha-\alpha^2)\mu_\vare \\
         \end{bmatrix}.
 \]
\end{Thm}

\noindent{\bf Proof.}
Using \eqref{SEGED45}, \eqref{SEGED43} and \eqref{SEGED44}, by the very same
 arguments as in the proof of \eqref{CONVERGENCE6}, one can obtain \eqref{CONVERGENCE10}.
Now we turn to prove \eqref{CONVERGENCE11}.
Using the notation
 \[
   B_n^{\dag}:=
          \begin{bmatrix}
             1+(\talpha_n^{\,\dag})^2 & 0 \\
             0 & 1+(\talpha_n^{\,\dag})^2 \\
           \end{bmatrix},
 \]
 by \eqref{SEGED46} and \eqref{SEGED47}, we have
  \begin{align*}
    \begin{bmatrix}
      \ttheta_{1,n}^{\,\dag} \\
      \ttheta_{2,n}^{\,\dag} \\
    \end{bmatrix}
     = (B_n^{\dag})^{-1}
        \begin{bmatrix}
         (1+(\talpha_n^{\,\dag})^2)Y_{s_1} - \talpha_n^{\,\dag}(Y_{s_1-1}+Y_{s_1+1}) - (1-\talpha_n^{\,\dag})\mu_\vare\\
         (1+(\talpha_n^{\,\dag})^2)Y_{s_2} - \talpha_n^{\,\dag}(Y_{s_2-1}+Y_{s_2+1}) - (1-\talpha_n^{\,\dag})\mu_\vare \\
        \end{bmatrix}
  \end{align*}
 holds asymptotically as \ $n\to\infty$ \ with probability one.
Theorem \ref{THEOREM2} yields that
 \[
     \PP\left(\lim_{n\to\infty}B_n^{\dag}
                 = \begin{bmatrix}
                    1+\alpha^2 & 0 \\
                    0 & 1+\alpha^2 \\
                   \end{bmatrix}=:B^{\dag}
         \right)=1.
 \]
By \eqref{Strong_consistency11} and \eqref{Strong_consistency14}, we have
 \begin{align*}
   &\begin{bmatrix}
      \sqrt{n}\big(\ttheta_{1,n}^{\,\dag} - \lim_{k\to\infty}\ttheta_{1,k}^{\,\dag}\big) \\
      \sqrt{n}\big(\ttheta_{2,n}^{\,\dag} - \lim_{k\to\infty}\ttheta_{2,k}^{\,\dag}\big) \\
    \end{bmatrix} \\
   &\qquad
    =\sqrt{n}(B_n^{\dag})^{-1}
    \left(
      \begin{bmatrix}
         (1+(\talpha_n^{\,\dag})^2)Y_{s_1} - \talpha_n^{\,\dag}(Y_{s_1-1}+Y_{s_1+1}) - (1-\talpha_n^{\,\dag})\mu_\vare\\
         (1+(\talpha_n^{\,\dag})^2)Y_{s_2} - \talpha_n^{\,\dag}(Y_{s_2-1}+Y_{s_2+1}) - (1-\talpha_n^{\,\dag})\mu_\vare \\
      \end{bmatrix}\right.\\
   &\phantom{\qquad = \sqrt{n}(B_n^{\dag})^{-1}\Big(\;\;}\left.
     - B_n^{\dag}(B^{\dag})^{-1}
      \begin{bmatrix}
         (1+\alpha^2)Y_{s_1} - \alpha(Y_{s_1-1}+Y_{s_1+1}) - (1-\alpha)\mu_\vare\\
         (1+\alpha^2)Y_{s_2} - \alpha(Y_{s_2-1}+Y_{s_2+1}) - (1-\alpha)\mu_\vare \\
      \end{bmatrix}
     \right)\\
   &\qquad=\sqrt{n}(B_n^{\dag})^{-1}
      \left(
       \begin{bmatrix}
         (1+(\talpha_n^{\,\dag})^2)Y_{s_1} - \talpha_n^{\,\dag}(Y_{s_1-1}+Y_{s_1+1}) - (1-\talpha_n^{\,\dag})\mu_\vare\\
         (1+(\talpha_n^{\,\dag})^2)Y_{s_2} - \talpha_n^{\,\dag}(Y_{s_2-1}+Y_{s_2+1}) - (1-\talpha_n^{\,\dag})\mu_\vare \\
       \end{bmatrix}\right.\\
  &\phantom{\qquad =\sqrt{n}(B_n^{\dag})^{-1}\Big(\;}\left.
    -\begin{bmatrix}
         (1+\alpha^2)Y_{s_1} - \alpha(Y_{s_1-1}+Y_{s_1+1}) - (1-\alpha)\mu_\vare\\
         (1+\alpha^2)Y_{s_2} - \alpha(Y_{s_2-1}+Y_{s_2+1}) - (1-\alpha)\mu_\vare \\
      \end{bmatrix}
       \right)\\
 &\phantom{\qquad=\;}
      +\sqrt{n}
      \left((B_n^{\dag})^{-1} -  (B^{\dag})^{-1}
       \right)
       \begin{bmatrix}
         (1+\alpha^2)Y_{s_1} - \alpha(Y_{s_1-1}+Y_{s_1+1}) - (1-\alpha)\mu_\vare\\
         (1+\alpha^2)Y_{s_2} - \alpha(Y_{s_2-1}+Y_{s_2+1}) - (1-\alpha)\mu_\vare \\
      \end{bmatrix}, \\
 \end{align*}
  and hence
 \begin{align*}
    &\begin{bmatrix}
      \sqrt{n}\big(\ttheta_{1,n}^{\,\dag} - \lim_{k\to\infty}\ttheta_{1,k}^{\,\dag}\big) \\
      \sqrt{n}\big(\ttheta_{2,n}^{\,\dag} - \lim_{k\to\infty}\ttheta_{2,k}^{\,\dag}\big) \\
   \end{bmatrix}\\
  &=\sqrt{n}(B_n^{\dag})^{-1}
    \begin{bmatrix}
      (\talpha_n^{\,\dag}-\alpha)\big((\talpha_n^{\,\dag}+\alpha)Y_{s_1} - Y_{s_1-1}- Y_{s_1+1} + \mu_\vare\big) \\
      (\talpha_n^{\,\dag}-\alpha)\big((\talpha_n^{\,\dag}+\alpha)Y_{s_2} - Y_{s_2-1}- Y_{s_2+1} + \mu_\vare\big) \\
    \end{bmatrix}\\
 &\phantom{=\;}
    +\sqrt{n}(B_n^{\dag})^{-1}
       \big( B^{\dag}- B_n^{\dag} \big) (B^{\dag})^{-1}
        \begin{bmatrix}
         (1+\alpha^2)Y_{s_1} - \alpha(Y_{s_1-1}+Y_{s_1+1}) - (1-\alpha)\mu_\vare\\
         (1+\alpha^2)Y_{s_2} - \alpha(Y_{s_2-1}+Y_{s_2+1}) - (1-\alpha)\mu_\vare \\
         \end{bmatrix}.
 \end{align*}
Then
 \begin{align}\label{SEGED48}
   &\begin{bmatrix}
      \sqrt{n}\big(\ttheta_{1,n}^{\,\dag} - \lim_{k\to\infty}\ttheta_{1,k}^{\,\dag}\big) \\
      \sqrt{n}\big(\ttheta_{2,n}^{\,\dag} - \lim_{k\to\infty}\ttheta_{2,k}^{\,\dag}\big) \\
    \end{bmatrix}
   = \sqrt{n}(\talpha_n^{\,\dag}-\alpha)
       \begin{bmatrix}
         K_n^{\dag} \\
         L_n^{\dag} \\
       \end{bmatrix}
  \end{align}
 holds asymptotically as \ $n\to\infty$ \ with probability one, where
 \begin{align*}
    \begin{bmatrix}
         K_n^{\dag} \\
         L_n^{\dag} \\
    \end{bmatrix}
    &:=  (B_n^{\dag})^{-1}
       \begin{bmatrix}
         -(\talpha_n^{\,\dag}+\alpha) & 0 \\
         0 & -(\talpha_n^{\,\dag}+\alpha) \\
       \end{bmatrix}
        (B^{\dag})^{-1}
        \begin{bmatrix}
         (1+\alpha^2)Y_{s_1} - \alpha(Y_{s_1-1}+Y_{s_1+1}) - (1-\alpha)\mu_\vare\\
         (1+\alpha^2)Y_{s_2} - \alpha(Y_{s_2-1}+Y_{s_2+1}) - (1-\alpha)\mu_\vare \\
        \end{bmatrix}\\
    &\phantom{=\;\;\,}
      + (B_n^{\dag})^{-1}
     \begin{bmatrix}
         (\talpha_n^{\,\dag}+\alpha)Y_{s_1} - Y_{s_1-1}- Y_{s_1+1} + \mu_\vare \\
          (\talpha_n^{\,\dag}+\alpha)Y_{s_2} - Y_{s_2-1}- Y_{s_2+1} + \mu_\vare \\
       \end{bmatrix}.
  \end{align*}
By \eqref{Strong_consistency10}, we have
 \ $\begin{bmatrix}
     K_n^{\dag} & L_n^{\dag} \\
    \end{bmatrix}^\top$ \
 converges almost surely as \ $n\to\infty$ \ to
 \begin{align*}
   & (B^{\dag})^{-1}
     \begin{bmatrix}
       2\alpha Y_{s_1} - Y_{s_1-1} - Y_{s_1+1} +\mu_\vare \\
       2\alpha Y_{s_2} - Y_{s_2-1} - Y_{s_2+1} +\mu_\vare  \\
     \end{bmatrix} \\
   & +(B^{\dag})^{-1}
    \begin{bmatrix}
      -2\alpha & 0 \\
      0 & -2\alpha \\
    \end{bmatrix}
    (B^{\dag})^{-1}
      \begin{bmatrix}
        (1+\alpha^2)Y_{s_1} - \alpha(Y_{s_1-1}+Y_{s_1+1}) - (1-\alpha)\mu_\vare\\
        (1+\alpha^2)Y_{s_2} - \alpha(Y_{s_2-1}+Y_{s_2+1}) - (1-\alpha)\mu_\vare \\
      \end{bmatrix} \\
   & = \frac{1}{(1+\alpha^2)^2}
       \begin{bmatrix}
         (\alpha^2-1)(Y_{s_1-1}+Y_{s_1+1}) + (1+2\alpha -\alpha^2)\mu_\vare \\
         (\alpha^2-1)(Y_{s_2-1}+Y_{s_2+1}) + (1+2\alpha -\alpha^2)\mu_\vare \\
       \end{bmatrix}
      =e_{\alpha,\vare}.
 \end{align*}
By \eqref{SEGED48}, \eqref{CONVERGENCE10} and Slutsky's lemma, we have \eqref{CONVERGENCE11}.
\proofend

\subsection{Two neighbouring outliers, estimation of the mean of the offspring distribution
            and the outliers' sizes}

In this section we assume that \ $I=2$ \ and that the relevant time points
 \ $s_1$, $s_2\in\NN$ \ are known.
We also suppose that \ $s_1:=s$ \ and \ $s_2:=s+1$, \ i.e., the time points \ $s_1$ \ and
 \ $s_2$ \ are neighbouring.
We concentrate on the CLS estimation of \ $\alpha$, \ $\theta_1$ \ and \ $\theta_2$.
Then, by \eqref{SEGED94},
 \begin{align*}
   \EE(Y_k\mid\cF^Y_{k-1})
     = \begin{cases}
        \alpha Y_{k-1}+\mu_\vare & \text{\quad if \ $1\leq k\leq s_1-1=s-1$,}\\
        \alpha Y_{k-1}+\mu_\vare + \theta_1 & \text{\quad if \ $k=s_1=s$,}\\
        \alpha Y_{k-1}+\mu_\vare - \alpha\theta_1 + \theta_2 & \text{\quad if \ $k=s+1=s_1+1=s_2$,}\\
        \alpha Y_{k-1}+\mu_\vare - \alpha\theta_2 & \text{\quad if \ $k=s+2=s_1+2=s_2+1$,}\\
        \alpha Y_{k-1}+\mu_\vare & \text{\quad if \ $k\geq s+2=s_2+2$.}
       \end{cases}
 \end{align*}
Hence
 \begin{align}\label{SEGED95}
  \begin{split}
   &\sum_{k=1}^n\big(Y_k-\EE(Y_k\mid \cF^Y_{k-1})\big)^2\\
         & =\DS\sum_{\substack{k=1 \\ k\not\in \{s,s+1,s+2\}}}^n
            \big(Y_k-\alpha Y_{k-1}-\mu_\vare\big)^2
           + \big(Y_{s}-\alpha Y_{s-1}-\mu_\vare-\theta_1\big)^2
           + \big(Y_{s+1}-\alpha Y_{s}-\mu_\vare+\alpha\theta_1-\theta_2\big)^2 \\
         & \phantom{=\;}
            + \big(Y_{s+2}-\alpha Y_{s+1}-\mu_\vare+\alpha\theta_2\big)^2,
    \qquad n\geq s+2,\;\; n\in\NN.
  \end{split}
 \end{align}
For all \ $n\geq s+2$, \ $n\in\NN$, \ we define the function
  \ $Q_n^{\dag\dag}:\RR^{n+1}\times\RR^3\to\RR$, \ as
 \begin{align*}
    &Q_n^{\dag\dag}({\bf y}_n;\alpha',\theta_1',\theta_2')\\
     &:=\DS\sum_{\substack{k=1 \\ k\not\in \{s,s+1,s+2\}}}^n
            \big(y_k-\alpha' y_{k-1}-\mu_\vare\big)^2
           + \big(y_{s}-\alpha' y_{s-1}-\mu_\vare-\theta_1'\big)^2
           + \big(y_{s+1}-\alpha' y_{s}-\mu_\vare+\alpha'\theta_1'-\theta_2'\big)^2 \\
         & \phantom{=\;}
            + \big(y_{s+2}-\alpha' y_{s+1}-\mu_\vare+\alpha'\theta_2'\big)^2,
    \qquad {\bf y}_n\in\RR^{n+1},\;\alpha',\theta_1',\theta_2'\in\RR.
 \end{align*}
By definition, for all \ $n\geq s+2$, \ a CLS estimator for
 the parameter \ $(\alpha,\theta_1,\theta_2)\in(0,1)\times\NN^2$ \ is
 a measurable function
 \ $(\talpha_n^{\,\dag\dag},\ttheta_{1,n}^{\,\dag\dag},\ttheta_{2,n}^{\,\dag\dag}):S_n\to\RR^3$ \ such that
  \begin{align*}
   Q_n^{\dag\dag}({\bf y}_n;\,&\talpha_n^{\,\dag\dag}({\bf y}_n),
    \ttheta_{1,n}^{\,\dag\dag}({\bf y}_n),
             \ttheta_{2,n}^{\,\dag\dag}({\bf y}_n))
       = \inf_{(\alpha',\theta_1',\theta_2')\in\RR^3}Q_n^{\dag\dag}({\bf y}_n;\alpha',\theta_1',\theta_2')
       \qquad \forall\;\;  {\bf y}_n\in S_n,
 \end{align*}
 where \ $S_n$ \ is suitable subset of \ $\RR^{n+1}$ \ (defined in the proof of Lemma \ref{LEMMA9}).
We note that we do not define the CLS estimator
 \ $(\talpha_n^{\,\dag\dag},\ttheta_{1,n}^{\,\dag\dag},\ttheta_{2,n}^{\,\dag\dag})$ \
 for all samples \ ${\bf y}_n\in \RR^{n+1}$.
We have
 \begin{align*}
    &\frac{\partial Q_n^{\dag\dag}}{\partial \alpha'}({\bf y}_n;\alpha',\theta_1',\theta_2')\\
        &= \DS\sum_{\substack{k=1 \\ k\not\in \{s,s+1,s+2\}}}^n
            \big(y_k-\alpha' y_{k-1}-\mu_\vare\big)(-2y_{k-1})
           -2\big(y_s-\alpha' y_{s-1}-\mu_\vare-\theta_1'\big)y_{s-1} \\
        &\;\;\; + 2\big(y_{s+1}-\alpha' y_s-\mu_\vare+\alpha'\theta_1'-\theta_2'\big)(-y_s+\theta_1')
                + 2\big(y_{s+2}-\alpha' y_{s+1}-\mu_\vare+\alpha'\theta_2'\big)(-y_{s+1}+\theta_2'),
 \end{align*}
 and
 \begin{align*}
   & \frac{\partial Q_n^{\dag\dag}}{\partial \theta_1'}({\bf y}_n;\alpha',\theta_1',\theta_2')
       = -2(y_{s}-\alpha' y_{s-1}-\mu_\vare-\theta_1')
         +2\alpha'(y_{s+1}-\alpha' y_s-\mu_\vare+\alpha'\theta_1'-\theta_2'),\\
   & \frac{\partial Q_n^{\dag\dag}}{\partial \theta_2'}({\bf y}_n;\alpha',\theta_1',\theta_2')
       =-2(y_{s+1}-\alpha' y_s-\mu_\vare+\alpha'\theta_1'-\theta_2')
         +2\alpha'(y_{s+2}-\alpha' y_{s+1}-\mu_\vare+\alpha'\theta_2').
 \end{align*}

The next lemma is about the existence and uniqueness of the CLS estimator of \ $(\alpha,\theta_1,\theta_2)$.

\begin{Lem}\label{LEMMA9}
There exist subsets \ $S_n\subset\RR^{n+1}$, $n\geq s+2$ \ with the following properties:
 \begin{enumerate}
  \item[\upshape{(i)}]
   there exists a unique CLS estimator
   \ $(\talpha_n^{\,\dag\dag},\ttheta_{1,n}^{\,\dag\dag},\ttheta_{2,n}^{\,\dag\dag}):S_n\to\RR^3$,
  \item[\upshape{(ii)}]
   for all \ ${\bf y}_n\in S_n$,
   \ $(\talpha_n^{\,\dag\dag}({\bf y}_n),
       \ttheta_{1,n}^{\,\dag\dag}({\bf y}_n),\ttheta_{2,n}^{\,\dag\dag}({\bf y}_n))$
   \ is the unique solution of the system of equations
   \begin{align}\label{Additive_CLSE_EQ6}
    \begin{split}
    \frac{\partial Q_n^{\dag\dag}}{\partial \alpha'}({\bf y}_n;\alpha',\theta_1',\theta_2')=0,\\
    \frac{\partial Q_n^{\dag\dag}}{\partial \theta_1'}({\bf y}_n;\alpha',\theta_1',\theta_2')=0,\\
    \frac{\partial Q_n^{\dag\dag}}{\partial \theta_2'}({\bf y}_n;\alpha',\theta_1',\theta_2')=0,
   \end{split}
  \end{align}  \item [\upshape{(iii)}]
  ${\bf Y}_n\in S_n$ \ holds asymptotically as \ $n\to\infty$ \ with probability one.
 \end{enumerate}
\end{Lem}

\noindent{\bf Proof.}
For any fixed \ ${\bf y}_n \in \RR^{n+1}$ \ and \ $\alpha' \in \RR$, \ the quadratic function
 \ $\RR^2 \ni (\theta_1',\theta_2') \mapsto Q_n^{\dag\dag}({\bf y}_n;\alpha',\theta_1',\theta_2')$
 \ can be written in the form
 \begin{align*}
  &Q_n^{\dag\dag}({\bf y}_n;\alpha',\theta_1',\theta_2') \\[2mm]
  &=  \!\!\left(\! \begin{bmatrix}
             \theta_1' \smallskip \\
             \theta_2'
            \end{bmatrix}
           - A_n(\alpha')^{-1} t_n({\bf y}_n;\alpha') \!\! \right)^{\hspace*{-1mm}\top}
      \!\!\! A_n(\alpha') \! \left( \! \begin{bmatrix}
                          \theta_1' \smallskip \\
                          \theta_2'
                         \end{bmatrix}
                         - A_n(\alpha')^{-1} t_n({\bf y}_n;\alpha') \! \right)
    \!\! + \widetilde{Q}_n^{\dag\dag}({\bf y}_n;\alpha'),
 \end{align*}
 where
 \begin{align*}
  t_n({\bf y}_n;\alpha')
  &:= \begin{bmatrix}
       ( 1 + (\alpha')^2 ) y_s - \alpha' ( y_{s-1} + y_{s+1} )  - (1-\alpha')\mu_\vare \smallskip \\
       ( 1 + (\alpha')^2 ) y_{s+1} - \alpha' ( y_s + y_{s+2}) - (1-\alpha')\mu_\vare
      \end{bmatrix} , \\[2mm]
  \widetilde{Q}_n^{\dag\dag}({\bf y}_n;\alpha')
  &:= \sum_{k=1}^n \big( y_k-\alpha' y_{k-1} \big)^2
      - t_n({\bf y}_n;\alpha')^\top A_n(\alpha')^{-1} t_n({\bf y}_n;\alpha') ,\\[2mm]
 A_n(\alpha')
   & := \begin{bmatrix}
       1+(\alpha')^2 & -\alpha' \smallskip \\
                    -\alpha' &  1+(\alpha')^2 \\
     \end{bmatrix}.
 \end{align*}
Then \ $\widetilde{Q}_n^{\dag\dag}({\bf y}_n;\alpha') = R_n({\bf y}_n;\alpha') / D_n(\alpha')$, \ where
 \ $D_n(\alpha'):=(1+(\alpha')^2)^2 - (\alpha')^2 = (\alpha')^4 + (\alpha')^2 +1>0$ \ and
 \ $\RR\ni\alpha' \mapsto R_n({\bf y}_n;\alpha')$ \ is a polynomial of order 6 with leading coefficient
 \begin{align*}
   c_n({\bf y}_n) &:= \sum_{k=1}^n y_{k-1}^2  - (y_s^2 + y_{s+1}^2).
 \end{align*}
Let
 \[
   \widetilde{S}^{\dag\dag}_n := \left\{{\bf y}_n\in\RR^{n+1} : c_n({\bf y}_n) > 0 \right\}.
 \]
For \ ${\bf y}_n \in \widetilde{S}^{\dag\dag}_n$, \ we have
 \ $\lim_{|\alpha'|\to\infty} \widetilde{Q}_n^{\,\dag\dag}({\bf y}_n;\alpha') = \infty$ \ and
 the continuous function
 \ $\RR \ni \alpha' \mapsto \widetilde{Q}_n^{\,\dag\dag}({\bf y}_n;\alpha')$ \ attains its infimum.
Consequently, for all \ $n\geq s+2$ \ there exists a CLS estimator
 \ $(\talpha_n^{\,\dag\dag},
    \ttheta_{1,n}^{\,\dag\dag},
    \ttheta_{2,n}^{\,\dag\dag}) : \widetilde{S}^{\dag\dag}_n\to\RR^3$,
\ where
 \begin{align}\nonumber
   \widetilde{Q}_n^{\,\dag\dag}({\bf y}_n;\talpha_n^{\,\dag\dag}({\bf y}_n))
   &=\inf_{\alpha'\in\RR} \widetilde{Q}_n^{\,\dag\dag}({\bf y}_n;\alpha')
   \qquad \forall\;{\bf y}_n\in \widetilde{S}^{\dag\dag}_n ,\\[2mm] \label{SEGED_UJ8}
  \begin{bmatrix}
   \ttheta_{1,n}^{\,\dag\dag}({\bf y}_n) \smallskip \\
   \ttheta_{2,n}^{\,\dag\dag}({\bf y}_n)
  \end{bmatrix}
  &= A_n(\talpha_n^{\,\dag\dag}({\bf y}_n))^{-1} t_n({\bf y}_n;\talpha_n^{\,\dag\dag}({\bf y}_n)) ,
  \qquad {\bf y}_n\in \widetilde{S}^{\dag\dag}_n,
 \end{align}
 and for all \ ${\bf y}_n\in \widetilde{S}^{\dag\dag}_n$,
   \ $(\talpha_n^{\,\dag\dag}({\bf y}_n),
    \ttheta_{1,n}^{\,\dag\dag}({\bf y}_n),
    \ttheta_{2,n}^{\,\dag\dag}({\bf y}_n))$
   \ is a  solution of the system of equations \eqref{Additive_CLSE_EQ6}.

By \eqref{Ergodic1} and \eqref{Ergodic2}, we get
 \ $\PP\left(\lim_{n\to\infty} n^{-1} c_n({\bf Y}_n) = \EE\widetilde X^2 \right)=1$, \ where \ $\widetilde X$
 \ denotes a random variable with the unique stationary distribution of the INAR(1) model in \eqref{INAR1}.
Hence \ ${\bf Y}_n\in \widetilde{S}^{\dag\dag}_n$ \ holds asymptotically as
 \ $n\to\infty$ \ with probability one.

Now we turn to find sets \ $S_n \subset \widetilde{S}^{\dag\dag}_n$, $n\geq s+2$ \ such that
 the system of equations \eqref{Additive_CLSE_EQ6} has a unique solution with respect to
 \ $(\alpha',\theta_1',\theta_2')$ \ for all \ ${\bf y}_n\in S_n$.
\ Let us introduce the \ $(3\times 3)$ \ Hessian matrix
 \[
    H_n({\bf y}_n;\alpha',\theta_1',\theta_2'):=
     \begin{bmatrix}
       \frac{\partial^2 Q_n^{\dag\dag}}{\partial (\alpha')^2}
       & \frac{\partial^2 Q_n^{\dag\dag}}{\partial \theta_1' \partial \alpha'}
       & \frac{\partial^2 Q_n^{\dag\dag}}{\partial \theta_2' \partial \alpha'}\\
          \frac{\partial^2 Q_n^{\dag\dag}}{\partial \alpha' \partial \theta_1'}
        & \frac{\partial^2 Q_n^{\dag\dag}}{ \partial (\theta_1')^2}
        &  \frac{\partial^2 Q_n^{\dag\dag}}{\partial \theta_2' \partial \theta_1'} \\
           \frac{\partial^2 Q_n^{\dag\dag}}{\partial \alpha' \partial \theta_2'}
        & \frac{\partial^2 Q_n^{\dag\dag}}{ \partial \theta_1' \partial \theta_2'}
        &  \frac{\partial^2 Q_n^{\dag\dag}}{ \partial (\theta_2')^2 } \\
     \end{bmatrix}
       ({\bf y}_n;\alpha',\theta_1',\theta_2'),
 \]
 and let us denote by \ $\Delta_{i,n}({\bf y}_n;\alpha',\theta_1',\theta_2')$ \
 its \ $i$-th order leading principal minor, \ $i=1,2,3$.
\ Further, for all \ $n\geq s+2$, \ let
 \[
   S_n:=\Big\{{\bf y}_n\in\widetilde{S}^{\dag\dag} : \Delta_{i,n}({\bf y}_n;\alpha',\theta_1',\theta_2')>0,
                                   \;\, i=1,2,3,\, \forall\;(\alpha',\theta_1',\theta_2')
                                   \in\RR^3 \Big\}.
 \]
By Berkovitz \cite[Theorem 3.3, Chapter III]{Ber}, the function
 \ $\RR^3 \ni (\alpha',\theta_1',\theta_2')
     \mapsto Q_n^{\dag\dag}({\bf y}_n;\alpha',\theta_1',\theta_2')$
 \ is strictly convex for all \ ${\bf y}_n\in S_n$.
\ Since it was already proved that the system of equations \eqref{Additive_CLSE_EQ6} has a solution for all
 \ ${\bf y}_n\in \widetilde{S}^{\dag\dag}_n$, \ we obtain that this solution is unique for all
 \ ${\bf y}_n\in S_n$.

Next we check that \ ${\bf Y}_n\in S_n$ \ holds asymptotically as \ $n\to\infty$ \ with probability one.
For all \ $(\alpha',\theta_1',\theta_2')\in\RR^3$,
 \begin{align*}
   \frac{\partial^2 Q_n^{\dag\dag}}{\partial (\alpha')^2}({\bf Y}_n;\alpha',\theta_1',\theta_2')
   &= 2 \DS\sum_{\substack{k=1 \\ k\not\in \{s,s+1,s+2\}}}^n Y_{k-1}^2
            +2Y_{s-1}^2+2(Y_s-\theta_1')^2
            +2(Y_{s+1}-\theta_2')^2 \\
   & =2\DS\sum_{\substack{k=1 \\ k\not\in \{s+1,s+2\}}}^n X_{k-1}^2
            + 2(X_{s}+\theta_1 - \theta_1')^2
            + 2(X_{s+1}+\theta_2 - \theta_2')^2 ,
 \end{align*}
 and
 \begin{align*}
  &\frac{\partial^2 Q_n^{\dag\dag}}{\partial\theta_1'\partial\alpha'}({\bf Y}_n;\alpha',\theta_1',\theta_2')
    = \frac{\partial^2 Q_n^{\dag\dag}}
     {\partial\alpha'\partial\theta_1'}({\bf Y}_n;\alpha',\theta_1',\theta_2')\\
  &\phantom{\frac{\partial^2 Q_n^{\dag\dag}}
   {\partial\alpha'\partial\theta_1'}({\bf Y}_n;\alpha',\theta_1',\theta_2')}
    = 2\big(Y_{s-1}+Y_{s+1}-2\alpha' Y_{s}-\mu_\vare+2\alpha'\theta_1'-\theta_2'\big)\\
  &\phantom{\frac{\partial^2 Q_n^{\dag\dag}}
   {\partial\alpha'\partial\theta_1'}({\bf Y}_n;\alpha',\theta_1',\theta_2')}
    =2\big(X_{s-1}+X_{s+1}-2\alpha' X_{s}-\mu_\vare-2\alpha'(\theta_1-\theta_1')+(\theta_2-\theta_2')\big),\\
  &\frac{\partial^2 Q_n^{\dag\dag}}{\partial\theta_2'\partial\alpha'}({\bf Y}_n;\alpha',\theta_1',\theta_2')
    = \frac{\partial^2 Q_n^{\dag\dag}}
     {\partial\alpha'\partial\theta_2'}({\bf Y}_n;\alpha',\theta_1',\theta_2')\\
  &\phantom{\frac{\partial^2 Q_n^{\dag\dag}}
   {\partial\alpha'\partial\theta_2'}({\bf Y}_n;\alpha',\theta_1',\theta_2')}
    = 2\big(Y_{s}+Y_{s+2}-2\alpha' Y_{s+1}-\mu_\vare-\theta_1'+2\alpha'\theta_2'\big)\\
  &\phantom{\frac{\partial^2 Q_n^{\dag\dag}}
   {\partial\alpha'\partial\theta_2'}({\bf Y}_n;\alpha',\theta_1',\theta_2')}
    =2\big(X_{s}+X_{s+2}-2\alpha' X_{s+1}-\mu_\vare+(\theta_1-\theta_1')-2\alpha'(\theta_2-\theta_2')\big),\\
  &\frac{\partial^2 Q_n^{\dag\dag}}{\partial(\theta_1')^2}({\bf Y}_n;\alpha',\theta_1',\theta_2')
    =\frac{\partial^2 Q_n^{\dag\dag}}{\partial(\theta_2')^2}({\bf Y}_n;\alpha',\theta_1',\theta_2')
    = 2((\alpha')^2+1),\\
  &\frac{\partial^2 Q_n^{\dag\dag}}
   {\partial\theta_1'\partial\theta_2'}({\bf Y}_n;\alpha',\theta_1',\theta_2')
    =\frac{\partial^2 Q_n^{\dag\dag}}
     {\partial\theta_2'\partial\theta_1'}({\bf Y}_n;\alpha',\theta_1',\theta_2')
     = -2 \alpha'.
 \end{align*}
Then \ $H_n({\bf Y}_n;\alpha',\theta_1',\theta_2')$ \ has the following leading principal minors
 \begin{align*}
    &\Delta_{1,n}({\bf Y}_n;\alpha',\theta_1',\theta_2')
       =2\DS\sum_{\substack{k=1 \\ k\not\in \{s+1,s+2\}}}^n X_{k-1}^2
            + 2(X_{s}+\theta_1 - \theta_1')^2
            + 2(X_{s+1}+\theta_2 - \theta_2')^2,\\[1mm]
    &\Delta_{2,n}({\bf Y}_n;\alpha',\theta_1',\theta_2')
      = 4\Bigg( ((\alpha')^2+1)\left(
                \DS\sum_{\substack{k=1 \\ k\not\in \{s+1,s+2\}}}^n X_{k-1}^2
                + (X_{s}+\theta_1 - \theta_1')^2
                + (X_{s+1}+\theta_2 - \theta_2')^2
                                \right) \\
    & \phantom{\Delta_{2,n}({\bf Y}_n;\alpha',\theta_1',\theta_2')=\Bigg(}
       -\big(X_{s-1}+X_{s+1}-2\alpha' X_{s}-\mu_\vare-2\alpha'(\theta_1-\theta_1')+(\theta_2-\theta_2')\big)^2
       \Bigg),
 \end{align*}
 and
 \begin{align*}
   &\Delta_{3,n}({\bf Y}_n;\alpha',\theta_1',\theta_2') = \det H_n({\bf Y}_n;\alpha',\theta_1',\theta_2')\\
   & =8\Bigg[
          ((\alpha')^4+(\alpha')^2+1)\left(
                              \DS\sum_{\substack{k=1 \\ k\not\in \{s+1,s+2\}}}^n X_{k-1}^2
            + (X_{s}+\theta_1 - \theta_1')^2
            + (X_{s+1}+\theta_2 - \theta_2')^2
                            \right) \\
         &\phantom{=8\Bigg( }
         - 2\alpha' ab -((\alpha')^2+1)b^2 -((\alpha')^2+1) a^2
         \Bigg],
 \end{align*}
 where
 \begin{align*}
   &a:= X_{s-1}+X_{s+1}-2\alpha' X_{s}-\mu_\vare-2\alpha'(\theta_1-\theta_1')+(\theta_2-\theta_2'),\\
   &b:= X_{s}+X_{s+2}-2\alpha' X_{s+1}-\mu_\vare+(\theta_1-\theta_1')-2\alpha'(\theta_2-\theta_2').
 \end{align*}
By \eqref{Ergodic2},
 \begin{align*}
    &\PP\left(\lim_{n\to\infty}\frac{1}{n}\Delta_{1,n}({\bf Y}_n;\alpha',\theta_1',\theta_2')
               = 2 \EE\widetilde X^2,
               \;\; \forall\;\; (\alpha',\theta_1',\theta_2')\in\RR^3\right)=1, \\[1mm]
    & \PP\left(\lim_{n\to\infty}\frac{1}{n}\Delta_{2,n}({\bf Y}_n;\alpha',\theta_1',\theta_2')
               = 4((\alpha')^2+1)\EE\widetilde X^2,
                \;\; \forall\;\; (\alpha',\theta_1',\theta_2')\in\RR^3 \right)=1,\\[1mm]
    &\PP\left(\lim_{n\to\infty}\frac{1}{n}\Delta_{3,n}({\bf Y}_n;\alpha',\theta_1',\theta_2')
               = 8((\alpha')^4+(\alpha')^2+1)\EE\widetilde X^2,
                \;\; \forall\;\; (\alpha',\theta_1',\theta_2')\in\RR^3 \right)=1,
 \end{align*}
 where \ $\widetilde X$ \ denotes a random variable with the unique stationary
 distribution of the INAR(1) model in \eqref{INAR1}.
Hence
 \begin{align*}
   & \PP\big(\lim_{n\to\infty}\Delta_{i,n}({\bf Y}_n;\alpha',\theta_1',\theta_2')=\infty,
                \;\; \forall\;\; (\alpha',\theta_1',\theta_2')\in\RR^3   \big)=1,
      \qquad i=1,2,3,
 \end{align*}
  which yields that \ ${\bf Y}_n\in S_n$ \ asymptotically as \ $n\to\infty$
 \ with probability one, since we have already proved that \ ${\bf Y}_n\in \widetilde{S}^{\dag\dag}_n$
 \ asymptotically as \ $n\to\infty$ \ with probability one.
\proofend

By Lemma \ref{LEMMA9},
   \ $(\talpha_n^{\,\dag\dag}({\bf Y}_n),
    \ttheta_{1,n}^{\,\dag\dag}({\bf Y}_n),
    \ttheta_{2,n}^{\,\dag\dag}({\bf Y}_n))$
 \ exists uniquely asymptotically as \ $n\to\infty$ \ with probability one.
In the sequel we will simply denote it by
 \ $(\talpha_n^{\,\dag\dag},\ttheta_{1,n}^{\,\dag\dag},\ttheta_{2,n}^{\,\dag\dag})$.

An easy calculation shows that
 \begin{align}\label{SEGED49}
   \talpha_n^{\,\dag\dag}
       = \frac{\sum_{k=1}^n(Y_k-\mu_\vare)Y_{k-1} - \ttheta_{1,n}^{\,\dag\dag}(Y_{s-1}+Y_{s+1}-\mu_\vare)
               - \ttheta_{2,n}^{\,\dag\dag}(Y_{s}+Y_{s+2}-\mu_\vare)+\ttheta_{1,n}^{\,\dag\dag}\ttheta_{2,n}^{\,\dag\dag}}
              {\sum_{k=1}^nY_{k-1}^2 - 2\ttheta_{1,n}^{\,\dag\dag} Y_{s} + (\ttheta_{1,n}^{\,\dag\dag})^2
               - 2 \ttheta_{2,n}^{\,\dag\dag} Y_{s+1} + (\ttheta_{2,n}^{\,\dag\dag})^2},
  \end{align}
 and
  \begin{align}\label{SEGED96}
   \begin{split}
   \begin{bmatrix}
     1+(\talpha_n^{\,\dag\dag})^2 &  -\talpha_n^{\,\dag\dag} \\
     -\talpha_n^{\,\dag\dag} & 1+(\talpha_n^{\,\dag\dag})^2 \\
   \end{bmatrix}
   &\begin{bmatrix}
      \ttheta_{1,n}^{\,\dag\dag} \\
      \ttheta_{2,n}^{\,\dag\dag} \\
   \end{bmatrix}\\
   & = \begin{bmatrix}
      Y_s - \talpha_n^{\,\dag\dag}Y_{s-1} -\mu_\vare - \talpha_n^{\,\dag\dag}(Y_{s+1}
          - \talpha_n^{\,\dag\dag}Y_s -\mu_\vare) \\
       Y_{s+1} - \talpha_n^{\,\dag\dag}Y_s -\mu_\vare - \talpha_n^{\,\dag\dag}(Y_{s+2}
           - \talpha_n^{\,\dag\dag}Y_{s+1} -\mu_\vare) \\
     \end{bmatrix}
   \end{split}
 \end{align}
 hold asymptotically as \ $n\to\infty$ \ with probability one.
Recalling that
 \ $D_n(\talpha_n^{\,\dag\dag})=(\talpha_n^{\,\dag\dag})^4+(\talpha_n^{\,\dag\dag})^2+1>0$, \ we have
 \begin{align}\label{SEGED50}
  \begin{split}
   \ttheta_{1,n}^{\,\dag\dag}
      =\frac{1}{D_n(\talpha_n^{\,\dag\dag})}
        \Big(&(1+(\talpha_n^{\,\dag\dag})^2)
              \Big[Y_s-\talpha_n^{\,\dag\dag} Y_{s-1}-\mu_\vare
                   - \talpha_n^{\,\dag\dag}(Y_{s+1}-\talpha_n^{\,\dag\dag}Y_s-\mu_\vare)
              \Big] \\
             & + \talpha_n^{\,\dag\dag}\Big[Y_{s+1}-\talpha_n^{\,\dag\dag}Y_s-\mu_\vare
                  - \talpha_n^{\,\dag\dag}(Y_{s+2}- \talpha_n^{\,\dag\dag}Y_{s+1}-\mu_\vare)\Big]
            \Big),
  \end{split}
 \end{align}
 and
 \begin{align}\label{SEGED51}
  \begin{split}
   \ttheta_{2,n}^{\,\dag\dag}
      =\frac{1}{D_n(\talpha_n^{\,\dag\dag})}
        \Big(&\talpha_n^{\,\dag\dag}\Big[Y_s-\talpha_n^{\,\dag\dag}Y_{s-1}-\mu_\vare
                  - \talpha_n^{\,\dag\dag}(Y_{s+1}- \talpha_n^{\,\dag\dag}Y_s-\mu_\vare)\Big] \\
              & + (1+(\talpha_n^{\,\dag\dag})^2)
                  \Big[Y_{s+1}-\talpha_n^{\,\dag\dag} Y_s-\mu_\vare
                   - \talpha_n^{\,\dag\dag}(Y_{s+2}-\talpha_n^{\,\dag\dag}Y_{s+1}-\mu_\vare)
                  \Big]
         \Big)
   \end{split}
 \end{align}
 hold asymptotically as \ $n\to\infty$ \ with probability one.

The next result shows that \ $\talpha_n^{\,\dag\dag}$ \ is a strongly consistent estimator
 of \ $\alpha$, \ whereas \ $\ttheta_{1,n}^{\,\dag\dag}$ \ and  \ $\ttheta_{2,n}^{\,\dag\dag}$ \ fail to be
 strongly consistent estimators of \ $\theta_1$ \ and \ $\theta_2$, \ respectively.

\begin{Thm}\label{THEOREM3}
For the CLS estimators
 $(\talpha_n^{\,\dag\dag},\ttheta_{1,n}^{\,\dag\dag},\ttheta_{2,n}^{\,\dag\dag})_{n\in\NN}$ of
 the parameter $(\alpha,\theta_1,\theta_2)\in(0,1)\times\NN^2$,
 \ the sequence \ $(\talpha_n^{\,\dag\dag})_{n\in\NN}$ \ is strongly consistent
 for all \ $(\alpha,\theta_1,\theta_2)\in(0,1)\times\NN^2$, \ i.e.,
 \begin{align} \label{Strong_consistency15}
   \PP(\lim_{n\to\infty}\talpha_n^{\,\dag\dag}=\alpha)=1,
     \qquad \forall\;(\alpha,\theta_1,\theta_2)\in(0,1)\times\NN^2,
 \end{align}
 whereas the sequences \ $(\ttheta_{1,n}^{\,\dag\dag})_{n\in\NN}$ \ and
  \ $(\ttheta_{2,n}^{\,\dag\dag})_{n\in\NN}$ \ are not strongly consistent
  for any \ $(\alpha,\theta_1,\theta_2)\in(0,1)\times\NN^2$, \ namely,
 \begin{align}\label{Strong_consistency16}
   \PP\left(\lim_{n\to\infty}
         \begin{bmatrix}
           \ttheta_{1,n}^{\,\dag\dag} \\[2mm]
           \ttheta_{2,n}^{\,\dag\dag} \\
         \end{bmatrix}
         =
          \begin{bmatrix}
            Y_s \\
            Y_{s+1} \\
          \end{bmatrix}
           +
           \begin{bmatrix}
             \frac{-\alpha(1+\alpha^2)Y_{s-1}
                      -\alpha^2Y_{s+2}
                      -(1-\alpha^3)\mu_\vare}{1+\alpha^2+\alpha^4} \\[2mm]
              \frac{-\alpha^2Y_{s-1}
                      -\alpha(1+\alpha^2)Y_{s+2}
                      -(1-\alpha^3)\mu_\vare}{1+\alpha^2+\alpha^4} \\
           \end{bmatrix}
          \right)=1
 \end{align}
 for all \ $(\alpha,\theta_1,\theta_2)\in(0,1)\times\NN^2$.
\end{Thm}

\noindent{\bf Proof.}
Using that for all \ $p_i\in\RR$, \ $i=0,1,\ldots,4$, \
 \[
   \sup_{x\in\RR}\frac{p_0+p_1x+p_2x^2+p_3x^3+p_4x^4}{1+x^2+x^4}
            <\infty,
 \]
 by \eqref{SEGED50} and \eqref{SEGED51}, we get the sequences
 \ $(\ttheta_{1,n}^{\,\dag\dag})_{n\in\NN}$
 \ and \ $(\ttheta_{2,n}^{\,\dag\dag})_{n\in\NN}$ \ are bounded with probability one.
 Hence using \eqref{SEGED49}, by the same arguments as in the proof of Theorem \ref{THEOREM1}, one can
 derive \eqref{Strong_consistency15}.
Then \eqref{Strong_consistency15}, \eqref{SEGED50} and \eqref{SEGED51} yield \eqref{Strong_consistency16}.
\proofend

The asymptotic distribution of the CLS estimation is given in the next theorem.

\begin{Thm}
Under the additional assumptions \ $\EE X_0^3<\infty$ \ and \ $\EE\vare_1^3<\infty$, \ we have
 \begin{align}\label{CONVERGENCE12}
   \sqrt{n}(\talpha_n^{\,\dag\dag}-\alpha)\distr \cN(0,\sigma_{\alpha,\,\vare}^2)
      \qquad \text{as \ $n\to\infty$,}
  \end{align}
 where \ $\sigma_{\alpha,\,\vare}^2$ \ is defined in \eqref{SEGED_SZIGMA_ALPHA}.
Moreover, conditionally on the values \ $Y_{s-1}$ \  and \ $Y_{s+2}$,
 \begin{align}\label{CONVERGENCE13}
    \begin{bmatrix}
      \sqrt{n}\big(\ttheta_{1,n}^{\,\dag\dag} - \lim_{k\to\infty}\ttheta_{1,k}^{\,\dag\dag}\big) \\
      \sqrt{n}\big(\ttheta_{2,n}^{\,\dag\dag} - \lim_{k\to\infty}\ttheta_{2,k}^{\,\dag\dag}\big) \\
    \end{bmatrix}
      \distr \cN\left(\begin{bmatrix}
                        0 \\
                        0 \\
                      \end{bmatrix},
                      f_{\alpha,\vare} \sigma_{\alpha,\,\vare}^2 f_{\alpha,\vare}^\top
      \right)
      \qquad \text{as \ $n\to\infty$,}
 \end{align}
 where \ $f_{\alpha,\vare}$ \ defined by
 \[
   \frac{1}{(1+\alpha^2+\alpha^4)^2}
         \begin{bmatrix}
           (\alpha^2-1)(\alpha^4+3\alpha^2+1)Y_{s-1} + 2\alpha(\alpha^4-1)Y_{s+2}
                    + \alpha(2-\alpha)(1+\alpha+\alpha^2)^2\mu_\vare \\
            2\alpha(\alpha^4-1)Y_{s-1} + (\alpha^2-1)(\alpha^4+3\alpha^2+1)Y_{s+2}
                    + \alpha(2-\alpha)(1+\alpha+\alpha^2)^2\mu_\vare  \\
         \end{bmatrix}.
 \]
\end{Thm}

\noindent{\bf Proof.}
Using \eqref{SEGED49} and that the sequences
 \ $(\ttheta_{1,n}^{\,\dag\dag})_{n\in\NN}$ \ and \ $(\ttheta_{2,n}^{\,\dag\dag})_{n\in\NN}$ \ are
 bounded with probability one, by the very same arguments as in the proof of \eqref{CONVERGENCE6},
 one can obtain \eqref{CONVERGENCE12}.
Now we turn to prove \eqref{CONVERGENCE13}.
Using the notation
 \[
   B_n^{\dag\dag}:=
          \begin{bmatrix}
             1+(\talpha_n^{\,\dag\dag})^2 & -\talpha_n^{\,\dag\dag} \\
             -\talpha_n^{\,\dag\dag} & 1+(\talpha_n^{\,\dag\dag})^2 \\
           \end{bmatrix},
 \]
 by \eqref{SEGED96}, we have
  \begin{align*}
    \begin{bmatrix}
      \ttheta_{1,n}^{\,\dag\dag} \\
      \ttheta_{2,n}^{\,\dag\dag} \\
    \end{bmatrix}
     = (B_n^{\dag\dag})^{-1}
        \begin{bmatrix}
         (1+(\talpha_n^{\,\dag\dag})^2)Y_{s} - \talpha_n^{\,\dag\dag}(Y_{s-1}+Y_{s+1})
            - (1-\talpha_n^{\,\dag\dag})\mu_\vare\\
         (1+(\talpha_n^{\,\dag\dag})^2)Y_{s+1} - \talpha_n^{\,\dag\dag}(Y_{s}+Y_{s+2}) - (1-\talpha_n^{\,\dag\dag})\mu_\vare \\
        \end{bmatrix}
  \end{align*}
 holds asymptotically as \ $n\to\infty$ \ with probability one.
Theorem \ref{THEOREM3} yields that
 \[
     \PP\left(\lim_{n\to\infty}B_n^{\dag\dag}
                 = \begin{bmatrix}
                    1+\alpha^2 & -\alpha \\
                    -\alpha & 1+\alpha^2 \\
                   \end{bmatrix}=: B^{\dag\dag}
         \right)=1.
 \]
By \eqref{Strong_consistency16}, we have
 \begin{align*}
   &\begin{bmatrix}
      \sqrt{n}\big(\ttheta_{1,n}^{\,\dag\dag} - \lim_{k\to\infty}\ttheta_{1,k}^{\,\dag\dag}\big) \\
      \sqrt{n}\big(\ttheta_{2,n}^{\,\dag\dag} - \lim_{k\to\infty}\ttheta_{2,k}^{\,\dag\dag}\big) \\
   \end{bmatrix}=\\
   &\qquad
    =\sqrt{n}(B_n^{\dag\dag})^{-1}
      \left(
         \begin{bmatrix}
          (1+(\talpha_n^{\,\dag\dag})^2)Y_{s} - \talpha_n^{\,\dag\dag}(Y_{s-1}+Y_{s+1}) - (1-\talpha_n^{\,\dag\dag})\mu_\vare\\
          (1+(\talpha_n^{\,\dag\dag})^2)Y_{s+1} - \talpha_n^{\,\dag\dag}(Y_{s}+Y_{s+2}) - (1-\talpha_n^{\,\dag\dag})\mu_\vare \\
         \end{bmatrix}\right.\\
   &\phantom{\qquad = \sqrt{n}(B_n^{\dag\dag})^{-1}\Big(\;\;}\left.
     - B_n^{\dag\dag} (B^{\dag\dag})^{-1}
      \begin{bmatrix}
         (1+\alpha^2)Y_{s} - \alpha(Y_{s-1}+Y_{s+1}) - (1-\alpha)\mu_\vare\\
         (1+\alpha^2)Y_{s+1} - \alpha(Y_{s}+Y_{s+2}) - (1-\alpha)\mu_\vare \\
      \end{bmatrix}
     \right),
 \end{align*}
 and hence
 \begin{align*}
    &\begin{bmatrix}
      \sqrt{n}\big(\ttheta_{1,n}^{\,\dag\dag} - \lim_{k\to\infty}\ttheta_{1,k}^{\,\dag\dag}\big) \\
      \sqrt{n}\big(\ttheta_{2,n}^{\,\dag\dag} - \lim_{k\to\infty}\ttheta_{2,k}^{\,\dag\dag}\big) \\
    \end{bmatrix}=\\
    &\qquad=\sqrt{n}(B_n^{\dag\dag})^{-1}
      \left(
       \begin{bmatrix}
         (1+(\talpha_n^{\,\dag\dag})^2)Y_{s} - \talpha_n^{\,\dag\dag}(Y_{s-1}+Y_{s+1}) - (1-\talpha_n^{\,\dag\dag})\mu_\vare\\
         (1+(\talpha_n^{\,\dag\dag})^2)Y_{s+1} - \talpha_n^{\,\dag\dag}(Y_{s}+Y_{s+2}) - (1-\talpha_n^{\,\dag\dag})\mu_\vare \\
       \end{bmatrix}\right.\\
  &\phantom{\qquad=\sqrt{n}(B_n^{\dag\dag})^{-1}\Big(\;}\left.
    -\begin{bmatrix}
         (1+\alpha^2)Y_{s} - \alpha(Y_{s-1}+Y_{s+1}) - (1-\alpha)\mu_\vare\\
         (1+\alpha^2)Y_{s+1} - \alpha(Y_{s}+Y_{s+2}) - (1-\alpha)\mu_\vare \\
      \end{bmatrix}
       \right)\\
 &\phantom{\qquad=\;}
      +\sqrt{n}
      \left((B_n^{\dag\dag})^{-1} -  (B^{\dag\dag})^{-1}
       \right)
       \begin{bmatrix}
         (1+\alpha^2)Y_{s} - \alpha(Y_{s-1}+Y_{s+1}) - (1-\alpha)\mu_\vare\\
         (1+\alpha^2)Y_{s+1} - \alpha(Y_{s}+Y_{s+2}) - (1-\alpha)\mu_\vare \\
      \end{bmatrix} \\
 &\qquad=\sqrt{n}(B_n^{\dag\dag})^{-1}
    \begin{bmatrix}
      (\talpha_n^{\,\dag\dag}-\alpha)\big((\talpha_n^{\,\dag\dag}+\alpha)Y_{s} - Y_{s-1}- Y_{s+1} + \mu_\vare\big) \\
      (\talpha_n^{\,\dag\dag}-\alpha)\big((\talpha_n^{\,\dag\dag}+\alpha)Y_{s+1} - Y_{s}- Y_{s+2} + \mu_\vare\big) \\
    \end{bmatrix}\\
 &\phantom{\qquad=\;}
    +\sqrt{n}(B_n^{\dag\dag})^{-1}
       \big( B^{\dag\dag}- B_n^{\dag\dag} \big) (B^{\dag\dag})^{-1}
        \begin{bmatrix}
         (1+\alpha^2)Y_{s} - \alpha(Y_{s-1}+Y_{s+1}) - (1-\alpha)\mu_\vare\\
         (1+\alpha^2)Y_{s+1} - \alpha(Y_{s}+Y_{s+2}) - (1-\alpha)\mu_\vare \\
         \end{bmatrix}.
 \end{align*}
Then
 \begin{align}\label{SEGED52}
   &\begin{bmatrix}
      \sqrt{n}\big(\ttheta_{1,n}^{\,\dag\dag} - \lim_{k\to\infty}\ttheta_{1,k}^{\,\dag\dag}\big) \\
      \sqrt{n}\big(\ttheta_{2,n}^{\,\dag\dag} - \lim_{k\to\infty}\ttheta_{2,k}^{\,\dag\dag}\big) \\
    \end{bmatrix}
   = \sqrt{n}(\talpha_n^{\,\dag\dag}-\alpha)
       \begin{bmatrix}
         K_n^{\dag\dag} \\
         L_n^{\dag\dag} \\
       \end{bmatrix}
  \end{align}
 holds asymptotically as \ $n\to\infty$ \ with probability one, where
 \begin{align*}
    \begin{bmatrix}
         K_n^{\dag\dag} \\
         L_n^{\dag\dag} \\
    \end{bmatrix}
    &:=      (B_n^{\dag\dag})^{-1}
       \begin{bmatrix}
         -(\talpha_n^{\,\dag\dag}+\alpha) & 1 \\
         1 & -(\talpha_n^{\,\dag\dag}+\alpha) \\
       \end{bmatrix}
        (B^{\dag\dag})^{-1}
        \begin{bmatrix}
         (1+\alpha^2)Y_{s} - \alpha(Y_{s-1}+Y_{s+1}) - (1-\alpha)\mu_\vare\\
         (1+\alpha^2)Y_{s+1} - \alpha(Y_{s}+Y_{s+2}) - (1-\alpha)\mu_\vare \\
        \end{bmatrix} \\
    &\phantom{=\;\;\,}
     + (B_n^{\dag\dag})^{-1}
     \begin{bmatrix}
         (\talpha_n^{\,\dag\dag}+\alpha)Y_{s} - Y_{s-1}- Y_{s+1} + \mu_\vare \\
          (\talpha_n^{\,\dag\dag}+\alpha)Y_{s+1} - Y_{s}- Y_{s+2} + \mu_\vare \\
       \end{bmatrix}.
  \end{align*}
By \eqref{Strong_consistency15}, we have
 \ $\begin{bmatrix}
     K_n^{\dag\dag} & L_n^{\dag\dag} \\
    \end{bmatrix}^\top$ \
 converges almost surely as \ $n\to\infty$ \ to
 \begin{align*}
   & (B^{\dag\dag})^{-1}
     \begin{bmatrix}
       2\alpha Y_{s} - Y_{s-1} - Y_{s+1} +\mu_\vare \\
       2\alpha Y_{s+1} - Y_{s} - Y_{s+2} +\mu_\vare  \\
     \end{bmatrix} \\
   &\phantom{(B^{\dag\dag})^{-1}}
     + (B^{\dag\dag})^{-1}
    \begin{bmatrix}
      -2\alpha & 1 \\
      1 & -2\alpha \\
    \end{bmatrix}
    (B^{\dag\dag})^{-1}
      \begin{bmatrix}
        (1+\alpha^2)Y_{s} - \alpha(Y_{s-1}+Y_{s+1}) - (1-\alpha)\mu_\vare\\
        (1+\alpha^2)Y_{s+1} - \alpha(Y_{s}+Y_{s+2}) - (1-\alpha)\mu_\vare \\
      \end{bmatrix}\\
   &\;
     =\frac{1}{1+\alpha^2+\alpha^4}
       \begin{bmatrix}
         1+\alpha^2 & \alpha \\
         \alpha & 1+\alpha^2 \\
       \end{bmatrix}
       \left(
         \begin{bmatrix}
           2\alpha Y_{s} - Y_{s-1} - Y_{s+1} +\mu_\vare \\
           2\alpha Y_{s+1} - Y_{s} - Y_{s+2} +\mu_\vare  \\
         \end{bmatrix}\right. \\
   &\left.
     + \frac{1}{(1+\alpha^2+\alpha^4)^2}
          \begin{bmatrix}
            -2\alpha^5-4\alpha^3 & 1-\alpha^2-3\alpha^4 \\
            1-\alpha^2-3\alpha^4 & -2\alpha^5-4\alpha^3 \\
          \end{bmatrix}
          \begin{bmatrix}
         (1+\alpha^2)Y_{s} - \alpha(Y_{s-1}+Y_{s+1}) - (1-\alpha)\mu_\vare\\
         (1+\alpha^2)Y_{s+1} - \alpha(Y_{s}+Y_{s+2}) - (1-\alpha)\mu_\vare \\
         \end{bmatrix}
         \right),
  \end{align*}
 which is equal to \ $f_{\alpha,\vare}$, \ by an easy, but tedious calculation.
Hence, by \eqref{SEGED52}, \eqref{CONVERGENCE12} and Slutsky's lemma, we have \eqref{CONVERGENCE13}.
\proofend

\subsection{Two not neighbouring outliers, estimation of the mean of the offspring and
            innovation distributions and the outliers' sizes}

In this section we assume that \ $I=2$ \ and that the relevant time points
 \ $s_1$, $s_2\in\NN$ \ are known.
We also suppose that \ $s_1<s_2-1$, \ i.e., the time points \ $s_1$ \ and
 \ $s_2$ \ are not neighbouring.
We concentrate on the CLS estimation of \ $\alpha$, \ $\mu_\vare$, \ $\theta_1$
 \ and \ $\theta_2$.

Motivated by \eqref{SEGED53}, for all \ $n\geq s_2+1$, \ $n\in\NN$, \ we define the function
 \ $Q_n^{\dag}:\RR^{n+1}\times\RR^4\to\RR$, \ as
 \begin{align*}
    &Q_n^{\dag}({\bf y}_n;\alpha',\mu_\vare',\theta_1',\theta_2')\\
      & :=\DS\sum_{\substack{k=1 \\ k\not\in \{s_1,s_1+1,s_2,s_2+1\}}}^n
            \big(y_k-\alpha' y_{k-1}-\mu_\vare'\big)^2
           + \big(y_{s_1}-\alpha' y_{s_1-1}-\mu_\vare'-\theta_1'\big)^2 \\
         &\phantom{=\;}
           + \big(y_{s_1+1}-\alpha' y_{s_1}-\mu_\vare'+\alpha'\theta_1'\big)^2
           + \big(y_{s_2}-\alpha' y_{s_2-1}-\mu_\vare'-\theta_2'\big)^2
           + \big(y_{s_2+1}-\alpha' y_{s_2}-\mu_\vare'+\alpha'\theta_2'\big)^2,
 \end{align*}
 for all \ ${\bf y}_n\in\RR^{n+1}$, $\alpha',\mu_\vare',\theta_1',\theta_2'\in\RR$.
\ By definition, for all \ $n\geq s_2+1$, \ a CLS estimator for
 the parameter \ $(\alpha,\mu_\vare,\theta_1,\theta_2)\in(0,1)\times(0,\infty)\times\NN^2$ \ is
 a measurable function
 \[
   (\halpha_n^{\,\dag},\hmuen^{\,\dag},\htheta_{1,n}^{\,\dag},\htheta_{2,n}^{\,\dag}):S_n\to\RR^4
 \]
 such that
 \begin{align*}
   Q_n^{\dag}({\bf y}_n;\,&\halpha_n^{\,\dag}({\bf y}_n),\hmuen^{\,\dag}({\bf y}_n),
           \htheta_{1,n}^{\,\dag}({\bf y}_n),
             \htheta_{2,n}^{\,\dag}({\bf y}_n))\\
   &=\inf_{(\alpha',\mu_\vare',\theta_1',\theta_2')\in\RR^4}
       Q_n^{\dag}({\bf y}_n;\alpha',\mu_\vare',\theta_1',\theta_2')
         \qquad \forall\;{\bf y}_n\in S_n,
 \end{align*}
 where \ $S_n$ \ is suitable subset of \ $\RR^{n+1}$ \ (defined in the proof of
 Lemma \ref{LEMMA13}).
\ We note that we do not define the CLS estimator
 \ $(\halpha_n^{\,\dag},\hmuen^{\,\dag},\htheta_{1,n}^{\,\dag},\htheta_{2,n}^{\,\dag})$ \
 for all samples \ ${\bf y}_n\in \RR^{n+1}$.

The next result is about the existence and uniqueness of
 \ $(\halpha_n^{\,\dag},\hmuen^{\,\dag},\htheta_{1,n}^{\,\dag},\htheta_{2,n}^{\,\dag})$.

\begin{Lem}\label{LEMMA13}
There exist subsets \ $S_n\subset\RR^{n+1}$, $n\geq \max(5,s_2+1)$ \ with the following properties:
 \begin{enumerate}
  \item[\upshape{(i)}]
   there exists a unique CLS estimator
   \ $(\halpha_n^{\,\dag},\hmuen^{\,\dag},\htheta_{1,n}^{\,\dag},\htheta_{2,n}^{\,\dag}):S_n\to\RR^4$,
  \item[\upshape{(ii)}]
   for all \ ${\bf y}_n\in S_n$,
   \ $(\halpha_n^{\,\dag}({\bf y}_n),\hmuen^{\,\dag}({\bf y}_n),
       \htheta_{1,n}^{\,\dag}({\bf y}_n),\htheta_{2,n}^{\,\dag}({\bf y}_n))$
   \ is the unique solution of the system of equations
  \begin{align}\label{Additive_CLSE_EQ7}
   \begin{split}
    \frac{\partial Q_n^{\dag}}{\partial \alpha'}
     ({\bf y}_n;\alpha',\mu_\vare',\theta_1',\theta_2')=0,
   \quad\quad
    \frac{\partial Q_n^{\dag}}{\partial \mu_\vare'}
     ({\bf y}_n;\alpha',\mu_\vare',\theta_1',\theta_2')=0,\\[2mm]
    \frac{\partial Q_n^{\dag}}{\partial \theta_1'}
     ({\bf y}_n;\alpha',\mu_\vare',\theta_1',\theta_2')=0,
       \quad\quad
    \frac{\partial Q_n^{\dag}}{\partial \theta_2'}
      ({\bf y}_n;\alpha',\mu_\vare',\theta_1',\theta_2')=0,
   \end{split}
 \end{align}
  \item [\upshape{(iii)}]
  ${\bf Y}_n\in S_n$ \ holds asymptotically as \ $n\to\infty$ \ with probability one.
 \end{enumerate}
\end{Lem}

\noindent{\bf Proof.}
For any fixed \ ${\bf y}_n \in \RR^{n+1}$, \ $n\geq \max(5,s_2+1)$ \ and \ $\alpha' \in \RR$,
 \ the quadratic function
 \ $\RR^3 \ni (\mu_\vare',\theta_1',\theta_2')
      \mapsto Q_n^{\dag}({\bf y}_n;\alpha',\mu_\vare',\theta_1',\theta_2')$
 \ can be written in the form
 \begin{align*}
  &Q_n^{\dag}({\bf y}_n;\alpha',\mu_\vare',\theta_1',\theta_2') \\[2mm]
  &=  \!\!\left(\! \begin{bmatrix}
             \mu_\vare' \smallskip \\
             \theta_1' \smallskip \\
             \theta_2'
            \end{bmatrix}
           - A_n(\alpha')^{-1} t_n({\bf y}_n;\alpha') \!\! \right)^{\hspace*{-1mm}\top}
      \!\!\! A_n(\alpha') \! \left( \! \begin{bmatrix}
                          \mu_\vare' \smallskip \\
                          \theta_1' \smallskip \\
                          \theta_2'
                         \end{bmatrix}
                         - A_n(\alpha')^{-1} t_n({\bf y}_n;\alpha') \!\!\! \right)
    \!\! + \widehat{Q}_n^{\dag}({\bf y}_n;\alpha'),
 \end{align*}
 where
 \begin{align*}
  t_n({\bf y}_n;\alpha')
  &:= \begin{bmatrix}
       \sum_{k=1}^n ( y_k - \alpha' y_{k-1}) \smallskip \\
       ( 1 + (\alpha')^2 ) y_{s_1} - \alpha' ( y_{s_1-1} + y_{s_1+1} ) \smallskip \\
       ( 1 + (\alpha')^2 ) y_{s_2} - \alpha' ( y_{s_2-1} + y_{s_2+1} )
      \end{bmatrix} , \\[2mm]
  \widehat{Q}_n^{\dag}({\bf y}_n;\alpha')
  &:= \sum_{k=1}^n \big( y_k-\alpha' y_{k-1} \big)^2
      - t_n({\bf y}_n;\alpha')^\top A_n(\alpha')^{-1} t_n({\bf y}_n;\alpha') ,
 \end{align*}
 and the matrix
 \begin{align*}
  A_n(\alpha')
  := \begin{bmatrix}
      n & 1-\alpha' & 1-\alpha' \smallskip \\
      1-\alpha' & 1+(\alpha')^2 & 0 \smallskip \\
      1-\alpha' & 0 &  1+(\alpha')^2 \\
     \end{bmatrix}
 \end{align*}
 is strictly positive definite for all \ $n\geq5$ \ and \ $\alpha' \in \RR$.
\ Indeed, the leading principal minors of \ $A_n(\alpha')$ \ take the following forms:
 \ $n,$
 \begin{align*}
  & n(1+(\alpha')^2) - (1-\alpha')^2 = (n-1)(\alpha')^2 + 2\alpha' + n-1,\\
  & D_n(\alpha'):=(1+(\alpha')^2)\big( (n-2)(\alpha')^2 + 4\alpha' + n-2 \big),
 \end{align*}
 and for all \ $n\geq 5$, \ the discriminant \ $16-4(n-2)^2$ \ of the equation
 \ $(n-2)x^2 + 4x +n-2 = 0$ \ is negative.

The inverse matrix \ $A_n(\alpha')^{-1}$ \ takes the form
 \begin{align*}
  \frac{1}{D_n(\alpha')}\!\!
   \begin{bmatrix}
    (1+(\alpha')^2)^2 \, & \, -(1-\alpha')(1+(\alpha')^2)
     \, & \, -(1-\alpha')(1+(\alpha')^2)  \\
    -(1-\alpha')(1+(\alpha')^2) \, & \, n(1+(\alpha')^2)-(1-\alpha')^2 \, & \, (1-\alpha')^2 \!\! \\
    -(1-\alpha')(1+(\alpha')^2) \, & \, (1-\alpha')^2 \, & \, n(1+(\alpha')^2)-(1-\alpha')^2
  \end{bmatrix}.
 \end{align*}
The polynomial \ $\RR\ni\alpha' \mapsto D_n(\alpha')$ \ is of order 4 with leading coefficient \ $n-2$.
\ We have \ $\widehat{Q}_n^{\dag}({\bf y}_n;\alpha') = R_n({\bf y}_n;\alpha') / D_n(\alpha')$, \ where
 \ $\RR\ni\alpha' \mapsto R_n({\bf y}_n;\alpha')$ \ is a polynomial of order 6 with leading coefficient
 \begin{align*}
   c_n({\bf y}_n) &:= (n-2) \sum_{k=1}^n y_{k-1}^2 - \left( \sum_{k=1}^n y_{k-1} \right)^2
                      - (n-1)(y_{s_1}^2 + y_{s_1}^2) \\[2mm]
                  &\;\quad + 2 (y_{s_1} + y_{s_1}) \sum_{k=1}^n y_{k-1} - 2 y_{s_1} y_{s_1} .
 \end{align*}
Let
 \[
   \widehat{S}^{\dag}_n := \left\{{\bf y}_n\in\RR^{n+1} : c_n({\bf y}_n) > 0 \right\}.
 \]
For \ ${\bf y}_n \in \widehat{S}^{\dag}_n$, \ we have
 \ $\lim_{|\alpha'|\to\infty} \widehat{Q}_n^{\,\dag}({\bf y}_n;\alpha') = \infty$ \ and the continuous function
 \ $\RR \ni \alpha' \mapsto \widehat{Q}_n^{\,\dag}({\bf y}_n;\alpha')$ \ attains its infimum.
Consequently, for all \ $n\geq\max(5,s_2+1)$ \ there exists a CLS estimator
 \ $(\halpha_n^{\,\dag},
    \hmuen^{\,\dag},
    \htheta_{1,n}^{\,\dag},
    \htheta_{2,n}^{\,\dag}) : \widehat{S}_n^{\dag} \to \RR^4$,
\ where
 \begin{align}\nonumber
   \widehat{Q}_n^{\,\dag}({\bf y}_n;\halpha_n^{\,\dag}({\bf y}_n))
   & =\inf_{\alpha'\in\RR} \widehat{Q}_n^{\,\dag}({\bf y}_n;\alpha')
   \qquad \forall\;{\bf y}_n\in  \widehat{S}^{\dag}_n ,\\[2mm] \label{CLSE_mu_theta}
  \begin{bmatrix}
   \hmuen^{\,\dag}({\bf y}_n) \smallskip \\
   \htheta_{1,n}^{\,\dag}({\bf y}_n) \smallskip \\
   \htheta_{2,n}^{\,\dag}({\bf y}_n)
  \end{bmatrix}
  & = A_n(\halpha_n^{\,\dag}({\bf y}_n))^{-1} t_n({\bf y}_n;\halpha_n^{\,\dag}({\bf y}_n)) ,
  \qquad {\bf y}_n\in  \widehat{S}^{\dag}_n,
 \end{align}
 and for all \ ${\bf y}_n\in \widehat{S}^{\dag}_n$,
   \ $(\halpha_n^{\,\dag}({\bf y}_n),\hmuen^{\,\dag}({\bf y}_n),\htheta_{1,n}^{\,\dag}({\bf y}_n),
       \htheta_{2,n}^{\,\dag}({\bf y}_n))$
   \ is a  solution of the system of equations \eqref{Additive_CLSE_EQ7}.

By \eqref{Ergodic1} and \eqref{Ergodic2}, we get
 \ $\PP\left(\lim_{n\to\infty} n^{-2} c_n({\bf Y}_n) = \var\widetilde X \right)=1$, \ where \ $\widetilde X$
 \ denotes a random variable with the unique stationary distribution of the INAR(1) model in \eqref{INAR1}.
Hence \ ${\bf Y}_n\in \widehat{S}_n^{\dag}$ \ holds asymptotically as \ $n\to\infty$ \ with probability one.

Now we turn to find sets \ $S_n \subset \widehat{S}^{\dag}_n$, $n\geq \max(5,s_2+1)$ \ such that
 the system of equations
 \eqref{Additive_CLSE_EQ7} has a unique solution with respect to
 \ $(\alpha',\mu_\vare',\theta_1',\theta_2')$ \ for all \ ${\bf y}_n\in S_n$.
\ Let us introduce the \ $(4\times 4)$ \ Hessian matrix
 \[
    H_n({\bf y}_n;\alpha',\mu_\vare',\theta_1',\theta_2'):=
     \begin{bmatrix}
       \frac{\partial^2 Q_n^{\dag}}{\partial (\alpha')^2}
      \hspace{1mm} & \frac{\partial^2 Q_n^{\dag}}{\partial \mu_\vare' \partial \alpha'}
      \hspace{1mm} & \frac{\partial^2 Q_n^{\dag}}{\partial \theta_1' \partial \alpha'}
     \hspace{1mm}  & \frac{\partial^2 Q_n^{\dag}}{\partial \theta_2' \partial \alpha'} \smallskip \\
                     \frac{\partial^2 Q_n^{\dag}}{\partial \alpha' \partial \mu_\vare'}
     \hspace{1mm}   & \frac{\partial^2 Q_n^{\dag}}{\partial (\mu_\vare')^2}
     \hspace{1mm}   & \frac{\partial^2 Q_n^{\dag}}{ \partial \theta_1' \partial \mu_\vare' }
     \hspace{1mm}   &  \frac{\partial^2 Q_n^{\dag}}{\partial \theta_2' \partial \mu_\vare'} \smallskip \\
                      \frac{\partial^2 Q_n^{\dag}}{\partial \alpha' \partial \theta_1'}
    \hspace{1mm}    & \frac{\partial^2 Q_n^{\dag}}{\partial \mu_\vare' \partial \theta_1'}
    \hspace{1mm}    & \frac{\partial^2 Q_n^{\dag}}{ \partial (\theta_1')^2}
    \hspace{1mm}    &  \frac{\partial^2 Q_n^{\dag}}{\partial \theta_2' \partial \theta_1'} \smallskip \\
                       \frac{\partial^2 Q_n^{\dag}}{\partial \alpha' \partial \theta_2'}
   \hspace{1mm}     & \frac{\partial^2 Q_n^{\dag}}{\partial \mu_\vare' \partial \theta_2'}
   \hspace{1mm}     & \frac{\partial^2 Q_n^{\dag}}{ \partial \theta_1' \partial \theta_2'}
   \hspace{1mm}     &  \frac{\partial^2 Q_n^{\dag}}{ \partial (\theta_2')^2 } \\
     \end{bmatrix}
       ({\bf y}_n;\alpha',\mu_\vare',\theta_1',\theta_2'),
 \]
 and let us denote by \ $\Delta_{i,n}({\bf y}_n;\alpha',\mu_\vare',\theta_1',\theta_2')$ \
 its \ $i$-th order leading principal minor, \ $i=1,2,3,4$.
\ Further, for all \ $n\geq \max(5,s_2+1)$, \ let
 \[
   S_n:=\Big\{{\bf y}_n\in\widehat{S}^{\dag}_n
          : \Delta_{i,n}({\bf y}_n;\alpha',\mu_\vare',\theta_1',\theta_2')>0,
                                   \;\, i=1,2,3,4,\, \forall\;(\alpha',\mu_\vare',\theta_1',\theta_2')
                                   \in\RR^4 \Big\}.
 \]
By Berkovitz \cite[Theorem 3.3, Chapter III]{Ber}, the function
 \ $\RR^4 \ni (\alpha',\mu_\vare',\theta_1',\theta_2')
     \mapsto Q_n^{\dag}({\bf y}_n;\alpha',\mu_\vare',\theta_1',\theta_2')$
 \ is strictly convex for all \ ${\bf y}_n\in S_n$.
\ Since it was already proved that the system of equations \eqref{Additive_CLSE_EQ7} has a solution for all
 \ ${\bf y}_n\in \widehat{S}^{\dag}_n$,
 \ we obtain that this solution is unique for all \ ${\bf y}_n\in S_n$.

For all \ ${\bf y}_n\in\RR^{n+1}$ \ and
 \ $(\alpha',\mu_\vare',\theta_1',\theta_2')\in\RR^4$, \ we have
 \begin{align*}
    &\frac{\partial Q_n^{\dag}}{\partial \alpha'}({\bf y}_n;\alpha',\mu_\vare',\theta_1',\theta_2')\\
    &\phantom{\;\;}
       = \DS\sum_{\substack{k=1 \\ k\not\in \{s_1,s_1+1,s_2,s_2+1\}}}^n\!
            \big(y_k-\alpha' y_{k-1}-\mu_\vare'\big)(-2y_{k-1})
           -2\big(y_{s_1}-\alpha' y_{s_1-1}-\mu_\vare'-\theta_1'\big)y_{s_1-1}\\
    &\phantom{\;\; = \;\;}
       + 2\big(y_{s_1+1}-\alpha' y_{s_1}-\mu_\vare'+\alpha'\theta_1'\big)(-y_{s_1}+\theta_1')
                - 2\big(y_{s_2}-\alpha' y_{s_2-1}-\mu_\vare'-\theta_2'\big)y_{s_2-1} \\
    &\phantom{\;\; = \;\;}
      + 2\big(y_{s_2+1}-\alpha' y_{s_2}-\mu_\vare'+\alpha'\theta_2'\big)(-y_{s_2}+\theta_2'),\\[2mm]
   &\frac{\partial Q_n^{\dag}}{\partial \mu_\vare'}({\bf y}_n;\alpha',\mu_\vare',\theta_1',\theta_2')\\
   &\phantom{\quad }
     = \DS\sum_{\substack{k=1 \\ k\not\in \{s_1,s_1+1,s_2,s_2+1\}}}^n
            (-2)\big(y_k-\alpha' y_{k-1}-\mu_\vare'\big)
           -2\big(y_{s_1}-\alpha' y_{s_1-1}-\mu_\vare'-\theta_1'\big) \\
   &\phantom{\;\; = \;\;}
      - 2\big(y_{s_1+1}-\alpha' y_{s_1}-\mu_\vare'+\alpha'\theta_1'\big)
                - 2\big(y_{s_2}-\alpha' y_{s_2-1}-\mu_\vare'-\theta_2'\big)\\
   &\phantom{\;\; = \;\;}
       - 2\big(y_{s_2+1}-\alpha' y_{s_2}-\mu_\vare'+\alpha'\theta_2'\big),
 \end{align*}
 and
 \begin{align}\label{SEGED_pluszban}
  \begin{split}
   & \frac{\partial Q_n^{\dag}}{\partial \theta_i'}({\bf y}_n;\alpha',\mu_\vare',\theta_1',\theta_2')\\
   & \phantom{\;\; =}
       = -2(y_{s_i}-\alpha' y_{s_i-1}-\mu_\vare'-\theta_i')
         +2\alpha'(y_{s_i+1}-\alpha' y_{s_i}-\mu_\vare'+\alpha'\theta_i'),
       \qquad i=1,2.
 \end{split}
 \end{align}
We also get for all \ $(\alpha',\mu_\vare',\theta_1',\theta_2')\in\RR^4$,
 \begin{align*}
   &\frac{\partial^2 Q_n^{\dag}}{\partial (\alpha')^2} ({\bf Y}_n;\alpha',\mu_\vare',\theta_1',\theta_2')\\
   &\phantom{\frac{\partial^2 Q_n^{\dag}}{\partial (\alpha')^2}\;}
    = 2 \DS\sum_{\substack{k=1 \\ k\not\in \{s_1,s_1+1,s_2,s_2+1\}}}^n Y_{k-1}^2
            +2Y_{s_1-1}^2+2(Y_{s_1}-\theta_1')^2
            +2Y_{s_2-1}^2+2(Y_{s_2}-\theta_2')^2 \\
   &\phantom{\frac{\partial^2 Q_n^{\dag}}{\partial (\alpha')^2}\;}
      =2\DS\sum_{\substack{k=1 \\ k\not\in \{s_1+1,s_2+1\}}}^n \!\!\!X_{k-1}^2
            + 2(X_{s_1}+\theta_1 - \theta_1')^2
            + 2(X_{s_2}+\theta_2 - \theta_2')^2,
  \end{align*}
 and
 \begin{align*}
  &\frac{\partial^2 Q_n^{\dag}}
    {\partial\mu_\vare'\partial\alpha'}({\bf Y}_n;\alpha',\mu_\vare',\theta_1',\theta_2')
    = \frac{\partial^2 Q_n^{\dag}}
     {\partial\alpha'\partial\mu_\vare'}({\bf Y}_n;\alpha',\mu_\vare',\theta_1',\theta_2')\\
   &\phantom{\frac{\partial^2 Q_n^{\dag}}
      {\partial\alpha'\partial\mu_\vare'}(}
    = 2\sum_{k=1}^nY_{k-1}-2\theta_1'-2\theta_2'
    = 2\sum_{k=1}^nX_{k-1}+2(\theta_1-\theta_1')+2(\theta_2-\theta_2'),\\
  &\frac{\partial^2 Q_n^{\dag}}
    {\partial\theta_i'\partial\alpha'}({\bf Y}_n;\alpha',\mu_\vare',\theta_1',\theta_2')
    = \frac{\partial^2 Q_n^{\dag}}
       {\partial\alpha'\partial\theta_i'}({\bf Y}_n;\alpha',\mu_\vare',\theta_1',\theta_2')\\
  &\phantom{\qquad\qquad\qquad\qquad}
    = 2(Y_{s_i-1}+Y_{s_i+1}-2\alpha' Y_{s_i}-\mu_\vare'+2\alpha'\theta_2')\\
  &\phantom{\qquad\qquad\qquad\qquad}
    =2(X_{s_i-1}+X_{s_i+1}-2\alpha' X_{s_i}-\mu_\vare'-2\alpha'(\theta_2-\theta_2')),\qquad i=1,2, \\
  & \frac{\partial^2 Q_n^{\dag}}{\partial(\mu_\vare')^2}({\bf Y}_n;\alpha',\mu_\vare',\theta_1',\theta_2')
     = 2n,\\[1mm]
  &\frac{\partial^2 Q_n^{\dag}}{\partial(\theta_1')^2}({\bf Y}_n;\alpha',\mu_\vare',\theta_1',\theta_2')
    =\frac{\partial^2 Q_n^{\dag}}{\partial(\theta_2')^2}({\bf Y}_n;\alpha',\mu_\vare',\theta_1',\theta_2')
    = 2((\alpha')^2+1),\\[1mm]
  &\frac{\partial^2 Q_n^{\dag}}{\partial\theta_1'\partial\theta_2'}
     ({\bf Y}_n;\alpha',\mu_\vare',\theta_1',\theta_2')
    =\frac{\partial^2 Q_n^{\dag}}
      {\partial\theta_2'\partial\theta_1'}({\bf Y}_n;\alpha',\mu_\vare',\theta_1',\theta_2')
     = 0,\\
 &\frac{\partial^2 Q_n^{\dag}}
   {\partial\theta_i'\partial\mu_\vare'}({\bf Y}_n;\alpha',\mu_\vare',\theta_1',\theta_2')
    =\frac{\partial^2 Q_n^{\dag}}
       {\partial\mu_\vare'\partial\theta_i'}({\bf Y}_n;\alpha',\mu_\vare',\theta_1',\theta_2')
     = 2(1-\alpha'),\qquad i=1,2.
 \end{align*}
The matrix \ $H_n({\bf Y}_n;\alpha',\mu_\vare',\theta_1',\theta_2')$ \ has the following
 leading principal minors
 \begin{align*}
    &\Delta_{1,n}({\bf Y}_n;\alpha',\mu_\vare',\theta_1',\theta_2')
       =2\!\!\DS\sum_{\substack{k=1 \\ k\not\in \{s_1+1,s_2+1\}}}^n \!\!\!\!\!X_{k-1}^2
            + 2(X_{s_1}+\theta_1 - \theta_1')^2
            + 2(X_{s_2}+\theta_2 - \theta_2')^2,\\
    &\Delta_{2,n}({\bf Y}_n;\alpha',\mu_\vare',\theta_1',\theta_2')
      = 4n\!\left(\!\DS\sum_{\substack{k=1 \\ k\not\in \{s_1+1,s_2+1\}}}^n \!\!\!\!\!\!\!\!\!X_{k-1}^2
                                   + (X_{s_1}+\theta_1 - \theta_1')^2
                                   + (X_{s_2}+\theta_2 - \theta_2')^2\!\!\right) \\
    & \phantom{\Delta_{2,n}({\bf Y}_n;\alpha',\mu_\vare',\theta_1',\theta_2'):=\;}
       -4\left(\sum_{k=1}^nX_{k-1}+(\theta_1-\theta_1')+(\theta_2-\theta_2')\right)^2,\\
    &\Delta_{3,n}({\bf Y}_n;\alpha',\mu_\vare',\theta_1',\theta_2')
     =8\Big(((\alpha')^2+1)n - (1-\alpha')^2\Big)\\
     & \phantom{\Delta_{3,n}({\bf Y}_n;\alpha',\mu_\vare',\theta_1',\theta_2') :=\;}\times
          \!\left(\!\DS\sum_{\substack{k=1 \\ k\not\in \{s_1+1,s_2+1\}}}^n \!\!\!\!\!\!\!\!\!\!X_{k-1}^2
                                   + (X_{s_1}+\theta_1 - \theta_1')^2
                                   + (X_{s_2}+\theta_2 - \theta_2')^2\!\right) \\
    &\phantom{\Delta_{3,n}({\bf Y}_n;\alpha',\mu_\vare',\theta_1',\theta_2') :=\;}
     +16(1-\alpha')L\left(\sum_{k=1}^nX_{k-1}+(\theta_1-\theta_1')+(\theta_2-\theta_2')\right)
              \!\! -8nL^2 \\
    & \phantom{\Delta_{3,n}({\bf Y}_n;\alpha',\mu_\vare',\theta_1',\theta_2') :=\;}
         -8((\alpha')^2+1)\left(\sum_{k=1}^nX_{k-1}+(\theta_1-\theta_1')+(\theta_2-\theta_2')\right)^2
  \end{align*}
  and
  \begin{align*}
    \Delta_{4,n}({\bf Y}_n;\alpha',\mu_\vare',\theta_1',\theta_2')
       = \det H_n({\bf Y}_n;\alpha',\mu_\vare',\theta_1',\theta_2'),
 \end{align*}
 where \ $L:=X_{s_1-1}+X_{s_1+1}-2\alpha' X_{s_1}-\mu_\vare'-2\alpha'(\theta_1-\theta_1')$.
\ By \eqref{Ergodic1} and \eqref{Ergodic2}, we get the following events
 have probability one
 \begin{align*}
    &\left\{\lim_{n\to\infty}\frac{1}{n}\Delta_{1,n}({\bf Y}_n;\alpha',\mu_\vare',\theta_1',\theta_2')
               = 2 \EE\widetilde X^2, \,\;\forall \;(\alpha',\mu_\vare',\theta_1',\theta_2')
               \in\RR^4 \right\}, \\[2mm]
    &\left\{\lim_{n\to\infty}\frac{1}{n^2}\Delta_{2,n}({\bf Y}_n;\alpha',\mu_\vare',\theta_1',\theta_2')
               = 4(\EE\widetilde X^2-(\EE\widetilde X)^2)\right.\\
    &\phantom{\left\{\lim_{n\to\infty}\frac{1}{n^2}\Delta_{2,n}
              ({\bf Y}_n;\alpha',\mu_\vare',\theta_1',\theta_2') \right.}
              \left.=4\var\widetilde X,
               \,\;\forall \;(\alpha',\mu_\vare',\theta_1',\theta_2') \in\RR^4 \right\},  \\[2mm]
    &\left\{\lim_{n\to\infty}\frac{1}{n^2}\Delta_{3,n}({\bf Y}_n;\alpha',\mu_\vare',\theta_1',\theta_2')
               = 8((\alpha')^2+1)\var\widetilde X,
               \,\;\forall \;(\alpha',\mu_\vare',\theta_1',\theta_2') \in\RR^4 \right\},  \\[2mm]
    &\left\{\lim_{n\to\infty}\frac{1}{n^2}\Delta_{4,n}({\bf Y}_n;\alpha',\mu_\vare',\theta_1',\theta_2')
               =16((\alpha')^2+1)^2\var\widetilde X,
               \,\;\forall \;(\alpha',\mu_\vare',\theta_1',\theta_2') \in\RR^4 \right\},
 \end{align*}
 where \ $\widetilde X$ \ denotes a random variable with the unique stationary
 distribution of the INAR(1) model in \eqref{INAR1}.
Hence
 \begin{align*}
    \PP\Big(\lim_{n\to\infty}\Delta_{i,n}({\bf Y}_n;\alpha',\mu_\vare',\theta_1',\theta_2')=\infty,
           \quad \forall\;(\alpha',\mu_\vare',\theta_1',\theta_2')\in\RR^4,
           \;\; i=1,2,3,4\Big)=1,
 \end{align*}
 which yields that \ ${\bf Y}_n\in S_n$ \ asymptotically as \ $n\to\infty$
 \ with probability one, since we have already proved that \ ${\bf Y}_n\in \widehat{S}^{\dag}_n$
 \ asymptotically as \ $n\to\infty$ \ with probability one.
\proofend

By Lemma \ref{LEMMA13}, \ $(\halpha_n^{\,\dag}({\bf Y}_n),
    \hmuen^{\,\dag}({\bf Y}_n),
    \htheta_{1,n}^{\,\dag}({\bf Y}_n),
    \htheta_{2,n}^{\,\dag}({\bf Y}_n))$
 \ exists uniquely asymptotically as \ $n\to\infty$ \ with probability one.
In the sequel we will simply denote it
 by \ $(\halpha_n^{\,\dag},\hmuen^{\,\dag},\htheta_{1,n}^{\,\dag},
            \htheta_{2,n}^{\,\dag})$.

The next result shows that \ $\halpha_n^{\,\dag}$ \ is a strongly consistent estimator
 of \ $\alpha$, \ $\hmuen^{\,\dag}$ \ is a strongly consistent estimator of \ $\mu_\vare$,
 \ whereas \ $\htheta_{1,n}^{\,\dag}$ \ and  \ $\htheta_{2,n}^{\,\dag}$ \ fail to be
 strongly consistent estimators of \ $\theta_1$ \ and \ $\theta_2$, \ respectively.

\begin{Thm}\label{THEOREM4}
Consider the CLS estimators
 \ $(\halpha_n^{\,\dag},\hmuen^{\,\dag},\htheta_{1,n}^{\,\dag},\htheta_{2,n}^{\,\dag})_{n\in\NN}$
 \ of the parameter
 \ $(\alpha,\mu_\vare,\theta_1,\theta_2) \in(0,1)\times(0,\infty)\times\NN^2$.
\ Then the sequences \ $(\halpha_n^{\,\dag})_{n\in\NN}$ \ and \ $(\hmuen^{\,\dag})_{n\in\NN}$ \ are strongly
 consistent for all
 \ $(\alpha,\mu_\vare,\theta_1,\theta_2)\in(0,1)\times(0,\infty)\times\NN^2$, \ i.e.,
 \begin{align} \label{Strong_consistency18}
   &\PP(\lim_{n\to\infty}\halpha_n^{\,\dag}=\alpha)=1,
     \qquad \forall\;(\alpha,\mu_\vare,\theta_1,\theta_2)\in(0,1)\times(0,\infty)\times\NN^2,\\ \label{Strong_consistency19}
   &\PP(\lim_{n\to\infty}\hmuen^{\,\dag}=\mu_\vare)=1,
     \qquad \forall\;(\alpha,\mu_\vare,\theta_1,\theta_2)\in(0,1)\times(0,\infty)\times\NN^2,
 \end{align}
  whereas the sequences \ $(\htheta_{1,n}^{\,\dag})_{n\in\NN}$ \ and \ $(\htheta_{2,n}^{\,\dag})_{n\in\NN}$
 \ are not strongly consistent for any
 \ $(\alpha,\mu_\vare,\theta_1,\theta_2)\in(0,1)\times(0,\infty)\times\NN^2$, \ namely,
 \begin{align}\label{Strong_consistency20}
    \PP\left(\lim_{n\to\infty}\htheta_{i,n}^{\,\dag}
         = Y_{s_i}-\frac{\alpha}{1+\alpha^2}(Y_{s_i-1}+Y_{s_i+1})
                    - \frac{1-\alpha}{1+\alpha^2}\mu_\vare \right)=1,
              \qquad i=1,2,
 \end{align}
 for all \ $(\alpha,\mu_\vare,\theta_1,\theta_2)\in(0,1)\times(0,\infty)\times\NN^2$.
\end{Thm}

\noindent{\bf Proof.}
The aim of the following discussion is to show that the sequences
 \ $(\htheta_{1,n}^{\,\dag}-\theta_1)_{n\in\NN}$
 \ and \ $(\htheta_{2,n}^{\,\dag}-\theta_2)_{n\in\NN}$ \ are bounded with probability one.
By \eqref{Additive_CLSE_EQ7}, \eqref{SEGED_pluszban} and Lemma \ref{LEMMA13}, we get
 \begin{align}\label{SEGED56}
  \htheta_{i,n}^{\,\dag} = Y_{s_i}-\frac{\halpha_n^{\,\dag}}{1+(\halpha_n^{\,\dag})^2}(Y_{s_i-1}+Y_{s_i+1})
                    - \frac{1-\halpha_n^{\,\dag}}{1+(\halpha_n^{\,\dag})^2}\hmuen^{\,\dag} \qquad i=1,2.
 \end{align}
By \eqref{CLSE_mu_theta} and the explicit form of the inverse matrix \ $A_n(\alpha')^{-1}$, \ we obtain
 \begin{align*}
   \begin{bmatrix}
     \hmuen^{\,\dag} \\
     \htheta_{1,n}^{\,\dag} \\
     \htheta_{2,n}^{\,\dag}
   \end{bmatrix}
   = \frac{1}{D_n(\halpha_n^{\,\dag})}
       \begin{bmatrix}
         G_n \\
         H_n \\
         J_n \\
       \end{bmatrix} ,
 \end{align*}
 where
 \begin{align*}
   & G_n:= - (1-\halpha_n^{\,\dag})(1+(\halpha_n^{\,\dag})^2)\\
   &\phantom{\;G_n:=}
              \times\Big((1+(\halpha_n^{\,\dag})^2)(Y_{s_1}+Y_{s_2})
                    - \halpha_n^{\,\dag}(Y_{s_1-1}+Y_{s_1+1}+Y_{s_2-1}+Y_{s_2+1})\Big) \\
         &\phantom{ G_n:=\;} + (1+(\halpha_n^{\,\dag})^2)^2 \sum_{k=1}^n (Y_k - \halpha_n^{\,\dag} Y_{k-1}) ,\\
   & H_n:= \big( n(1+(\halpha_n^{\,\dag})^2)-(1-\halpha_n^{\,\dag})^2\big)
              \Big((1+(\halpha_n^{\,\dag})^2) Y_{s_1} - \halpha_n^{\,\dag}(Y_{s_1-1}+Y_{s_1+1}) \Big) \\
   &\phantom{H_n:=\;} + (1-\halpha_n^{\,\dag})^2
      \Big((1+(\halpha_n^{\,\dag})^2) Y_{s_2} - \halpha_n^{\,\dag}(Y_{s_2-1}+Y_{s_2+1})\Big)\\
   &\phantom{H_n:=\;}
           -(1-\halpha_n^{\,\dag})(1+(\halpha_n^{\,\dag})^2) \sum_{k=1}^n (Y_k - \halpha_n^{\,\dag} Y_{k-1}),\\
  & J_n:= (1-\halpha_n^{\,\dag})^2
        \Big((1+(\halpha_n^{\,\dag})^2) Y_{s_1} - \halpha_n^{\,\dag}(Y_{s_1-1}+Y_{s_1+1})\Big)\\
  &\phantom{H_n:=\;}
          + \big( n(1+(\halpha_n^{\,\dag})^2)-(1-\halpha_n^{\,\dag})^2\big)
     \Big((1+(\halpha_n^{\,\dag})^2) Y_{s_2} - \halpha_n^{\,\dag}(Y_{s_2-1}+Y_{s_2+1}) \Big) \\
  &\phantom{H_n:=\;}
      -(1-\halpha_n^{\,\dag})(1+(\halpha_n^{\,\dag})^2) \sum_{k=1}^n (Y_k - \halpha_n^{\,\dag} Y_{k-1})    .
 \end{align*}
Using \eqref{Ergodic1} and that for all \ $p_i\in\RR$, \ $i=0,\ldots,4$,
 \[
    \sup_{x\in\RR,\;n\geq 5}
        \frac{n(p_4x^4+p_3x^3+p_2x^2+p_1x+p_0)}{(1+x^2)((n-2)x^2+4x+n-2)}
       <\infty,
 \]
 one can think it over that \ $H_n/D_n(\halpha_n^{\,\dag})$, \ $n\in\NN$, \ and
 \ $J_n/D_n(\halpha_n^{\,\dag})$, \ $n\in\NN$, \ are bounded with probability one, which yields also that
 the sequences \ $(\htheta_{1,n}^{\,\dag}-\theta_1)_{n\in\NN}$
 \ and \ $(\htheta_{2,n}^{\,\dag}-\theta_2)_{n\in\NN}$ \ are bounded
 with probability one.

Again by Lemma \ref{LEMMA13} and equations \eqref{Additive_CLSE_EQ7} we get that
 \begin{align*}
    \begin{bmatrix}
      \halpha_n^{\,\dag} \\
      \hmuen^{\,\dag} \\
     \end{bmatrix}
    = \begin{bmatrix}
         a_n & b_n \\
         b_n &  n \\
        \end{bmatrix}^{-1}
       \begin{bmatrix}
           c_n \\
           d_n \\
       \end{bmatrix}
 \end{align*}
 holds asymptotically as \ $n\to\infty$ \ with probability one, where
 \begin{align*}
  & a_n:= \sum_{k=1}^n X_{k-1}^2 + (\theta_1-\htheta_{1,n}^{\,\dag})(\theta_1-\htheta_{1,n}^{\,\dag}+2X_{s_1})
                               + (\theta_2-\htheta_{2,n}^{\,\dag})(\theta_2-\htheta_{2,n}^{\,\dag}+2X_{s_2}),\\
  & b_n:= \sum_{k=1}^n X_{k-1} + \theta_1-\htheta_{1,n}^{\,\dag} + \theta_2-\htheta_{2,n}^{\,\dag},\\
  & c_n:= \sum_{k=1}^n X_{k-1}X_k + (\theta_1-\htheta_{1,n}^{\,\dag})(X_{s_1-1}+X_{s_1+1})
                                 + (\theta_2-\htheta_{2,n}^{\,\dag})(X_{s_2-1}+X_{s_2+1}),\\
  & d_n:= \sum_{k=1}^n X_k + \theta_1-\htheta_{1,n}^{\,\dag}
                                 + \theta_2-\htheta_{2,n}^{\,\dag}.
 \end{align*}
Here we emphasize that the matrix
 \[
   \begin{bmatrix}
     a_n & b_n \\
     b_n &  n \\
   \end{bmatrix}
 \]
 is invertible asymptotically as \ $n\to\infty$ \ with probability one, since
 using \eqref{Ergodic1}, \eqref{Ergodic2} and that the sequences
 \ $(\htheta_{1,n}^{\,\dag}-\theta_1)_{n\in\NN}$ \ and
 \ $(\htheta_{2,n}^{\,\dag}-\theta_2)_{n\in\NN}$ \ are bounded with probability one
 we get
 \begin{align}\label{SEGED96_new}
   \PP\left(\lim_{n\to\infty}\frac{a_n}{n}=\EE\widetilde X^2\right)=1,
     \qquad\qquad \PP\left(\lim_{n\to\infty}\frac{b_n}{n}=\EE\widetilde X\right)=1,
 \end{align}
 and hence
 \[
    \PP\left(\lim_{n\to\infty}\frac{1}{n^2}(na_n - b_n^2)
                = \EE\widetilde X^2 - (\EE\widetilde X)^2
                =\var\widetilde X\right)=1.
 \]
This yields that
 \[
    \PP\left(\lim_{n\to\infty}(na_n - b_n^2)
                = \infty \right)=1.
 \]

Further
 \begin{align}\label{SEGED55}
   \begin{bmatrix}
     \halpha_n^{\,\dag} - \alpha \\
     \hmuen^{\,\dag} - \mu_\vare\\
   \end{bmatrix}
    = \begin{bmatrix}
         a_n & b_n \\
         b_n &  n \\
        \end{bmatrix}^{-1}
       \begin{bmatrix}
          e_n \\
         f_n \\
       \end{bmatrix}
 \end{align}
 holds asymptotically as \ $n\to\infty$ \ with probability one, where
 \begin{align*}
  &e_n:=
     \sum_{k=1}^n \!X_{k-1}(X_k-\alpha X_{k-1}-\mu_\vare)\\
  &\phantom{e_n:=\;}
              + (\theta_1-\htheta_{1,n}^{\,\dag})
                \big(X_{s_1-1}+X_{s_1+1}-2\alpha X_{s_1}
                  - \mu_\vare - \alpha(\theta_1-\htheta_{1,n}^{\,\dag})\!\big)\\
  &\phantom{e_n:=\;}
              + (\theta_2-\htheta_{2,n}^{\,\dag})
               \big(X_{s_2-1}+X_{s_2+1}-2\alpha X_{s_2}
               - \mu_\vare - \alpha(\theta_2-\htheta_{2,n}^{\,\dag})\big),\\
  &f_n:= \sum_{k=1}^n (X_k -\alpha X_{k-1}-\mu_\vare)
              + (1-\alpha)(\theta_1-\htheta_{1,n}^{\,\dag}
                                 + \theta_2-\htheta_{2,n}^{\,\dag}).
 \end{align*}

Then, using again \eqref{STAC_MOMENT1}, \eqref{Ergodic1}, \eqref{Ergodic2}, \eqref{Ergodic3}
 and that the sequences \ $(\htheta_{1,n}^{\,\dag}-\theta_1)_{n\in\NN}$ \ and
 \ $(\htheta_{2,n}^{\,\dag}-\theta_2)_{n\in\NN}$ \ are bounded with probability one,
 we get
 \begin{align*}
   &\PP\left(\lim_{n\to\infty}\frac{e_n}{n}
      =\alpha\EE\widetilde X^2+\mu_\vare\EE\widetilde X
       -\alpha\EE\widetilde X^2-\mu_\vare\EE\widetilde X =0 \right)=1,\\[1mm]
   &\PP\left(\lim_{n\to\infty}\frac{f_n}{n}
      =\EE\widetilde X - \alpha\EE\widetilde X-\mu_\vare = 0\right)=1.
 \end{align*}
Hence, by \eqref{SEGED55}, we obtain
 \begin{align*}
   \PP\left(
     \lim_{n\to\infty}
          \begin{bmatrix}
            \halpha_n^{\,\dag} - \alpha \\
            \hmuen^{\,\dag} - \mu_\vare\\
           \end{bmatrix}
       = \begin{bmatrix}
           \EE\widetilde X^2 & \EE\widetilde X \\
           \EE\widetilde X &  1 \\
         \end{bmatrix}^{-1}
          \begin{bmatrix}
                  0 \\
                  0 \\
          \end{bmatrix}
      =  \begin{bmatrix}
            0 \\
            0 \\
          \end{bmatrix}
   \right)=1,
 \end{align*}
 which yields \eqref{Strong_consistency18} and \eqref{Strong_consistency19}.
Then \eqref{Strong_consistency18}, \eqref{Strong_consistency19} and \eqref{SEGED56}
 imply \eqref{Strong_consistency20}. 
\proofend

The asymptotic distribution of the CLS estimation is given in the next theorem.

\begin{Thm}
Under the additional assumptions \ $\EE X_0^3<\infty$ \ and \ $\EE\vare_1^3<\infty$, \ we have
 \begin{align}\label{CONVERGENCE14}
   \begin{bmatrix}
     \sqrt{n}(\halpha_n^{\,\dag}-\alpha) \\
     \sqrt{n}(\hmuen^{\,\dag}-\mu_\vare) \\
   \end{bmatrix}
   \distr \cN\left(\begin{bmatrix}
                     0 \\
                     0 \\
                   \end{bmatrix}
     ,B_{\alpha,\,\vare}\right)
      \qquad \text{as \ $n\to\infty$,}
  \end{align}
 where \ $B_{\alpha,\vare}$ \ is defined in \eqref{SEGED_BALPHA}.
Moreover, conditionally on the values \ $Y_{s_1-1}$, \ $Y_{s_2-1}$ \ and \ $Y_{s_1+1}$,
 \ $Y_{s_2+1}$,
 \begin{align}\label{CONVERGENCE15}
    \begin{bmatrix}
      \sqrt{n}\big(\htheta_{1,n}^{\,\dag} - \lim_{k\to\infty}\htheta_{1,k}^{\,\dag}\big) \\
      \sqrt{n}\big(\htheta_{2,n}^{\,\dag} - \lim_{k\to\infty}\htheta_{2,k}^{\,\dag}\big) \\
    \end{bmatrix}
      \distr \cN\left(\begin{bmatrix}
                        0 \\
                        0 \\
                      \end{bmatrix},
                      C_{\alpha,\vare} B_{\alpha,\vare} C_{\alpha,\vare}^\top
      \right)
      \qquad \text{as \ $n\to\infty$,}
 \end{align}
 where
 \[
   C_{\alpha,\vare}
      := \frac{1}{(1+\alpha^2)^2}
         \begin{bmatrix}
           (\alpha^2-1)(Y_{s_1-1}+Y_{s_1+1}) + (1+2\alpha-\alpha^2)\mu_\vare & (\alpha-1)(1+\alpha^2)\\
           (\alpha^2-1)(Y_{s_2-1}+Y_{s_2+1}) + (1+2\alpha-\alpha^2)\mu_\vare & (\alpha-1)(1+\alpha^2) \\
         \end{bmatrix}.
 \]
\end{Thm}

\noindent{\bf Proof.}
By \eqref{SEGED_BALPHA} and \eqref{SEGED55}, to prove \eqref{CONVERGENCE14} it is enough to show that
 \begin{align}\label{SEGED41_new}
   &\PP\left(\lim_{n\to\infty}\frac{1}{n}
              \begin{bmatrix}
                  a_n & b_n \\
                  b_n & n \\
                \end{bmatrix}
              = \begin{bmatrix}
                  \EE\widetilde X^2 & \EE\widetilde X \\
                  \EE\widetilde X & 1 \\
                \end{bmatrix}
              \right)=1,\\[2mm] \label{SEGED42_new}
   &\frac{1}{\sqrt n}
          \begin{bmatrix}
            e_n \\
            f_n \\
          \end{bmatrix}
         \distr \cN\left(
           \begin{bmatrix}
             0 \\
             0 \\
          \end{bmatrix}\!,\,A_{\alpha,\vare}\right)
         \qquad \text{as \ $n\to\infty$, }
 \end{align}
 where \ $\widetilde X$ \ is a random variable having the unique stationary distribution
 of the INAR(1) model in \eqref{INAR1} and the \ $(2\times 2)$-matrix \ $A_{\alpha,\vare}$
 \ is defined in \eqref{SEGED_AALPHA}.
By \eqref{SEGED96_new}, we have \eqref{SEGED41_new}.
By formula (6.43) in Hall and Heyde \cite[Section 6.3]{HalHey},
 \begin{align*}
   \begin{bmatrix}
     \frac{1}{\sqrt n} \sum_{k=1}^n X_{k-1}(X_k-\alpha X_{k-1}-\mu_\vare)  \\
     \frac{1}{\sqrt n} \sum_{k=1}^n (X_k-\alpha X_{k-1}-\mu_\vare) \\
   \end{bmatrix}
         \distr
   \cN\left(\begin{bmatrix}
              0 \\
              0 \\
            \end{bmatrix}
             \!,\,A_{\alpha,\vare}\right)
     \qquad \text{as \ $n\to\infty$.}
 \end{align*}
Hence using that the sequences \ $(\htheta_{1,n}^{\,\dag}-\theta_1)_{n\in\NN}$
 \ and \ $(\htheta_{2,n}^{\,\dag}-\theta_2)_{n\in\NN}$ \ are bounded with probability one,
 by Slutsky's lemma (see, e.g., Lemma 2.8 in van der Vaart \cite{Vaart}), we get \eqref{SEGED42_new}.

Now we turn to prove \eqref{CONVERGENCE15}.
Using the notation
 \[
   B_n^{\dag}:=
          \begin{bmatrix}
             1+(\halpha_n^{\,\dag})^2 & 0 \\
             0 & 1+(\halpha_n^{\,\dag})^2 \\
           \end{bmatrix},
 \]
 by \eqref{SEGED56}, we have 
  \begin{align*}
    \begin{bmatrix}
      \htheta_{1,n}^{\,\dag} \smallskip  \\
      \htheta_{2,n}^{\,\dag} \\
    \end{bmatrix}
     = (B_n^{\dag})^{-1}
        \begin{bmatrix}
         (1+(\halpha_n^{\,\dag})^2)Y_{s_1} - \halpha_n^{\,\dag}(Y_{s_1-1}+Y_{s_1+1})
          - (1-\halpha_n^{\,\dag})\hmuen^{\,\dag} \smallskip  \\
         (1+(\halpha_n^{\,\dag})^2)Y_{s_2} - \halpha_n^{\,\dag}(Y_{s_2-1}+Y_{s_2+1})
            - (1-\halpha_n^{\,\dag})\hmuen^{\,\dag} \\
        \end{bmatrix}
  \end{align*}
 holds asymptotically as \ $n\to\infty$ \ with probability one.
Theorem \ref{THEOREM4} yields that
 \[
     \PP\left(\lim_{n\to\infty}B_n^{\dag}
                 = \begin{bmatrix}
                    1+\alpha^2 & 0 \\
                    0 & 1+\alpha^2 \\
                   \end{bmatrix}=:B^{\dag}
         \right)=1.
 \]
By \eqref{Strong_consistency20}, we have 
 \begin{align*}
   &\begin{bmatrix}
      \sqrt{n}\big(\htheta_{1,n}^{\,\dag} - \lim_{k\to\infty}\htheta_{1,k}^{\,\dag}\big) \smallskip \\
      \sqrt{n}\big(\htheta_{2,n}^{\,\dag} - \lim_{k\to\infty}\htheta_{2,k}^{\,\dag}\big) \\
   \end{bmatrix}\\
   &\qquad =\sqrt{n}(B_n^{\dag})^{-1}
    \left(
      \begin{bmatrix}
         (1+(\halpha_n^{\,\dag})^2)Y_{s_1} - \halpha_n^{\,\dag}(Y_{s_1-1}+Y_{s_1+1})
             - (1-\halpha_n^{\,\dag})\hmuen^{\,\dag} \smallskip  \\
         (1+(\halpha_n^{\,\dag})^2)Y_{s_2} - \halpha_n^{\,\dag}(Y_{s_2-1}+Y_{s_2+1})
             - (1-\halpha_n^{\,\dag})\hmuen^{\,\dag} \\
      \end{bmatrix}\right.\\
   &\phantom{\qquad = \sqrt{n}(B_n^{\dag})^{-1}\Big(\;\;}\left.
     - B_n^{\dag}(B^{\dag})^{-1}
      \begin{bmatrix}
         (1+\alpha^2)Y_{s_1} - \alpha(Y_{s_1-1}+Y_{s_1+1}) - (1-\alpha)\mu_\vare \smallskip  \\
         (1+\alpha^2)Y_{s_2} - \alpha(Y_{s_2-1}+Y_{s_2+1}) - (1-\alpha)\mu_\vare \\
      \end{bmatrix}
     \right)\\
  &\qquad\;=\sqrt{n}(B_n^{\dag})^{-1}
      \left(
       \begin{bmatrix}
         (1+(\halpha_n^{\,\dag})^2)Y_{s_1} - \halpha_n^{\,\dag}(Y_{s_1-1}+Y_{s_1+1})
            - (1-\halpha_n^{\,\dag})\hmuen^{\,\dag} \smallskip  \\
         (1+(\halpha_n^{\,\dag})^2)Y_{s_2} - \halpha_n^{\,\dag}(Y_{s_2-1}+Y_{s_2+1})
            - (1-\halpha_n^{\,\dag})\hmuen^{\,\dag} \\
       \end{bmatrix}\right.\\
  &\phantom{\qquad=\sqrt{n}(B_n^{\dag})^{-1}\Big(\;}\left.
    -\begin{bmatrix}
         (1+\alpha^2)Y_{s_1} - \alpha(Y_{s_1-1}+Y_{s_1+1}) - (1-\alpha)\mu_\vare \smallskip \\
         (1+\alpha^2)Y_{s_2} - \alpha(Y_{s_2-1}+Y_{s_2+1}) - (1-\alpha)\mu_\vare \\
      \end{bmatrix}
       \right)\\
 &\phantom{\qquad\;=\;}
      +\sqrt{n}
      \left((B_n^{\dag})^{-1} -  (B^{\dag})^{-1}
       \right)
       \begin{bmatrix}
         (1+\alpha^2)Y_{s_1} - \alpha(Y_{s_1-1}+Y_{s_1+1}) - (1-\alpha)\mu_\vare \smallskip \\
         (1+\alpha^2)Y_{s_2} - \alpha(Y_{s_2-1}+Y_{s_2+1}) - (1-\alpha)\mu_\vare \\
      \end{bmatrix}\\[1mm]
  &\qquad=\sqrt{n}(B_n^{\dag})^{-1}
    \begin{bmatrix}
      (\halpha_n^{\,\dag}+\alpha)Y_{s_1}
         - (Y_{s_1-1}+ Y_{s_1+1})
         + \hmuen^{\,\dag}
        \hspace{1mm} & \alpha-1 \smallskip \\
       (\halpha_n^{\,\dag}+\alpha)Y_{s_2}
         - (Y_{s_2-1}+ Y_{s_2+1})
         + \hmuen^{\,\dag}
        \hspace{1mm} & \alpha-1 \\
    \end{bmatrix}
    \begin{bmatrix}
      \halpha_n^{\,\dag}-\alpha \smallskip \\
      \hmuen^{\,\dag}-\mu_\vare \\
    \end{bmatrix} \\
 &\phantom{\qquad=\;}
    +\sqrt{n}(B_n^{\dag})^{-1}
       \big( B^{\dag}- B_n^{\dag} \big) (B^{\dag})^{-1}
        \begin{bmatrix}
         (1+\alpha^2)Y_{s_1} - \alpha(Y_{s_1-1}+Y_{s_1+1}) - (1-\alpha)\mu_\vare \smallskip \\
         (1+\alpha^2)Y_{s_2} - \alpha(Y_{s_2-1}+Y_{s_2+1}) - (1-\alpha)\mu_\vare \\
         \end{bmatrix}.
 \end{align*}
Then
 \begin{align}\label{SEGED58}
   &\begin{bmatrix}
      \sqrt{n}\big(\htheta_{1,n}^{\,\dag} - \lim_{k\to\infty}\htheta_{1,k}^{\,\dag}\big) \smallskip \\
      \sqrt{n}\big(\htheta_{2,n}^{\,\dag} - \lim_{k\to\infty}\htheta_{2,k}^{\,\dag}\big) \\
    \end{bmatrix}
   =  C_{n,\alpha,\vare}
       \begin{bmatrix}
        \sqrt{n}(\halpha_n^{\,\dag}-\alpha)  \smallskip \\
        \sqrt{n}(\hmuen^{\,\dag}-\mu_\vare)  \\
       \end{bmatrix}
  \end{align}
 holds asymptotically as \ $n\to\infty$ \ with probability one, where \ $C_{n,\alpha,\vare}$
 \ is defined by
 \begin{align*}
    (B_n^{\dag})^{-1}
      &\begin{bmatrix}
      (\halpha_n^{\,\dag}+\alpha)Y_{s_1}
         - Y_{s_1-1} -Y_{s_1+1}
         + \hmuen^{\,\dag}
         \hspace{1mm} & \alpha-1 \smallskip \\
       (\halpha_n^{\,\dag}+\alpha)Y_{s_2}
         - Y_{s_2-1} - Y_{s_2+1}
         + \hmuen^{\,\dag}
        \hspace{1mm}  & \alpha-1 \\
    \end{bmatrix}\\
    &-(\halpha_n^{\,\dag}+\alpha)(B_n^{\dag})^{-1}
        (B^{\dag})^{-1}
        \begin{bmatrix}
         (1+\alpha^2)Y_{s_1} - \alpha(Y_{s_1-1}+Y_{s_1+1}) - (1-\alpha)\mu_\vare \hspace{1mm} & 0 \smallskip \\
         (1+\alpha^2)Y_{s_2} - \alpha(Y_{s_2-1}+Y_{s_2+1}) - (1-\alpha)\mu_\vare \hspace{1mm} & 0 \\
        \end{bmatrix}.
  \end{align*}
By \eqref{Strong_consistency18} and \eqref{Strong_consistency19}, we have \ $C_{n,\alpha,\vare}$
 \ converges almost surely as \ $n\to\infty$ \ to
 \begin{align*}
   & (B^{\dag})^{-1}
     \begin{bmatrix}
       2\alpha Y_{s_1} - Y_{s_1-1} - Y_{s_1+1} +\mu_\vare \hspace{1mm} & \alpha-1 \smallskip \\
       2\alpha Y_{s_2} - Y_{s_2-1} - Y_{s_2+1} +\mu_\vare \hspace{1mm} & \alpha-1  \\
     \end{bmatrix} \\
   & +(B^{\dag})^{-1}
    \begin{bmatrix}
      -2\alpha & 0 \\
      0 & -2\alpha \\
    \end{bmatrix}
    (B^{\dag})^{-1}
      \begin{bmatrix}
        (1+\alpha^2)Y_{s_1} - \alpha(Y_{s_1-1}+Y_{s_1+1}) - (1-\alpha)\mu_\vare \hspace{1mm} & 0 \smallskip \\
        (1+\alpha^2)Y_{s_2} - \alpha(Y_{s_2-1}+Y_{s_2+1}) - (1-\alpha)\mu_\vare \hspace{1mm} & 0 \\
      \end{bmatrix} \\
   & = \frac{1}{(1+\alpha^2)^2}
       \begin{bmatrix}
         (\alpha^2-1)(Y_{s_1-1}+Y_{s_1+1}) + (1+2\alpha -\alpha^2)\mu_\vare \hspace{1mm}
              & (\alpha-1)(1+\alpha^2) \smallskip \\
         (\alpha^2-1)(Y_{s_2-1}+Y_{s_2+1}) + (1+2\alpha -\alpha^2)\mu_\vare \hspace{1mm}
              & (\alpha-1)(1+\alpha^2) \\
       \end{bmatrix}\\
   & =C_{\alpha,\vare}.
 \end{align*}
By \eqref{SEGED58}, \eqref{CONVERGENCE14} and Slutsky's lemma, we have \eqref{CONVERGENCE15}.
\proofend

\subsection{Two neighbouring outliers, estimation of the mean of the offspring and
            innovation distributions and the outliers' sizes}

In this section we assume that \ $I=2$ \ and that the relevant time points
 \ $s_1$, $s_2\in\NN$ \ are known.
We also suppose that \ $s_1:=s$ \ and \ $s_2:=s+1$, \ i.e., the time points
 \ $s_1$ \ and \ $s_2$ \ are neighbouring.
We concentrate on the CLS estimation of \ $\alpha$, \ $\mu_\vare$, \ $\theta_1$
 \ and \ $\theta_2$.

Motivated by \eqref{SEGED95}, for all \ $n\geq s+2$, \ $n\in\NN$, \ we define the function
 \ $Q_n^{\dag\dag}:\RR^{n+1}\times\RR^4\to\RR$, \ as
 \begin{align*}
    &Q_n^{\dag\dag}({\bf y}_n;\alpha',\mu_\vare',\theta_1',\theta_2')\\
      & :=\DS\sum_{\substack{k=1 \\ k\not\in \{s,s+1,s+2\}}}^n
            \big(y_k-\alpha' y_{k-1}-\mu_\vare'\big)^2
           + \big(y_{s}-\alpha' y_{s-1}-\mu_\vare'-\theta_1'\big)^2
           + \big(y_{s+1}-\alpha' y_{s}-\mu_\vare'+\alpha'\theta_1'-\theta_2'\big)^2 \\
         &\phantom{=\;}
           + \big(y_{s+2}-\alpha' y_{s+1}-\mu_\vare'+\alpha'\theta_2'\big)^2,
       \qquad {\bf y}_n\in\RR^{n+1},\;\alpha',\mu_\vare',\theta_1',\theta_2'\in\RR.
 \end{align*}
By definition, for all \ $n\geq s+2$, \ a CLS estimator for
 the parameter \ $(\alpha,\mu_\vare,\theta_1,\theta_2)\in(0,1)\times(0,\infty)\times\NN^2$ \ is
 a measurable function
 \ $(\halpha_n^{\,\dag\dag},\hmuen^{\,\dag\dag},\htheta_{1,n}^{\,\dag\dag},\htheta_{2,n}^{\,\dag\dag}):
                                 S_n\to\RR^4$
 \ such that
 \begin{align*}
   Q_n^{\dag\dag}({\bf y}_n;\,\halpha_n^{\,\dag\dag}({\bf y}_n),
    &\hmuen^{\,\dag\dag}({\bf y}_n),
           \htheta_{1,n}^{\,\dag\dag}({\bf y}_n),
             \htheta_{2,n}^{\,\dag\dag}({\bf y}_n)) \\
   &=\inf_{(\alpha',\mu_\vare',\theta_1',\theta_2')\in\RR^4}
       Q_n^{\dag\dag}({\bf y}_n;\alpha',\mu_\vare',\theta_1',\theta_2')
         \qquad \forall\;\;  {\bf y}_n\in S_n,
 \end{align*}
 where \ $S_n$ \ is suitable subset of \ $\RR^{n+1}$ \ (defined in the proof of
 Lemma \ref{LEMMA14}).
We note that we do not define the CLS estimator
 \ $(\halpha_n^{\,\dag\dag},\hmuen^{\,\dag\dag},\htheta_{1,n}^{\,\dag\dag},\htheta_{2,n}^{\,\dag\dag})$ \
 for all samples \ ${\bf y}_n\in \RR^{n+1}$.
\ For all \ ${\bf y}_n\in\RR^{n+1}$ \ and \ $(\alpha',\mu_\vare',\theta_1',\theta_2')\in\RR^4$,
 \begin{align*}
    &\frac{\partial Q_n^{\dag\dag}}{\partial \alpha'}({\bf y}_n;\alpha',\mu_\vare',\theta_1',\theta_2')\\
    &\phantom{\qquad}
       = \DS\sum_{\substack{k=1 \\ k\not\in \{s,s+1,s+2\}}}^n
            \big(y_k-\alpha' y_{k-1}-\mu_\vare'\big)(-2y_{k-1})
           -2\big(y_{s}-\alpha' y_{s-1}-\mu_\vare'-\theta_1'\big)y_{s-1} \\
   &\phantom{\qquad = \;\;}
      + 2\big(y_{s+1}-\alpha' y_{s}-\mu_\vare'+\alpha'\theta_1'-\theta_2'\big)(-y_{s}+\theta_1')
      + 2\big(y_{s+2}-\alpha' y_{s+1}-\mu_\vare'+\alpha'\theta_2'\big)(-y_{s+1}+\theta_2'),
 \end{align*}
 and
 \begin{align*}
   &\frac{\partial Q_n^{\dag\dag}}{\partial \mu_\vare'}({\bf y}_n;\alpha',\mu_\vare',\theta_1',\theta_2')\\
   &\phantom{\qquad }
     = \DS\sum_{\substack{k=1 \\ k\not\in \{s,s+1,s+2\}}}^n
            (-2)\big(y_k-\alpha' y_{k-1}-\mu_\vare'\big)
           -2\big(y_{s}-\alpha' y_{s-1}-\mu_\vare'-\theta_1'\big) \\
   &\phantom{\qquad = \;\;}
      - 2\big(y_{s+1}-\alpha' y_{s}-\mu_\vare'+\alpha'\theta_1'-\theta_2'\big)
       - 2\big(y_{s+2}-\alpha' y_{s+1}-\mu_\vare'+\alpha'\theta_2'\big), \\
   & \frac{\partial Q_n^{\dag\dag}}{\partial \theta_1'}({\bf y}_n;\alpha',\mu_\vare',\theta_1',\theta_2') \\
   &\qquad\qquad
      = -2(y_{s}-\alpha' y_{s-1}-\mu_\vare'-\theta_1')
         +2\alpha'(y_{s+1}-\alpha' y_{s}-\mu_\vare'+\alpha'\theta_1'-\theta_2'),\\
   & \frac{\partial Q_n^{\dag\dag}}{\partial \theta_2'}({\bf y}_n;\alpha',\mu_\vare',\theta_1',\theta_2') \\
   &\qquad\qquad
        = -2(y_{s+1}-\alpha' y_{s}-\mu_\vare'+\alpha'\theta_1'-\theta_2')
         +2\alpha'(y_{s+2}-\alpha' y_{s+1}-\mu_\vare'+\alpha'\theta_2').
 \end{align*}

The next lemma is about the existence and uniqueness of the CLS estimator of
 \ $(\alpha,\mu_\vare,\theta_1,\theta_2)$.

\begin{Lem}\label{LEMMA14}
There exist subsets \ $S_n\subset\RR^{n+1}$, $n\geq \max(3,s+2)$ \ with the following properties:
 \begin{enumerate}
  \item[\upshape{(i)}]
   there exists a unique CLS estimator
   \ $(\halpha_n^{\,\dag\dag},\hmuen^{\,\dag\dag},\htheta_{1,n}^{\,\dag\dag},
       \htheta_{2,n}^{\,\dag\dag}):S_n\to\RR^4$,
  \item[\upshape{(ii)}]
   for all \ ${\bf y}_n\in S_n$,
   \ $(\halpha_n^{\,\dag\dag}({\bf y}_n),\hmuen^{\,\dag\dag}({\bf y}_n),
       \htheta_{1,n}^{\,\dag\dag}({\bf y}_n),\htheta_{2,n}^{\,\dag\dag}({\bf y}_n))$
   \ is the unique solution of the system of equations
  \begin{align}\label{Additive_CLSE_EQ8}
   \begin{split}
    & \frac{\partial Q_n^{\dag\dag}}{\partial \alpha'}
     ({\bf y}_n;\alpha',\mu_\vare',\theta_1',\theta_2')=0,\qquad
    \frac{\partial Q_n^{\dag\dag}}{\partial \mu_\vare'}
     ({\bf y}_n;\alpha',\mu_\vare',\theta_1',\theta_2')=0,\\
   & \frac{\partial Q_n^{\dag\dag}}{\partial \theta_1'}
     ({\bf y}_n;\alpha',\mu_\vare',\theta_1',\theta_2')=0,\qquad
    \frac{\partial Q_n^{\dag\dag}}{\partial \theta_2'}
      ({\bf y}_n;\alpha',\mu_\vare',\theta_1',\theta_2')=0,
    \end{split}
   \end{align}
  \item [\upshape{(iii)}]
  ${\bf Y}_n\in S_n$ \ holds asymptotically as \ $n\to\infty$ \ with probability one.
 \end{enumerate}
\end{Lem}

\noindent{\bf Proof.}
For any fixed \ ${\bf y}_n \in \RR^{n+1}$, \ $n\geq \max(3,s+2)$ \ and \ $\alpha' \in \RR$,
 \ the quadratic function
 \ $\RR^3 \ni (\mu_\vare',\theta_1',\theta_2')
      \mapsto Q_n^{\dag\dag}({\bf y}_n;\alpha',\mu_\vare',\theta_1',\theta_2')$
 \ can be written in the form
 \begin{align*}
  &Q_n^{\dag\dag}({\bf y}_n;\alpha',\mu_\vare',\theta_1',\theta_2') \\[2mm]
  &=  \!\!\left(\! \begin{bmatrix}
             \mu_\vare' \smallskip \\
             \theta_1' \smallskip \\
             \theta_2'
            \end{bmatrix}
           - A_n(\alpha')^{-1} t_n({\bf y}_n;\alpha') \!\! \right)^{\hspace*{-1mm}\top}
      \!\!\! A_n(\alpha') \! \left( \! \begin{bmatrix}
                          \mu_\vare' \smallskip \\
                          \theta_1' \smallskip \\
                          \theta_2'
                         \end{bmatrix}
                         - A_n(\alpha')^{-1} t_n({\bf y}_n;\alpha') \! \right)
    \!\! + \widehat{Q}_n^{\dag\dag}({\bf y}_n;\alpha'),
 \end{align*}
 where
 \begin{align*}
  t_n({\bf y}_n;\alpha')
  &:= \begin{bmatrix}
       \sum_{k=1}^n ( y_k - \alpha' y_{k-1}) \smallskip \\
       ( 1 + (\alpha')^2 ) y_s - \alpha' ( y_{s-1} + y_{s+1} ) \smallskip \\
       ( 1 + (\alpha')^2 ) y_{s+1} - \alpha' ( y_s + y_{s+2} )
      \end{bmatrix} , \\[2mm]
  \widehat{Q}_n^{\dag\dag}({\bf y}_n;\alpha')
  &:= \sum_{k=1}^n \big( y_k-\alpha' y_{k-1} \big)^2
      - t_n({\bf y}_n;\alpha')^\top A_n(\alpha')^{-1} t_n({\bf y}_n;\alpha') ,
 \end{align*}
 and the matrix
 \begin{align*}
  A_n(\alpha')
  := \begin{bmatrix}
      n & 1-\alpha' & 1-\alpha' \smallskip \\
      1-\alpha' & 1+(\alpha')^2 &  -\alpha' \smallskip \\
      1-\alpha' & -\alpha' &  1+(\alpha')^2 \\
     \end{bmatrix}
 \end{align*}
 is strictly positive definite for all \ $n\geq3$ \ and \ $\alpha' \in \RR$.
\ Indeed, the leading principal minors of \ $A_n(\alpha')$ \ take the following forms:
 \ $n,$
 \begin{align*}
  & n(1+(\alpha')^2) - (1-\alpha')^2 = (n-1)(\alpha')^2 + 2\alpha' + n-1,\\
  & D_n(\alpha')
   := n(1 + (\alpha')^2)^2 - 2(1-\alpha')^2(1+(\alpha')^2)
       - 2\alpha'(1-\alpha')^2-n(\alpha')^2\\
  &\phantom{D_n(\alpha') :=}
    =   n(1+\alpha'+(\alpha')^2)(1-\alpha'+(\alpha')^2)
         - 2(1-\alpha')^2(1+\alpha'+(\alpha')^2)\\
  &\phantom{D_n(\alpha') :=}
   = (1+\alpha'+(\alpha')^2)
          \big( (n-2)(\alpha')^2 - (n-4)\alpha' + n-2 \big),
 \end{align*}
 and for all \ $n\geq 3$, \ the discriminant \ $(n-4)^2-4(n-2)^2=-3n^2+8n$ \ of the equation
  \ $(n-2)x^2-(n-4)x+n-2=0$ \ is negative.
The inverse matrix \ $A_n(\alpha')^{-1}$ \ takes the form
 \begin{align*}
  \frac{1}{D_n(\alpha')}\!\!
   \begin{bmatrix}
    1+(\alpha')^2+(\alpha')^4 \, & \, -(1-\alpha')(1+\alpha'+(\alpha')^2)
     \, & \, -(1-\alpha')(1+\alpha'+(\alpha')^2)  \\
    -(1-\alpha')(1+\alpha'+(\alpha')^2) \, & \, n(1+(\alpha')^2)-(1-\alpha')^2
     \, & \, (1-\alpha')^2 + n\alpha' \!\! \\
    -(1-\alpha')(1+\alpha'+(\alpha')^2) \, & \, (1-\alpha')^2 + n\alpha' \, & \, n(1+(\alpha')^2)-(1-\alpha')^2
  \end{bmatrix}.
 \end{align*}
The polynomial \ $\RR\ni\alpha' \mapsto D_n(\alpha')$ \ is of order 4 with leading coefficient \ $n-2$.
\ We have \ $\widehat{Q}_n^{\dag\dag}({\bf y}_n;\alpha') = R_n({\bf y}_n;\alpha') / D_n(\alpha')$, \ where
 \ $\RR\ni\alpha' \mapsto R_n({\bf y}_n;\alpha')$ \ is a polynomial of order 6 with leading coefficient
 \begin{align*}
   c_n({\bf y}_n) &:= (n-2) \sum_{k=1}^n y_{k-1}^2 - \left( \sum_{k=1}^n y_{k-1} \right)^2
                      - (n-1)(y_s^2 + y_{s+1}^2) \\[2mm]
                  &\;\quad + 2 (y_s + y_{s+1}) \sum_{k=1}^n y_{k-1} - 2 y_s y_{s+1}.
 \end{align*}
Let
 \[
   \widehat{S}^{\,\dag\dag}_n := \left\{{\bf y}_n\in\RR^{n+1} : c_n({\bf y}_n) > 0 \right\}.
 \]
For \ ${\bf y}_n \in \widehat{S}^{\,\dag\dag}_n$, \ we have
 \ $\lim_{|\alpha'|\to\infty} \widehat{Q}_n^{\,\dag\dag}({\bf y}_n;\alpha') = \infty$ \
  and the continuous function
 \ $\RR \ni \alpha' \mapsto \widehat{Q}_n^{\,\dag\dag}({\bf y}_n;\alpha')$ \ attains its infimum.
Consequently, for all \ $n\geq\max(3,s+2)$ \ there exists a CLS estimator
 \ $(\halpha_n^{\,\dag\dag},
    \hmuen^{\,\dag\dag},
    \htheta_{1,n}^{\,\dag\dag},
    \htheta_{2,n}^{\,\dag\dag}): \widehat{S}^{\,\dag\dag}_n\to\RR^4$,
\ where
 \begin{align}\nonumber
  \widehat{Q}_n^{\,\dag\dag}({\bf y}_n;\halpha_n^{\,\dag\dag}({\bf y}_n))
   & =\inf_{\alpha'\in\RR} \widehat{Q}_n^{\,\dag\dag}({\bf y}_n;\alpha')
   \qquad \forall\;{\bf y}_n\in \widehat{S}^{\,\dag\dag}_n ,\\[2mm]\label{SEGED_UJ11}
  \begin{bmatrix}
   \hmuen^{\,\dag\dag}({\bf y}_n) \smallskip \\
   \htheta_{1,n}^{\,\dag\dag}({\bf y}_n) \smallskip \\
   \htheta_{2,n}^{\,\dag\dag}({\bf y}_n)
  \end{bmatrix}
  & = A_n(\halpha_n^{\,\dag\dag}({\bf y}_n))^{-1} t_n({\bf y}_n;\halpha_n^{\,\dag\dag}({\bf y}_n)) ,
  \qquad {\bf y}_n\in \widehat{S}^{\,\dag\dag}_n,
 \end{align}
 and for all \ ${\bf y}_n\in \widehat{S}^{\dag\dag}_n$,
 \ $(\halpha_n^{\,\dag\dag}({\bf y}_n),\hmuen^{\,\dag\dag}({\bf y}_n),\htheta_{1,n}^{\,\dag\dag}({\bf y}_n),
       \htheta_{2,n}^{\,\dag\dag}({\bf y}_n))$
  \ is a  solution of the system of equations \eqref{Additive_CLSE_EQ8}.

By \eqref{Ergodic1} and \eqref{Ergodic2}, we get
 \ $\PP\left(\lim_{n\to\infty} n^{-2} c_n({\bf Y}_n) = \var\widetilde X \right)=1$, \ where \ $\widetilde X$
 \ denotes a random variable with the unique stationary distribution of the INAR(1) model in \eqref{INAR1}.
Hence \ ${\bf Y}_n\in \widehat{S}^{\,\dag\dag}_n$ \ holds asymptotically as \ $n\to\infty$ \ with
 probability one.

Now we turn to find sets \ $S_n \subset \widehat{S}^{\dag\dag}_n$, $n\geq \max(3,s+2)$
 \ such that the system of equations
 \eqref{Additive_CLSE_EQ8} has a unique solution with respect to
 \ $(\alpha',\mu_\vare',\theta_1',\theta_2')$ \ for all \ ${\bf y}_n\in S_n$.
Let us introduce the \ $(4\times 4)$ \ Hessian matrix
 \[
    H_n({\bf y}_n;\alpha',\mu_\vare',\theta_1',\theta_2'):=
     \begin{bmatrix}
       \frac{\partial^2 Q_n^{\dag}}{\partial (\alpha')^2}
      \hspace{1mm} & \frac{\partial^2 Q_n^{\dag}}{\partial \mu_\vare' \partial \alpha'}
      \hspace{1mm} & \frac{\partial^2 Q_n^{\dag}}{\partial \theta_1' \partial \alpha'}
     \hspace{1mm}  & \frac{\partial^2 Q_n^{\dag}}{\partial \theta_2' \partial \alpha'} \smallskip \\
                     \frac{\partial^2 Q_n^{\dag}}{\partial \alpha' \partial \mu_\vare'}
     \hspace{1mm}   & \frac{\partial^2 Q_n^{\dag}}{\partial (\mu_\vare')^2}
     \hspace{1mm}   & \frac{\partial^2 Q_n^{\dag}}{ \partial \theta_1' \partial \mu_\vare' }
     \hspace{1mm}   &  \frac{\partial^2 Q_n^{\dag}}{\partial \theta_2' \partial \mu_\vare'} \smallskip \\
                      \frac{\partial^2 Q_n^{\dag}}{\partial \alpha' \partial \theta_1'}
    \hspace{1mm}    & \frac{\partial^2 Q_n^{\dag}}{\partial \mu_\vare' \partial \theta_1'}
    \hspace{1mm}    & \frac{\partial^2 Q_n^{\dag}}{ \partial (\theta_1')^2}
    \hspace{1mm}    &  \frac{\partial^2 Q_n^{\dag}}{\partial \theta_2' \partial \theta_1'} \smallskip \\
                       \frac{\partial^2 Q_n^{\dag}}{\partial \alpha' \partial \theta_2'}
   \hspace{1mm}     & \frac{\partial^2 Q_n^{\dag}}{\partial \mu_\vare' \partial \theta_2'}
   \hspace{1mm}     & \frac{\partial^2 Q_n^{\dag}}{ \partial \theta_1' \partial \theta_2'}
   \hspace{1mm}     &  \frac{\partial^2 Q_n^{\dag}}{ \partial (\theta_2')^2 } \\
     \end{bmatrix}
       ({\bf y}_n;\alpha',\mu_\vare',\theta_1',\theta_2'),
 \]
 and let us denote by \ $\Delta_{i,n}({\bf y}_n;\alpha',\mu_\vare',\theta_1',\theta_2')$ \
 its \ $i$-th order leading principal minor, \ $i=1,2,3,4$.
\ Further, for all \ $n\geq \max(3,s+2)$, \ let
 \[
   S_n:=\Big\{{\bf y}_n\in \widehat{S}^{\dag\dag}_n
         : \Delta_{i,n}({\bf y}_n;\alpha',\mu_\vare',\theta_1',\theta_2')>0,
                                   \;\, i=1,2,3,4,\, \forall\;(\alpha',\mu_\vare',\theta_1',\theta_2')
                                   \in\RR^4 \Big\}.
 \]
By Berkovitz \cite[Theorem 3.3, Chapter III]{Ber}, the function
 \ $\RR^4 \ni (\alpha',\mu_\vare',\theta_1',\theta_2')
     \mapsto Q_n^{\dag\dag}({\bf y}_n;\alpha',\mu_\vare',\theta_1',\theta_2')$
 \ is strictly convex for all \ ${\bf y}_n\in S_n$.
\ Since it was already proved that the system of equations \eqref{Additive_CLSE_EQ8} has a solution for all
 \ ${\bf y}_n\in \widehat{S}^{\dag\dag}_n$,
 \ we obtain that this solution is unique for all \ ${\bf y}_n\in S_n$.

Next we check that \ ${\bf Y}_n\in S_n$ \ holds asymptotically as \ $n\to\infty$ \ with probability one.
For all \ $(\alpha',\mu_\vare',\theta_1',\theta_2')\in\RR^4$,
 \begin{align*}
   &\frac{\partial^2 Q_n^{\dag\dag}}{\partial (\alpha')^2}
      ({\bf Y}_n;\alpha',\mu_\vare',\theta_1',\theta_2')\\
   &\phantom{\frac{\partial^2 Q_n^{\dag\dag}}{\partial (\alpha')^2}\;}
    = 2 \DS\sum_{\substack{k=1 \\ k\not\in \{s,s+1,s+2\}}}^n Y_{k-1}^2
            +2Y_{s-1}^2+2(Y_{s}-\theta_1')^2
            +2(Y_{s+1}-\theta_2')^2 \\
   &\phantom{\frac{\partial^2 Q_n^{\dag\dag}}{\partial (\alpha')^2}\;}
      =2\DS\sum_{\substack{k=1 \\ k\not\in \{s+1,s+2\}}}^n X_{k-1}^2
            + 2(X_{s}+\theta_1 - \theta_1')^2
            + 2(X_{s+1}+\theta_2 - \theta_2')^2 ,
 \end{align*}
 and
 \begin{align*}
   &\frac{\partial^2 Q_n^{\dag\dag}}
    {\partial\mu_\vare'\partial\alpha'}({\bf Y}_n;\alpha',\mu_\vare',\theta_1',\theta_2')
    = \frac{\partial^2 Q_n^{\dag\dag}}
     {\partial\alpha'\partial\mu_\vare'}({\bf Y}_n;\alpha',\mu_\vare',\theta_1',\theta_2')\\
   &\phantom{\frac{\partial^2 Q_n^{\dag\dag}}
      {\partial\alpha'\partial\mu_\vare'}({\bf Y}_n;\alpha',\mu_\vare',\theta_1',\theta_2')}
    = 2\sum_{k=1}^nY_{k-1}-2\theta_1'-2\theta_2'
    = 2\sum_{k=1}^nX_{k-1}+2(\theta_1-\theta_1')+2(\theta_2-\theta_2'),  \\
  &\frac{\partial^2 Q_n^{\dag\dag}}
    {\partial\theta_1'\partial\alpha'}({\bf Y}_n;\alpha',\mu_\vare',\theta_1',\theta_2')
    = \frac{\partial^2 Q_n^{\dag\dag}}
     {\partial\alpha'\partial\theta_1'}({\bf Y}_n;\alpha',\mu_\vare',\theta_1',\theta_2')\\
  &\phantom{\frac{\partial^2 Q_n^{\dag\dag}}
     {\partial\alpha'\partial\theta_1'}({\bf Y}_n;\alpha',\mu_\vare',\theta_1',\theta_2')}
    = 2(Y_{s-1}+Y_{s+1}-2\alpha' Y_{s}-\mu_\vare'+2\alpha'\theta_1'-\theta_2')\\
  &\phantom{\frac{\partial^2 Q_n^{\dag\dag}}
     {\partial\alpha'\partial\theta_1'}({\bf Y}_n;\alpha',\mu_\vare',\theta_1',\theta_2')}
    =2(X_{s-1}+X_{s+1}-2\alpha' X_{s}-\mu_\vare'-2\alpha'(\theta_1-\theta_1') + (\theta_2-\theta_2')),\\
  &\frac{\partial^2 Q_n^{\dag\dag}}
    {\partial\theta_2'\partial\alpha'}({\bf Y}_n;\alpha',\mu_\vare',\theta_1',\theta_2')
    = \frac{\partial^2 Q_n^{\dag\dag}}
       {\partial\alpha'\partial\theta_2'}({\bf Y}_n;\alpha',\mu_\vare',\theta_1',\theta_2')\\
  &\phantom{\frac{\partial^2 Q_n^{\dag\dag}}
    {\partial\alpha'\partial\theta_2'}({\bf Y}_n;\alpha',\mu_\vare',\theta_1',\theta_2')}
    = 2(Y_{s}+Y_{s+2}-2\alpha' Y_{s+1}-\mu_\vare'-\theta_1'+2\alpha'\theta_2')\\
  &\phantom{\frac{\partial^2 Q_n^{\dag\dag}}
     {\partial\alpha'\partial\theta_2'}({\bf Y}_n;\alpha',\mu_\vare',\theta_1',\theta_2')}
    =2(X_{s}+X_{s+2}-2\alpha' X_{s+1}-\mu_\vare'+(\theta_1-\theta_1')-2\alpha'(\theta_2-\theta_2')),\\
  & \frac{\partial^2 Q_n^{\dag\dag}}
      {\partial(\mu_\vare')^2}({\bf Y}_n;\alpha',\mu_\vare',\theta_1',\theta_2')
     = 2n,
 \end{align*}
 and
 \begin{align*}
  &\frac{\partial^2 Q_n^{\dag\dag}}{\partial(\theta_1')^2}({\bf Y}_n;\alpha',\mu_\vare',\theta_1',\theta_2')
    =\frac{\partial^2 Q_n^{\dag\dag}}
     {\partial(\theta_2')^2}({\bf Y}_n;\alpha',\mu_\vare',\theta_1',\theta_2')
    = 2((\alpha')^2+1),\\
  &\frac{\partial^2 Q_n^{\dag\dag}}
    {\partial\theta_1'\partial\theta_2'}({\bf Y}_n;\alpha',\mu_\vare',\theta_1',\theta_2')
    =\frac{\partial^2 Q_n^{\dag\dag}}
      {\partial\theta_2'\partial\theta_1'}({\bf Y}_n;\alpha',\mu_\vare',\theta_1',\theta_2')
     = -2\alpha',\\
 &\frac{\partial^2 Q_n^{\dag\dag}}
   {\partial\theta_1'\partial\mu_\vare'}({\bf Y}_n;\alpha',\mu_\vare',\theta_1',\theta_2')
    =\frac{\partial^2 Q_n^{\dag\dag}}
       {\partial\mu_\vare'\partial\theta_1'}({\bf Y}_n;\alpha',\mu_\vare',\theta_1',\theta_2')
     = 2(1-\alpha'),\\
 &\frac{\partial^2 Q_n^{\dag\dag}}
   {\partial\theta_2'\partial\mu_\vare'}({\bf Y}_n;\alpha',\mu_\vare',\theta_1',\theta_2')
    =\frac{\partial^2 Q_n^{\dag\dag}}
       {\partial\mu_\vare'\partial\theta_2'}({\bf Y}_n;\alpha',\mu_\vare',\theta_1',\theta_2')
     = 2(1-\alpha').
 \end{align*}
Then \ $H_n({\bf Y}_n;\alpha',\mu_\vare',\theta_1',\theta_2')$ \  has the following leading principal minors
 \begin{align*}
    &\Delta_{1,n}({\bf Y}_n;\alpha',\mu_\vare',\theta_1',\theta_2')
       =2\DS\sum_{\substack{k=1 \\ k\not\in \{s+1,s+2\}}}^n X_{k-1}^2
            + 2(X_{s}+\theta_1 - \theta_1')^2
            + 2(X_{s+1}+\theta_2 - \theta_2')^2, \\
    &\Delta_{2,n}({\bf Y}_n;\alpha',\mu_\vare',\theta_1',\theta_2')
      = 4n\left(\DS\sum_{\substack{k=1 \\ k\not\in \{s+1,s+2\}}}^n X_{k-1}^2
                                   + (X_{s}+\theta_1 - \theta_1')^2
                                   + (X_{s+1}+\theta_2 - \theta_2')^2\right) \\
    & \phantom{\Delta_{2,n}({\bf Y}_n;\alpha',\mu_\vare',\theta_1',\theta_2'):=\;}
       -4\left(\sum_{k=1}^nX_{k-1}+(\theta_1-\theta_1')+(\theta_2-\theta_2')\right)^2,
 \end{align*}
  and
  \begin{align*}
    &\Delta_{3,n}({\bf Y}_n;\alpha',\mu_\vare',\theta_1',\theta_2') \\
    &\quad =8\big(((\alpha')^2+1)n - (1-\alpha')^2\big)
           \left(\DS\sum_{\substack{k=1 \\ k\not\in \{s+1,s+2\}}}^n X_{k-1}^2
                                   + (X_{s}+\theta_1 - \theta_1')^2
                                   + (X_{s+1}+\theta_2 - \theta_2')^2\right) \\
    &\quad \qquad+16L\left(\sum_{k=1}^nX_{k-1}+(\theta_1-\theta_1')+(\theta_2-\theta_2')\right)
               -8nL^2 \\
    &\quad \qquad -8((\alpha')^2+1)\left(\sum_{k=1}^nX_{k-1}+(\theta_1-\theta_1')+(\theta_2-\theta_2')\right)^2,\\
    &\Delta_{4,n}({\bf Y}_n;\alpha',\mu_\vare',\theta_1',\theta_2')
      = \det H_n({\bf Y}_n;\alpha',\mu_\vare',\theta_1',\theta_2'),
 \end{align*}
 where \ $L:=X_{s-1}+X_{s+1}-2\alpha' X_{s}-\mu_\vare'-2\alpha'(\theta_1-\theta_1')+\theta_2-\theta_2'$.
\ By \eqref{Ergodic1} and \eqref{Ergodic2}, we get
 \begin{align*}
    &\PP\left(\lim_{n\to\infty}\frac{1}{n}\Delta_{1,n}({\bf Y}_n;\alpha',\mu_\vare',\theta_1',\theta_2')
               = 2 \EE\widetilde X^2,
                \;\;\forall\;\; (\alpha',\mu_\vare',\theta_1',\theta_2')\in\RR^4 \right)=1, \\
    & \PP\left(\lim_{n\to\infty}\frac{1}{n^2}\Delta_{2,n}({\bf Y}_n;\alpha',\mu_\vare',\theta_1',\theta_2')
               = 4(\EE\widetilde X^2-(\EE\widetilde X)^2)=4\var\widetilde X,
                \;\;\forall\;\; (\alpha',\mu_\vare',\theta_1',\theta_2')\in\RR^4 \right)=1,
  \end{align*}
  and
  \begin{align*}
    &  \PP\left(\lim_{n\to\infty}\frac{1}{n^2}\Delta_{3,n}({\bf Y}_n;\alpha',\mu_\vare',\theta_1',\theta_2')
               = 8((\alpha')^2+1)\var\widetilde X,
                \;\;\forall\;\; (\alpha',\mu_\vare',\theta_1',\theta_2')\in\RR^4 \right)=1,  \\
    & \PP\left(\lim_{n\to\infty}\frac{1}{n^2}\Delta_{4,n}({\bf Y}_n;\alpha',\mu_\vare',\theta_1',\theta_2')
               =16((\alpha')^4+(\alpha')^2+1)\var\widetilde X,
                \;\;\forall\;\; (\alpha',\mu_\vare',\theta_1',\theta_2')\in\RR^4 \right)=1,
 \end{align*}
 where \ $\widetilde X$ \ denotes a random variable with the unique stationary
 distribution of the INAR(1) model in \eqref{INAR1}.
Hence
 \begin{align*}
   & \PP\big(\lim_{n\to\infty}\Delta_{i,n}({\bf Y}_n;\alpha',\mu_\vare',\theta_1',\theta_2')=\infty,
             \;\;\forall\;\; (\alpha',\mu_\vare',\theta_1',\theta_2')\in\RR^4 \big)=1,
      \qquad i=1,2,3,4,
 \end{align*}
 which yields that \ ${\bf Y}_n\in S_n$ \ asymptotically as \ $n\to\infty$
 \ with probability one, since we have already proved that \ ${\bf Y}_n\in \widehat{S}^{\dag\dag}_n$
 \ asymptotically as \ $n\to\infty$ \ with probability one.
\proofend

By Lemma \ref{LEMMA14}, \ $(\halpha_n^{\,\dag\dag}({\bf Y}_n),
    \hmuen^{\,\dag\dag}({\bf Y}_n),
    \htheta_{1,n}^{\,\dag\dag}({\bf Y}_n),
    \htheta_{2,n}^{\,\dag\dag}({\bf Y}_n))$
 \ exists uniquely asymptotically as \ $n\to\infty$ \ with probability one.
In the sequel we will simply denote it by
 $(\halpha_n^{\,\dag\dag},\hmuen^{\,\dag\dag},\htheta_{1,n}^{\,\dag\dag},\htheta_{2,n}^{\,\dag\dag})$.

The next result shows that \ $\halpha_n^{\,\dag\dag}$ \ is a strongly consistent estimator
 of \ $\alpha$, \ $\hmuen^{\,\dag\dag}$ \ is a strongly consistent estimator of \ $\mu_\vare$,
 \ whereas \ $\htheta_{1,n}^{\,\dag\dag}$ \ and  \ $\htheta_{2,n}^{\,\dag\dag}$ \ fail to be
 strongly consistent estimators of \ $\theta_1$ \ and \ $\theta_2$, \ respectively.

\begin{Thm}\label{THEOREM6}
Consider the CLS estimators
 \ $(\halpha_n^{\,\dag\dag},\hmuen^{\,\dag\dag},\htheta_{1,n}^{\,\dag\dag},\htheta_{2,n}^{\,\dag\dag})_{n\in\NN}$ \ of the parameter
 \ $(\alpha,\mu_\vare,\theta_1,\theta_2)\in(0,1)\times(0,\infty)\times\NN^2$.
\ The sequences \ $(\halpha_n^{\,\dag\dag})_{n\in\NN}$ \ and \ $(\hmuen^{\,\dag\dag})_{n\in\NN}$
 \ are strongly consistent for all
 \ $(\alpha,\mu_\vare,\theta_1,\theta_2)\in(0,1)\times(0,\infty)\times\NN^2$, \ i.e.,
 \begin{align} \label{Strong_consistency22}
   &\PP(\lim_{n\to\infty}\halpha_n^{\,\dag\dag}=\alpha)=1,
     \qquad \forall\;(\alpha,\mu_\vare,\theta_1,\theta_2)\in(0,1)\times(0,\infty)\times\NN^2,\\
       \label{Strong_consistency23}
  &\PP(\lim_{n\to\infty}\hmuen^{\,\dag\dag}=\mu_\vare)=1,
     \qquad \forall\;(\alpha,\mu_\vare,\theta_1,\theta_2)\in(0,1)\times(0,\infty)\times\NN^2,
 \end{align}
 whereas the sequences \ $(\htheta_{1,n}^{\,\dag\dag})_{n\in\NN}$ \ and
  \ $(\htheta_{2,n}^{\,\dag\dag})_{n\in\NN}$ \ are not strongly consistent for any
  \ $(\alpha,\mu_\vare,\theta_1,\theta_2)\in(0,1)\times(0,\infty)\times\NN^2$, \ namely,
 \begin{align}\label{Strong_consistency24}
   \PP\left(\lim_{n\to\infty}
             \begin{bmatrix}
               \ttheta_{1,n}^{\,\dag\dag} \\[2mm]
               \ttheta_{2,n}^{\,\dag\dag} \\
             \end{bmatrix}
            = \begin{bmatrix}
                Y_s \\
                Y_{s+1} \\
              \end{bmatrix}
             +
              \begin{bmatrix}
                \frac{-\alpha(1+\alpha^2)Y_{s-1}
                   -\alpha^2Y_{s+2}
                   -(1-\alpha^3)\mu_\vare}{1+\alpha^2+\alpha^4} \\[2mm]
                \frac{-\alpha^2Y_{s-1}
                   -\alpha(1+\alpha^2)Y_{s+2}
                   -(1-\alpha^3)\mu_\vare}{1+\alpha^2+\alpha^4} \\
              \end{bmatrix}
          \right)=1
 \end{align}
 for all \ $(\alpha,\mu_\vare,\theta_1,\theta_2)\in(0,1)\times(0,\infty)\times\NN^2$.
\end{Thm}

\noindent{\bf Proof.}
An easy calculation shows that
 \begin{align*}
   &\left(\DS\sum_{\substack{k=1 \\ k\not\in \{s+1,s+2\}}}^n Y_{k-1}^2
        + (Y_{s}-\htheta_{1,n}^{\,\dag\dag})^2 + (Y_{s+1}-\htheta_{2,n}^{\,\dag\dag})^2 \right)\halpha_n^{\,\dag\dag}
        +\left(\sum_{k=1}^n Y_{k-1} - \htheta_{1,n}^{\,\dag\dag} - \htheta_{2,n}^{\,\dag\dag} \right)\hmuen^{\,\dag\dag} \\
   &\qquad\qquad
      = \sum_{k=1}^n Y_{k-1}Y_k - \htheta_{1,n}^{\,\dag\dag}(Y_{s-1}+Y_{s+1})
        - \htheta_{2,n}^{\,\dag\dag}(Y_{s}+Y_{s+2})+\htheta_{1,n}^{\,\dag\dag}\htheta_{2,n}^{\,\dag\dag} , \\
   &\left(\sum_{k=1}^n Y_{k-1} - \htheta_{1,n}^{\,\dag\dag} - \htheta_{2,n}^{\,\dag\dag} \right)\halpha_n^{\,\dag\dag}
        + n\hmuen^{\,\dag\dag}
      = \sum_{k=1}^n Y_k - \htheta_{1,n}^{\,\dag\dag} - \htheta_{2,n}^{\,\dag\dag} ,
 \end{align*}
 hold asymptotically as \ $n\to\infty$ \ with probability one, and hence
 \begin{align}\label{SEGED63}
    \begin{bmatrix}
      \halpha_n^{\,\dag\dag} \\
      \hmuen^{\,\dag\dag} \\
     \end{bmatrix}
    = \begin{bmatrix}
         a_n & b_n \\
         b_n &  n \\
        \end{bmatrix}^{-1}
       \begin{bmatrix}
        k_n \\
        \ell_n \\
       \end{bmatrix},
 \end{align}
 holds asymptotically as \ $n\to\infty$ \ with probability one, where
 \begin{align*}
  & a_n:= \sum_{k=1}^n X_{k-1}^2 + (\theta_1-\htheta_{1,n}^{\,\dag\dag})(\theta_1-\htheta_{1,n}^{\,\dag\dag}+2X_{s})
                               + (\theta_2-\htheta_{2,n}^{\,\dag\dag})(\theta_2-\htheta_{2,n}^{\,\dag\dag}+2X_{s+1}),\\
  & b_n:= \sum_{k=1}^n X_{k-1} + \theta_1-\htheta_{1,n}^{\,\dag\dag} + \theta_2-\htheta_{2,n}^{\,\dag\dag},
  \end{align*}
  and
  \begin{align*}
  & k_n:= \sum_{k=1}^n X_{k-1}X_k + (\theta_1-\htheta_{1,n}^{\,\dag\dag})(X_{s-1}+X_{s+1})
                                 + (\theta_2-\htheta_{2,n}^{\,\dag\dag})(X_{s}+X_{s+2})
                                 + (\theta_1-\htheta_{1,n}^{\,\dag\dag})(\theta_2-\htheta_{2,n}^{\,\dag\dag}),\\
  & \ell_n:= \sum_{k=1}^n X_k + \theta_1-\htheta_{1,n}^{\,\dag\dag}
                                 + \theta_2-\htheta_{2,n}^{\,\dag\dag}.
 \end{align*}
 Furthermore,
 \begin{align}\label{SEGED61}
   \begin{bmatrix}
     \halpha_n^{\,\dag\dag} - \alpha \\
     \hmuen^{\,\dag\dag} - \mu_\vare\\
   \end{bmatrix}
    = \begin{bmatrix}
         a_n & b_n \\
         b_n &  n \\
        \end{bmatrix}^{-1}
       \begin{bmatrix}
          c_n \\
         d_n \\
       \end{bmatrix}
 \end{align}
 hold asymptotically as \ $n\to\infty$ \ with probability one,
 where
 \begin{align*}
  &c_n:=
     \sum_{k=1}^n X_{k-1}(X_k-\alpha X_{k-1}-\mu_\vare)
              + (\theta_1-\htheta_{1,n}^{\,\dag\dag})
                \big(X_{s-1}+X_{s+1}-2\alpha X_{s} - \mu_\vare - \alpha(\theta_1-\htheta_{1,n}^{\,\dag\dag})\big)\\
  &\phantom{c_n:=\;}
              + (\theta_2-\htheta_{2,n}^{\,\dag\dag})
               \big(X_{s}+X_{s+2}-2\alpha X_{s+1} - \mu_\vare - \alpha(\theta_2-\htheta_{2,n}^{\,\dag\dag})\big)
              + (\theta_1-\htheta_{1,n}^{\,\dag\dag})(\theta_2-\htheta_{2,n}^{\,\dag\dag}),\\
  &d_n:= \sum_{k=1}^n (X_k -\alpha X_{k-1}-\mu_\vare)
              + (1-\alpha)(\theta_1-\htheta_{1,n}^{\,\dag\dag}
                                 + \theta_2-\htheta_{2,n}^{\,\dag\dag}).
 \end{align*}

We show that the sequences
 \ $(\htheta_{1,n}^{\,\dag\dag}-\theta_1)_{n\in\NN}$
 \ and \ $(\htheta_{2,n}^{\,\dag\dag}-\theta_2)_{n\in\NN}$ \ are bounded with probability one.
\ An easy calculation shows that
 \begin{align*}
   &n \hmuen^{\,\dag\dag} + (1-\halpha_n^{\,\dag\dag})\htheta_{1,n}^{\,\dag\dag} + (1-\halpha_n^{\,\dag\dag})\htheta_{2,n}^{\,\dag\dag}
       = \sum_{k=1}^n (Y_k - \halpha_n^{\,\dag\dag} Y_{k-1}),\\
   & (1-\halpha_n^{\,\dag\dag})\hmuen^{\,\dag\dag}  + (1+(\halpha_n^{\,\dag\dag})^2)\htheta_{1,n}^{\,\dag\dag}
       - \halpha_n^{\,\dag\dag}\htheta_{2,n}^{\,\dag\dag}
       = (1+(\halpha_n^{\,\dag\dag})^2) Y_{s} - \halpha_n^{\,\dag\dag}(Y_{s-1}+Y_{s+1}),\\
   & (1-\halpha_n^{\,\dag\dag})\hmuen^{\,\dag\dag}  - \halpha_n^{\,\dag\dag}\htheta_{1,n}^{\,\dag\dag}
      + (1+(\halpha_n^{\,\dag\dag})^2)\htheta_{2,n}^{\,\dag\dag}
       = (1+(\halpha_n^{\,\dag\dag})^2) Y_{s+1} - \halpha_n^{\,\dag\dag}(Y_{s}+Y_{s+2}),
 \end{align*}
 hold asymptotically as \ $n\to\infty$ \ with probability one, or equivalently
 \begin{align}\label{SEGED97}
   \begin{bmatrix}
     n & 1-\halpha_n^{\,\dag\dag} & 1-\halpha_n^{\,\dag\dag} \\
     1-\halpha_n^{\,\dag\dag} & 1+(\halpha_n^{\,\dag\dag})^2 & -\halpha_n^{\,\dag\dag} \\
     1-\halpha_n^{\,\dag\dag} & -\halpha_n^{\,\dag\dag} &  1+(\halpha_n^{\,\dag\dag})^2 \\
   \end{bmatrix}
   \begin{bmatrix}
     \hmuen^{\,\dag\dag} \\
     \htheta_{1,n}^{\,\dag\dag} \\
     \htheta_{2,n}^{\,\dag\dag} \\
    \end{bmatrix}
   = \begin{bmatrix}
       \sum_{k=1}^n (Y_k - \halpha_n^{\,\dag\dag} Y_{k-1}) \\
       (1+(\halpha_n^{\,\dag\dag})^2) Y_{s} - \halpha_n^{\,\dag\dag}(Y_{s-1}+Y_{s+1}) \\
       (1+(\halpha_n^{\,\dag\dag})^2) Y_{s+1} - \halpha_n^{\,\dag\dag}(Y_{s}+Y_{s+2}) \\
     \end{bmatrix}
 \end{align}
 holds asymptotically as \ $n\to\infty$ \ with probability one.
Since for all \ $n\geq 3$
 \[
   D_n(\halpha_n^{\,\dag\dag})
     = (1+\halpha_n^{\,\dag\dag}+(\halpha_n^{\,\dag\dag})^2)
          \big( (n-2)(\halpha_n^{\,\dag\dag})^2 - (n-4)\halpha_n^{\,\dag\dag} + n-2 \big)
     >0,
 \]
 we get asymptotically as \ $n\to\infty$ \ with probability one we have
 \begin{align}\label{SEGED59}
   \begin{split}
    \begin{bmatrix}
       \hmuen^{\,\dag\dag} \\
       \htheta_{1,n}^{\,\dag\dag} \\
       \htheta_{2,n}^{\,\dag\dag} \\
     \end{bmatrix}
   &= \frac{1}{D_n(\halpha_n^{\,\dag\dag}) }
      \begin{bmatrix}
        1+(\halpha_n^{\,\dag\dag})^2+(\halpha_n^{\,\dag\dag})^4
              & u_n  & u_n \\
        u_n  & w_n
             &  v_n \\
        u_n  & v_n & w_n  \\
      \end{bmatrix}
       \begin{bmatrix}
         \sum_{k=1}^n (Y_k - \halpha_n^{\,\dag\dag} Y_{k-1}) \\
         (1+(\halpha_n^{\,\dag\dag})^2) Y_{s} - \halpha_n^{\,\dag\dag}(Y_{s-1}+Y_{s+1}) \\
         (1+(\halpha_n^{\,\dag\dag})^2) Y_{s+1} - \halpha_n^{\,\dag\dag}(Y_{s}+Y_{s+2}) \\
       \end{bmatrix}\\
  & \;\; =: \frac{1}{D_n(\halpha_n^{\,\dag\dag}) }
       \begin{bmatrix}
         G_n \\
         H_n \\
         J_n \\
       \end{bmatrix},
  \end{split}
 \end{align}
 where
 \begin{align*}
    & u_n:=-(1-\halpha_n^{\,\dag\dag})(1+\halpha_n^{\,\dag\dag}+(\halpha_n^{\,\dag\dag})^2), \\
    & v_n:=(1-\halpha_n^{\,\dag\dag})^2 + n \halpha_n^{\,\dag\dag}, \\
    & w_n:=n(1+(\halpha_n^{\,\dag\dag})^2)-(1-\halpha_n^{\,\dag\dag})^2,
 \end{align*}
 and
 \begin{align*}
   & G_n:= - (1-\halpha_n^{\,\dag\dag})(1+\halpha_n^{\,\dag\dag}+(\halpha_n^{\,\dag\dag})^2)
              \Big((1+(\halpha_n^{\,\dag\dag})^2)(Y_{s}+Y_{s+1})
                    - \halpha_n^{\,\dag\dag}(Y_{s-1}+Y_{s+1}+Y_{s}+Y_{s+2})\Big) \\
         &\phantom{ G_n:=\;} + (1+(\halpha_n^{\,\dag\dag})^2+(\halpha_n^{\,\dag\dag})^4)
                      \sum_{k=1}^n (Y_k - \halpha_n^{\,\dag\dag} Y_{k-1}) ,\\
   & H_n:= \big( n(1+(\halpha_n^{\,\dag\dag})^2)-(1-\halpha_n^{\,\dag\dag})^2\big)
              \Big((1+(\halpha_n^{\,\dag\dag})^2) Y_{s} - \halpha_n^{\,\dag\dag}(Y_{s-1}+Y_{s+1}) \Big) \\
   &\phantom{H_n:=\;} + ((1-\halpha_n^{\,\dag\dag})^2+n\halpha_n^{\,\dag\dag})
      \Big((1+(\halpha_n^{\,\dag\dag})^2) Y_{s+1} - \halpha_n^{\,\dag\dag}(Y_{s}+Y_{s+2})\Big)\\
   &\phantom{H_n:=\;}
           -(1-\halpha_n^{\,\dag\dag})(1+\halpha_n^{\,\dag\dag}+(\halpha_n^{\,\dag\dag})^2)
                \sum_{k=1}^n (Y_k - \halpha_n^{\,\dag\dag} Y_{k-1}),\\
  & J_n:= \big((1-\halpha_n^{\,\dag\dag})^2+n\halpha_n^{\,\dag\dag}\big)
      \Big((1+(\halpha_n^{\,\dag\dag})^2) Y_{s} - \halpha_n^{\,\dag\dag}(Y_{s-1}+Y_{s+1})\Big)\\
  &\phantom{H_n:=\;}
          + \big( n(1+(\halpha_n^{\,\dag\dag})^2)-(1-\halpha_n^{\,\dag\dag})^2\big)
     \Big((1+(\halpha_n^{\,\dag\dag})^2) Y_{s+1} - \halpha_n^{\,\dag\dag}(Y_{s}+Y_{s+2}) \Big) \\
  &\phantom{H_n:=\;}
      -(1-\halpha_n^{\,\dag\dag})(1+\halpha_n^{\,\dag\dag}+(\halpha_n^{\,\dag\dag})^2)
       \sum_{k=1}^n (Y_k - \halpha_n^{\,\dag\dag} Y_{k-1})    .
 \end{align*}
Using \eqref{Ergodic1} and that for all \ $p_i\in\RR$, \ $i=0,\ldots,4$,
 \[
    \sup_{x\in\RR,\;n\in\NN}
        \frac{n(p_4x^4+p_3x^3+p_2x^2+p_1x+p_0)}{(1+x+x^2)((n-2)x^2-(n-4)x+n-2)}
       <\infty,
 \]
 one can think it over that \ $H_n/D_n(\halpha_n^{\,\dag\dag})$, \ $n\in\NN$, \ and
 \ $J_n/D_n(\halpha_n^{\,\dag\dag})$, \ $n\in\NN$, \ are bounded with probability one, which yields also that
 the sequences \ $(\htheta_{1,n}^{\,\dag\dag}-\theta_1)_{n\in\NN}$
 \ and \ $(\htheta_{2,n}^{\,\dag\dag}-\theta_2)_{n\in\NN}$ \ are bounded
 with probability one.

By the same arguments as in the proof of Theorem \ref{THEOREM4}, one can
 derive \eqref{Strong_consistency22} and \eqref{Strong_consistency23}.
Indeed, using \eqref{STAC_MOMENT1}, \eqref{Ergodic1}, \eqref{Ergodic2} and \eqref{Ergodic3}, we get
 \begin{align*}
   &\PP\left(\lim_{n\to\infty}\frac{a_n}{n}=\EE\widetilde X^2\right)=1,
     \qquad\qquad \PP\left(\lim_{n\to\infty}\frac{b_n}{n}=\EE\widetilde X\right)=1,\\
   &\PP\left(\lim_{n\to\infty}\frac{k_n}{n}=\alpha\EE\widetilde X^2+\mu_\vare\EE\widetilde X\right)=1,
     \qquad\qquad \PP\left(\lim_{n\to\infty}\frac{\ell_n}{n}=\EE\widetilde X\right)=1,\\
   &\PP\left(\lim_{n\to\infty}\frac{c_n}{n}
      =\alpha\EE\widetilde X^2+\mu_\vare\EE\widetilde X
       -\alpha\EE\widetilde X^2-\mu_\vare\EE\widetilde X = 0 \right)=1,\\
   &\PP\left(\lim_{n\to\infty}\frac{d_n}{n}
      =\EE\widetilde X - \alpha\EE\widetilde X-\mu_\vare = 0 \right)=1.
 \end{align*}
Hence, by \eqref{SEGED61}, we obtain
 \begin{align*}
   \PP\left(
     \lim_{n\to\infty}
          \begin{bmatrix}
            \halpha_n^{\,\dag\dag} - \alpha \\
            \hmuen^{\,\dag\dag} - \mu_\vare\\
           \end{bmatrix}
       = \begin{bmatrix}
           \EE\widetilde X^2 & \EE\widetilde X \\
           \EE\widetilde X &  1 \\
         \end{bmatrix}^{-1}
          \begin{bmatrix}
              0 \\
              0 \\
          \end{bmatrix}
      =  \begin{bmatrix}
            0 \\
            0 \\
          \end{bmatrix}
   \right)=1,
 \end{align*}
 which yields \eqref{Strong_consistency22} and \eqref{Strong_consistency23}.
Then \eqref{Strong_consistency22}, \eqref{Strong_consistency23} and \eqref{SEGED59}
  imply \eqref{Strong_consistency24}. 
Indeed,
 \begin{align*}
   \PP\left(\lim_{n\to\infty}\frac{D_n(\halpha_n^{\,\dag\dag})}{n} = 1+\alpha^2+\alpha^4 \right)=1,
    \qquad \forall\;\; (\alpha,\mu_\vare,\theta_1,\theta_2)\in(0,1)\times(0,\infty)\times\NN^2,
 \end{align*}
 and \ $\frac{H_n}{n}$ \ converges almost surely as \ $n\to\infty$ \ to
 \begin{align*}
   &(1+\alpha^2)\Big((1+\alpha^2)Y_s-\alpha(Y_{s-1}+Y_{s+1})\Big)
     + \alpha \Big((1+\alpha^2)Y_{s+1}-\alpha(Y_s+Y_{s+2})\Big)\\
   &\qquad - (1-\alpha)(1+\alpha+\alpha^2)(1-\alpha)\EE\widetilde X \\
   & =-\alpha(1+\alpha^2)Y_{s-1} + (1+\alpha^2+\alpha^4)Y_s
      - \alpha^2 Y_{s+2} - (1-\alpha^3)\mu_\vare,
 \end{align*}
 and \ $\frac{J_n}{n}$ \ converges almost surely as \ $n\to\infty$ \ to
 \begin{align*}
   &\alpha\Big((1+\alpha^2)Y_s-\alpha(Y_{s-1}+Y_{s+1})\Big)
     + (1+\alpha^2) \Big((1+\alpha^2)Y_{s+1}-\alpha(Y_s+Y_{s+2})\Big) \\
   &\qquad - (1-\alpha)(1+\alpha+\alpha^2)(1-\alpha)\EE\widetilde X \\
   & =-\alpha^2 Y_{s-1} + (1+\alpha^2+\alpha^4)Y_{s+1}
      - \alpha(1+\alpha^2) Y_{s+2} - (1-\alpha^3)\mu_\vare.
 \end{align*}
\proofend

The asymptotic distribution of the CLS estimation is given in the next theorem.

\begin{Thm}
Under the additional assumptions \ $\EE X_0^3<\infty$ \ and \ $\EE\vare_1^3<\infty$, \ we have
 \begin{align}\label{CONVERGENCE16}
   \begin{bmatrix}
     \sqrt{n}(\halpha_n^{\,\dag\dag}-\alpha) \\
     \sqrt{n}(\hmuen^{\,\dag\dag}-\mu_\vare) \\
   \end{bmatrix}
   \distr \cN\left(\begin{bmatrix}
                     0 \\
                     0 \\
                   \end{bmatrix}
     ,B_{\alpha,\,\vare}\right)
      \qquad \text{as \ $n\to\infty$,}
  \end{align}
 where \ $B_{\alpha,\vare}$ \ is defined in \eqref{SEGED_BALPHA}.
Moreover, conditionally on the values \ $Y_{s-1}$ \ and \ $Y_{s+2}$,
 \begin{align}\label{CONVERGENCE17}
    \begin{bmatrix}
      \sqrt{n}\big(\htheta_{1,n}^{\,\dag\dag} - \lim_{k\to\infty}\htheta_{1,k}^{\,\dag\dag}\big) \\
      \sqrt{n}\big(\htheta_{2,n}^{\,\dag\dag} - \lim_{k\to\infty}\htheta_{2,k}^{\,\dag\dag}\big) \\
    \end{bmatrix}
      \distr \cN\left(\begin{bmatrix}
                        0 \\
                        0 \\
                      \end{bmatrix},
                      D_{\alpha,\vare} B_{\alpha,\vare} D_{\alpha,\vare}^\top
      \right)
      \qquad \text{as \ $n\to\infty$,}
 \end{align}
 where the $(2\times 2)$-matrix  \ $D_{\alpha,\vare}$ \ is defined by
 \begin{align*}
    D_{\alpha,\vare}
       := \begin{bmatrix}
            \frac{(\alpha^2-1)(\alpha^4+3\alpha^2+1)Y_{s-1} + 2\alpha(\alpha^4-1)Y_{s+2}
                     + \alpha(2-\alpha)(1+\alpha+\alpha^2)^2\mu_\vare} {(1+\alpha^2+\alpha^4)^2}
              \hspace{1mm}  & \frac{\alpha^3-1}{1+\alpha^2+\alpha^4} \medskip \\
            \frac{2\alpha(\alpha^4-1)Y_{s-1} + (\alpha^2-1)(\alpha^4+3\alpha^2+1)Y_{s+2}
                     + \alpha(2-\alpha)(1+\alpha+\alpha^2)^2\mu_\vare} {(1+\alpha^2+\alpha^4)^2}
              \hspace{1mm} & \frac{\alpha^3-1}{1+\alpha^2+\alpha^4} \\
           \end{bmatrix}.
 \end{align*}
\end{Thm}

\noindent{\bf Proof.}
Using \eqref{SEGED61} and that the sequences
 \ $(\htheta_{1,n}^{\,\dag\dag}-\theta_1)_{n\in\NN}$
 \ and \ $(\htheta_{2,n}^{\,\dag\dag}-\theta_2)_{n\in\NN}$ \ are bounded with probability one,
 by the very same arguments as in the proof of \eqref{CONVERGENCE8},
 one can obtain \eqref{CONVERGENCE16}.
Now we turn to prove \eqref{CONVERGENCE17}.
Using the notation
 \[
   B_n^{\dag\dag}:=
          \begin{bmatrix}
             1+(\halpha_n^{\,\dag\dag})^2 & -\halpha_n^{\,\dag\dag} \\
             -\halpha_n^{\,\dag\dag} & 1+(\halpha_n^{\,\dag\dag})^2 \\
           \end{bmatrix},
 \]
 by \eqref{SEGED97}, we have
  \begin{align*}
    \begin{bmatrix}
      \htheta_{1,n}^{\,\dag\dag} \\
      \htheta_{2,n}^{\,\dag\dag} \\
    \end{bmatrix}
     = (B_n^{\dag\dag})^{-1}
        \begin{bmatrix}
         (1+(\halpha_n^{\,\dag\dag})^2)Y_{s} - \halpha_n^{\,\dag\dag}(Y_{s-1}+Y_{s+1})
          - (1-\halpha_n^{\,\dag\dag})\hmuen^{\,\dag\dag}\\
         (1+(\halpha_n^{\,\dag\dag})^2)Y_{s+1} - \halpha_n^{\,\dag\dag}(Y_{s}+Y_{s+2})
            - (1-\halpha_n^{\,\dag\dag})\hmuen^{\,\dag\dag} \\
        \end{bmatrix}
  \end{align*}
 holds asymptotically as \ $n\to\infty$ \ with probability one.
Theorem \ref{THEOREM6} yields that
 \[
     \PP\left(\lim_{n\to\infty}B_n^{\dag\dag}
                 = \begin{bmatrix}
                    1+\alpha^2 & -\alpha \\
                    -\alpha & 1+\alpha^2 \\
                   \end{bmatrix}=:  B^{\dag\dag}
         \right)=1.
 \]
Again by Theorem \ref{THEOREM6}, we have
 \begin{align*}
   &\begin{bmatrix}
      \sqrt{n}\big(\htheta_{1,n}^{\,\dag\dag} - \lim_{k\to\infty}\htheta_{1,k}^{\,\dag\dag}\big) \\
      \sqrt{n}\big(\htheta_{2,n}^{\,\dag\dag} - \lim_{k\to\infty}\htheta_{2,k}^{\,\dag\dag}\big) \\
   \end{bmatrix}\\
   &=\sqrt{n}(B_n^{\dag\dag})^{-1}
    \left(
      \begin{bmatrix}
         (1+(\halpha_n^{\,\dag\dag})^2)Y_{s} - \halpha_n^{\,\dag\dag}(Y_{s-1}+Y_{s+1})
             - (1-\halpha_n^{\,\dag\dag})\hmuen^{\,\dag\dag}\\
         (1+(\halpha_n^{\,\dag\dag})^2)Y_{s+1} - \halpha_n^{\,\dag\dag}(Y_{s}+Y_{s+2})
             - (1-\halpha_n^{\,\dag\dag})\hmuen^{\,\dag\dag} \\
      \end{bmatrix}\right.\\
   &\phantom{= \sqrt{n}(B_n^{\dag\dag})^{-1}\Big(\;\;}\left.
     - B_n^{\dag\dag} (B^{\dag\dag})^{-1}
      \begin{bmatrix}
         (1+\alpha^2)Y_{s} - \alpha(Y_{s-1}+Y_{s+1}) - (1-\alpha)\mu_\vare\\
         (1+\alpha^2)Y_{s+1} - \alpha(Y_{s}+Y_{s+2}) - (1-\alpha)\mu_\vare \\
      \end{bmatrix}
     \right)\\
    &=\sqrt{n}(B_n^{\dag\dag})^{-1}
      \left(
       \begin{bmatrix}
         (1+(\halpha_n^{\,\dag\dag})^2)Y_{s} - \halpha_n^{\,\dag\dag}(Y_{s-1}+Y_{s+1})
            - (1-\halpha_n^{\,\dag\dag})\hmuen^{\,\dag\dag} \\
         (1+(\halpha_n^{\,\dag\dag})^2)Y_{s+1} - \halpha_n^{\,\dag\dag}(Y_{s}+Y_{s+2})
            - (1-\halpha_n^{\,\dag\dag})\hmuen^{\,\dag\dag} \\
       \end{bmatrix}\right.\\
  &\phantom{=\sqrt{n}(B_n^{\dag\dag})^{-1}\Big(\;}\left.
    -\begin{bmatrix}
         (1+\alpha^2)Y_{s} - \alpha(Y_{s-1}+Y_{s+1}) - (1-\alpha)\mu_\vare\\
         (1+\alpha^2)Y_{s+1} - \alpha(Y_{s}+Y_{s+2}) - (1-\alpha)\mu_\vare \\
      \end{bmatrix}
       \right)\\
 &\phantom{=\;}
      +\sqrt{n}
      \left((B_n^{\dag\dag})^{-1} - (B^{\dag\dag})^{-1}
       \right)
       \begin{bmatrix}
         (1+\alpha^2)Y_{s} - \alpha(Y_{s-1}+Y_{s+1}) - (1-\alpha)\mu_\vare\\
         (1+\alpha^2)Y_{s+1} - \alpha(Y_{s}+Y_{s+2}) - (1-\alpha)\mu_\vare \\
      \end{bmatrix} \\
 &=\sqrt{n}(B_n^{\dag\dag})^{-1}
    \begin{bmatrix}
      (\halpha_n^{\,\dag\dag}+\alpha)Y_{s}
         - (Y_{s-1}+ Y_{s+1})
         + \hmuen^{\,\dag\dag}
         & \alpha-1\\
       (\halpha_n^{\,\dag\dag}+\alpha)Y_{s+1}
         - (Y_{s}+ Y_{s+2})
         + \hmuen^{\,\dag\dag}
         & \alpha-1 \\
    \end{bmatrix}
    \begin{bmatrix}
      \halpha_n^{\,\dag\dag}-\alpha \\
      \hmuen^{\,\dag\dag}-\mu_\vare \\
    \end{bmatrix} \\
 &\phantom{=\;}
    +\sqrt{n}(B_n^{\dag\dag})^{-1}
       \big( B^{\dag\dag}- B_n^{\dag\dag} \big) (B^{\dag\dag})^{-1}
        \begin{bmatrix}
         (1+\alpha^2)Y_{s} - \alpha(Y_{s-1}+Y_{s+1}) - (1-\alpha)\mu_\vare\\
         (1+\alpha^2)Y_{s+1} - \alpha(Y_{s}+Y_{s+2}) - (1-\alpha)\mu_\vare \\
         \end{bmatrix}.
 \end{align*}
Hence
 \begin{align}\label{SEGED62}
   &\begin{bmatrix}
      \sqrt{n}\big(\htheta_{1,n}^{\,\dag\dag} - \lim_{k\to\infty}\htheta_{1,k}^{\,\dag\dag}\big) \\
      \sqrt{n}\big(\htheta_{2,n}^{\,\dag\dag} - \lim_{k\to\infty}\htheta_{2,k}^{\,\dag\dag}\big) \\
    \end{bmatrix}
   =  D_{n,\alpha,\vare}
       \begin{bmatrix}
        \sqrt{n}(\halpha_n^{\,\dag\dag}-\alpha)  \\
        \sqrt{n}(\hmuen^{\,\dag\dag}-\mu_\vare)  \\
       \end{bmatrix}
  \end{align}
 holds asymptotically as \ $n\to\infty$ \ with probability one, where
 \begin{align*}
  &D_{n,\alpha,\vare}
    :=
    (B_n^{\dag\dag})^{-1}
      \begin{bmatrix}
      (\halpha_n^{\,\dag\dag}+\alpha)Y_{s}
         - Y_{s-1} -Y_{s+1}
         + \hmuen^{\,\dag\dag}
         & \alpha-1\\
       (\halpha_n^{\,\dag\dag}+\alpha)Y_{s+1}
         - Y_{s} - Y_{s+2}
         + \hmuen^{\,\dag\dag}
         & \alpha-1 \\
    \end{bmatrix}\\
    &+ (B_n^{\dag\dag})^{-1}
       \begin{bmatrix}
         -(\halpha_n^{\,\dag\dag}+\alpha) & 1 \\
         1 & -(\halpha_n^{\,\dag\dag}+\alpha) \\
       \end{bmatrix}
        (B^{\dag\dag})^{-1}
        \begin{bmatrix}
         (1+\alpha^2)Y_{s} - \alpha(Y_{s-1}+Y_{s+1}) - (1-\alpha)\mu_\vare & 0\\
         (1+\alpha^2)Y_{s+1} - \alpha(Y_{s}+Y_{s+2}) - (1-\alpha)\mu_\vare & 0 \\
        \end{bmatrix}.
  \end{align*}
By \eqref{Strong_consistency22} and \eqref{Strong_consistency23}, we have \ $D_{n,\alpha,\vare}$
 \ converges almost surely as \ $n\to\infty$ \ to
 \begin{align*}
   & (B^{\dag\dag})^{-1}
     \begin{bmatrix}
       2\alpha Y_{s} - Y_{s-1} - Y_{s+1} +\mu_\vare & \alpha-1 \\
       2\alpha Y_{s+1} - Y_{s} - Y_{s+2} +\mu_\vare & \alpha-1  \\
     \end{bmatrix} \\
   & + (B^{\dag\dag})^{-1}
    \begin{bmatrix}
      -2\alpha & 1 \\
      1 & -2\alpha \\
    \end{bmatrix}
     (B^{\dag\dag})^{-1}
      \begin{bmatrix}
        (1+\alpha^2)Y_{s} - \alpha(Y_{s-1}+Y_{s+1}) - (1-\alpha)\mu_\vare & 0 \\
        (1+\alpha^2)Y_{s+1} - \alpha(Y_{s}+Y_{s+2}) - (1-\alpha)\mu_\vare & 0 \\
      \end{bmatrix}
     =D_{\alpha,\vare}.
 \end{align*}
Hence, by \eqref{SEGED62}, \eqref{CONVERGENCE16} and Slutsky's lemma, we have \eqref{CONVERGENCE17}.
\proofend

\section{The INAR(1) model with innovational outliers}\label{INAR1_innovational}

\subsection{The model and some preliminaries}

In this section we introduce INAR(1) models contaminated with innovational outliers and we also
 give some preliminaries.

\begin{Def}
Let \ $(\vare_\ell)_{\ell\in\NN}$ \ be an i.i.d.\ sequence of non-negative integer-valued random
 variables.
A stochastic process \ $(Y_k)_{k\in\ZZ_+}$ \ is called an INAR(1) model with finitely many
 innovational outliers if
 \[
   Y_k = \sum_{j=1}^{Y_{k-1}}\xi_{k,j}+\eta_k,  \qquad  k\in\NN,
 \]
 where for all \ $k\in\NN$, \ $(\xi_{k,j})_{j\in\NN}$ \ is a sequence of i.i.d.\ Bernoulli random
 variables with mean \ $\alpha\in(0,1)$ \ such that these sequences are mutually independent and
 independent of the sequence \ $(\vare_\ell)_{\ell\in\NN}$, \ and \ $Y_0$ \ is a non-negative
 integer-valued random variable independent of the sequences \ $(\xi_{k,j})_{j\in\NN}$, \ $k\in\NN$,
 \ and \ $(\vare_\ell)_{\ell\in\NN}$, \ and
 \[
   \eta_k:=\vare_k+\sum_{i=1}^I\delta_{k,s_i}\theta_i,
         \qquad k\in\ZZ_+,
 \]
 where \ $I\in\NN$ \ and \ $s_i,\;\theta_i\in\NN$, \ $i=1,\ldots,I$.
\ We assume that \ $\EE Y_0^2<\infty$ and that \ $\EE\vare_1^2<\infty$, \ $\PP(\vare_1\ne 0)>0$.
\end{Def}

In case of one (innovational) outlier a more suitable representation of \ $Y$ \ is given in the
 following proposition.

\begin{Pro}\label{LEMMA3_Decomposition}
Let \ $(Y_k)_{k\in\ZZ_+}$ \ be an INAR(1) model with one innovational outlier
 \ $\theta_1:=\theta$ \ at time point \ $s_1:=s$.
\ Then for all \ $\omega\in\Omega$ \ and \ $k\in\ZZ_+$, \ $Y_k(\omega)=X_k(\omega)+Z_k(\omega)$,
 \ where \ $(X_k)_{k\in\ZZ_+}$ \ is an INAR(1) model given by
  \[
   X_k:=
        \sum_{j=1}^{X_{k-1}}\xi_{k,j} + \vare_k,  \qquad k\in\NN,
 \]
 with \ $X_0:=Y_0$, \ and
 \begin{align}\label{Decomposition}
     Z_k:=
         \begin{cases}
           0 & \quad \text{if \ $k=0,1,\ldots,s-1$,}\\
           \theta & \quad \text{if \ $k=s$,}\\
            \sum_{j=X_{k-1}+1}^{X_{k-1}+Z_{k-1}}\xi_{k,j} & \quad \text{if \ $k\geq s+1$.}
         \end{cases}
 \end{align}
Moreover, the processes \ $X$ \ and \ $Z$ \ are independent, and
 \ $\PP(\lim_{k\to\infty}Z_k=0)=1$ \ and \ $Z_k\stackrel{L_p}\longrightarrow 0$ \ as \ $k\to\infty$
 \ for all \ $p\in\NN$, \ where \ $\stackrel{L_p}\longrightarrow$ \ denotes convergence in \ $L_p$.
\end{Pro}

\noindent{\bf Proof.}
Clearly, \ $Y_j=X_j+Z_j$ \ for \ $j=0,1,\ldots,s-1$, \ and
 \begin{align*}
   &Y_s=\sum_{j=1}^{Y_{s-1}}\xi_{s,j}+\eta_s
      =\sum_{j=1}^{X_{s-1}}\xi_{s,j}+\vare_s+\theta
      =X_{s}+\theta=X_s+Z_s,\\
   &Y_{s+1}=\sum_{j=1}^{Y_s}\xi_{s+1,j}+\eta_{s+1}
      =\sum_{j=1}^{X_s+Z_s}\xi_{s+1,j}+\vare_{s+1}
      =\sum_{j=1}^{X_s}\xi_{s+1,j}
       +\sum_{j=X_s+1}^{X_s+Z_s}\xi_{s+1,j}
       +\vare_{s+1}\\
   &\phantom{Y_{s+1}=}
       = X_{s+1} + \sum_{j=X_s+1}^{X_s+Z_s}\xi_{s+1,j}
       = X_{s+1} + Z_{s+1}.
 \end{align*}
By induction, we easily conclude that \ $Y_k(\omega)=X_k(\omega)+Z_k(\omega)$ \ for all \ $\omega\in\Omega$
 \ and \ $k\in\ZZ_+$.

In proving the independence of the processes \ $X$ \ and \ $Z$, \ it is enough to check
 that the conditions of Lemma \ref{LEMMA_INNOVATIONAL_INDEPENDENCE} (see Appendix) are satisfied.
For all \ $n>s$, \ $i_{n-1},i_n,j_{n-1},j_n\in\ZZ_+$ \ and for all
 \ $B\in\sigma(\xi_{i,j}:i=1,\ldots,n-2,\;j\in\NN)$ \ with the property that the event
 \ $A:=\{X_{n-1}=i_{n-1},Z_{n-1}=j_{n-1}\}\cap B$ \ has positive probability, we get
 \begin{align}\label{SEGED105_elozetes}
  \begin{split}
   \PP(X_n=i_n,Z_n=j_n\mid A)
    & =\PP\left(\sum_{j=1}^{i_{n-1}}\xi_{n,j}+\vare_n=i_n,
                \sum_{j=i_{n-1}+1}^{i_{n-1}+j_{n-1}}\xi_{n,j}=j_n \right)\\
    & =\PP\left(\sum_{j=1}^{i_{n-1}}\xi_{n,j}+\vare_n=i_n\right)
      \PP\left(\sum_{j=i_{n-1}+1}^{i_{n-1}+j_{n-1}}\xi_{n,j}=j_n \right),
   \end{split}
 \end{align}
 where we used the measurability of \ $(X_{n-1},Z_{n-1})$ \ with respect to the
 $\sigma$--algebra \ $\sigma(\xi_{i,j}:i=1,\ldots,n-1,\;j\in\NN)$
 \ and that the random variables \ $\vare_n$, \ $(\xi_{n,1},\ldots,\xi_{n,i_{n-1}})$ \ and
 \ $(\xi_{n,i_{n-1}+1},\ldots,\xi_{n,i_{n-1}+j_{n-1}})$ \ are independent of this \ $\sigma$--algebra
 and also from each other.
Hence, for all \ $n>s$,
 \begin{align}\label{SEGED105}
   \PP(X_n=i_n,Z_n=j_n\mid A)
        =\PP(X_n=i_n,Z_n=j_n\mid X_{n-1}=i_{n-1},Z_{n-1}=j_{n-1}).
 \end{align}
Since \ $Z_0=Z_1=\cdots=Z_{s-1}=0$, \ $Z_s=\theta$, \ and \ $(X_n)_{n\in\ZZ_+}$ \ is a Markov chain,
 we have \eqref{SEGED105} is satisfied also for \ $n=1,2,\ldots,s$, \ which yields that
 \ $(X_n,Z_n)_{n\in\ZZ_+}$ \ is a Markov chain.
Since \ $Z_0=0$, \ $X_0$ \ and \ $Z_0$ \ are independent.
Similar arguments along with the result in \eqref{SEGED105_elozetes}, with the special
 choice \ $B:=\Omega$ \ lead to
 \begin{align*}
    \PP&(X_n=i_n,Z_n=j_n\mid X_{n-1}=i_{n-1},Z_{n-1}=j_{n-1}) \\
       &=\PP\left(\sum_{j=1}^{i_{n-1}}\xi_{n,j}+\vare_n=i_n\,\Bigg\vert\; X_{n-1}=i_{n-1} \right)
         \PP\left(\sum_{j=1}^{j_{n-1}}\xi_{n,j+i_{n-1}}=j_n \,\Bigg\vert\; Z_{n-1}=j_{n-1}\right) \\
       & = \PP(X_n=i_n\mid X_{n-1}=i_{n-1}) \PP(Z_n=j_n\mid Z_{n-1}=j_{n-1}),
 \end{align*}
 which yields that the conditions of Lemma
 \ref{LEMMA_INNOVATIONAL_INDEPENDENCE} are satisfied.

Since
 \[
   Z_{k+1}=\sum_{j=X_{k}+1}^{X_{k}+Z_{k}}\xi_{k+1,j}
           \leq \sum_{j=X_{k}+1}^{X_{k}+Z_{k}}1
           = Z_k,\qquad k\geq s,
 \]
 the sequence \ $(Z_k(\omega))_{k\geq s+1}$ \ is monotone decreasing for all \ $\omega\in\Omega$.
\ Using the fact that \ $Z_k\geq0$, \ $k\in\NN$, \ we have \ $(Z_k(\omega))_{k\in\ZZ_+}$ \ converges
 for all \ $\omega\in\Omega$.
Hence, if we check that \ $Z_k$ \ converges in probability to \ $0$ \ as
 \ $k\to\infty$, \ then, by Riesz's theorem, we get \ $\PP(\lim_{k\to\infty}Z_k=0)=1$.
\ Let \ $\cF_k^{X,Z}$ \ be the $\sigma$--algebra generated by the random variables
 \ $X_0,X_1,\ldots,X_k$ \ and \ $Z_0,Z_1,\ldots,Z_k$.
\ Using that \ $\EE(Z_k\mid \cF^{X,Z}_{k-1}) = \alpha Z_{k-1}$, \ $k\geq s+1$, \ we get
 \ $\EE Z_k=\alpha\EE Z_{k-1}$, \ $k\geq s+1$, \ and hence
 \ $\EE Z_{s+k}=\alpha^k \EE Z_s=\theta\alpha^k$, \ $k\geq 0$.
\ For all \ $\vare>0$, \ by Markov's inequality,
 \begin{align*}
     \PP(Z_{s+k}\geq\vare) \leq \frac{\EE Z_{s+k}}{\vare}
                            =\frac{\theta\alpha^k}{\vare}
                            \to0
     \qquad \text{as \ $k\to\infty$},
 \end{align*}
 as desired.
We note that the fact that \ $\PP(\lim_{k\to\infty}Z_k=0)=1$ \ is in accordance with
 Theorem 2 in Chapter I in Athreya and Ney \cite{AthNey}.

Since the sequence \ $(Z_k(\omega))_{k\geq s+1}$ \ is monotone decreasing  for all \ $\omega\in\Omega$,
 \ we get for all \ $p\in\NN$ \ and for any constant \ $M>0$, \ the sequence
 \ $\left(\vert Z_k\vert^p\bone_{\{\vert Z_k\vert\geq M\}}\right)_{k\geq s+1}$ \ is monotone
 decreasing.
\ Hence
  \[
    \sup_{k\geq s+1} \EE\left(\vert Z_k\vert^p\bone_{\{\vert Z_k\vert\geq M\}}\right)
       = \EE\left(\vert Z_{s+1}\vert^p\bone_{\{\vert Z_{s+1}\vert\geq M\}}\right)
       \to 0
       \qquad \text{as \ $M\to\infty$,}
  \]
 which yields the uniformly integrability of \ $(Z_k^p)_{k\in\NN}$.
\ By Lemma \ref{LEMMA4_LP} (see Appendix), we conclude that \ $Z_k\lpmean 0$ \ as
 \ $k\to\infty$, \ i.e., \ $\lim_{k\to\infty}\EE Z_k^p=0$.
\proofend

For our later purposes we need the following lemma about the explicit forms of the first
 and second moments of the process \ $Z$.

\begin{Lem}\label{LEMMA_Z}
We have
 \begin{align}\label{SEGED_Z1}
   &\EE Z_{s+k}=\theta\alpha^k,\qquad k\in\ZZ_+,\\ \label{SEGED_Z2}
   &\EE Z_{s+k}^2=\theta^2\alpha^{2k}-\theta\alpha^k(\alpha^k-1),\qquad k\in\ZZ_+,\\ \label{SEGED_Z3}
   &\EE(Z_{s+k-1}Z_{s+k})
                 = \alpha\EE Z_{s+k-1}^2
                 =\theta^2\alpha^{2k-1} - \theta\alpha^k(\alpha^{k-1}-1),
                 \qquad k\in\NN.
 \end{align}
\end{Lem}

\noindent{\bf Proof.}
Recall that \ $\cF_k^{X,Z}$ \ denotes the $\sigma$--algebra generated by the random variables
 \ $X_0,X_1,\ldots,X_k$ \ and \ $Z_0,Z_1,\ldots,Z_k$.
\ Using that \ $\EE(Z_k\mid \cF^{X,Z}_{k-1}) = \alpha Z_{k-1}$, \ $k\geq s+1$, \ we get
 \ $\EE Z_k=\alpha\EE Z_{k-1}$, \ $k\geq s+1$, \ and hence
 \ $\EE Z_{s+k}=\alpha^k \EE Z_s=\theta\alpha^k$, \ $k\in\ZZ_+$.
\ Since \ $\alpha\in(0,1)$, \ we have \ $\lim_{k\to\infty}\EE Z_k=0$.
\ Moreover, using that
 \begin{align}\label{SEGED_Z4}
  \EE((Z_k-\alpha Z_{k-1})^2\mid \cF^{X,Z}_{k-1})
      = \EE\left(\left(\sum_{j=X_{k-1}+1}^{X_{k-1}+Z_{k-1}}
           (\xi_{k,j}-\alpha)\right)^2\,\Big\vert\,\cF^{X,Z}_{k-1}\right)
      =\alpha(1-\alpha) Z_{k-1},
      \qquad k\geq s+1,
 \end{align}
 we get
 \[
   \EE(Z_k^2\mid \cF^{X,Z}_{k-1})
      = \EE\Big(\big((Z_k-\alpha Z_{k-1})+\alpha Z_{k-1}\big)^2\mid \cF^{X,Z}_{k-1}\Big)
      =\alpha(1-\alpha) Z_{k-1}+\alpha^2 Z_{k-1}^2,
     \qquad k\geq s+1,
 \]
 and hence \ $\EE Z_k^2=\alpha^2\EE Z_{k-1}^2+\alpha(1-\alpha)\EE Z_{k-1}$, \ $k\geq s+1$.
\ Then
 \begin{align*}
    \begin{bmatrix}
      \EE Z_k \\
      \EE Z_k^2 \\
    \end{bmatrix}
     =\begin{bmatrix}
        \alpha & 0 \\
        \alpha(1-\alpha) & \alpha^2 \\
      \end{bmatrix}
      \begin{bmatrix}
      \EE Z_{k-1} \\
      \EE Z_{k-1}^2 \\
    \end{bmatrix},
    \qquad k\geq s+1,
 \end{align*}
 and hence, by an easy calculation, for all \ $k\geq 0$,
 \begin{align*}
    \begin{bmatrix}
      \EE Z_{s+k} \\
      \EE Z_{s+k}^2 \\
    \end{bmatrix}
    & =\begin{bmatrix}
        \alpha^k & 0 \\
        (1-\alpha)\alpha^k\sum_{\ell=0}^{k-1}\alpha^\ell & \alpha^{2k} \\
      \end{bmatrix}
      \begin{bmatrix}
       \EE Z_s \\
       \EE Z_s^2 \\
      \end{bmatrix}
     = \begin{bmatrix}
        \alpha^k & 0 \\
        (1-\alpha)\alpha^k\frac{\alpha^k-1}{\alpha-1} & \alpha^{2k} \\
      \end{bmatrix}
      \begin{bmatrix}
       \theta \\
       \theta^2 \\
      \end{bmatrix}  \\
   & =  \begin{bmatrix}
       \theta\alpha^k \\
       \theta^2\alpha^{2k}-\theta\alpha^k(\alpha^k-1) \\
      \end{bmatrix}.
 \end{align*}
Finally, for all \ $k\in\NN$,
 \begin{align*}
   \EE(Z_{s+k-1}Z_{s+k})
      & = \EE\big(\EE(Z_{s+k-1}Z_{s+k}\mid \cF^Z_{s+k-1})\big)
        = \EE\big(Z_{s+k-1}\EE(Z_{s+k}\mid \cF^Z_{s+k-1})\big)
        = \EE\big(Z_{s+k-1}\alpha Z_{s+k-1}\big) \\
      & = \alpha\EE Z_{s+k-1}^2,
 \end{align*}
 which yields \eqref{SEGED_Z3}.
\proofend

In case of two (innovational) outliers a similar representation of \ $Y$ \ is given in the
 following proposition.

\begin{Pro}\label{LEMMA3_Decomposition2}
Let \ $(Y_k)_{k\in\ZZ_+}$ \ be an INAR(1) model with two innovational outliers
 \ $\theta_1$ \ and \ $\theta_2$ \ at time points \ $s_1$ \ and \ $s_2$, \ $s_1 < s_2$,
 \[
   Y_k=
        \sum_{j=1}^{Y_{k-1}}\xi_{k,j} + \eta_k, \qquad k\in\NN,
 \]
 where for all \ $k\in\NN$, \ $(\xi_{k,j})_{j\in\NN}$ \ is a sequence of i.i.d.\ Bernoulli random
 variables with mean \ $\alpha\in(0,1)$ \ such that these sequences are mutually independent and
 independent of the sequence \ $(\vare_\ell)_{\ell\in\NN}$, \ and \ $Y_0$ \ is a non-negative
 integer-valued random variable independent of the sequences \ $(\xi_{k,j})_{j\in\NN}$, \ $k\in\NN$,
 \ and \ $(\vare_\ell)_{\ell\in\NN}$, \ and \ $\eta_k:=\vare_k+\delta_{k,s_1}\theta_1+\delta_{k,s_2}\theta_2$,
 \ $k\in\ZZ_+$.
\ Then for all \ $\omega\in\Omega$ \ and \ $k\in\ZZ_+$,
 \ $Y_k(\omega)=X_k(\omega)+Z_k^{(1)}(\omega)+Z_k^{(2)}(\omega)$,
 \ where \ $(X_k)_{k\in\ZZ_+}$ \ is an INAR(1) model given by
  \[
   X_k:=
        \sum_{j=1}^{X_{k-1}}\xi_{k,j} + \vare_k,  \qquad k\in\NN,
 \]
 with \ $X_0:=Y_0$, \ and
 \begin{align}\label{Decomposition2}
     Z_k^{(1)}:=
         \begin{cases}
           0 & \quad \text{if \ $k=0,1,\ldots,s_1-1$,}\\
           \theta_1 & \quad \text{if \ $k=s_1$,}\\
            \sum_{j=X_{k-1}+1}^{X_{k-1}+Z_{k-1}^{(1)}}\xi_{k,j} & \quad \text{if \ $k\geq s_1+1$,}
         \end{cases}
 \end{align}
 and
 \begin{align}\label{Decomposition3}
     Z_k^{(2)}:=
         \begin{cases}
           0 & \quad \text{if \ $k=0,1,\ldots,s_2-1$,}\\
           \theta_2 & \quad \text{if \ $k=s_2$,}\\
            \sum_{j=X_{k-1}+Z_{k-1}^{(1)}+1}^{X_{k-1}+Z_{k-1}^{(1)}+Z_{k-1}^{(2)}}\xi_{k,j}
               & \quad \text{if \ $k\geq s_2+1$.}
         \end{cases}
 \end{align}
Moreover, the processes \ $X$, \ $Z^{(1)}$ \ and \ $Z^{(2)}$ \ are (pairwise) independent, and
 \ $\PP(\lim_{k\to\infty}Z_k^{(i)}=0)=1$, \ $i=1,2$, \ and
 \ $Z_k^{(i)}\stackrel{L_p}\longrightarrow 0$ \ as \ $k\to\infty$ \ for all \ $p\in\NN$, \ $i=1,2$.
\end{Pro}

\noindent{\bf Proof.}
The proof is the very same as the proof of Proposition \ref{LEMMA3_Decomposition}.
We only note that the independence of \ $Z^{(1)}$ \ and \ $Z^{(2)}$ \ follows by the definitions
 of the processes \ $Z^{(1)}$ \ and \ $Z^{(2)}$.
\proofend

In the sequel we denote by \ $\cF_k^Y$ \ the $\sigma$--algebra generated by the random
 variables \ $Y_0,Y_1,\ldots,Y_k$.
\ For all \ $n\in\NN$, \ $y_0,\ldots,y_n\in\RR$ \ and \ $\omega\in\Omega$, \ let us introduce
 the notations
 \begin{align*}
    {\bf Y}_n(\omega):=(Y_0(\omega),Y_1(\omega),\ldots,Y_n(\omega)),\qquad\;
    {\bf Y}_n:=(Y_0,Y_1,\ldots,Y_n), \qquad\;
    {\bf y}_n:=(y_0,y_1,\ldots,y_n).
 \end{align*}

\subsection{One outlier, estimation of the mean of the offspring distribution
            and the outlier's size}

First we suppose that \ $I=1$ \ and that \ $s_1:=s$ \ is known.
We concentrate on the CLS estimation of the parameter \ $(\alpha,\theta)$,
 \ where \ $\theta:=\theta_1$.
\ An easy calculation shows that
 \begin{align}\label{SEGED_COND1}
   \EE(Y_k\mid\cF^Y_{k-1})
       = \alpha Y_{k-1} + \EE\eta_k
       = \alpha Y_{k-1} + \mu_\vare + \delta_{k,s}\theta,
      \qquad k\in\NN.
 \end{align}
Hence for \ $n\geq s$, \ $n\in\NN$,
 \begin{align}\label{Innovation_CLSE}
   \begin{split}
    \sum_{k=1}^n&\big(Y_k-\EE(Y_k\mid \cF^Y_{k-1})\big)^2\\
        & =\sums \big(Y_k-\alpha Y_{k-1}-\mu_\vare\big)^2
           + \big(Y_s-\alpha Y_{s-1}-\mu_\vare-\theta\big)^2.
    \end{split}
 \end{align}
For all \ $n\geq s$, \ $n\in\NN$, \ we define the function \ $Q_n:\RR^{n+1}\times\RR^2\to\RR$, \ as
 \begin{align*}
    Q_n({\bf y}_n;\alpha',\theta')
        :=\sums \big(y_k-\alpha' y_{k-1}-\mu_\vare\big)^2
           + \big(y_s-\alpha' y_{s-1}-\mu_\vare-\theta'\big)^2,
 \end{align*}
 for all \ ${\bf y}_n\in\RR^{n+1}$ \ and \ $\alpha',\theta'\in\RR$.
\ By definition, for all \ $n\geq s$, \ a CLS estimator for
 the parameter \ $(\alpha,\theta)\in(0,1)\times\NN$ \ is a measurable function
 \ $(\talpha_n,\ttheta_n):S_n\to\RR^2$ \ such that
 \begin{align*}
   Q_n({\bf y}_n;\,&\talpha_n({\bf y}_n),\ttheta_n({\bf y}_n))
       = \inf_{(\alpha',\theta')\in\RR^2}Q_n({\bf y}_n;\alpha',\theta')
       \qquad \forall\;\;  {\bf y}_n\in S_n,
 \end{align*}
 where \ $S_n$ \ is suitable subset of \ $\RR^{n+1}$ \ (defined in the proof of
 Lemma \ref{LEMMA15}).
We note that we do not define the CLS estimator
 \ $(\talpha_n,\ttheta_n)$ \ for all samples \ ${\bf y}_n\in \RR^{n+1}$.
\ We have
 \begin{align*}
   &\frac{\partial Q_n}{\partial \alpha'}({\bf y}_n;\alpha',\theta')
     =-2\sums \big(y_k-\alpha' y_{k-1}-\mu_\vare\big)y_{k-1}
             - 2\big(y_s-\alpha' y_{s-1}-\mu_\vare-\theta'\big)y_{s-1}, \\
   &\frac{\partial Q_n}{\partial \theta'}({\bf y}_n;\alpha',\theta')
      =-2\big(y_s-\alpha' y_{s-1}-\mu_\vare-\theta'\big).
 \end{align*}

The next lemma is about the existence and uniqueness of the CLS estimator of \ $(\alpha,\theta)$.

\begin{Lem}\label{LEMMA15}
There exist subsets \ $S_n\subset\RR^{n+1}$, $n\geq \max(3,s+1)$ \ with the following properties:
 \begin{enumerate}
  \item[\upshape{(i)}]
   there exists a unique CLS estimator
   \ $(\talpha_n,\ttheta_n):S_n\to\RR^2$,
  \item[\upshape{(ii)}]
   for all \ ${\bf y}_n\in S_n$, \ the system of equations
  \begin{align}\label{Innovation_CLSE_EQ}
    \frac{\partial Q_n}{\partial \alpha'}({\bf y}_n;\alpha',\theta')=0,\qquad
    \frac{\partial Q_n}{\partial \theta'}({\bf y}_n;\alpha',\theta')=0,
 \end{align}
 has the unique solution
   \begin{align}\label{SEGED3}
   &\talpha_n({\bf y}_n)
     =\frac{\DS\sums(y_k-\mu_\vare)y_{k-1}}
           {\DS\sums y_{k-1}^2},\\[1mm] \label{SEGED4}
   &\ttheta_n({\bf y}_n)
      = y_s-\talpha_n({\bf y}_n) y_{s-1} -\mu_\vare,
 \end{align}
  \item[\upshape{(iii)}]
  ${\bf Y}_n\in S_n$ \ holds asymptotically as \ $n\to\infty$ \ with probability one.
 \end{enumerate}
\end{Lem}

\noindent{\bf Proof.}
One can easily check that the unique solution of the system of equations \eqref{Innovation_CLSE_EQ}
 takes the form \eqref{SEGED3} and \eqref{SEGED4} whenever \ $\DS\sums y_{k-1}^2>0$.

Next we prove that the function
  \ $\RR^2 \ni (\alpha',\theta')
     \mapsto Q_n({\bf y}_n;\alpha',\theta')$
  \ is strictly convex for all \ ${\bf y}_n\in S_n$, \ where
 \[
   S_n := \left\{{\bf y}_n\in\RR^{n+1} : \DS\sums y_{k-1}^2 > 0 \right\}.
 \]
For this it is enough to check that the \ $(2\times 2)$-matrix
 \[
  H_n({\bf y}_n;\alpha',\theta')
   :=\begin{bmatrix}
    \frac{\partial^2 Q_n}{\partial(\alpha')^2}
    & \frac{\partial^2 Q_n}{\partial\theta'\partial\alpha'}  \\
    \frac{\partial^2 Q_n}{\partial\alpha'\partial\theta'}
    & \frac{\partial^2 Q_n}{\partial(\theta')^2} \\
   \end{bmatrix}
     ({\bf y}_n;\alpha',\theta')
 \]
 is (strictly) positive definite for all \ ${\bf y}_n\in S_n$, \ see, e.g.,
 Berkovitz \cite[Theorem 3.3, Chapter III]{Ber}.
For all \ ${\bf y}_n\in\RR^{n+1}$ \ and \ $(\alpha',\theta')\in\RR^2$,
 \begin{align*}
   &\frac{\partial^2 Q_n}{\partial(\alpha')^2}({\bf y}_n;\alpha',\theta')
       = 2\sums y_{k-1}^2 + 2y_{s-1}^2
       =2\sum_{k=1}^n y_{k-1}^2,\\
   &\frac{\partial^2 Q_n}{\partial\alpha'\partial\theta'}({\bf y}_n;\alpha',\theta')
       =\frac{\partial^2 Q_n}{\partial\theta'\partial\alpha'}({\bf y}_n;\alpha',\theta')
       = 2y_{s-1},\\
   & \frac{\partial^2 Q_n}{\partial(\theta')^2}({\bf y}_n;\alpha',\theta')=2.
 \end{align*}
Then \ $H_n({\bf y}_n;\alpha',\theta')$ \ has leading principal minors
 \[
   2\sum_{k=1}^n y_{k-1}^2 \qquad \text{and} \qquad
    4\DS\sums y_{k-1}^2,
 \]
 which are positive for all \ ${\bf y}_n\in S_n$.
\ Hence \ $H_n({\bf y}_n;\alpha',\theta')$ \ is (strictly) positive definite
 for all \ ${\bf y}_n\in S_n$.

Since the function \ $\RR^2 \ni (\alpha',\theta')
     \mapsto Q_n({\bf y}_n;\alpha',\theta')$
  \ is strictly convex for all \ ${\bf y}_n\in S_n$ \ and the system of equations
 \eqref{Innovation_CLSE_EQ} has a unique solution for all \ ${\bf y}_n\in S_n$,
 we get the function in question attains its (global) minimum at this unique solution,
 which yields (i) and (ii).

Next we check that \ ${\bf Y}_n\in S_n$ \ holds asymptotically as \ $n\to\infty$ \ with probability one.
By Proposition \ref{LEMMA3_Decomposition}, we get
 \[
    \sums Y_{k-1}^2
      =     \sums X_{k-1}^2 + 2 \sums X_{k-1}Z_{k-1} + \sums Z_{k-1}^2,
      \qquad n\geq s+1.
 \]
 Using again Proposition \ref{LEMMA3_Decomposition} and \eqref{Ergodic2}, we have
 \begin{align*}
   \PP\left(\lim_{n\to\infty}\frac{1}{n}\sums Z_{k-1}^2=0\right)=1, \qquad\qquad
   \PP\left(\lim_{n\to\infty}\frac{1}{n}\sums X_{k-1}^2=\EE \widetilde X^2\right)=1.
 \end{align*}
By Cauchy-Schwartz's inequality,
 \begin{align*}
   \frac{1}{n}\left\vert \sum_{k=s+1}^nX_{k-1}Z_{k-1} \right\vert
      \leq \sqrt{\frac{1}{n}\sum_{k=s+1}^nX_{k-1}^2\frac{1}{n}\sum_{k=s+1}^nZ_{k-1}^2}
      \to \sqrt{\EE\widetilde X^2}\sqrt{\lim_{n\to\infty}\frac{1}{n}\sum_{k=s+1}^nZ_{k-1}^2}=0,
 \end{align*}
 and hence
 \begin{align}\label{SEGED98}
   \PP\left(\lim_{n\to\infty}\frac{1}{n}\sums X_{k-1}Z_{k-1}=0\right)=1.
 \end{align}
Then
 \[
    \PP\left(\lim_{n\to\infty}\frac{1}{n}\sum_{k=1}^n Y_{k-1}^2 = \EE\widetilde X^2\right)=1,
 \]
 which implies that
 \[
    \PP\left(\lim_{n\to\infty}\sum_{k=1}^n Y_{k-1}^2 = \infty\right)=1.
 \]
Hence \ ${\bf Y}_n\in S_n$ \ holds asymptotically as \ $n\to\infty$ \ with probability one.
\proofend

By Lemma \ref{LEMMA15}, \ $(\talpha_n({\bf Y}_n),
                           \ttheta_n({\bf Y}_n))$
 \ exists uniquely asymptotically as \ $n\to\infty$ \ with probability one.
In the sequel we will simply denote it by \ $(\talpha_n,\ttheta_n)$.

The next result shows that \ $\talpha_n$ \ is a strongly consistent estimator
 of \ $\alpha$, \ whereas \ $\ttheta_n$ \ fails to be a strongly consistent estimator
 of \ $\theta$.

\begin{Thm}\label{THEOREM7}
Consider the CLS estimators
 \ $(\talpha_n,\ttheta_n)_{n\in\NN}$ \ of the parameter
 \ $(\alpha,\theta)\in(0,1)\times\NN$.
\ The sequence \ $(\talpha_n)_{n\in\NN}$ \ is strongly consistent for all
 \ $(\alpha,\theta)\in(0,1)\times\NN$, \ i.e.,
 \begin{align}\label{Strong_consistency3}
   \PP(\lim_{n\to\infty}\talpha_n=\alpha)=1, \qquad \forall\;(\alpha,\theta)\in(0,1)\times\NN,
 \end{align}
 whereas the sequence \ $(\ttheta_n)_{n\in\NN}$ \ is not strongly consistent for any
  \ $(\alpha,\theta)\in(0,1)\times\NN$, \ namely,
 \begin{align}\label{Strong_consistency4}
   \PP\left(
        \lim_{n\to\infty}\ttheta_n
         = Y_s-\alpha Y_{s-1}-\mu_\vare
        \right)=1,\qquad \forall\;(\alpha,\theta)\in(0,1)\times\NN.
 \end{align}
\end{Thm}

\noindent{\bf Proof.}
By \eqref{SEGED3} and Proposition \ref{LEMMA3_Decomposition}, we have asymptotically as \ $n\to\infty$
 \ with probability one,
 \begin{align*}
   \talpha_n
      &=\frac{\DS\sums (X_k-\mu_\vare+Z_k)(X_{k-1}+Z_{k-1})}
            {\DS\sums (X_{k-1}+Z_{k-1})^2} \\
      &=\frac{\DS\sums (X_k-\mu_\vare)X_{k-1}
              +\DS\sum_{k=s+1}^n(X_k-\mu_\vare)Z_{k-1}
              +\DS\sum_{k=s+1}^nX_{k-1}Z_k
              +\DS\sum_{k=s+1}^nZ_{k-1}Z_k}
            {\DS\sums X_{k-1}^2
              +2\DS\sum_{k=s+1}^nX_{k-1}Z_{k-1}
              +\DS\sum_{k=s+1}^nZ_{k-1}^2}.
 \end{align*}
By \eqref{Strong_consistency_INAR}, to prove \eqref{Strong_consistency3},
 it is enough to check that
  \begin{align} \label{SEGED23}
   \PP\!\left(\!
     \lim_{n\to\infty}\!
              \frac{-(X_s-\mu_\vare)X_{s-1}
                + \sum_{k=s+1}^n(X_k-\mu_\vare)Z_{k-1}
                + \sum_{k=s+1}^nX_{k-1}Z_k
                + \sum_{k=s+1}^nZ_{k-1}Z_k }{n} =0
       \!\right)\!=1,
  \end{align}
 and
  \begin{align} \label{SEGED24}
     \PP\left(\lim_{n\to\infty}
              \frac{-X_{s-1}^2 + 2 \sum_{k=s+1}^nX_{k-1}Z_{k-1}
                                + \sum_{k=s+1}^nZ_{k-1}^2}{n} =0 \right)=1.
  \end{align}
By Proposition \ref{LEMMA3_Decomposition} and Cauchy-Schwartz's inequality, we have
 \[
   \PP\left( \lim_{n\to\infty} \frac{1}{n}\sum_{k=s+1}^nZ_{k-1}^2 = 0 \right)
      = \PP\left( \lim_{n\to\infty} \frac{1}{n}\sum_{k=s+1}^n Z_{k-1}Z_k = 0 \right)
      =1.
 \]
Hence, using also \eqref{SEGED98}, we get \eqref{SEGED23} and \eqref{SEGED24}.
By \eqref{SEGED4} and \eqref{Strong_consistency3}, we get \eqref{Strong_consistency4}.
\proofend

The asymptotic distribution of the CLS estimation is given in the next theorem.

\begin{Thm}\label{Proposition1}
Under the additional assumptions \ $\EE Y_0^3<\infty$ \ and \ $\EE\vare_1^3<\infty$, \ we have
 \begin{align}\label{CONVERGENCE3}
   \sqrt{n}(\talpha_n-\alpha)\distr \cN(0,\sigma_{\alpha,\,\vare}^2)
      \qquad \text{as \ $n\to\infty$,}
  \end{align}
 where \ $\sigma_{\alpha,\,\vare}^2$ \ is defined in \eqref{SEGED_SZIGMA_ALPHA}.
Moreover, conditionally on the value \ $Y_{s-1}$,
 \begin{align}\label{CONVERGENCE4}
  \sqrt{n}\big(\ttheta_n - \lim_{k\to\infty}\ttheta_k\big)
         \distr \cN(0,Y_{s-1}^2\sigma_{\alpha,\,\vare}^2)
     \qquad \text{as \ $n\to\infty$.}
 \end{align}
\end{Thm}

\noindent{\bf Proof.}
By \eqref{SEGED3}, we have
 \begin{align*}
   \talpha_n-\alpha
      =\frac{\DS\sums (Y_k-\alpha Y_{k-1}-\mu_\vare)Y_{k-1}}
            {\DS\sums Y_{k-1}^2}
 \end{align*}
 holds asymptotically as \ $n\to\infty$ \ with probability one.
For all \ $n\geq s+1$, \ by Proposition \ref{LEMMA3_Decomposition}, we have
 \begin{align*}
   &\DS\sums (Y_k-\alpha Y_{k-1}-\mu_\vare)Y_{k-1}\\
   &\qquad = \DS\sums \big[(X_k-\alpha X_{k-1}-\mu_\vare)+(Z_k-\alpha Z_{k-1})\big]
                                            (X_{k-1}+Z_{k-1}) \\
   &\qquad
     =\DS\sums (X_k-\alpha X_{k-1}-\mu_\vare)X_{k-1}
      +\DS\sum_{k=s+1}^n(X_k-\alpha X_{k-1}-\mu_\vare)Z_{k-1}
      +\DS\sum_{k=s+1}^n(Z_k-\alpha Z_{k-1})X_{k-1}\\
   &\qquad \;\;\;\;\;
      +\DS\sum_{k=s+1}^n(Z_k-\alpha Z_{k-1})Z_{k-1},
 \end{align*}
 and
 \begin{align*}
   \DS\sums Y_{k-1}^2
    = \DS\sums X_{k-1}^2
        +2\DS\sum_{k=s+1}^nX_{k-1}Z_{k-1}
        +\DS\sum_{k=s+1}^nZ_{k-1}^2.
 \end{align*}
By \eqref{SEGED91} and \eqref{SEGED11}, we have
 \begin{align}\label{CONVERGENCE5}
  \sqrt{n}\;
    \frac{\DS\sum_{k=1}^n(X_k-\alpha X_{k-1}-\mu_\vare)X_{k-1}}
         {\DS\sum_{k=1}^nX_{k-1}^2}
     \distr \cN(0,\sigma_{\alpha,\,\vare}^2)
    \qquad \text{as \ $n\to\infty$.}
 \end{align}
In what follows we show that
 \begin{align}\label{SEGED12}
   & \frac{1}{\sqrt n}
      \DS\sum_{k=s+1}^n(X_k-\alpha X_{k-1}-\mu_\vare)Z_{k-1}
      \mean 0  \qquad \text{as \ $n\to\infty$,}\\ \label{SEGED13}
   & \frac{1}{\sqrt n}
     \DS\sum_{k=s+1}^n(Z_k-\alpha Z_{k-1})X_{k-1}
      \mean 0  \qquad \text{as \ $n\to\infty$,}\\  \label{SEGED14}
   & \frac{1}{\sqrt n}
     \DS\sum_{k=s+1}^n(Z_k-\alpha Z_{k-1})Z_{k-1}
           \mean 0  \qquad \text{as \ $n\to\infty$,}\\ \label{SEGED15}
   & \PP\left(
        \lim_{n\to\infty}
        \frac{1}{n}\DS\sum_{k=s+1}^nX_{k-1}Z_{k-1}
          =0\right)=1, \\ \label{SEGED16}
   & \PP\left(
       \lim_{n\to\infty}
       \frac{1}{n}
       \DS\sum_{k=s+1}^nZ_{k-1}^2
       = 0\right)=1,
 \end{align}
 where \ $\mean$ \ denotes convergence in \ $L_1$.
\ We recall that if \ $(\eta_n)_{n\in\NN}$ \ is a sequence of square integrable random
 variables such that
 \ $\lim_{n\to\infty}\EE\eta_n=0$ \ and \ $\lim_{n\to\infty}\EE\eta_n^2=0$, \ then
 \ $\eta_n$ \ converges in \ $L_2$ \ and hence in \ $L_1$ \ to \ $0$ \ as \ $n\to\infty$.
\ Hence to prove \eqref{SEGED12} it is enough to check that
 \begin{align} \label{SEGED17}
   &\lim_{n\to\infty}
      \frac{1}{\sqrt n}\sum_{k=s+1}^n\EE\big[(X_k-\alpha X_{k-1}-\mu_\vare)Z_{k-1}\big] =0,\\ \label{SEGED18}
   &\lim_{n\to\infty}
     \frac{1}{n} \EE\left(\sum_{k=s+1}^n(X_k-\alpha X_{k-1}-\mu_\vare)Z_{k-1}\right)^2=0.
 \end{align}
Since \ $\EE X_k=\alpha\EE X_{k-1}+\mu_\vare$, \ $k\in\NN$, \ and the processes \ $X$ \ and \ $Z$
 \ are independent, we have \eqref{SEGED17}.
Similarly, we get
 \begin{align*}
   \EE\left(\sum_{k=s+1}^n(X_k-\alpha X_{k-1}-\mu_\vare)Z_{k-1}\right)^2
     =\sum_{k=s+1}^n\EE(X_k-\alpha X_{k-1}-\mu_\vare)^2\EE Z_{k-1}^2,
     \quad n\geq s+1.
 \end{align*}
By \eqref{SEGED64},
 \begin{align*}
    \EE(X_k & - \alpha X_{k-1}-\mu_\vare)^2
       = \alpha^k(1-\alpha)\EE X_0 + \alpha(1-\alpha^{k-1})\mu_\vare + \sigma_\vare^2,
     \qquad k\in\NN,
 \end{align*}
 and then
 \begin{align*}
  \lim_{k\to\infty}\EE(X_k-\alpha X_{k-1}-\mu_\vare)^2
     = \lim_{k\to\infty} \big(\alpha^k(1-\alpha) \EE X_0
           + \alpha(1-\alpha^{k-1})\mu_\vare+\sigma_\vare^2\big)
     = \alpha\mu_\vare+\sigma_\vare^2.
 \end{align*}
Hence there exists some \ $L>0$ \ such that \ $\EE(X_k-\alpha X_{k-1}-\mu_\vare)^2<L$
 \ for all \ $k\in\NN$.
\ By Proposition \ref{LEMMA3_Decomposition}, \ $\lim_{k\to\infty}\EE Z_k^2=0$, \ and hence
 \ $\lim_{n\to\infty}\frac{1}{n}\sum_{k=1}^n\EE Z_k^2=0$, \ which yields that
 \[
    \frac{1}{n}\sum_{k=s+1}^n\EE(X_k-\alpha X_{k-1}-\mu_\vare)^2\EE Z_{k-1}^2
               \leq \frac{L}{n} \sum_{k=s+1}^n \EE Z_{k-1}^2\to 0
      \qquad \text{as \ $n\to\infty$.}
 \]
To prove \eqref{SEGED13}, it is enough to check that
 \begin{align} \label{SEGED19}
   &\lim_{n\to\infty}\frac{1}{\sqrt{n}}\sum_{k=s+1}^n\EE\big[(Z_k-\alpha Z_{k-1})X_{k-1}\big]
         =0,\\ \label{SEGED20}
   &\lim_{n\to\infty}\frac{1}{n}
         \EE\left(\sum_{k=s+1}^n(Z_k-\alpha Z_{k-1})X_{k-1}\right)^2
        =0.
 \end{align}
Since \ $\EE Z_k=\alpha \EE Z_{k-1}$, \ $k\geq s+1$, \ and the processes \ $X$ \ and \ $Z$
 \ are independent, we have \eqref{SEGED19}.
Similarly, we get
 \begin{align*}
     \EE\left(\sum_{k=s+1}^n(Z_k-\alpha Z_{k-1})X_{k-1}\right)^2
       =\sum_{k=s+1}^n\EE(Z_k-\alpha Z_{k-1})^2\EE X_{k-1}^2, \qquad n\geq s+1.
 \end{align*}
Using that \ $\lim_{k\to\infty}\EE X_k^2 = \EE \widetilde X^2$ \ (see, \eqref{Seged_stationary1} and
 \eqref{STAC_MOMENT2}), there exists some \ $L>0$ \ such that \ $\EE X_k^2<L$ \ for all \ $k\in\NN$.
\ By Proposition \ref{LEMMA3_Decomposition},
 \ $\lim_{k\to\infty}\EE (Z_k-\alpha Z_{k-1})^2\leq \lim_{k\to\infty}2\EE (Z_k^2+\alpha^2 Z_{k-1}^2)=0$,
 \ and hence \ $\lim_{n\to\infty}\frac{1}{n}\sum_{k=1}^n\EE (Z_k-\alpha Z_{k-1})^2=0$,
 \ which yields that
 \[
    \frac{1}{n}
      \sum_{k=s+1}^n\EE(Z_k-\alpha Z_{k-1})^2\EE X_{k-1}^2
      \leq \frac{L}{n}  \sum_{k=s+1}^n\EE (Z_k-\alpha Z_{k-1})^2
       \to 0 \qquad \text{as \ $n\to\infty$.}
 \]
To prove \eqref{SEGED14}, it is enough to check that
 \begin{align} \label{SEGED21}
   &\lim_{n\to\infty}\frac{1}{\sqrt n}
     \DS\sum_{k=s+1}^n\EE\big[(Z_k-\alpha Z_{k-1})Z_{k-1}\big]=0,\\ \label{SEGED22}
   &\lim_{n\to\infty}\frac{1}{n}
      \EE\left(\DS\sum_{k=s+1}^n(Z_k-\alpha Z_{k-1})Z_{k-1}\right)^2=0.
 \end{align}
Using that \ $\EE\big[(Z_k-\alpha Z_{k-1})Z_{k-1}\big]=\EE(Z_{k-1}\EE(Z_k-\alpha Z_{k-1} \mid Z_{k-1}))=0$,
 \ $k\in\NN$, \ we get \eqref{SEGED21}.
For all \ $k>\ell$, \ $k,\ell\in\NN$, \ we get
 \begin{align*}
    \EE[(Z_k-\alpha Z_{k-1})Z_{k-1}(Z_\ell-\alpha Z_{\ell-1})Z_{\ell-1}]
      &= \EE\big[\EE[(Z_k-\alpha Z_{k-1})Z_{k-1}(Z_\ell-\alpha Z_{\ell-1})Z_{\ell-1}\mid\cF_{k-1}^Z]\big]  \\
      & =\EE\big[Z_{k-1}(Z_\ell-\alpha Z_{\ell-1})Z_{\ell-1}\EE(Z_k-\alpha Z_{k-1}\mid\cF_{k-1}^Z)\big]
        =0,
 \end{align*}
 and hence, by \eqref{SEGED_Z4}, we obtain
 \begin{align*}
   \frac{1}{n}\EE\left(\DS\sum_{k=s+1}^n(Z_k-\alpha Z_{k-1})Z_{k-1}\right)^2
      & =\frac{1}{n}\sum_{k=s+1}^n \EE[(Z_k-\alpha Z_{k-1})^2Z_{k-1}^2]\\
      & =\frac{1}{n}\sum_{k=s+1}^n \EE\big[Z_{k-1}^2\EE((Z_k-\alpha Z_{k-1})^2\mid\cF_{k-1}^Z)\big]\\
      &  = \frac{1}{n}\sum_{k=s+1}^n \EE\big[Z_{k-1}^2\alpha(1-\alpha)Z_{k-1}\big]
         = \frac{\alpha(1-\alpha)}{n}\sum_{k=s+1}^n \EE Z_{k-1}^3.
 \end{align*}
By Proposition \ref{LEMMA3_Decomposition}, this implies \eqref{SEGED22}.
Condition \eqref{SEGED15} was already proved, see \eqref{SEGED98}.
Finally, Proposition \ref{LEMMA3_Decomposition} easily yields \eqref{SEGED16}.
Using \eqref{CONVERGENCE5} -- \eqref{SEGED16}, Slutsky's lemma yields \eqref{CONVERGENCE3}.

By \eqref{SEGED4} and \eqref{Strong_consistency4},
 \[
  \sqrt{n}\big(\ttheta_n - \lim_{k\to\infty}\ttheta_k\big)
   = \sqrt{n}\big(\ttheta_n - (Y_s-\alpha Y_{s-1}-\mu_\vare)\big)
   = -\sqrt{n}(\talpha_n-\alpha)Y_{s-1},
  \]
 holds asymptotically as \ $n\to\infty$ \ with probability one, and hence by \eqref{CONVERGENCE3},
 we get \eqref{CONVERGENCE4}.
\proofend

\begin{Rem}\label{REMARK1}
By \eqref{SEGED_Z1} and \eqref{SEGED_POT},
 \begin{align*}
   \EE Y_k
      =\begin{cases}
         \alpha^k\EE Y_0 + \mu_\vare\frac{1-\alpha^k}{1-\alpha}
              & \quad \text{if \ $k=1,\ldots,s-1$,}\\[1mm]
          \alpha^k\EE Y_0 + \theta\alpha^{k-s}+\mu_\vare\frac{1-\alpha^k}{1-\alpha}
              & \quad \text{if \ $k\geq s$}.
        \end{cases}
 \end{align*}
Hence \ $\EE(Y_s-\alpha Y_{s-1}-\mu_\vare)=\theta$, \ $\theta\in\NN$.
\ Moreover, by \eqref{SEGED64},
 \begin{align*}
  \var(Y_s-\alpha Y_{s-1}-\mu_\vare)
     & = \var(X_s-\alpha X_{s-1}-\mu_\vare+\theta)
       = \var(X_s - \alpha X_{s-1}-\mu_\vare)\\
     & = \alpha^s(1-\alpha)\EE X_0 + \alpha\mu_\vare (1-\alpha^{s-1}) + \sigma_\vare^2 \\
     & = \alpha^s(1-\alpha)\EE Y_0 + \alpha\mu_\vare (1-\alpha^{s-1}) + \sigma_\vare^2.
 \end{align*}
If \ $k\geq s+1$, \ then one can derive a more complicated formula for
 \ $ \var(Y_k-\alpha Y_{k-1}-\mu_\vare)$ \ containing the moments of \ $Z$, \ too.

We also check that \ $\ttheta_n$ \ is an asymptotically unbiased estimator of \ $\theta$ \
 as \ $n\to\infty$ \ for all \ $(\alpha,\theta)\in(0,1)\times\NN$.
\ By \eqref{Strong_consistency4}, the sequence \ $\ttheta_n-\theta$, \ $n\in\NN$,
 \ converges with probability one and hence bounded with probability one, and then the dominated
 convergence theorem yields that
 \ $\lim_{n\to\infty}\EE(\ttheta_n-\theta)=0$.
\proofend
\end{Rem}

\subsection{One outlier, estimation of the mean of the offspring and innovation distributions
            and the outlier's size}

We suppose that \ $I=1$ \ and that \ $s_1:=s$ \ is known.
We concentrate on the CLS estimation of \ $(\alpha,\mu_\vare,\theta)$, \ where \ $\theta:=\theta_1$.
Motivated by \eqref{Innovation_CLSE}, for all \ $n\geq s$, \ $n\in\NN$, \ we define the function
 \ $Q_n:\RR^{n+1}\times\RR^3\to\RR$, \ as
 \begin{align*}
    Q_n({\bf y}_n;\alpha',\mu_\vare',\theta')
        :=\sums \big(y_k-\alpha' y_{k-1}-\mu_\vare'\big)^2
           + \big(y_s-\alpha' y_{s-1}-\mu_\vare'-\theta'\big)^2  ,
 \end{align*}
 for all \ ${\bf y}_n\in\RR^{n+1}$ \ and \ $\alpha',\mu_\vare',\theta'\in\RR$.
\ By definition, for all \ $n\geq s$, \ a CLS estimator for
 the parameter \ $(\alpha,\mu_\vare,\theta)\in(0,1)\times(0,\infty)\times\NN$ \ is
 a measurable function \ $(\halpha_n,\hmuen,\htheta_n):S_n\to\RR^3$ \ such that
 \begin{align*}
   Q_n({\bf y}_n;\,&\halpha_n({\bf y}_n),\hmuen({\bf y}_n),\htheta_n({\bf y}_n))
       = \inf_{(\alpha',\mu_\vare',\theta')\in\RR^3}Q_n({\bf y}_n;\alpha',\mu_\vare',\theta')
       \qquad \forall\;\;  {\bf y}_n\in S_n,
 \end{align*}
 where \ $S_n$ \ is suitable subset of \ $\RR^{n+1}$ \ (defined in the proof of
 Lemma \ref{LEMMA17}).
We note that we do not define the CLS estimator
 \ $(\halpha_n,\hmuen,\htheta_n)$ \ for all samples \ ${\bf y}_n\in \RR^{n+1}$.
\ We get
 \begin{align*}
   &\frac{\partial Q_n}{\partial \alpha'}({\bf y}_n;\alpha',\mu_\vare',\theta')
     = -2\sums \big(y_k-\alpha' y_{k-1}-\mu_\vare'\big)y_{k-1}
             - 2\big(y_s-\alpha' y_{s-1}-\mu_\vare'-\theta'\big)y_{s-1}, \\
   &\frac{\partial Q_n}{\partial \mu_\vare'}({\bf y}_n;\alpha',\mu_\vare',\theta')
     = -2\sums \big(y_k-\alpha' y_{k-1}-\mu_\vare'\big)
       -2\big(y_s-\alpha' y_{s-1}-\mu_\vare'-\theta'\big),\\
   &\frac{\partial Q_n}{\partial \theta'}({\bf y}_n;\alpha',\mu_\vare',\theta')
     = -2\big(y_s-\alpha' y_{s-1}-\mu_\vare'-\theta'\big).
 \end{align*}

The next lemma is about the existence and uniqueness of the CLS estimator of \ $(\alpha,\mu_\vare,\theta)$.

\begin{Lem}\label{LEMMA17}
There exist subsets \ $S_n\subset\RR^{n+1}$, $n\geq s$ \ with the following properties:
 \begin{enumerate}
  \item[\upshape{(i)}]
   there exists a unique CLS estimator
   \ $(\halpha_n,\hmuen,\htheta_n):S_n\to\RR^3$,
  \item[\upshape{(ii)}]
   for all \ ${\bf y}_n\in S_n$, \ the system of equations
  \begin{align}\label{Innovation_CLSE_EQ3}
   \begin{split}
    &\frac{\partial Q_n}{\partial \alpha'}({\bf y}_n;\alpha',\mu_\vare',\theta')=0,\\
    &\frac{\partial Q_n}{\partial \mu_\vare'}({\bf y}_n;\alpha',\mu_\vare',\theta')=0,\\
    &\frac{\partial Q_n}{\partial \theta'}({\bf y}_n;\alpha',\mu_\vare',\theta')=0,
   \end{split}
 \end{align}
 has the unique solution
 \begin{align}\label{SEGED7}
   &\halpha_n({\bf y}_n)
     =\frac{ (n-1)\DS\sums y_{k-1} y_k
            - \DS\sums y_k \DS\sums y_{k-1}}
           {D_n({\bf y}_n)},\\[2mm]\label{SEGED8}
   &\hmuen({\bf y}_n)
       = \frac{\DS\sums y_{k-1}^2 \DS\sums y_k
          - \DS\sums y_{k-1} \DS\sums y_{k-1}y_k} {D_n({\bf y}_n)},
             \\[2mm]\label{SEGED9}
   &\htheta_n({\bf y}_n)
      = y_s - \halpha_n({\bf y}_n) y_{s-1}-\hmuen({\bf y}_n),
 \end{align}
 where
 \[
    D_n({\bf y}_n):=(n-1)\sums y_{k-1}^2
              - \left(\sums y_{k-1} \right)^2,
       \qquad n\geq s,
 \]
  \item[\upshape{(iii)}]
  ${\bf Y}_n\in S_n$ \ holds asymptotically as \ $n\to\infty$ \ with probability one.
 \end{enumerate}
\end{Lem}

\noindent{\bf Proof.}
One can easily check that the unique solution of the system of equations \eqref{Innovation_CLSE_EQ3}
 takes the form \eqref{SEGED7}-\eqref{SEGED8}-\eqref{SEGED9} whenever \ $D_n({\bf y}_n)>0$.

For all \ $n\geq s+1$, \ let
 \[
   S_n := \left\{{\bf y}_n\in\RR^{n+1} : D_n({\bf y}_n)>0,
                  \;\Delta_{i,n}({\bf y}_n;\alpha',\mu_\vare',\theta') > 0, \;\; i=1,2,3,
                  \; \forall\, (\alpha',\mu_\vare',\theta')\in\RR^3  \right\},
 \]
 where \ $\Delta_{i,n}({\bf y}_n;\alpha',\mu_\vare',\theta')$, $i=1,2,3$, \ denotes the $i$-th order
 leading principal minor of the \ $3\times 3$ \ matrix
 \begin{align*}
  H_n({\bf y}_n;\alpha',\mu_\vare',\theta')
  :=\begin{bmatrix}
     \frac{\partial^2 Q_n}{\partial(\alpha')^2}
     & \frac{\partial^2 Q_n}{\partial\mu_\vare'\partial\alpha'}
     & \frac{\partial^2 Q_n}{\partial\theta'\partial\alpha'}  \\
     \frac{\partial^2 Q_n}{\partial\alpha'\partial\mu_\vare'}
     & \frac{\partial^2 Q_n}{\partial(\mu_\vare')^2}
     & \frac{\partial^2 Q_n}{\partial\theta'\partial \mu_\vare} \\
       \frac{\partial^2 Q_n}{\partial\alpha'\partial\theta'}
     & \frac{\partial^2 Q_n}{\partial\mu_\vare'\partial\theta'}
     & \frac{\partial^2 Q_n}{\partial(\theta')^2} \\
   \end{bmatrix}
      ({\bf y}_n;\alpha',\mu_\vare',\theta').
 \end{align*}
Then the function
  \ $\RR^3 \ni (\alpha',\mu_\vare',\theta')
     \mapsto Q_n({\bf y}_n;\alpha',\mu_\vare',\theta')$
  \ is strictly convex for all \ ${\bf y}_n\in S_n$,
 \ see, e.g., Berkovitz \cite[Theorem 3.3, Chapter III]{Ber}.

Since the function \ $\RR^3 \ni (\alpha',\mu_\vare',\theta')
     \mapsto Q_n({\bf y}_n;\alpha',\mu_\vare',\theta')$
  \ is strictly convex for all \ ${\bf y}_n\in S_n$ \ and the system of equations
 \eqref{Innovation_CLSE_EQ3} has a unique solution for all \ ${\bf y}_n\in S_n$,
 we get the function in question attains its (global) minimum at this unique solution,
 which yields (i) and (ii).

Further, for all \ ${\bf y}_n\in\RR^{n+1}$ \ and
 \ $(\alpha',\mu_\vare',\theta')\in\RR^3$, \ we have
 \begin{align*}
   &\frac{\partial^2 Q_n}{\partial(\alpha')^2}({\bf y}_n;\alpha',\mu_\vare',\theta')
       = 2\sums y_{k-1}^2 + 2y_{s-1}^2
       =2\sum_{k=1}^n y_{k-1}^2,\\
   &\frac{\partial^2 Q_n}{\partial\alpha'\partial\theta'}({\bf y}_n;\alpha',\mu_\vare',\theta')
       =\frac{\partial^2 Q_n}{\partial\theta'\partial\alpha'}({\bf y}_n;\alpha',\mu_\vare',\theta')
       = 2y_{s-1},
  \end{align*}
  and
  \begin{align*}
   &\frac{\partial^2 Q_n}{\partial\alpha'\partial\mu_\vare'}({\bf y}_n;\alpha',\mu_\vare',\theta')
       =\frac{\partial^2 Q_n}{\partial\mu_\vare'\partial\alpha'}({\bf y}_n;\alpha',\mu_\vare',\theta')
       = 2\sum_{k=1}^ny_{k-1},\\
   &\frac{\partial^2 Q_n}{\partial\theta'\partial\mu_\vare'}({\bf y}_n;\alpha',\mu_\vare',\theta')
       =\frac{\partial^2 Q_n}{\partial\mu_\vare'\partial\theta'}({\bf y}_n;\alpha',\mu_\vare',\theta')
       = 2,\\
   &\frac{\partial^2 Q_n}{\partial(\theta')^2}({\bf y}_n;\alpha',\mu_\vare',\theta')=2,
     \qquad
     \frac{\partial^2 Q_n}{\partial(\mu_\vare')^2}({\bf y}_n;\alpha',\mu_\vare',\theta')=2n.
 \end{align*}
Then \ $H_n({\bf y}_n;\alpha',\mu_\vare',\theta')$ \ has the following leading principal minors
 \[
  \Delta_{1,n}({\bf y}_n;\alpha',\mu_\vare',\theta') = 2\sum_{k=1}^n y_{k-1}^2, \quad
   \Delta_{2,n}({\bf y}_n;\alpha',\mu_\vare',\theta')
         = 4\left(n\sum_{k=1}^n y_{k-1}^2 - \left(\sum_{k=1}^n y_{k-1}\right)^2\right),
 \]
 and
 \begin{align*}
    \Delta_{3,n}({\bf y}_n;\alpha',\mu_\vare',\theta')
      & = \det H_n({\bf y}_n;\alpha',\mu_\vare',\theta')\\
      &  = 8 \left( (n-1)\sum_{k=1}^n y_{k-1}^2 + 2Y_{s-1}\sum_{k=1}^n y_{k-1}
                  - n (y_{s-1})^2 - \left( \sum_{k=1}^n y_{k-1}\right)^2\right).
 \end{align*}
Note that \ $\Delta_{i,n}({\bf y}_n;\alpha',\mu_\vare',\theta')$, $i=1,2,3$, \ do not depend
 on  \ $(\alpha',\mu_\vare',\theta')$, \ and hence we will simply denote
 \ $\Delta_{i,n}({\bf y}_n;\alpha',\mu_\vare',\theta')$ \ by
 \ $\Delta_{i,n}({\bf y}_n)$.

Next we check that \ ${\bf Y}_n\in S_n$ \ holds asymptotically as \ $n\to\infty$ \ with probability one.
By \eqref{Ergodic1} and \eqref{Ergodic2}, using the very same arguments as in the proof of
 Lemma \ref{LEMMA15}, one can get
 \begin{align*}
   &\PP\left(\lim_{n\to\infty}\frac{\Delta_{1,n}({\bf Y}_n)}{n} = 2
              \EE\widetilde X^2 \right)=1,\\
   &\PP\left(\lim_{n\to\infty}\frac{\Delta_{2,n}({\bf Y}_n)}{n^2} = 4
             \var\widetilde X \right)=1,\\
   &\PP\left(\lim_{n\to\infty}\frac{\Delta_{3,n}({\bf Y}_n)}{n^2}
              = 8 \var\widetilde X
              \right)=1,
 \end{align*}
 where \ $\widetilde X$ \ denotes a random variable with the unique stationary distribution of
 the INAR(1) model in \eqref{INAR1}.
Hence
 \begin{align*}
    \PP\left(\lim_{n\to\infty}\Delta_{i,n}({\bf Y}_n) = \infty \right)=1,
        \qquad i=1,2,3.
 \end{align*}
By \eqref{Ergodic1} and \eqref{Ergodic2}, we also get
 \begin{align}\label{LEMMA7_allitasa}
   \PP\left(\lim_{n\to\infty}\frac{D_n({\bf Y}_n)}{n^2}=\var\widetilde X\right)=1
 \end{align}
 and hence \ $\PP(\lim_{n\to\infty}D_n({\bf Y}_n)=\infty)=1$.
\proofend

By Lemma \ref{LEMMA17}, \ $(\halpha_n({\bf Y}_n),
    \hmuen({\bf Y}_n),
    \htheta_n({\bf Y}_n))$
 \ exists uniquely asymptotically as \ $n\to\infty$ \ with probability one.
In the sequel we will simply denote it
  by \ $(\halpha_n,\hmuen,\htheta_n)$, \ we will also denote \ $D_n({\bf Y}_n)$
  \ by \ $D_n$.

By \eqref{SEGED7} and \eqref{SEGED8}, we also get
 \[
    \hmuen = \frac{1}{n-1}\left(\sums Y_k -\halpha_n\sums Y_{k-1} \right)
 \]
 holds asymptotically as \ $n\to\infty$ \ with probability one.

The next result shows that \ $\halpha_n$ \ and \ $\hmuen$ \ are strongly consistent estimators
 of \ $\alpha$ \ and \ $\mu_\vare$, \ whereas \ $\htheta_n$ \ fails to be a
 strongly consistent estimator of \ $\theta$.

\begin{Thm}\label{THEOREM8}
Consider the CLS estimators
 \ $(\halpha_n,\hmuen,\htheta_n)_{n\in\NN}$ \ of the parameter
 \ $(\alpha,\mu_\vare,\theta)\in(0,1)\times(0,\infty)\times\NN$.
\ The sequences \ $(\halpha_n)_{n\in\NN}$ \ and \ $(\hmuen)_{n\in\NN}$
 \ are strongly consistent for all
 \ $(\alpha,\mu_\vare,\theta)\in(0,1)\times(0,\infty)\times\NN$, \ i.e.,
 \begin{align}\label{Strong_consistency8}
   &\PP(\lim_{n\to\infty}\halpha_n=\alpha)=1,
     \qquad \forall\;(\alpha,\mu_\vare,\theta)\in(0,1)\times(0,\infty)\times\NN,\\
     \label{Strong_consistency9}
   &\PP(\lim_{n\to\infty}\hmuen=\mu_\vare)=1,
     \qquad \forall\; (\alpha,\mu_\vare,\theta)\in(0,1)\times(0,\infty)\times\NN,
 \end{align}
 whereas the sequence \ $(\htheta_n)_{n\in\NN}$ \ is not strongly consistent for any
  \ $(\alpha,\mu_\vare,\theta)\in(0,1)\times(0,\infty)\times\NN$, \ namely,
  \begin{align}\label{Strong_consistency28}
  \PP\left(\lim_{n\to\infty}\htheta_n = Y_s-\alpha Y_{s-1}-\mu_\vare\right)=1,
    \qquad \forall\; (\alpha,\mu_\vare,\theta)\in(0,1)\times(0,\infty)\times\NN.
 \end{align}
\end{Thm}

\noindent{\bf Proof.}
By \eqref{SEGED7}, \eqref{SEGED8} and Proposition \ref{LEMMA3_Decomposition}, we get
 \begin{align*}
    \begin{bmatrix}
       \halpha_n  \\
       \hmuen  \\
     \end{bmatrix}
    = \frac{1}{D_n}
      \begin{bmatrix}
        K_n  \\
        L_n  \\
      \end{bmatrix},
  \end{align*}
 where
 \begin{align*}
   &K_n :=
           (n-1)\DS\sums (X_{k-1}+Z_{k-1})(X_k+Z_k)
             - \DS\sums (X_k+Z_k)\!\DS\sums (X_{k-1}+Z_{k-1}),\\
   &L_n:=\!\!\DS\sums \!(X_{k-1}+Z_{k-1})^2\!\DS\sums (X_k+Z_k)
           - \!\DS\sums \!(X_{k-1}+Z_{k-1})\!\DS\sums (X_{k-1}+Z_{k-1})(X_k+Z_k).
 \end{align*}
Using the very same arguments as in the proof of Theorem \ref{THEOREM7}, we obtain
 \begin{align*}
   &\PP\left(\lim_{n\to\infty}\frac{K_n}{n^2}
             = \alpha\EE\widetilde X^2 + \mu_\vare\EE\widetilde X - (\EE\widetilde X)^2 \right)=1,\\
   &\PP\left(\lim_{n\to\infty}\frac{L_n}{n^2}
             = \EE\widetilde X^2 \EE\widetilde X
                - \EE\widetilde X(\alpha\EE\widetilde X^2+\mu_\vare\EE\widetilde X) \right)=1.
 \end{align*}
By \eqref{LEMMA7_allitasa} and \eqref{STAC_MOMENT1}, \eqref{STAC_MOMENT2}, we obtain
 \begin{align*}
  \PP\left(\lim_{n\to\infty}\halpha_n
            =\lim_{n\to\infty}\frac{K_n}{D_n}
            = \frac{\alpha \var\widetilde X + (\alpha-1)(\EE\widetilde X)^2 + \mu_\vare \EE\widetilde X}
               {\var\widetilde X}
            =\alpha
     \right)=1,
 \end{align*}
 and
 \begin{align*}
  \PP\left(\lim_{n\to\infty}\hmuen
            =\lim_{n\to\infty}\frac{L_n}{D_n}
            = \frac{(1-\alpha)\EE\widetilde X \EE\widetilde X^2 - \mu_\vare (\EE\widetilde X)^2}
               {\var\widetilde X}
            =\mu_\vare
     \right)=1,
 \end{align*}
 where we used that
 \begin{align*}
   \frac{(1-\alpha)\EE\widetilde X \EE\widetilde X^2 - \mu_\vare (\EE\widetilde X)^2}{\var\widetilde X}
      & = \frac{1}{\var\widetilde X}
         \left[(1-\alpha)\frac{\mu_\vare}{1-\alpha}
         \left(\frac{\sigma_\vare^2+\alpha\mu_\vare}{1-\alpha^2} + \frac{\mu_\vare^2}{(1-\alpha)^2}\right)
         -\mu_\vare\frac{\mu_\vare^2}{(1-\alpha)^2}\right]\\
      & = \frac{\mu_\vare}{\var\widetilde X}\frac{\sigma_\vare^2+\alpha\mu_\vare}{1-\alpha^2}
        =\mu_\vare.
 \end{align*}
Finally, using \eqref{SEGED9}, \eqref{Strong_consistency8} and \eqref{Strong_consistency9}
 we get \eqref{Strong_consistency28}.
\proofend

The asymptotic distribution of the CLS estimation is given in the next theorem.

\begin{Thm}\label{Proposition2}
Under the additional assumptions \ $\EE Y_0^3<\infty$ \ and \ $\EE\vare_1^3<\infty$, \ we have
 \begin{align}\label{CONVERGENCE18}
      \begin{bmatrix}
        \sqrt{n}(\halpha_n-\alpha) \\
        \sqrt{n}(\hmuen-\mu_\vare) \\
      \end{bmatrix}
     \distr
      \cN\left( \begin{bmatrix}
                  0 \\
                  0 \\
                \end{bmatrix}
               , \;B_{\alpha,\vare}
      \right)
   \qquad \text{as \ $n\to\infty$,}
 \end{align}
 where \ $B_{\alpha,\vare}$ \ is  defined in \eqref{SEGED_BALPHA}.
Moreover, conditionally on the value \ $Y_{s-1}$,
 \begin{align}\label{CONVERGENCE19}
   \sqrt{n}(\htheta_n-\lim_{k\to\infty}\htheta_k) \distr \cN(0,[Y_{s-1} \;0]B_{\alpha,\vare}[Y_{s-1}\; 0]^\top)
   \qquad \text{as \ $n\to\infty$.}
 \end{align}
\end{Thm}

\noindent{\bf Proof.}
By \eqref{SEGED7} and \eqref{SEGED8}, we have
 \begin{align*}
    \begin{bmatrix}
      \halpha_n - \alpha \\
      \hmuen - \mu_\vare \\
    \end{bmatrix}
 =
    \frac{1}{D_n}
     \begin{bmatrix}
       (n-1)\DS\sums (Y_k-\alpha Y_{k-1})Y_{k-1}
            - \DS\sums (Y_k-\alpha Y_{k-1})\DS\sums Y_{k-1} \\
       \DS\sums Y_{k-1}^2 \DS\sums (Y_k-\mu_\vare)
         - \DS\sums Y_{k-1} \DS\sums (Y_k-\mu_\vare)Y_{k-1} \\
    \end{bmatrix},
 \end{align*}
 holds asymptotically as \ $n\to\infty$ \ with probability one.
By Proposition \ref{LEMMA3_Decomposition}, we get
 \begin{align*}
    \begin{bmatrix}
      \halpha_n - \alpha \\
      \hmuen - \mu_\vare \\
    \end{bmatrix}
    = \frac{1}{D_n}
      \begin{bmatrix}
        (n-1)\DS\sums (X_k-\alpha X_{k-1})X_{k-1}
            - \DS\sums (X_k-\alpha X_{k-1})\DS\sums X_{k-1} +R_n  \\
      \DS\sums X_{k-1}^2 \DS\sums (X_k-\mu_\vare)
         - \DS\sums X_{k-1} \DS\sums (X_k-\mu_\vare)X_{k-1} + Q_n \\
      \end{bmatrix},
 \end{align*}
 holds asymptotically as \ $n\to\infty$ \ with probability one, where
 \begin{align*}
   R_n:= & (n-1)\sums (Z_k-\alpha Z_{k-1})(X_{k-1}+Z_{k-1})
           + (n-1)\sums (X_k-\alpha X_{k-1})Z_{k-1} \\
         &  - \sums (Z_k-\alpha Z_{k-1})\sums (X_{k-1}+Z_{k-1})
            - \sums (X_k-\alpha X_{k-1})\sums Z_{k-1},
 \end{align*}
 and
 \begin{align*}
   Q_n:=& \sums (2X_{k-1}Z_{k-1}+Z_{k-1}^2)\sums (X_k+Z_k-\mu_\vare)
          + \sums X_{k-1}^2 \sums Z_k \\
        & - \sums Z_{k-1}\sums (X_k+Z_k-\mu_\vare)(X_{k-1}+Z_{k-1})
          - \sums X_{k-1} \sums Z_k(X_{k-1}+Z_{k-1})\\
        & - \sums X_{k-1}\sums (X_k-\mu_\vare)Z_{k-1}.
 \end{align*}
By \eqref{LEMMA7_allitasa}, \eqref{ALPHA_MU_CONVERGENCE_INAR} and Slutsky's lemma, we have
 \begin{align*}
   \frac{\sqrt n}{D_n}
      \begin{bmatrix}
        n\DS\sums (X_k-\alpha X_{k-1})X_{k-1}
            - \DS\sums (X_k-\alpha X_{k-1})\DS\sums X_{k-1} \\
      \DS\sums X_{k-1}^2 \DS\sums (X_k-\mu_\vare)
         - \DS\sums X_{k-1} \DS\sums (X_k-\mu_\vare)X_{k-1} \\
      \end{bmatrix}
     \distr
       \cN\left( \begin{bmatrix}
                  0 \\
                  0 \\
                \end{bmatrix}
               , \;B_{\alpha,\vare}
      \right),
 \end{align*}
 as \ $n\to\infty$, \ and hence to prove \eqref{CONVERGENCE18}, by \eqref{LEMMA7_allitasa} and Slutsky's lemma,
 it is enough to check that
 \begin{align}\label{SEGED69}
    \frac{R_n}{n^{3/2}}\stoch 0
      \qquad \text{as \ $n\to\infty$,}
 \end{align}
 and
 \begin{align}\label{SEGED70}
    \frac{Q_n}{n^{3/2}}\stoch 0
      \qquad \text{as \ $n\to\infty$,}
 \end{align}
 where \ $\stoch$ \ denotes convergence in probability.
By \eqref{SEGED13} and \eqref{SEGED14}, to prove \eqref{SEGED69}
 it remains to check that
 \begin{align}\label{SEGED71}
   &\frac{1}{\sqrt n}\sums (X_k-\alpha X_{k-1})Z_{k-1}
      \stoch 0
       \qquad \text{as \ $n\to\infty$,} \\  \label{SEGED72}
   &\frac{1}{\sqrt n} \sums (Z_k-\alpha Z_{k-1})
     \cdot\frac{1}{n} \sums (X_{k-1}+Z_{k-1})
       \stoch 0
      \qquad \text{as \ $n\to\infty$,} \\  \label{SEGED73}
   & \frac{1}{n} \sums (X_k-\alpha X_{k-1})
     \cdot\frac{1}{\sqrt n} \sums Z_{k-1}
       \stoch 0
      \qquad \text{as \ $n\to\infty$.}
 \end{align}
Using \eqref{SEGED12} and that
 \begin{align*}
   \frac{1}{\sqrt n} \sums(X_k-\alpha X_{k-1}) Z_{k-1}
      = \frac{1}{\sqrt n} \sums(X_k-\alpha X_{k-1}-\mu_\vare) Z_{k-1}
         + \mu_\vare\frac{1}{\sqrt n}\sums Z_{k-1},
 \end{align*}
 to prove \eqref{SEGED71}, it is enough to check that
 \begin{align}\label{SEGED74}
    \frac{1}{\sqrt n}\sums Z_{k-1}
      \stoch 0
      \qquad \text{as \ $n\to\infty$.}
 \end{align}
Using that \ $Z_k\geq 0$, \ $k\in\NN$, \ by Markov's inequality, it is enough to check that
 \begin{align}\label{SEGED75}
    \lim_{n\to\infty}
      \frac{1}{\sqrt n}\sums \EE Z_{k-1} = 0.
 \end{align}
Since, by \eqref{SEGED_Z1}, \ $\EE Z_{s+k}=\theta\alpha^k$, \ $k\geq 0$, \ we have
 \begin{align}\label{SEGED82}
   \frac{1}{\sqrt n}\sums \EE Z_{k-1}
       \leq \frac{\theta}{\sqrt n}\sum_{k=0}^n \alpha^k
        = \frac{\theta}{\sqrt n}\frac{\alpha^{n+1}-1}{\alpha-1}
         \to 0
      \qquad \text{as \ $n\to\infty$.}
 \end{align}

\noindent Using that
 \[
    \PP\left( \lim_{n\to\infty} \frac{1}{n}\sums(X_{k-1}+Z_{k-1})=\EE\widetilde X\right)=1,
 \]
 to prove \eqref{SEGED72} it is enough to check that
 \begin{align}\label{SEGED76}
   \frac{1}{\sqrt n} \sums (Z_k-\alpha Z_{k-1}) \mean 0
   \qquad \text{as \ $n\to\infty$.}
 \end{align}
To verify \eqref{SEGED76} it is enough to show that
 \begin{align}\label{SEGED77}
  & \lim_{n\to\infty} \frac{1}{\sqrt n} \sums \EE(Z_k-\alpha Z_{k-1})=0,\\ \label{SEGED78}
  & \lim_{n\to\infty} \frac{1}{n} \EE\left(\sums (Z_k-\alpha Z_{k-1})\right)^2=0.
 \end{align}
Using that \ $\EE Z_k=\alpha\EE Z_{k-1}$, \ $k\geq s+1$, \ we get \eqref{SEGED77} is satisfied.
Using that \ $\EE[(Z_k-\alpha Z_{k-1})(Z_\ell-\alpha Z_{\ell-1})]=0$ \ for all \ $k\ne\ell$,
 $k,\ell\geq s+1$, \ we have
 \[
    \frac{1}{n} \EE \left(\sums (Z_k-\alpha Z_{k-1})\right)^2
       =  \frac{1}{n} \sums \EE(Z_k-\alpha Z_{k-1})^2\to 0
    \qquad \text{as \ $n\to\infty$,}
 \]
 as we showed in the proof of Theorem \ref{Proposition1}.
Using that
 \[
   \PP\left(\lim_{n\to\infty}\frac{1}{n}\sums (X_k-\alpha X_{k-1})
             = (1-\alpha)\EE\widetilde X\right)=1,
 \]
 to prove \eqref{SEGED73} it is enough to verify \eqref{SEGED74} which was already done.

Now we turn to prove \eqref{SEGED70}.
Using \eqref{SEGED74} and that
 \begin{align*}
   &\PP\left(\lim_{n\to\infty} \frac{1}{n}\sums(X_k+Z_k-\mu_\vare)
               =\EE\widetilde X-\mu_\vare\right)=1, \\
   &\PP\left(\lim_{n\to\infty} \frac{1}{n}\sums X_{k-1}^2
                 =\EE\widetilde X^2\right)=1, \\
   &\PP\left(\lim_{n\to\infty} \frac{1}{n}\sums(X_k+Z_k-\mu_\vare)(X_{k-1}+Z_{k-1})
               =\alpha\EE\widetilde X^2+\mu_\vare\EE\widetilde X - \mu_\vare\EE\widetilde X
               = \alpha\EE\widetilde X^2 \right)=1,
 \end{align*}
 it is enough to verify that
 \begin{align}\label{SEGED79}
   &\frac{1}{\sqrt n}\sums X_{k-1}Z_{k-1}
       \stoch 0
      \qquad \text{as \ $n\to\infty$,} \\ \label{SEGED80}
   &\frac{1}{\sqrt n}\sums Z_{k-1}^2
       \stoch 0
      \qquad \text{as \ $n\to\infty$,} \\ \label{SEGED81}
   &\frac{1}{\sqrt n}\sums Z_{k-1}Z_k
       \stoch 0
     \qquad \text{as \ $n\to\infty$.}
 \end{align}
To prove \eqref{SEGED79}, using that the processes \ $X$ \ and \ $Z$ \ are non-negative,
 by Markov's inequality, it is enough to verify that
  \[
    \lim_{n\to\infty} \frac{1}{\sqrt n} \sums \EE(X_{k-1}Z_{k-1})=0.
  \]
Using that the processes \ $X$ \ and \ $Z$ \ are independent and
 \ $\lim_{k\to\infty}\EE X_{k-1}=\EE\widetilde X$, \ as in the proof of
 \eqref{SEGED13}, we get it is enough to check that
 \[
     \lim_{n\to\infty} \frac{1}{\sqrt n} \sums \EE Z_{k-1}=0,
 \]
 which follows by \eqref{SEGED82}.
Similarly, to prove \eqref{SEGED80}, it is enough to show that
 \[
     \lim_{n\to\infty} \frac{1}{\sqrt n} \sums \EE Z_{k-1}^2=0.
 \]
By \eqref{SEGED_Z2}, we have for all \ $n\geq s+1$,
 \begin{align*}
   \frac{1}{\sqrt n} \sums \EE Z_{k-1}^2
     \leq \frac{1}{\sqrt n}
           \sum_{k=0}^n(\theta^2\alpha^{2k}-\theta\alpha^k(\alpha^k-1))
     \leq \frac{\theta^2}{\sqrt n}\frac{\alpha^{2(n+1)}-1}{\alpha^2-1}
          + \frac{\theta}{\sqrt n}\frac{\alpha^{n+1}-1}{\alpha-1}
     \to 0,
 \end{align*}
 as \ $n\to\infty$.
\ Similarly, to prove \eqref{SEGED81}, it is enough to check that
 \[
     \lim_{n\to\infty} \frac{1}{\sqrt n} \sums \EE (Z_{k-1}Z_k)=0.
 \]
By \eqref{SEGED_Z3}, we have for all \ $n\geq s+1$,
 \begin{align*}
   \frac{1}{\sqrt n} \sums \EE (Z_{k-1}Z_k)
     \leq \frac{1}{\sqrt n}
           \sum_{k=1}^n(\theta^2\alpha^{2k-1}+\theta\alpha^k(1-\alpha^k))
     \leq \frac{\theta^2}{\alpha\sqrt n}\frac{\alpha^{2n}-1}{\alpha^2-1}
          + \frac{\theta}{\sqrt n}\frac{\alpha^{n}-1}{\alpha-1}
     \to 0,
 \end{align*}
 as \ $n\to\infty$.

Finally, using \eqref{SEGED9} and \eqref{Strong_consistency28}, we get
 \begin{align*}
   \sqrt{n}(\htheta_n-\lim_{k\to\infty}\htheta_k)
     = -\sqrt{n}(\halpha_n - \alpha)Y_{s-1}
        - \sqrt{n}(\hmuen-\mu_\vare)
     =\begin{bmatrix}
        - Y_{s-1} & -1 \\
      \end{bmatrix}
       \begin{bmatrix}
         \sqrt{n}(\halpha_n - \alpha) \\
          \sqrt{n}(\hmuen-\mu_\vare) \\
       \end{bmatrix},
 \end{align*}
 and hence, by \eqref{CONVERGENCE18}, we have \eqref{CONVERGENCE19}.
\proofend

\subsection{Two outliers, estimation of the mean of the offspring distribution
            and the outliers' sizes}

We assume that \ $I=2$ \ and that the relevant time points \ $s_1$, $s_2\in\NN$,
 \ $s_1\ne s_2$, \ are known.
We concentrate on the CLS estimation of \ $(\alpha,\theta_1,\theta_2)$.
\ We have
 \begin{align*}
   \EE(Y_k\mid\cF^Y_{k-1})
       = \alpha Y_{k-1} + \mu_\vare + \delta_{k,s_1}\theta_1+\delta_{k,s_2}\theta_2,
      \qquad k\in\NN.
 \end{align*}
Hence for all \ $n\geq \max(s_1,s_2)$,
 \begin{align*}
   \begin{split}
    \sum_{k=1}^n&\big(Y_k-\EE(Y_k\mid \cF^Y_{k-1})\big)^2\\
        & = {\sum_{k=1}^n}^{(s_1,s_2)} \big(Y_k-\alpha Y_{k-1}-\mu_\vare\big)^2
           + \big(Y_{s_1}-\alpha Y_{s_1-1}-\mu_\vare-\theta_1\big)^2
           + \big(Y_{s_2}-\alpha Y_{s_2-1}-\mu_\vare-\theta_2\big)^2.
    \end{split}
 \end{align*}
For all \ $n\geq\max(s_1,s_2)$, \ $n\in\NN$, \ we define the function \ $Q_n:\RR^{n+1}\times\RR^3\to\RR$,
 \ as
 \begin{align*}
    &Q_n({\bf y}_n;\alpha',\theta_1',\theta_2')\\
      & :={\sum_{k=1}^n}^{(s_1,s_2)} \big(y_k-\alpha' y_{k-1}-\mu_\vare\big)^2
           + \big(y_{s_1}-\alpha' y_{s_1-1}-\mu_\vare-\theta_1'\big)^2
           + \big(y_{s_2}-\alpha' y_{s_2-1}-\mu_\vare-\theta_2'\big)^2,
 \end{align*}
  for all \ ${\bf y}_n\in\RR^{n+1}$, $\alpha',\theta_1',\theta_2'\in\RR$.
\ By definition, for all \ $n\geq \max(s_1,s_2)$, \ a CLS estimator for
 the parameter \ $(\alpha,\theta_1,\theta_2)\in(0,1)\times\NN^2$ \ is
 a measurable function \ $(\talpha_n,\ttheta_{1,n},\ttheta_{2,n}):S_n\to\RR^3$
 \ such that
 \begin{align*}
   Q_n({\bf y}_n;\,&\talpha_n({\bf y}_n),\ttheta_{1,n}({\bf y}_n),
             \ttheta_{2,n}({\bf y}_n))
       = \inf_{(\alpha',\theta_1',\theta_2')\in\RR^3}Q_n({\bf y}_n;\alpha',\theta_1',\theta_2')
       \qquad \forall\;\;  {\bf y}_n\in S_n,
 \end{align*}
  where \ $S_n$ \ is suitable subset of \ $\RR^{n+1}$ \ (defined in the proof of
 Lemma \ref{LEMMA11}).
We note that we do not define the CLS estimator
 \ $(\talpha_n,\ttheta_{1,n},\ttheta_{2,n})$ \ for all samples \ ${\bf y}_n\in \RR^{n+1}$.
\ For all \ ${\bf y}_n\in\RR^{n+1}$, $\alpha',\theta_1',\theta_2'\in\RR$,
 \begin{align*}
   & \frac{\partial Q_n}{\partial \alpha'}({\bf y}_n;\alpha',\theta_1',\theta_2')
      ={\sum_{k=1}^n}^{(s_1,s_2)}\!(-2)\big(y_k-\alpha' y_{k-1}-\mu_\vare\big)y_{k-1}
    - 2\big(y_{s_1}-\alpha' y_{s_1-1}-\mu_\vare-\theta_1'\big)y_{s_1-1}\\
   &\phantom{\frac{\partial Q_n}{\partial \alpha'}({\bf y}_n;\alpha',\theta_1',\theta_2')=\;}
             - 2\big(y_{s_2}-\alpha' y_{s_2-1}-\mu_\vare-\theta_2'\big)y_{s_2-1}, \\
   & \frac{\partial Q_n}{\partial \theta_1'}({\bf y}_n;\alpha',\theta_1',\theta_2')
      =-2\big(y_{s_1}-\alpha' y_{s_1-1}-\mu_\vare-\theta_1'\big),\\
   & \frac{\partial Q_n}{\partial \theta_2'}({\bf y}_n;\alpha',\theta_1',\theta_2')
      =-2\big(y_{s_2}-\alpha' y_{s_2-1}-\mu_\vare-\theta_2'\big).
 \end{align*}

The next lemma is about the existence and uniqueness of the CLS estimator of \ $(\alpha,\theta_1,\theta_2)$.

\begin{Lem}\label{LEMMA11}
There exist subsets \ $S_n\subset\RR^{n+1}$, $n\geq \max(3,s_1,s_2)$ \ with the following properties:
 \begin{enumerate}
  \item[\upshape{(i)}]
   there exists a unique CLS estimator
   \ $(\talpha_n,\ttheta_{1,n},\ttheta_{2,n}):S_n\to\RR^3$,
  \item[\upshape{(ii)}]
   for all \ ${\bf y}_n\in S_n$, \ the system of equations
   \begin{align}\label{Innovation_CLSE_EQ2}
   \begin{split}
    &\frac{\partial Q_n}{\partial \alpha'}({\bf y}_n;\alpha',\theta_1',\theta_2')=0,\\
    &\frac{\partial Q_n}{\partial \theta_1'}({\bf y}_n;\alpha',\theta_1',\theta_2')=0, \\
    &\frac{\partial Q_n}{\partial \theta_2'}({\bf y}_n;\alpha',\theta_1',\theta_2')=0,
   \end{split}
  \end{align}
  has the unique solution
     \begin{align}\label{SEGED83}
   &\talpha_n({\bf y}_n)
     =\frac{\DS{\sum_{k=1}^n}^{(s_1,s_2)} (y_k-\mu_\vare)y_{k-1}}
           {\DS{\sum_{k=1}^n}^{(s_1,s_2)} y_{k-1}^2},\\ \label{SEGED84}
   &\ttheta_{i,n}({\bf y}_n)
      = y_{s_i}-\talpha_n({\bf y}_n) y_{s_i-1}-\mu_\vare,
       \qquad i=1,2,
 \end{align}
  \item[\upshape{(iii)}]
  ${\bf Y}_n\in S_n$ \ holds asymptotically as \ $n\to\infty$ \ with probability one.
 \end{enumerate}
\end{Lem}

\noindent{\bf Proof.}
One can easily check that the unique solution of the system of equations \eqref{Innovation_CLSE_EQ2}
 takes the form \eqref{SEGED83} and \eqref{SEGED84} whenever
 \ $\DS{\sum_{k=1}^n}^{(s_1,s_2)} y_{k-1}^2>0$.

For all \ $n\geq\max(3,s_1,s_2)$, \ let
 \[
   S_n := \left\{{\bf y}_n\in\RR^{n+1} : \Delta_{i,n}({\bf y}_n;\alpha',\theta_1',\theta_2') > 0,
                       \;\; i=1,2,3,
                       \;\forall\, (\alpha',\theta_1',\theta_2')\in\RR^3 \right\},
 \]
 where \ $\Delta_{i,n}({\bf y}_n;\alpha',\mu_\vare',\theta')$, $i=1,2,3$, \ denotes the $i$-th order
 leading principal minor of the \ $3\times 3$ \ matrix
  \begin{align*}
  H_n({\bf y}_n;\alpha',\theta_1',\theta_2')
  := \begin{bmatrix}
    \frac{\partial^2 Q_n}{\partial(\alpha')^2}
    & \frac{\partial^2 Q_n}{\partial\theta_1'\partial\alpha'}
    & \frac{\partial^2 Q_n}{\partial\theta_2'\partial\alpha'} \\
    \frac{\partial^2 Q_n}{\partial\alpha'\partial\theta_1'}
    & \frac{\partial^2 Q_n}{\partial(\theta_1')^2}
    & \frac{\partial^2 Q_n}{\partial\theta_2'\partial \theta_1'} \\
    \frac{\partial^2 Q_n}{\partial\alpha'\partial\theta_2'}
    & \frac{\partial^2 Q_n}{\partial\theta_1'\partial\theta_2'}
    & \frac{\partial^2 Q_n}{\partial(\theta_2')^2} \\
   \end{bmatrix}
     ({\bf y}_n;\alpha',\theta_1',\theta_2')
 \end{align*}
Then the function
  \ $\RR^3 \ni (\alpha',\theta_1',\theta_2')
     \mapsto Q_n({\bf y}_n;\alpha',\theta_1',\theta_2')$
  \ is strictly convex for all \ ${\bf y}_n\in S_n$,
 \ see, e.g., Berkovitz \cite[Theorem 3.3, Chapter III]{Ber}.
Further, for all \ ${\bf y}_n\in\RR^{n+1}$ \ and
 \ $(\alpha',\theta_1',\theta_2')\in\RR^3$, \ we have
 \begin{align*}
   &\frac{\partial^2 Q_n}{\partial(\alpha')^2}({\bf y}_n;\alpha',\theta_1',\theta_2')
       = 2\sumssd y_{k-1}^2 + 2y_{s_1-1}^2 + 2y_{s_2-1}^2
       =2\sum_{k=1}^n y_{k-1}^2,\\
   &\frac{\partial^2 Q_n}{\partial\alpha'\partial\theta_1'}({\bf y}_n;\alpha',\theta_1',\theta_2')
       =\frac{\partial^2 Q_n}{\partial\theta_1'\partial\alpha'}({\bf y}_n;\alpha',\theta_1',\theta_2')
       = 2y_{s_1-1},\\
   &\frac{\partial^2 Q_n}{\partial\alpha'\partial\theta_2'}({\bf y}_n;\alpha',\theta_1',\theta_2')
       =\frac{\partial^2 Q_n}{\partial\theta_2'\partial\alpha'}({\bf y}_n;\alpha',\theta_1',\theta_2')
       = 2y_{s_2-1},\\
   &\frac{\partial^2 Q_n}{\partial\theta_1'\partial\theta_2'}({\bf y}_n;\alpha',\theta_1',\theta_2')
       =\frac{\partial^2 Q_n}{\partial\theta_2'\partial\theta_1'}({\bf y}_n;\alpha',\theta_1',\theta_2')
       = 0,\\
   &\frac{\partial^2 Q_n}{\partial(\theta_1')^2}({\bf y}_n;\alpha',\theta_1',\theta_2')=2,
     \qquad
     \frac{\partial^2 Q_n}{\partial(\theta_2')^2}({\bf y}_n;\alpha',\theta_1',\theta_2')=2.
 \end{align*}
This yields that the system of equations \eqref{Innovation_CLSE_EQ2} has a unique solution
 for all \ ${\bf y}_n\in S_n$.
Using also that the function \ $\RR^3 \ni (\alpha',\theta_1',\theta_2')
     \mapsto Q_n({\bf y}_n;\alpha',\theta_1',\theta_2')$
  \ is strictly convex for all \ ${\bf y}_n\in S_n$, \
 we get the function in question attains its (global) minimum at this unique solution,
 which yields (i) and (ii).

Then \ $H_n({\bf y}_n;\alpha',\theta_1',\theta_2')$ \ has the following leading principal minors
 \[
  \Delta_{1,n}({\bf y}_n;\alpha',\theta_1',\theta_2') = 2\sum_{k=1}^n y_{k-1}^2, \quad \quad
      \Delta_{2,n}({\bf y}_n;\alpha',\theta_1',\theta_2')
        = 4\left(\sum_{k=1}^n y_{k-1}^2  - (y_{s_1-1})^2\right),
 \]
 and
 \[
    \Delta_{3,n}({\bf y}_n;\alpha',\theta_1',\theta_2')
      =\det H_n({\bf y}_n;\alpha',\theta_1',\theta_2')
      = 8 \left( \sum_{k=1}^n y_{k-1}^2 - (y_{s_1-1})^2 - (y_{s_2-1})^2\right).
 \]
Note that \ $\Delta_{i,n}({\bf y}_n;\alpha',\theta_1',\theta_2')$, $i=1,2,3$, \ do not depend
 on  \ $(\alpha',\theta_1',\theta_2')$, \ and hence we will simply denote
 \ $\Delta_{i,n}({\bf y}_n;\alpha',\theta_1',\theta_2')$ \ by
 \ $\Delta_{i,n}({\bf y}_n)$.

 Next we check that \ ${\bf Y}_n\in S_n$ \ holds asymptotically as \ $n\to\infty$ \ with probability one.
\ By \eqref{Ergodic1} and \eqref{Ergodic2}, using the very same arguments as in the proof of
 Lemma \ref{LEMMA15}, one can get
 \begin{align*}
   &\PP\left(\lim_{n\to\infty}\frac{\Delta_{1,n}({\bf Y}_n)}{n}
                = 2 \EE\widetilde X^2 \right)=1,\\
   &\PP\left(\lim_{n\to\infty}\frac{\Delta_{2,n}({\bf Y}_n)}{n}
                = 4 \EE\widetilde X^2 \right)=1,\\
   &\PP\left(\lim_{n\to\infty}\frac{\Delta_{3,n}({\bf Y}_n)}{n}
              = 8 \EE\widetilde X^2 \right)=1,
 \end{align*}
 where \ $\widetilde X$ \ denotes a random variable with the unique stationary distribution of
 the INAR(1) model in \eqref{INAR1}.
Hence
 \begin{align*}
   \PP\left(\lim_{n\to\infty}\Delta_{i,n}({\bf Y}_n)
                = \infty \right)=1,\qquad i=1,2,3,
 \end{align*}
which yields that \ ${\bf Y}_n\in S_n$ \ holds asymptotically as \ $n\to\infty$ \ with probability one.
\proofend

By Lemma \ref{LEMMA11}, \ $(\talpha_n({\bf Y}_n),
    \ttheta_{1,n}({\bf Y}_n), \ttheta_{2,n}({\bf Y}_n))$
 \ exists uniquely asymptotically as \ $n\to\infty$ \ with probability one.
In the sequel we will simply denote it by
 \ $(\talpha_n,\ttheta_{1,n},\ttheta_{2,n})$.

The next result shows that \ $\talpha_n$ \ is a strongly consistent estimator
 of \ $\alpha$, \ whereas \ $\ttheta_{i,n}$, \ $i=1,2$, \ fail to be
 strongly consistent estimators of \ $\theta_{i,n}$, \ $i=1,2$, \  respectively.

\begin{Thm}
Consider the CLS estimators
 \ $(\talpha_n,\ttheta_{1,n},\ttheta_{2,n})_{n\in\NN}$ \ of the parameter
 \ $(\alpha,\theta_1,\theta_2)\in(0,1)\times\NN^2$.
\ The sequence \ $(\talpha_n)_{n\in\NN}$ \ is strongly consistent for all
 \ $(\alpha,\theta_1,\theta_2)\in(0,1)\times\NN^2$, \ i.e.,
 \begin{align}\label{Strong_consistency29}
   \PP(\lim_{n\to\infty}\talpha_n=\alpha)=1,
     \qquad \forall\;(\alpha,\theta_1,\theta_2)\in(0,1)\times\NN^2,
 \end{align}
 whereas the sequences \ $(\ttheta_{1,n})_{n\in\NN}$ \ and  \ $(\ttheta_{2,n})_{n\in\NN}$
 \ are not strongly consistent for any
  \ $(\alpha,\theta_1,\theta_2)\in(0,1)\times\NN^2$, \ namely,
 \begin{align}\label{Strong_consistency30}
  \PP\left(\lim_{n\to\infty}\ttheta_{i,n} = Y_{s_i}-\alpha Y_{s_i-1}-\mu_\vare\right)=1,
    \qquad \forall\; (\alpha,\theta_1,\theta_1)\in(0,1)\times\NN^2,\qquad i=1,2.
 \end{align}
\end{Thm}

\noindent{\bf Proof.}
By \eqref{SEGED83} and Proposition \ref{LEMMA3_Decomposition2}, we get
 \begin{align*}
    \talpha_n=\frac{\DS\sumssd(X_k-\mu_\vare)X_{k-1} + K_n}
               {\DS\sumssd X_{k-1}^2 + L_n},
 \end{align*}
 holds asymptotically as \ $n\to\infty$ \ with probability one, where
 \begin{align*}
   & K_n:= \sumssd (Z_k^{(1)}+Z_k^{(2)})(X_{k-1}+Z_{k-1}^{(1)}+Z_{k-1}^{(2)})
         + \sumssd (X_k-\mu_\vare)(Z_{k-1}^{(1)}+Z_{k-1}^{(2)}),\\
   & L_n:= \sumssd \big[(Z_{k-1}^{(1)})^2 + (Z_{k-1}^{(2)})^2
                         + 2X_{k-1}Z_{k-1}^{(1)}
                         + 2X_{k-1}Z_{k-1}^{(2)}
                         + 2Z_{k-1}^{(1)}Z_{k-1}^{(2)}
                       \big].
 \end{align*}
Using the very same arguments as in the proof of Theorem \ref{THEOREM7}, one can obtain
 \begin{align*}
   \PP\left(\lim_{n\to\infty}\frac{K_n}{n}=0\right)=1
   \qquad\qquad\text{and}\qquad\qquad
   \PP\left(\lim_{n\to\infty}\frac{L_n}{n}=0\right)=1.
 \end{align*}
Indeed, the only fact that was not verified in the proof of Theorem \ref{THEOREM7} is that
 \begin{align}\label{SEGED89}
    \PP\left(\lim_{n\to\infty}\frac{1}{n}\sumssd Z_k^{(1)}Z_k^{(2)}=0\right)=1.
 \end{align}
 By Cauchy-Schwartz's inequality and Proposition \ref{LEMMA3_Decomposition2},
 \begin{align*}
   \frac{1}{n}\sumssd Z_k^{(1)}Z_k^{(2)}
      \leq \sqrt{\frac{1}{n}\sum_{k=1}^n (Z_k^{(1)})^2\frac{1}{n}\sum_{k=1}^n(Z_k^{(2)})^2}\to0
      \qquad \text{as \ $n\to\infty$,}
 \end{align*}
 with probability one.

Finally, by \eqref{SEGED84} and \eqref{Strong_consistency29}, we have \eqref{Strong_consistency30}.
\proofend

\begin{Rem}
 Since
  \[
     \begin{bmatrix}
        Y_{s_1}-\alpha Y_{s_1-1}-\mu_\vare \\
        Y_{s_2}-\alpha Y_{s_2-1}-\mu_\vare \\
     \end{bmatrix}
      =
     \begin{bmatrix}
        X_{s_1}-\alpha X_{s_1-1}-\mu_\vare+\theta_1 \\
        X_{s_2}-\alpha X_{s_2-1}-\mu_\vare+\theta_2 \\
     \end{bmatrix},
  \]
 and
 \begin{align*}
  \cov(X_{s_1} & - \alpha X_{s_1-1} -\mu_\vare+\theta_1,X_{s_2}-\alpha X_{s_2-1}-\mu_\vare+\theta_2) \\
    & = \EE\big[(X_{s_1}-\alpha X_{s_1-1}-\mu_\vare)(X_{s_2}-\alpha X_{s_2-1}-\mu_\vare)\big]\\
    & = \EE\left[\left(\sum_{j=1}^{X_{s_1-1}}(\xi_{s_1,j}-\alpha)+(\vare_{s_1}-\mu_\vare)\right)
           \left(\sum_{j=1}^{X_{s_2-1}}(\xi_{s_2,j}-\alpha)+(\vare_{s_2}-\mu_\vare)\right)\right]
     =0,
 \end{align*}
 by Remark \ref{REMARK1}, we get
 \begin{align*}
   &\var\begin{bmatrix}
              Y_{s_1}-\alpha Y_{s_1-1}-\mu_\vare \\
              Y_{s_2}-\alpha Y_{s_2-1}-\mu_\vare \\
            \end{bmatrix}\\
   &\qquad
      = \begin{bmatrix}
        \alpha^{s_1}(1-\alpha)\EE Y_0
         + \alpha\mu_\vare(1-\alpha^{s_1-1}) + \sigma_\vare^2  & 0 \\
        0 & \alpha^{s_2}(1-\alpha)\EE Y_0 + \alpha\mu_\vare(1-\alpha^{s_2-1})  + \sigma_\vare^2 \\
      \end{bmatrix}.
 \end{align*}
\proofend
\end{Rem}

The asymptotic distribution of the CLS estimation is given in the next theorem.

\begin{Thm}
Under the additional assumptions \ $\EE Y_0^3<\infty$ \ and \ $\EE\vare_1^3<\infty$, \ we have
 \begin{align}\label{CONVERGENCE22}
    \sqrt{n}(\talpha_n-\alpha)
      \distr
        \cN(0,\sigma_{\alpha,\vare}^2)
   \qquad \text{as \ $n\to\infty$,}
 \end{align}
 where \ $\sigma_{\alpha,\,\vare}^2$ \ is  defined in \eqref{SEGED_SZIGMA_ALPHA}.
Moreover, conditionally on the values \ $Y_{s_1-1}$ \ and \ $Y_{s_2-1}$,
 \begin{align}\label{CONVERGENCE23}
   \begin{bmatrix}
     \sqrt{n}\big(\ttheta_{1,n}-\lim_{k\to\infty}\ttheta_{1,k}\big) \\
      \sqrt{n}\big(\ttheta_{2,n}-\lim_{k\to\infty}\ttheta_{2,k}\big) \\
   \end{bmatrix}
    \distr
    \cN\left(\begin{bmatrix}
               0 \\
               0 \\
             \end{bmatrix},
             \sigma_{\alpha,\vare}^2 \begin{bmatrix}
          Y_{s_1-1}^2 &  Y_{s_1-1}Y_{s_2-1} \\
          Y_{s_1-1}Y_{s_2-1} & Y_{s_2-1}^2 \\
        \end{bmatrix}
    \right)
   \qquad \text{as \ $n\to\infty$.}
 \end{align}
\end{Thm}

\noindent{\bf Proof.}
By \eqref{SEGED83}, we get
 \[
   \talpha_n-\alpha
      = \frac{\DS\sumssd(Y_k-\alpha Y_{k-1}-\mu_\vare)Y_{k-1}}
               {\DS\sumssd Y_{k-1}^2},
 \]
 holds asymptotically as \ $n\to\infty$ \ with probability one.
To prove \eqref{CONVERGENCE22}, using Proposition \ref{LEMMA3_Decomposition2} and
 \eqref{SEGED12}--\eqref{SEGED16}, it is enough to check that
 \begin{align}\label{SEGED85}
   &\frac{1}{\sqrt n}\sum_{k=1}^n(Z_k^{(1)}-\alpha Z_{k-1}^{(1)})Z_{k-1}^{(2)}
      \stoch 0
     \qquad \text{as \ $n\to\infty$,}\\ \label{SEGED86}
   &\frac{1}{\sqrt n}\sum_{k=1}^n(Z_k^{(2)}-\alpha Z_{k-1}^{(2)})Z_{k-1}^{(1)}
      \stoch 0
     \qquad \text{as \ $n\to\infty$.}
 \end{align}
To prove \eqref{SEGED85} it is enough to verify that
  \begin{align}\label{SEGED87}
   &\lim_{n\to\infty}
     \frac{1}{\sqrt n}\sum_{k=1}^n\EE[(Z_k^{(1)}-\alpha Z_{k-1}^{(1)})Z_{k-1}^{(2)}]
      = 0,\\ \label{SEGED88}
   &\lim_{n\to\infty}\frac{1}{n}\EE\left(\sum_{k=1}^n(Z_k^{(1)}-\alpha Z_{k-1}^{(1)})Z_{k-1}^{(2)}\right)^2
      =0.
 \end{align}
Since the processes \ $Z^{(1)}$ \ and \ $Z^{(2)}$ \ are independent, we have
 \[
   \EE[(Z_k^{(1)}-\alpha Z_{k-1}^{(1)})Z_{k-1}^{(2)}]
     = \EE(Z_k^{(1)}-\alpha Z_{k-1}^{(1)})\EE Z_{k-1}^{(2)}
     =0, \qquad k\in\NN,
 \]
 which yields \eqref{SEGED87}.
Using that for all \ $k,$ $\ell\in\NN$, \ $k>\ell$,
 \begin{align*}
   \EE\big[(Z_k^{(1)}-\alpha Z_{k-1}^{(1)})Z_{k-1}^{(2)}
           & (Z_\ell^{(1)}-\alpha Z_{\ell-1}^{(1)})Z_{\ell-1}^{(2)}\big]
      = \EE\big[(Z_k^{(1)}-\alpha Z_{k-1}^{(1)})
              (Z_\ell^{(1)}-\alpha Z_{\ell-1}^{(1)})\big]
            \EE(Z_{k-1}^{(2)}Z_{\ell-1}^{(2)}) \\
    & =  \EE\big[ (Z_\ell^{(1)}-\alpha Z_{\ell-1}^{(1)})
                  \EE(Z_k^{(1)}-\alpha Z_{k-1}^{(1)}\mid \cF_{k-1}^{Z^{(1)}})\big]
            \EE(Z_{k-1}^{(2)}Z_{\ell-1}^{(2)}) \\
    & = 0\cdot \EE(Z_{k-1}^{(2)}Z_{\ell-1}^{(2)})
      =0,
 \end{align*}
 we get
 \begin{align*}
    \frac{1}{n}\EE\left(\sum_{k=1}^n(Z_k^{(1)}-\alpha Z_{k-1}^{(1)})Z_{k-1}^{(2)}\right)^2
     &  = \frac{1}{n}\sum_{k=1}^n\EE[(Z_k^{(1)}-\alpha Z_{k-1}^{(1)})^2(Z_{k-1}^{(2)})^2] \\
     &  =  \frac{1}{n}\sum_{k=1}^n\EE(Z_k^{(1)}-\alpha Z_{k-1}^{(1)})^2\EE(Z_{k-1}^{(2)})^2 \\
     & \leq \frac{2}{n} \sum_{k=1}^n\EE((Z_k^{(1)})^2+\alpha^2(Z_{k-1}^{(1)})^2)\EE(Z_{k-1}^{(2)})^2.
 \end{align*}
Hence, by \eqref{SEGED_Z2},
 \begin{align*}
    \frac{1}{n}&\EE\left(\sum_{k=1}^n(Z_k^{(1)}-\alpha Z_{k-1}^{(1)})Z_{k-1}^{(2)}\right)^2 \\
      & \leq \frac{2}{n}
           \sum_{k=0}^n\Bigg[\Big(
              \theta_1^2\alpha^{2k}+\theta_1\alpha^k(1-\alpha^k)
              +\alpha^2\theta_1^2\alpha^{2(k-1)}
              +\alpha^2\theta_1\alpha^{k-1}(1-\alpha^{k-1})\Big)\\
       &\phantom{\leq \frac{2}{n}\sum_{k=1}^n\quad }\times
              \Big(\theta_2^2\alpha^{2(k-1)}+\theta_2\alpha^{k-1}(1-\alpha^{k-1})\Big) \Bigg] \\
       & \leq \frac{2(\theta_1^2+\theta_1)}{n}
             \Big(\theta_2^2\frac{\alpha^{2n}-1}{\alpha^2-1}
                  +\theta_2\frac{\alpha^{n}-1}{\alpha-1}
                  +\theta_2^2\frac{\alpha^{2n}-1}{\alpha^2-1}
                  +\alpha\theta_2\frac{\alpha^{n}-1}{\alpha-1}\Big)
    \to 0
   \qquad \text{as \ $n\to\infty$.}
 \end{align*}
Similarly one can check \eqref{SEGED86}.

Moreover, conditionally on the values \ $Y_{s_1-1}$ \ and \ $Y_{s_2-1}$,
 \ by \eqref{SEGED84}, \eqref{Strong_consistency30} and \eqref{CONVERGENCE22},
 \begin{align*}
      \begin{bmatrix}
           \sqrt{n}\big(\ttheta_{1,n} - \lim_{k\to\infty} \ttheta_{1,k}\big) \\
           \sqrt{n}\big(\ttheta_{2,n} - \lim_{k\to\infty} \ttheta_{2,k}\big) \\
      \end{bmatrix}
  &= \sqrt{n}
      \begin{bmatrix}
         \ttheta_{1,n} - (Y_{s_1}-\alpha Y_{s_1-1}-\mu_\vare) \\
         \ttheta_{2,n} - (Y_{s_2}-\alpha Y_{s_2-1}-\mu_\vare) \\
      \end{bmatrix}
   = \sqrt{n}
      \begin{bmatrix}
         -(\talpha_n-\alpha)Y_{s_1-1} \\
         -(\talpha_n-\alpha)Y_{s_2-1} \\
      \end{bmatrix}  \\
    & \distr
         \cN\left(
               \begin{bmatrix}
                 0 \\
                 0 \\
              \end{bmatrix}
               , \sigma_{\alpha,\,\vare}^2
                 \begin{bmatrix}
                   Y_{s_1-1}^2 & Y_{s_1-1}Y_{s_2-1} \\
                   Y_{s_1-1}Y_{s_2-1} & Y_{s_2-1}^2 \\
                 \end{bmatrix}
              \right)
   \qquad \text{as \ $n\to\infty$.}
 \end{align*}
\proofend

\subsection{Two outliers, estimation of the mean of the offspring and innovation distributions
            and the outliers' sizes}

We assume that \ $I=2$ \ and that the relevant time points \ $s_1$, $s_2\in\NN$,
 \ $s_1\ne s_2$, \ are known.
We concentrate on the CLS estimation of \ $(\alpha,\mu_\vare,\theta_1,\theta_2)$.
\ We have
 \begin{align*}
   \EE(Y_k\mid\cF^Y_{k-1})
       = \alpha Y_{k-1} + \mu_\vare + \delta_{k,s_1}\theta_1+\delta_{k,s_2}\theta_2,
      \qquad k\in\NN.
 \end{align*}
Hence for all \ $n\geq\max(s_1,s_2)$, $n\in\NN$,
 \begin{align*}
   \begin{split}
    \sum_{k=1}^n&\big(Y_k-\EE(Y_k\mid \cF^Y_{k-1})\big)^2\\
        & = {\sum_{k=1}^n}^{(s_1,s_2)} \big(Y_k-\alpha Y_{k-1}-\mu_\vare\big)^2
           + \big(Y_{s_1}-\alpha Y_{s_1-1}-\mu_\vare-\theta_1\big)^2
           + \big(Y_{s_2}-\alpha Y_{s_2-1}-\mu_\vare-\theta_2\big)^2.
    \end{split}
 \end{align*}
For all  \ $n\geq\max(s_1,s_2)$, $n\in\NN$, \ we define the function \ $Q_n:\RR^{n+1}\times\RR^4\to\RR$, \ as
 \begin{align*}
    Q_n&({\bf y}_n;\alpha',\mu_\vare',\theta_1',\theta_2')\\
      & :={\sum_{k=1}^n}^{(s_1,s_2)} \big(y_k-\alpha' y_{k-1}-\mu_\vare'\big)^2
           + \big(y_{s_1}-\alpha' y_{s_1-1}-\mu_\vare'-\theta_1'\big)^2
           + \big(y_{s_2}-\alpha' y_{s_2-1}-\mu_\vare'-\theta_2'\big)^2,
 \end{align*}
  for all \ ${\bf y}_n\in\RR^{n+1}$, $\alpha',\mu_\vare',\theta_1',\theta_2'\in\RR$.
\ By definition, for all \ $n\geq\max(s_1,s_2)$, \ a CLS estimator for
 the parameter \ $(\alpha,\mu_\vare,\theta_1,\theta_2)\in(0,1)\times(0,\infty)\times\NN^2$ \ is
 a measurable function \ $(\halpha_n,\hmuen,\htheta_{1,n},\htheta_{2,n}):S_n\to\RR^4$
 \ such that
 \begin{align*}
   Q_n({\bf y}_n;\,\halpha_n({\bf y}_n), &\hmuen({\bf y}_n),
             \htheta_{1,n}({\bf y}_n), \htheta_{2,n}({\bf y}_n))\\
       &= \inf_{(\alpha',\mu_\vare',\theta_1',\theta_2')\in\RR^4}
           Q_n({\bf y}_n;\alpha',\mu_\vare',\theta_1',\theta_2')
       \qquad \forall\;\;  {\bf y}_n\in S_n,
 \end{align*}
 where \ $S_n$ \ is suitable subset of \ $\RR^{n+1}$ \ (defined in the proof of
 Lemma \ref{LEMMA18}).
We note that we do not define the CLS estimator
 \ $(\halpha_n,\hmuen,\htheta_{1,n},\htheta_{2,n})$ \ for all samples \ ${\bf y}_n\in \RR^{n+1}$.
\ For all \ ${\bf y}_n\in\RR^{n+1}$, $\alpha',\mu_\vare',\theta_1',\theta_2'\in\RR$,
 \begin{align*}
   & \frac{\partial Q_n}{\partial \alpha'}({\bf y}_n;\alpha',\mu_\vare',\theta_1',\theta_2')\\
   &\qquad\qquad
      =-2{\sum_{k=1}^n}^{(s_1,s_2)} \big(y_k-\alpha' y_{k-1}-\mu_\vare'\big)y_{k-1}
       - 2\big(y_{s_1}-\alpha' y_{s_1-1}-\mu_\vare'-\theta_1'\big)y_{s_1-1} \\
   &\phantom{\qquad\qquad=\;\,}
       - 2\big(y_{s_2}-\alpha' y_{s_2-1}-\mu_\vare'-\theta_2'\big)y_{s_2-1}, \\
   & \frac{\partial Q_n}{\partial \mu_\vare'}({\bf y}_n;\alpha',\mu_\vare',\theta_1',\theta_2')
      =-2{\sum_{k=1}^n}^{(s_1,s_2)} \big(y_k-\alpha' y_{k-1}-\mu_\vare'\big)
             - 2\big(y_{s_1}-\alpha' y_{s_1-1}-\mu_\vare'-\theta_1'\big)\\
   &\phantom{\frac{\partial Q_n}{\partial \alpha'}({\bf y}_n;\alpha',\mu_\vare',\theta_1',\theta_2')=\;}
             - 2\big(y_{s_2}-\alpha' y_{s_2-1}-\mu_\vare'-\theta_2'\big), \\
   & \frac{\partial Q_n}{\partial \theta_1'}({\bf y}_n;\alpha',\mu_\vare',\theta_1',\theta_2')
      =-2\big(y_{s_1}-\alpha' y_{s_1-1}-\mu_\vare'-\theta_1'\big),\\
   & \frac{\partial Q_n}{\partial \theta_2'}({\bf y}_n;\alpha',\mu_\vare',\theta_1',\theta_2')
      =-2\big(y_{s_2}-\alpha' y_{s_2-1}-\mu_\vare'-\theta_2'\big).
 \end{align*}

The next lemma is about the existence and uniqueness of the CLS estimator
 of \ $(\alpha,\mu_\vare,\theta_1,\theta_2)$.

\begin{Lem}\label{LEMMA18}
There exist subsets \ $S_n\subset\RR^{n+1}$, $n\geq \max(s_1,s_2)$ \ with the following properties:
 \begin{enumerate}
  \item[\upshape{(i)}]
   there exists a unique CLS estimator
   \ $(\halpha_n,\hmuen,\htheta_{1,n},\htheta_{2,n}):S_n\to\RR^4$,
  \item[\upshape{(ii)}]
   for all \ ${\bf y}_n\in S_n$, \ the system of equations
  \begin{align}\label{Innovation_CLSE_EQ4}
   \begin{split}
    &\frac{\partial Q_n}{\partial \alpha'}
      ({\bf y}_n;\alpha',\mu_\vare',\theta_1',\theta_2')=0, \qquad\qquad
    \frac{\partial Q_n}{\partial \mu_\vare'}
      ({\bf y}_n;\alpha',\mu_\vare',\theta_1',\theta_2')=0,\\
    &\frac{\partial Q_n}{\partial \theta_1'}
       ({\bf y}_n;\alpha',\mu_\vare',\theta_1',\theta_2')=0, \qquad\qquad
    \frac{\partial Q_n}{\partial \theta_2'}
       ({\bf y}_n;\alpha',\mu_\vare',\theta_1',\theta_2')=0,
   \end{split}
 \end{align}
 has the unique solution
  \begin{align}\label{SEGED65}
    &\halpha_n({\bf y}_n)
     =\frac{ (n-2)\DS\sumssd y_{k-1} y_k
            - \DS\sumssd y_k \DS\sumssd y_{k-1}}
           {D_n({\bf y}_n)},\\[2mm]\label{SEGED66}
   &\hmuen({\bf y}_n)
       = \frac{\DS\sumssd y_{k-1}^2 \DS\sumssd y_k
          - \DS\sumssd y_{k-1} \DS\sumssd y_{k-1}y_k}
          {D_n({\bf y}_n)},
             \\[2mm]\label{SEGED67}
   &\htheta_{i,n}({\bf y}_n)
      = y_{s_i}-\halpha_n({\bf y}_n) y_{s_i-1}-\hmuen({\bf y}_n),\qquad i=1,2,
 \end{align}
   where
 \[
   D_n({\bf y}_n):=(n-2)\sumssd y_{k-1}^2 -\left(\sumssd y_{k-1}\right)^2,
     \qquad n\geq \max(s_1,s_2),
 \]
  \item[\upshape{(iii)}]
  ${\bf Y}_n\in S_n$ \ holds asymptotically as \ $n\to\infty$ \ with probability one.
 \end{enumerate}
\end{Lem}

\noindent{\bf Proof.}
One can easily check that the unique solution of the system of equations \eqref{Innovation_CLSE_EQ4}
 takes the form \eqref{SEGED65}-\eqref{SEGED66}-\eqref{SEGED67} whenever \ $D_n({\bf y}_n)>0$.

For all \ $n\geq\max(s_1,s_2)$, \ let
 \[
   S_n := \left\{{\bf y}_n\in\RR^{n+1} : D_n({\bf y}_n)>0,
                    \;\Delta_{i,n}({\bf y}_n;\alpha',\mu_\vare',\theta_1',\theta_2') > 0,
                   \;\; i=1,2,3,4,\; \forall\, (\alpha',\mu_\vare',\theta_1',\theta_2')\in\RR^4 \right\},
 \]
 where \ $\Delta_{i,n}({\bf y}_n;\alpha',\mu_\vare',\theta_1',\theta_2')$, $i=1,2,3,4$, \
 denotes the $i$-th order leading principal minor of the \ $4\times 4$ \ matrix
 \begin{align*}
  H_n({\bf y}_n;\alpha',\mu_\vare',\theta_1',\theta_2')
    : =
    \begin{bmatrix}
    \frac{\partial^2 Q_n}{\partial(\alpha')^2}
    & \frac{\partial^2 Q_n}{\partial\mu_\vare'\partial\alpha'}
    & \frac{\partial^2 Q_n}{\partial\theta_1'\partial\alpha'}
    & \frac{\partial^2 Q_n}{\partial\theta_2'\partial\alpha'} \\
    \frac{\partial^2 Q_n}{\partial\alpha'\partial\mu_\vare'}
    & \frac{\partial^2 Q_n}{\partial(\mu_\vare')^2}
    & \frac{\partial^2 Q_n}{\partial\theta_1'\partial\mu_\vare'}
    & \frac{\partial^2 Q_n}{\partial\theta_2'\partial\mu_\vare'} \\
    \frac{\partial^2 Q_n}{\partial\alpha'\partial\theta_1'}
    & \frac{\partial^2 Q_n}{\partial\mu_\vare'\partial\theta_1'}
    & \frac{\partial^2 Q_n}{\partial(\theta_1')^2}
    & \frac{\partial^2 Q_n}{\partial\theta_2'\partial \theta_1'} \\
    \frac{\partial^2 Q_n}{\partial\alpha'\partial\theta_2'}
    & \frac{\partial^2 Q_n}{\partial\mu_\vare'\partial\theta_2'}
    & \frac{\partial^2 Q_n}{\partial\theta_1'\partial\theta_2'}
    & \frac{\partial^2 Q_n}{\partial(\theta_2')^2} \\
   \end{bmatrix}
     ({\bf y}_n;\alpha',\mu_\vare',\theta_1',\theta_2').
 \end{align*}

Then the function
  \ $\RR^4 \ni (\alpha',\mu_\vare',\theta_1',\theta_2')
     \mapsto Q_n({\bf y}_n;\alpha',\mu_\vare',\theta_1',\theta_2')$
  \ is strictly convex for all \ ${\bf y}_n\in S_n$,
 \ see, e.g., Berkovitz \cite[Theorem 3.3, Chapter III]{Ber}.

Since the function \ $\RR^4 \ni (\alpha',\mu_\vare',\theta_1',\theta_2')
     \mapsto Q_n({\bf y}_n;\alpha',\mu_\vare',\theta_1',\theta_2')$
  \ is strictly convex for all \ ${\bf y}_n\in S_n$ \ and the system of equations
 \eqref{Innovation_CLSE_EQ4} has a unique solution for all \ ${\bf y}_n\in S_n$,
 we get the function in question attains its (global) minimum at this unique solution,
 which yields (i) and (ii).

Further, for all \ ${\bf y}_n\in\RR^{n+1}$ \ and
 \ $(\alpha',\mu_\vare',\theta_1',\theta_2')\in\RR^4$,
 \begin{align*}
   &\frac{\partial^2 Q_n}{\partial(\alpha')^2}({\bf y}_n;\alpha',\mu_\vare',\theta_1',\theta_2')
       = 2\sumssd y_{k-1}^2 + 2y_{s_1-1}^2 + 2y_{s_2-1}^2
       =2\sum_{k=1}^n y_{k-1}^2,\\
   &\frac{\partial^2 Q_n}{\partial(\mu_\vare')^2}({\bf y}_n;\alpha',\mu_\vare',\theta_1',\theta_2')
       =2n,\\
   &\frac{\partial^2 Q_n}{\partial\alpha'\partial\mu_\vare'}({\bf y}_n;\alpha',\mu_\vare',\theta_1',\theta_2')
       =\frac{\partial^2 Q_n}{\partial\mu_\vare'\partial\alpha'}({\bf y}_n;\alpha',\mu_\vare',\theta_1',\theta_2')
       = 2\sum_{k=1}^n y_{k-1},\\
    &\frac{\partial^2 Q_n}{\partial\theta_1'\partial\mu_\vare'}({\bf y}_n;\alpha',\mu_\vare',\theta_1',\theta_2')
       =\frac{\partial^2 Q_n}{\partial\mu_\vare'\partial\theta_1'}({\bf y}_n;\alpha',\mu_\vare',\theta_1',\theta_2')
       = 2,\\
  &\frac{\partial^2 Q_n}{\partial\theta_2'\partial\mu_\vare'}({\bf y}_n;\alpha',\mu_\vare',\theta_1',\theta_2')
       =\frac{\partial^2 Q_n}{\partial\mu_\vare'\partial\theta_2'}({\bf y}_n;\alpha',\mu_\vare',\theta_1',\theta_2')
       = 2,\\
   &\frac{\partial^2 Q_n}{\partial\alpha'\partial\theta_1'}({\bf y}_n;\alpha',\mu_\vare',\theta_1',\theta_2')
       =\frac{\partial^2 Q_n}{\partial\theta_1'\partial\alpha'}({\bf y}_n;\alpha',\mu_\vare',\theta_1',\theta_2')
       = 2y_{s_1-1},\\
   &\frac{\partial^2 Q_n}{\partial\alpha'\partial\theta_2'}({\bf y}_n;\alpha',\mu_\vare',\theta_1',\theta_2')
       =\frac{\partial^2 Q_n}{\partial\theta_2'\partial\alpha'}({\bf y}_n;\alpha',\mu_\vare',\theta_1',\theta_2')
       = 2y_{s_2-1},
 \end{align*}
 and
 \begin{align*}
   &\frac{\partial^2 Q_n}{\partial\theta_1'\partial\theta_2'}({\bf y}_n;\alpha',\mu_\vare',\theta_1',\theta_2')
       =\frac{\partial^2 Q_n}{\partial\theta_2'\partial\theta_1'}({\bf y}_n;\alpha',\mu_\vare',\theta_1',\theta_2')
       = 0,\\
   &\frac{\partial^2 Q_n}{\partial(\theta_1')^2}({\bf y}_n;\alpha',\mu_\vare',\theta_1',\theta_2')=2,
     \qquad
     \frac{\partial^2 Q_n}{\partial(\theta_2')^2}({\bf y}_n;\alpha',\mu_\vare',\theta_1',\theta_2')=2.
 \end{align*}
Then \ $H_n({\bf y}_n;\alpha',\mu_\vare',\theta_1',\theta_2')$ \ has the following leading principal minors
 \begin{align*}
  & \Delta_{1,n}({\bf y}_n;\alpha',\mu_\vare',\theta_1',\theta_2') = 2\sum_{k=1}^n y_{k-1}^2,\\
  & \Delta_{2,n}({\bf y}_n;\alpha',\mu_\vare',\theta_1',\theta_2')
          = 4\left(n\sum_{k=1}^n y_{k-1}^2 - \left(\sum_{k=1}^ny_{k-1}\right)^2\right),\\
  & \Delta_{3,n}({\bf y}_n;\alpha',\mu_\vare',\theta_1',\theta_2') =
      8\left((n-1)\sum_{k=1}^n y_{k-1}^2 - \left(\sum_{k=1}^ny_{k-1}\right)^2
              + 2y_{s_1-1}\sum_{k=1}^n y_{k-1} - n (y_{s_1-1})^2\right),\\
  & \Delta_{4,n}({\bf y}_n;\alpha',\mu_\vare',\theta_1',\theta_2') :=
     \det H_n({\bf y}_n;\alpha',\mu_\vare',\theta_1',\theta_2').
 \end{align*}
Note that \ $\Delta_{i,n}({\bf y}_n;\alpha',\mu_\vare',\theta_1',\theta_2')$, $i=1,2,3,4$, \ do not depend
 on  \ $(\alpha',\mu_\vare',\theta_1',\theta_2')$, \ and hence we will simply denote
 \ $\Delta_{i,n}({\bf y}_n;\alpha',\mu_\vare',\theta_1',\theta_2')$ \ by
 \ $\Delta_{i,n}({\bf y}_n)$.

Next we check that \ ${\bf Y}_n\in S_n$ \ holds asymptotically as \ $n\to\infty$ \ with probability one.
By \eqref{Ergodic1} and \eqref{Ergodic2}, using the very same arguments as in the proof of
 Lemma \ref{LEMMA15}, one can get
 \begin{align*}
   &\PP\left(\lim_{n\to\infty}\frac{\Delta_{1,n}({\bf Y}_n)}{n}
                 = 2 \EE\widetilde X^2 \right)=1,\\
   &\PP\left(\lim_{n\to\infty}\frac{\Delta_{2,n}({\bf Y}_n)}{n^2}
                 = 4 \var\widetilde X \right)=1,\\
   &\PP\left(\lim_{n\to\infty}\frac{\Delta_{3,n}({\bf Y}_n)}{n^2}
                 = 8 \var\widetilde X \right)=1,\\
   &\PP\left(\lim_{n\to\infty}\frac{\Delta_{4,n}({\bf Y}_n)}{n^2}
              = 16 \var\widetilde X \right)=1,
 \end{align*}
 where \ $\widetilde X$ \ denotes a random variable with the unique stationary distribution of
 the INAR(1) model in \eqref{INAR1}.
Hence
 \begin{align*}
   \PP\left(\lim_{n\to\infty}\Delta_{i,n}({\bf Y}_n)
                 = \infty \right)=1,\qquad
   i=1,2,3,4.
 \end{align*}
By \eqref{Ergodic1} and \eqref{Ergodic2}, we also get
 \begin{align}\label{LEMMA16_allitasa}
    \PP\left(\lim_{n\to\infty}\frac{D_n({\bf Y}_n)}{n^2}=\var\widetilde X\right)=1,
 \end{align}
 and hence \ $\PP(\lim_{n\to\infty}D_n({\bf Y}_n)=\infty)=1$.
\proofend

By Lemma \ref{LEMMA18},
\[
    (\halpha_n({\bf Y}_n),
    \hmuen({\bf Y}_n),
    \htheta_{1,n}({\bf Y}_n),
    \htheta_{2,n}({\bf Y}_n)
 \]
 exists uniquely asymptotically as \ $n\to\infty$ \ with probability one.
In the sequel we will simply denote it by
 \ $(\halpha_n,\hmuen,\htheta_{1,n},\htheta_{2,n})$, \ and we will also denote
 \ $D_n({\bf Y}_n)$  \ by \ $D_n$.

The next result shows that \ $\halpha_n$ \ and \ $\hmuen$ \ are strongly consistent estimators
 of \ $\alpha$ \ and \ $\mu_\vare$, \ respectively, whereas \ $\ttheta_{i,n}$, \ $i=1,2$, \ fail to be
 strongly consistent estimators of \ $\theta_{i,n}$, \ $i=1,2$, \  respectively.

\begin{Thm}\label{THEOREM9}
Consider the CLS estimators
 \ $(\halpha_n,\hmuen,\htheta_{1,n},\htheta_{1,n})_{n\in\NN}$ \ of the parameter
 \ $(\alpha,\mu_\vare,\theta_1,\theta_2)\in(0,1)\times(0,\infty)\times\NN^2$.
\ The sequences \ $(\halpha_n)_{n\in\NN}$ \ and \ $(\hmuen)_{n\in\NN}$
 \ are strongly consistent for all
 \ $(\alpha,\mu_\vare,\theta_1,\theta_2)\in(0,1)\times(0,\infty)\times\NN^2$, \ i.e.,
 \begin{align}\label{Strong_consistency26}
   &\PP(\lim_{n\to\infty}\halpha_n=\alpha)=1,
     \qquad \forall\;(\alpha,\mu_\vare,\theta_1,\theta_2)\in(0,1)\times(0,\infty)\times\NN^2,\\
     \label{Strong_consistency27}
   &\PP(\lim_{n\to\infty}\hmuen=\mu_\vare)=1,
     \qquad \forall\; (\alpha,\mu_\vare,\theta_1,\theta_2)\in(0,1)\times(0,\infty)\times\NN^2,
 \end{align}
 whereas the sequences \ $(\htheta_{1,n})_{n\in\NN}$ \ and \ $(\htheta_{2,n})_{n\in\NN}$ \
 are not strongly consistent for any
 \ $(\alpha,\mu_\vare,\theta_1,\theta_2)\in(0,1)\times(0,\infty)\times\NN^2$, \ namely,
 \begin{align}\label{Strong_consistency31}
  \PP\left(\lim_{n\to\infty}\htheta_{i,n} = Y_{s_i}-\alpha Y_{s_i-1}-\mu_\vare\right)=1
 \end{align}
 for all \ $(\alpha,\mu_\vare,\theta_1,\theta_2)\in(0,1)\times(0,\infty)\times\NN^2$
  \ and \ $i=1,2$.
\end{Thm}

\noindent{\bf Proof.}
To prove \eqref{Strong_consistency26} and \eqref{Strong_consistency27},
 using Proposition \ref{LEMMA3_Decomposition2} and the proof of Theorem \ref{THEOREM8},
 it is enough to check that
 \begin{align*}
   & \PP\left(\lim_{n\to\infty}
              \frac{1}{n}\sumssd(Z^{(1)}_{k-1}+Z^{(2)}_{k-1})(Z^{(1)}_k+Z^{(2)}_k)=0\right)=1,\\
   & \PP\left(\lim_{n\to\infty}
              \frac{1}{n^2}\sumssd(Z^{(1)}_{k-1}+Z^{(2)}_{k-1})^2=0\right)=1.
 \end{align*}
 The above relations follows by \eqref{SEGED89}.

By \eqref{SEGED67} and \eqref{Strong_consistency26}, \eqref{Strong_consistency27},
 we have \eqref{Strong_consistency31}.
\proofend

The asymptotic distribution of the CLS estimation is given in the next theorem.

\begin{Thm}
Under the additional assumptions \ $\EE Y_0^3<\infty$ \ and \ $\EE\vare_1^3<\infty$, \ we have
 \begin{align}\label{CONVERGENCE20}
      \begin{bmatrix}
        \sqrt{n}(\halpha_n-\alpha) \\
        \sqrt{n}(\hmuen-\mu_\vare) \\
      \end{bmatrix}
     \distr
      \cN\left( \begin{bmatrix}
                  0 \\
                  0 \\
                \end{bmatrix}
               , \;B_{\alpha,\vare}
      \right)
   \qquad \text{as \ $n\to\infty$,}
 \end{align}
 where \ $B_{\alpha,\vare}$ \ is  defined in \eqref{SEGED_BALPHA}.
Moreover, conditionally on the values \ $Y_{s_1-1}$ \ and \ $Y_{s_2-1}$,
 \begin{align}\label{CONVERGENCE21}
   \begin{bmatrix}
     \sqrt{n}\big(\htheta_{1,n}-\lim_{k\to\infty}\htheta_{1,k}\big) \\
      \sqrt{n}\big(\htheta_{2,n}-\lim_{k\to\infty}\htheta_{2,k}\big) \\
   \end{bmatrix}
    \distr
    \cN\left(\begin{bmatrix}
               0 \\
               0 \\
             \end{bmatrix},
             C_{\alpha,\vare} B_{\alpha,\vare} C_{\alpha,\vare}^\top
    \right)
   \qquad \text{as \ $n\to\infty$,}
 \end{align}
 where
 \[
   C_{\alpha,\vare}
      :=\begin{bmatrix}
          Y_{s_1-1} & 1 \\
          Y_{s_2-1} & 1 \\
        \end{bmatrix}.
 \]
\end{Thm}

\noindent{\bf Proof.}
Using Proposition \ref{LEMMA3_Decomposition2}, the proof of Theorem \ref{Proposition2},
and \eqref{SEGED85}, \eqref{SEGED86}, one can obtain \eqref{CONVERGENCE20}.
By \eqref{SEGED67} and \eqref{Strong_consistency31},
\begin{align*}
   \begin{bmatrix}
     \sqrt{n}\big(\htheta_{1,n}-\lim_{k\to\infty}\htheta_{1,k}\big) \\
     \sqrt{n}\big(\htheta_{2,n}-\lim_{k\to\infty}\htheta_{2,k}\big) \\
   \end{bmatrix}
   & = \sqrt{n}
      \begin{bmatrix}
       Y_{s_1}-\halpha_n Y_{s_1-1} -\hmuen - (Y_{s_1}-\alpha Y_{s_1-1}-\mu_\vare) \\
       Y_{s_2}-\halpha_n Y_{s_2-1} -\hmuen - (Y_{s_2}-\alpha Y_{s_2-1}-\mu_\vare) \\
      \end{bmatrix}\\
   & = \begin{bmatrix}
       -Y_{s_1-1} & -1 \\
       -Y_{s_2-1} & -1 \\
      \end{bmatrix}
      \begin{bmatrix}
        \sqrt{n}(\halpha_n-\alpha) \\
        \sqrt{n}(\hmuen-\mu_\vare) \\
      \end{bmatrix}
\end{align*}
 holds asymptotically as \ $n\to\infty$ \ with probability one.
Using \eqref{CONVERGENCE20} we obtain \eqref{CONVERGENCE21}.
\proofend

\section{Appendix}\label{Appendix}

\numberwithin{equation}{section}

\begin{Lem2}\label{LEMMA_UNIQUE_STATDISTR}
If \ $\alpha\in(0,1)$ \ and \ $\EE\vare_1<\infty$,
\ then the INAR(1) model in \eqref{INAR1} has a unique stationary distribution.
\end{Lem2}

\noindent{\bf Proof.}
We follow the train of thoughts given in Section 6.3
 in Hall and Heyde \cite{HalHey}, but we also complete the proof given there.
For all \ $n\in\ZZ_+$, \ let \ $P_n$ \ denote the probability generating function
 of \ $X_n$, \ i.e., \ $P_n(s):=\EE s^{X_n}$, $\vert s\vert\leq 1$, $s\in\CC$.
Let \ $A$ \ and \ $B$ \ be the probability generating function of the offspring \ ($\xi_{1,1}$) \
 and the innovation \ ($\vare_1$) \ distribution, respectively.
With the notation
 \[
   A^{(k)}(s):=(\underbrace{A\circ\cdots\circ A}_\text{$k$-times})(s),
      \qquad  \vert s\vert\leq 1,\; s\in\CC, \; k\in\NN,
 \]
 we get for all \ $\vert s\vert\leq 1$, $s\in\CC$, \ and \ $n\in\NN$,
  \begin{align*}
   P_n(s) & = \EE(\EE(s^{X_n}\mid \cF_{n-1}^X))
            = \EE\Big[\EE\big(s^{\sum_{j=1}^{X_{n-1}}\xi_{n,j}}\mid \cF_{n-1}^X\big)
                     \EE(s^{\vare_n}\mid \cF_{n-1}^X) \Big] \\
          & = \EE(A(s)^{X_{n-1}} B(s)) = P_{n-1}(A(s)) B(s).
  \end{align*}
By iteration, we have
 \begin{align}\label{SEGED101}
  \begin{split}
   P_n(s) & = P_{n-1}(A(s))B(s) = P_{n-2}((A\circ A)(s))B(A(s))B(s)
            = \cdots  \\
          & = P_0(A^{(n)}(s)) B(s) \prod_{k=1}^{n-1} B(A^{(k)}(s)),
              \qquad \vert s\vert\leq 1,\; s\in\CC, \; n\in\NN.
  \end{split}
 \end{align}
We check that \ $\lim_{n\to\infty} P_0(A^{(n)}(s))=P_0(1)=1$, \ $s\in\CC$, \ and
 verify that the sequence \ $\prod_{k=1}^n B(A^{(k)}(s))$, $n\in\NN$, \ is convergent
 for all \ $s\in[0,1]$.
\ By iteration, for all \ $n\in\NN$,
 \begin{align}\label{SEGED100}
  \begin{split}
   A^{(n)}(s) & = A^{(n-1)}(1-\alpha+\alpha s)
                = A^{(n-2)}(1-\alpha +\alpha(1-\alpha+\alpha s)) \\
              &= A^{(n-2)}(1-\alpha +\alpha(1-\alpha)+\alpha^2 s)
                = \cdots
                = (1-\alpha)\sum_{k=0}^{n-1} \alpha^k +\alpha^n s \\
               & =(s-1)\alpha^n + 1,
  \end{split}
 \end{align}
 and hence \ $\lim_{n\to\infty} A^{(n)}(s)=1$, \ $s\in\CC$.
\ Then \ $\lim_{n\to\infty} P_0(A^{(n)}(s))=P_0(1)=1$, \ $s\in\CC$.
Since \ $0\leq B(v)\leq 1$, \ $v\in[0,1]$, $v\in\RR$, \ we get
 for all \ $s\in[0,1]$, \ the sequence \ $\prod_{k=1}^n B(A^{(k)}(s))$, $n\in\NN$, \
 is nonnegative and monotone decreasing and hence convergent.

\noindent We will use the following known theorem
 (see, e.g., Chung \cite[Section I.6, Theorem 4 and Section I.7, Theorem 2]{Chu}.
Let \ $(\xi_n)_{n\in\ZZ_+}$ \ be a homogeneous Markov chain with state space \ $I$.
\ Let us suppose that there exists some subset \ $D$ \ of \ $I$ \ such that \ $D$ \ is an
 essential, aperiodic class and \ $I\setminus D$ \ is a subset of inessential states.
Then either
 \begin{enumerate}
    \item[(a)] for all \ $i\in I$, $j\in D$ \  we have \ $\lim_{n\to\infty}p_{i,j}^{(n)}=0$, \ and
               therefore, there does not exist any stationary distribution,
 \end{enumerate}
  or
  \begin{enumerate}
    \item[(b)] for all \ $i,j\in D$ \  we have \ $\lim_{n\to\infty}p_{i,j}^{(n)}:=\pi_j>0$, \ and
               in this case the unique stationary distribution is given by \ $(\widetilde\pi_j)_{j\in I}$ \
               where \ $\widetilde\pi_j:=\pi_j$ \ if \ $j\in D$ \ and
               \ $\widetilde\pi_j:= 0$ if \ $j\in I\setminus D$.
 \end{enumerate}
Here \ $p_{i,j}^{(n)}$ \ denotes the $n$-step transition probability from the state \ $i$ \
 to the state \ $j$.

\noindent
Let us introduce the notation
 \[
    i_{\min}:=\min\Big\{ i\in\ZZ_+ : \PP(\vare_1=i)>0\Big\}.
 \]
Using that the offspring distribution is Bernoulli, i.e., it can take values \ $0$ \ and \ $1$,
 \ and both of them with positive probability, since \ $\alpha\in(0,1)$, \ one can think it over that
 the set of states \ $D:=\Big\{ i\in\ZZ_+ : i\geq i_{\min}\Big\}$ \ is an essential class.
Note also that \ $I\setminus D$ \ is a finite set of inessential states.
The class \ $D$ \ is aperiodic, since
 \begin{align*}
   p_{i_{\min}, i_{\min}} = \PP(X_{n+1}=i_{\min} \mid X_n=i_{\min})
                     \geq \PP(\vare_{n+1}=i_{\min}) (1-\alpha)^{i_{\min}}>0.
 \end{align*}
Note that if the additional assumption \ $\PP(\vare_1=0)>0$ \ is satisfied, then
 the Markov chain is irreducible and aperiodic.

\noindent Let us assume that there is no stationary distribution.
With the notation
 \begin{align*}
  \widetilde P_n (s):=\frac{P_n(s)}{s^{i_{\min}}}
                     =\sum_{k=0}^\infty s^k \PP(X_n=k+i_{\min}),
       \qquad s\in[0,1],
 \end{align*}
 we get for all \ $n\in\NN$,
 \begin{align*}
   \widetilde\PP_n(0)= \PP(X_n=i_{\min})=\sum_{j=0}^\infty \PP(X_n=i_{\min}\mid X_0=j)\PP(X_0=j)
         =\sum_{j=0}^\infty p_{j,i_{\min}}^{(n)} \PP(X_0=j).
 \end{align*}
Hence, by part (a) of the above recalled theorem, we get \ $\lim_{n\to\infty}p_{j,i_{\min}}^{(n)}=0$
 \ for all \ $j\in\ZZ_+$.
\ Then the dominated convergence theorem yields that
 \[
   \lim_{n\to\infty} \widetilde P_n(0)= 0.
 \]
However, we show that \ $\lim_{n\to\infty} \widetilde P_n(0)>0$, \  which is a contradiction.
Using that \ $\PP(\vare_1=i_{\min})>0$ \ and that
 \[
   \widetilde\PP_n(0) = P_0(1-\alpha^n) \PP(\vare_1=i_{\min}) \prod_{k=1}^{n-1} B(1-\alpha^k),
 \]
 we have it is enough to prove that the limit of the
 sequence \ $\prod_{k=1}^n B(A^{(k)}(0))=\prod_{k=1}^n B(1-\alpha^k)$, $n\in\NN$, \ is positive.
It is known that for this it is enough to verify that
 \[
    \sum_{k=1}^\infty (1 - B(A^{(k)}(0)) )
     =\sum_{k=1}^\infty (1-B(1-\alpha^k))
      \qquad
      \text{ is convergent,}
 \]
 see, e.g., Br\'emaud \cite[Appendix, Theorem 1.9]{Bre}.
Just as in Section 6.3 in Hall and Heyde \cite{HalHey}, we show that for all \ $s\in[0,1)$, \
 \ $\sum_{k=1}^\infty (1 - B(A^{(k)}(s)) )$ \ is convergent.
For all \ $k\in\NN$, \ $s\in[0,1)$,
 \[
   1 - B(A^{(k)}(s))
     = \frac{1- B(A^{(k)}(s))}{1 - A^{(k)}(s)}(1- A^{(k)}(s)),
 \]
 and, by mean value theorem,
 \begin{align*}
  \frac{1- B(A^{(k)}(s))}{1 - A^{(k)}(s)}
    = \frac{B(A^{(k)}(1))- B(A^{(k)}(s))}{A^{(k)}(1) - A^{(k)}(s)}
    = \frac{B'(\theta(s))(A^{(k)}(1)- A^{(k)}(s))}{A^{(k)}(1) - A^{(k)}(s)}
    = B'(\theta(s)),
 \end{align*}
 with some \ $\theta(s)\in(s,1)$.
\ Since
 \begin{align}\label{SEGED104}
    B'(s)=\EE(\vare_1 s^{\vare_1-1})
         =\sum_{k=1}^\infty k s^{k-1}\PP(\vare_1=k)
         \leq \sum_{k=1}^\infty k \PP(\vare_1=k)
          = \EE\vare_1,
     \qquad s\in[0,1],
 \end{align}
 we have
 \[
    \frac{1- B(A^{(k)}(s))}{1 - A^{(k)}(s)}
      \leq \EE\vare_1=\mu_\vare, \qquad s\in[0,1).
 \]
Furthermore, by \eqref{SEGED100}, we get
 \[
   1- A^{(k)}(s) \leq 1 - A^{(k)}(0) = \alpha^k, \qquad k\in\NN, \;\; s\in[0,1),
 \]
 and hence \ $1- B(A^{(k)}(s)) \leq \mu_\vare \alpha^k$ \ for all \ $k\in\NN$, \ $s\in[0,1)$.
Then
 \[
   \sum_{k=1}^\infty (1-B(A^{(k)}(s)))
     \leq \mu_\vare \sum_{k=1}^\infty \alpha^k =\frac{\mu_\vare\alpha}{1-\alpha}<\infty,
       \qquad s\in[0,1).
 \]

Let us denote by \ $\widetilde X$ \ a random variable on \ $(\Omega,\cA,P)$ \ with a stationary
 distribution for the Markov chain \ $(X_n)_{n\in\ZZ_+}$.
\ The, by the dominated convergence theorem and part (b) of the above recalled theorem,
 we have for all \ $j\in\ZZ_+$,
 \begin{align*}
   \lim_{n\to\infty} \PP(X_n=j)
     & = \lim_{n\to\infty} \sum_{i=0}^\infty \PP(X_n=j\mid X_0=i)\PP(X_0=i) \\
     & = \sum_{i=0}^\infty \big(\lim_{n\to\infty} \PP(X_n=j\mid X_0=i)\big)\PP(X_0=i) \\
     & = \PP(\widetilde X=j) \sum_{i=0}^\infty \PP(X_0=i)
       = \PP(\widetilde X=j),
 \end{align*}
 which yields that \ $X_n$ \ converges in distribution to \ $\widetilde X$ \ as \ $n\to\infty$.
\ By the continuity theorem for probability generating functions (see, e.g., Feller \cite[Section 11]{Fel}),
 we also have \ $\widetilde X$ \ has the probability generating function
  \begin{align}\label{SEGED99}
    P(s):= B(s)\prod_{k=1}^{\infty} B(A^{(k)}(s)),
    \qquad s\in(0,1).
 \end{align}
The uniqueness of the stationary distribution follows by part (b) of the above recalled theorem.
\proofend

\vskip0.5cm

{\normalsize \bf Proofs of formulae \eqref{STAC_MOMENT1}, \eqref{STAC_MOMENT2} and \eqref{STAC_MOMENT3}.}
Let us introduce the probability generating functions
 \begin{align*}
    A(s):=\EE s^\xi=1-\alpha+\alpha s, \qquad s>0,
 \end{align*}
 where \ $\xi$ \ is a random variable with Bernoulli distribution having parameter
 \ $\alpha\in(0,1)$ \ and
 \begin{align*}
    B(s):=\EE s^\vare=\sum_{k=0}^\infty \PP(\vare=k)s^k,
                \qquad s>0,
 \end{align*}
 where \ $\vare$ \ is a non-negative integer-valued random variable with the same distribution
 as \ $\vare_1$.
\ In what follows we suppose that \ $\EE\vare^3<\infty$.
\ Since \ $\alpha\in(0,1)$ \ and \ $\EE\vare<\infty$, \ by Lemma \ref{LEMMA_UNIQUE_STATDISTR},
 there exists a uniquely determined stationary distribution of the INAR(1) model in \eqref{INAR1}.
Let us denote by \ $\widetilde X$ \ a random variable with this unique stationary distribution.
Due to Hall and Heyde \cite[formula (6.38)]{HalHey} or by the proof
 of Lemma \ref{LEMMA_UNIQUE_STATDISTR}, the probability generating function of \ $\widetilde X$ \ takes the form
 \begin{align}\label{SEGED_3RDMOMENT1}
   P(s):=\EE s^{\widetilde X}
        = B(s)B(A(s))B(A(A(s)))\cdots = B(s)\prod_{k=1}^\infty B(A^{(k)}(s)), \qquad s\in(0,1),
 \end{align}
 where for all \ $k\in\NN$,
 \[
   A^{(k)}(s)=(\underbrace{A\circ\cdots\circ A}_\text{$k$-times})(s), \qquad  s\in(0,1).
 \]
Hence for all \ $s\in(0,1)$,
 \begin{align}\label{SEGED38}
  \begin{split}
   \log P(s)
       & =\log \EE s^{\widetilde X}
         = \log B(s) + \log B(A(s)) +\log B(A(A(s)))+\cdots \\
       & =\log B(s) + \sum_{k=1}^\infty \log B(A^{(k)}(s)).
   \end{split}
 \end{align}
Using that \ $\EE\vare^3<\infty$, \ by Abel's theorem
 (see, e.g., Br\'emaud \cite[Appendix, Theorems 1.2 and 1.3]{Bre}), we get
 \begin{align*}
   & \lim_{s\uparrow 1}\left(\frac{\dd}{\dd s}\log B(s)\right)
      = \lim_{s\uparrow 1} \frac{\EE(\vare s^{\vare-1})}{\EE s^{\vare}}
      = \EE \vare,\\
   &\lim_{s\uparrow 1}\left( \frac{\dd^2}{\dd s^2} \log B(s)\right)
      = \lim_{s\uparrow 1}
             \frac{\EE(\vare(\vare-1)s^{\vare-2})\EE s^{\vare}
              - (\EE\vare s^{\vare-1})^2}
             {(\EE s^{\vare})^2}
     = \EE(\vare(\vare-1))- (\EE\vare)^2,
 \end{align*}
 and
 \[
   \lim_{s\uparrow 1}\left( \frac{\dd^3}{\dd s^3} \log B(s)\right)
     =\lim_{s\uparrow 1} \frac{N(s)}{(\EE s^{\vare})^4},
 \]
 where
 \begin{align*}
    N(s) := & \EE(\vare(\vare-1)(\vare-2)s^{\vare-3})(\EE s^{\vare})^3
              + \EE(\vare(\vare-1)s^{\vare-2}) \EE(\vare s^{\vare-1})
               (\EE s^{\vare})^2\\
            & - \Big[\EE(\vare(\vare-1)s^{\vare-2})\EE s^{\vare}
                   -(\EE\vare s^{\vare-1})^2\Big]
               2\EE s^{\vare}\EE(\vare s^{\vare-1})\\
             &  - 2\EE(\vare s^{\vare-1}) \EE(\vare(\vare-1)s^{\vare-2})
              (\EE s^{\vare})^2, \qquad s\in(0,1).
 \end{align*}
\ Hence
 \begin{align*}
   \lim_{s\uparrow 1}\left( \frac{\dd^3}{\dd s^3} \log B(s) \right)
     & = \EE(\vare(\vare-1)(\vare-2))
         - \EE(\vare(\vare-1))\EE\vare
         - 2[\EE(\vare(\vare-1))-(\EE\vare)^2]\EE\vare\\
     & = \EE\vare^3 -3\EE\vare^2+2\EE\vare
         -3(\EE\vare^2-\EE\vare)\EE\vare
         +2(\EE\vare)^3.
 \end{align*}
Then
 \begin{align}\label{SEGED_3RDMOMENT2}
  & \EE\vare = \lim_{s\uparrow 1}\left( \frac{\dd}{\dd s} \log B(s) \right),\\ \label{SEGED_3RDMOMENT3}
  & \EE\vare^2 = \lim_{s\uparrow 1}\left( \frac{\dd^2}{\dd s^2} \log B(s) \right)
                        + \EE \vare + (\EE \vare)^2,\\ \label{SEGED_3RDMOMENT4}
  & \EE\vare^3 = \lim_{s\uparrow 1}\left( \frac{\dd^3}{\dd s^3} \log B(s)\right)
                        + 3\EE\vare^2 - 2\EE\vare
                        + 3\EE\vare(\EE\vare^2-\EE\vare)
                        - 2(\EE\vare)^3.
 \end{align}
By \eqref{SEGED38},
 \[
   \log P(s)
      = \log \EE s^{\widetilde X}
      = b(s) + \sum_{k=1}^\infty b(A^{(k)}(s)),\qquad s\in(0,1],
 \]
 where \ $b(s):=\log B(s)$, \ $s\in(0,1]$.
\ We show that
 \begin{align}\label{SEGED102}
    \lim_{s\uparrow 1}\left( \frac{\dd}{\dd s} \log P(s) \right)
       = b'(1) + \sum_{k=1}^\infty \left[b'(A^{(k)}(1)) \prod_{\ell=0}^{k-1}A'(A^{(\ell)}(1))\right].
 \end{align}
First we note that
 $$
   \frac{\dd}{\dd s} \log P(s)
      = b'(s) + \sum_{k=1}^\infty \left[b'(A^{(k)}(s)) \prod_{\ell=0}^{k-1}A'(A^{(\ell)}(s))\right]  ,
     \qquad s\in(0,1),
 $$
 and, by \eqref{SEGED100}, we get for all \ $k\in\NN$, \ the functions \ $b'(A^{(k)}(s))$, $s\in[0,1]$ \
 are well-defined.
We check that the functions \ $b'(A^{(k)}(s))$, $s\in[0,1]$, \ $k\in\NN$, \ are bounded with a
 common bound.
By \eqref{SEGED100}, we have
 \[
   A^{(k)}(s) = (s-1)\alpha^k +1 \in[1-\alpha^k,1],\qquad s\in[0,1],\quad k\in\NN,
 \]
 and hence \ $A^{(k)}(s)\in[1-\alpha,1]$, $s\in[0,1]$, \ $k\in\NN$.
\ Then, using \eqref{SEGED104}, we get
 \[
    b'(A^{(k)}(s)) = \frac{B'(A^{(k)}(s))}{B(A^{(k)}(s)) }
                    \leq \frac{\EE\vare}{B(1-\alpha)}<\infty,
                    \qquad s\in[0,1],\quad k\in\NN.
 \]
Using that \ $A'(s) = \alpha = \EE\xi$, $\forall$ $s>0$, \ and that
 \[
     \sum_{k=1}^\infty\alpha^k
      =\frac{\alpha}{1-\alpha}<\infty,
 \]
 the dominated convergence theorem and \eqref{SEGED_3RDMOMENT2} yield \eqref{SEGED102}.
Hence, since \ $A(1)=1$,
 \begin{align}\label{SEGED103}
  \begin{split}
     \lim_{s\uparrow 1}\left(  \frac{\dd}{\dd s} \log P(s) \right)
      & = \EE\vare + (\EE\vare)\EE\xi + (\EE\vare)(\EE\xi)^2+\cdots
        =\sum_{k=0}^\infty(\EE\vare)(\EE\xi)^k \\
      &=\frac{\EE\vare}{1-\EE\xi}
        =\frac{\mu_\vare}{1-\alpha}
        <\infty.
   \end{split}
  \end{align}
Just as we derived \eqref{SEGED_3RDMOMENT2}, but without supposing \ $\EE\widetilde X<\infty$,
 \ Abel's theorem yields that
 \[
    \EE\widetilde X = \lim_{s\uparrow 1}\left( \frac{\dd}{\dd s} \log P(s) \right).
 \]
By \eqref{SEGED103}, we get\ $\EE\widetilde X = \frac{\mu_\vare}{1-\alpha}$, \ which also shows that
 \ $\EE\widetilde X$ \ is finite.

\noindent Using that
 \[
   b''(s)=\frac{B''(s)B(s)-(B'(s))^2}{(B(s))^2},\qquad s\in(0,1),
 \]
 we get
 \[
    b''(A^{(k)}(s)) \leq \frac{\EE(\vare(\vare-1))\EE\vare}{(B(1-\alpha))^2}
                    <\infty,
          \qquad s\in[0,1],\quad k\in\NN.
 \]
Using also that \ $b''(1)=\EE(\vare(\vare-1))-(\EE\vare)^2$ \ and \ $A''(s)=0$, \ $s>0$,
 \ by the dominated convergence theorem, one can check that
 \begin{align*}
    \lim_{s\uparrow 1}\left(  \frac{\dd^2}{\dd s^2} \log P(s) \right)
       &= b''(1) + \sum_{k=1}^\infty b''(A^{(k)}(1))\left(\prod_{\ell=0}^{k-1}A'(A^{\ell}(1))\right)^2
       = b''(1) \sum_{k=0}^\infty (\EE\xi)^{2k}\\
       & = \frac{\EE(\vare(\vare-1))-(\EE\vare)^2}{1-(\EE\xi)^2}
         = \frac{\var\vare-\EE\vare}{1-\alpha^2}
         = \frac{\sigma_\vare^2-\mu_\vare}{1-\alpha^2},
 \end{align*}
 which implies that \ $\EE\widetilde X^2$ \ is finite and
 \begin{align*}
    \EE\widetilde X^2
      = \frac{\sigma_\vare^2-\mu_\vare}{1-\alpha^2}
        + \frac{\mu_\vare}{1-\alpha}
        + \frac{\mu_\vare^2}{(1-\alpha)^2}
      = \frac{\sigma_\vare^2+\alpha\mu_\vare}{1-\alpha^2}
         + \frac{\mu_\vare^2}{(1-\alpha)^2}.
 \end{align*}
By a similar argument, using that \ $\EE\vare^3<\infty$ \ and
 \[
   b'''(1) = \EE(\vare(\vare-1)(\vare-2)) - 3(\EE\vare)(\EE\vare(\vare-1)) + 2(\EE\vare)^3,
 \]
 we get
 \begin{align*}
  \lim_{s\uparrow 1}&\left(\frac{\dd^3}{\dd s^3} \log P(s) \right)
      = b'''(1) + \sum_{k=1}^\infty b'''(A^{(k)}(1))\left(\prod_{\ell=0}^{k-1}A'(A^{\ell}(1))\right)^3\\
     & = \frac{\EE(\vare(\vare-1)(\vare-2))-3(\EE\vare)(\EE\vare(\vare-1))+2(\EE\vare)^3}{1-(\EE\xi)^3} \\
     & = \frac{\EE\vare^3-3\EE\vare^2+2\EE\vare - 3\EE\vare(\EE\vare^2-\EE\vare)+2(\EE\vare)^3}{1-\alpha^3} \\
     & = \frac{\EE\vare^3-3(\sigma_\vare^2+\mu_\vare^2)+2\mu_\vare
        - 3\mu_\vare(\sigma_\vare^2+\mu_\vare^2-\mu_\vare)+2\mu_\vare^3}{1-\alpha^3}
       = \frac{\EE\vare^3-3\sigma_\vare^2(1+\mu_\vare)-\mu_\vare^3+2\mu_\vare}{1-\alpha^3},
 \end{align*}
 which implies that \ $\EE\widetilde X^3$ \ is finite and
 \begin{align*}
  \EE \widetilde X^3
     = &\frac{\EE\vare^3-3\sigma_\vare^2(1+\mu_\vare)-\mu_\vare^3+2\mu_\vare}{1-\alpha^3}
         + 3 \frac{\sigma_\vare^2+\alpha\mu_\vare}{1-\alpha^2}
         - 2 \frac{\mu_\vare}{1-\alpha} \\
       & + 3 \frac{\mu_\vare(\sigma_\vare^2+\alpha\mu_\vare)}{(1-\alpha)(1-\alpha^2)}
         + \frac{\mu_\vare^3}{(1-\alpha)^3}.
 \end{align*}
This yields \eqref{STAC_MOMENT3}.
One can also write \eqref{STAC_MOMENT3} in the following form
 \begin{align*}
   \EE \widetilde X^3
     = \frac{1}{1-\alpha^3}
           & \left[
               3\alpha^2(1-\alpha)\EE\widetilde X^2
             + 3\alpha^2\mu_\vare\EE\widetilde X^2
             + 3\alpha\EE\widetilde X(\sigma_\vare^2+\mu_\vare^2)
             + \EE\vare^3 +3\alpha(1-\alpha)\mu_\vare\EE\widetilde X \right.\\
           & \left. \; + \alpha(1-\alpha)(1-2\alpha) \EE\widetilde X \right].
 \end{align*}
\proofend

\begin{Lem2}\label{LEMMA_INNOVATIONAL_INDEPENDENCE}
Let \ $(X_n)_{n\in\ZZ_+}$ \ and \ $(Z_n)_{n\in\ZZ_+}$ \ be two (not necessarily homogeneous)
 Markov chains with state space \ $\ZZ_+$.
\ Let us suppose that \ $(X_n,Z_n)_{n\in\ZZ_+}$ \ is a Markov chain,
 \ $X_0$ \ and \ $Z_0$ \ are independent, and that
 for all \ $n\in\NN$ \ and \ $i,j,k,\ell\in\ZZ_+$ \ such that
 \ $\PP(X_{n-1}=k, \, Z_{n-1}=\ell)>0$,
 \begin{align*}
   \PP(X_n=i,Z_n=j\mid X_{n-1}=k,Z_{n-1}=\ell)
             = \PP(X_n=i\mid X_{n-1}=k) \PP(Z_n=j\mid Z_{n-1}=\ell).
 \end{align*}
Then \ $(X_n)_{n\in\ZZ_+}$ \ and \ $(Z_n)_{n\in\ZZ_+}$ \ are independent.
\end{Lem2}

\noindent{\bf Proof.}
For all \ $n\in\NN$ \ and \ $i_0,i_1,\ldots,i_n,j_0,j_1,\ldots,j_n\in\ZZ_+$, \ we get
 \begin{align*}
   \PP&(X_n=i_n,\ldots,X_0=i_0,Z_n=j_n,\ldots,Z_0=j_0) \\
    & = \PP(X_n=i_n,Z_n=j_n\mid X_{n-1}=i_{n-1},Z_{n-1}=j_{n-1})
        \cdots \PP(X_1=i_1,Z_1=j_1\mid X_0=i_0,Z_0=j_0)\\
    &\phantom{=\;}
        \times \PP(X_0=i_0,Z_0=j_0) \\
    & = \PP(X_n=i_n \mid X_{n-1}=i_{n-1})\cdots \PP(X_1=i_1\mid X_0=i_0) \PP(X_0=i_0)\\
    &\phantom{=\;}
       \times \PP(Z_n=j_n\mid Z_{n-1}=j_{n-1}) \cdots \PP(Z_1=j_1\mid Z_0=j_0)\PP(Z_0=j_0) \\
    & = \PP(X_n=i_n,\ldots,X_0=i_0)\PP(Z_n=j_n,\ldots,Z_0=j_0),
 \end{align*}
 which yields that \ $X_n,\ldots,X_0$ \ and \ $Z_n,\ldots,Z_0$ \ are independent.
One can think it over that this implies the statement.
\proofend

The following result can be found in several textbooks, see, e.g.,
 Theorem 3.6 in Bhattacharya and Waymire \cite[Chapter 0]{BhaWay}.
For completeness we give a proof.

\begin{Lem2}\label{LEMMA4_LP}
Let \ $(\xi_n)_{n\in\NN}$ \ be a sequence of random variables such that \ $\PP(\lim_{n\to\infty}\xi_n=0)=1$
 \ and \ $\{\xi_n^p:n\in\NN\}$ \ is uniformly integrable for some \ $p\in\NN$, \ i.e.,
 \ $\lim_{M\to\infty}\sup_{n\in\NN}\EE\left(\vert\xi_n\vert^p\bone_{\{\vert\xi_n\vert>M\}}\right)=0.$
 \ Then \ $\xi_n\lpmean 0$ \ as \ $n\to\infty$, \ i.e., \ $\lim_{n\to\infty}\EE\vert\xi_n\vert^p=0$.
\end{Lem2}

\noindent{\bf Proof.}
 For all \ $n\in\NN$ \ and \ $M>0$, \ we get
 \begin{align*}
  \EE\vert\xi_n\vert^p
    = \EE\left(\vert\xi_n\vert^p\bone_{\{\vert\xi_n\vert^p>M\}}\right)
      + \EE\left(\vert\xi_n\vert^p\bone_{\{\vert\xi_n\vert^p\leq M\}}\right)
  \leq \sup_{n\in\NN} \EE\left(\vert\xi_n\vert^p\bone_{\{\vert\xi_n\vert^p>M\}}\right)
      + \EE\left(\vert\xi_n\vert^p\bone_{\{\vert\xi_n\vert^p\leq M\}}\right).
 \end{align*}
By \ $\PP(\lim_{n\to\infty}\xi_n=0)=1$,
 \[
  \lim_{n\to\infty}\vert\xi_n(\omega)\vert^p\bone_{\{\vert\xi_n(\omega)\vert^p\leq M\}}=0,
    \qquad \forall\;\; \omega\in\Omega,
  \]
 and \ $\EE\left(\vert\xi_n\vert^p\bone_{\{\vert\xi_n\vert^p\leq M\}}\right)\leq M^p<\infty$
 \ for all \ $n\in\NN$. \
Hence, by dominated convergence theorem, we have
 \[
  \lim_{n\to\infty}\EE\left(\vert\xi_n\vert^p\bone_{\{\vert\xi_n\vert^p\leq M\}}\right) =0
 \]
 for all \ $M>0$.
\ Then
 \[
  \limsup_{n\to\infty} \EE\vert\xi_n\vert^p
      \leq \sup_{n\in\NN} \EE\left(\vert\xi_n\vert^p\bone_{\{\vert\xi_n\vert^p>M\}}\right),
   \qquad \forall\;\; M>0.
 \]
By the uniformly integrability of \ $\{\xi_n^p:n\in\NN\}$, \ we have
 \ $\lim_{M\to\infty}\sup_{n\in\NN} \EE\left(\vert\xi_n\vert^p\bone_{\{\vert\xi_n\vert^p>M\}}\right)
  =0$, \ which yields the assertion.
\proofend

\end{document}